\documentclass[11pt,a4paper,reqno]{amsart}

\usepackage{times}
\usepackage{mathrsfs,mathtools}
\usepackage[inline]{enumitem}
\usepackage{stmaryrd}
\usepackage[usenames,dvipsnames]{xcolor}
\usepackage{amsmath}
\usepackage{amsthm}
\usepackage{amssymb}
\usepackage{amsfonts}
\usepackage{subfiles}
\usepackage{dsfont}
\usepackage{hyperref}
\usepackage[T1]{fontenc}
\usepackage[english]{babel}
\usepackage{todonotes}
\usepackage{geometry}
\usepackage{xspace,cancel}
\usepackage[capitalize]{cleveref}
\usepackage{etoolbox}

\usepackage[numbers,sort&compress,longnamesfirst]{natbib}

\allowdisplaybreaks

\setlist[enumerate]{label=\textup{(\roman*)},itemindent=.85cm,leftmargin=0.cm}

\newcounter{Conditions_counter}
\newcounter{Preliminary_counter}
\makeatletter \@addtoreset{equation}{section}
\def\theequation{\thesection.\arabic{equation}}
\newtheorem{theorem}{Theorem}[section]
\newtheorem*{assumption_named*}{\assumptionName}
	\providecommand{\assumptionName}{}
	
\newtheorem{corollary}[theorem]{Corollary}
\newtheorem{lemma}[theorem]{Lemma}
\newtheorem{proposition}[theorem]{Proposition}
\newtheorem{definition}[theorem]{Definition}
\newtheorem{assumption}[theorem]{Assumption}

\theoremstyle{plain}
\newtheorem{remark}[theorem]{Remark}

\DeclareFontFamily{U}{mathx}{\hyphenchar\font45}
\DeclareFontShape{U}{mathx}{m}{n}{
<5><6><7><8><9><10>
<10.95><12><14.4><17.28><20.74><24.88>
mathx10
}{}
\DeclareSymbolFont{mathx}{U}{mathx}{m}{n}
\DeclareFontSubstitution{U}{mathx}{m}{n}
\DeclareMathAccent{\widebar}{0}{mathx}{"73}

\makeatletter
\def\namedlabel#1#2{\begingroup
    #2%
    \def\@currentlabel{#2}%
    \phantomsection\label{#1}\endgroup
}

\numberwithin{equation}{section}

\newcommand\numberthis{\addtocounter{equation}{1}\tag{\theequation}}

\definecolor{red}{rgb}{0.7,0.15,0.15}
\definecolor{green}{rgb}{0,0.5,0}
\definecolor{blue}{rgb}{0,0,0.7}
\hypersetup{	
				colorlinks, linkcolor={red},
				citecolor={green}, urlcolor={blue}
			}

\newtoggle{full} 
 \togglefalse{full}

\makeatletter
\newcommand \Dotfill {\leavevmode \leaders \hb@xt@ 6pt{\hss \hss }\hfill \kern \z@}
\makeatother

\makeatletter
\def\@tocline#1#2#3#4#5#6#7{\relax
  \ifnum #1>\c@tocdepth 
  \else
    \par \addpenalty\@secpenalty\addvspace{#2}%
    \begingroup \hyphenpenalty\@M
    \@ifempty{#4}{%
      \@tempdima\csname r@tocindent\number#1\endcsname\relax
    }{%
      \@tempdima#4\relax
    }%
    \parindent\z@ \leftskip#3\relax \advance\leftskip\@tempdima\relax
    \rightskip\@pnumwidth plus4em \parfillskip-\@pnumwidth
    #5\leavevmode\hskip-\@tempdima
      \ifcase #1
       \or\or \hskip 1.65em \or \hskip 3.3em \else \hskip 4.95em \fi%
      #6\nobreak\relax
    \Dotfill
    \hbox to\@pnumwidth{\@tocpagenum{#7}}\par
    \nobreak
    \endgroup
  \fi}
\makeatother

\makeatletter
\def\l@section{\@tocline{1}{0pt}{1pc}{}{\scshape}}
\renewcommand{\tocsection}[3]{%
\indentlabel{\@ifnotempty{#2}{\ignorespaces#1 #2.\hskip 0.7em}}#3}
\def\l@subsection{\@tocline{2}{0pt}{1pc}{5pc}{}}

\def\l@subsubsection{\@tocline{3}{0pt}{1pc}{7pc}{}}

\makeatother

\begin{document}
\title[A general and sharp regularity condition for integro-differential equations]{A general and sharp regularity condition for integro-differential equations with non-dominated measures}

\author[A. Saplaouras]{Alexandros Saplaouras}

\address{Department of Mathematics, National Technical University of Athens, Zografou Campus, 15780 Athens, Greece}
\email{alsapl@mail.ntua.gr}

\thanks{Tha author gratefully acknowledges the financial support from the Hellenic Foundation for Research and Innovation (H.F.R.I.) under the “2nd Call for H.F.R.I. Research Projects to support Post-Doctoral Researchers”
(Project Number: 235).}


\subjclass[2010]{}

\date{}

\begin{abstract}
The aim of this work is to present the regularity condition (also known in the literature as structure condition) an integro-differential operator may satisfy in order for the domination principle to hold for (sub-, \mbox{super-)} solutions of polynomial growth.
More precisely, the framework presented in \citet{Hollender2017}, in which power functions are used in order to determine the integrability conditions, is weakened by substituting the power functions with Young functions.  
The use of Young functions allows for sharp integrability conditions, which are crucial when one deals with limit theorems.
As an immediate application, it is considered the case of parabolic Hamilton--Jacobi--Bellman (HJB) operators, for which the regularity condition is satisfied and, consequently, the comparison principle as well.
The parabolic HJB operator presented in this work can be associated to second-order (decoupled) forward-backward stochastic differential equations with jumps. 
\end{abstract}

\maketitle \frenchspacing

\iftoggle{full}{\tableofcontents}{}

\vspace{-2em}
\section{Introduction}
The notion of viscosity solutions seems to be the suitable generalization of the classical one when one seeks the establishment of a connection between stochastic differential equations and integro-differential equations.
Indeed, it has long been understood that the solution of a forward-backward stochastic differential equation provides a viscosity solution of an integro-differential equation; \emph{e.g.}, the celebrated Feynman--Kac formula, or more generally see among others 
\citet{barles1997backward,pardoux1998backward,el1997backward,pardoux1997probabilistic} for semilinear (integro-)differential equations, \citet{peng1991probabilistic,pardoux1992backward,sow2004} for quasilinear parabolic PDEs,
\citet{PhamOptimalStoppingControlled} for Hamilton--Jacobi--Bellman (HJB) integro-differential equations associated to optimal stopping of controlled jump-diffusion processes,
\citet{neufeld2016nonlinear} for the connection between (non-)linear stochastic processes and the associated (non-)linear generators of integro-differential equations and 
\citet{soner2012wellposedness,KaziTani2BSDE_PIDE} for parabolic second order fully nonlinear (integro-)differential equations.   

The probabilistic representation of solutions of integro-differential equations especially allows for an alternative approach of their numerical approximations. 
To this direction there are numerous works; indicatively \citet{BouchardElieTouziDiscreteTime,cheridito2007second,PossamaiTanWeakApprox,Tan2013Splitting,bouchard2016numerical,fahim2011probabilistic} for (fully nonlinear) local cases and \citet{bouchard2008discrete} for semilinear nonlocal cases.
In such approximating problems the integro-differential equation under interest has to be well defined, in the sense that it admits a unique solution. 
However, to the best of author's knowledge, a careful inspection of the relevant literature reveals that the case of solutions of arbitrary growth cannot be handled in the nonlocal case, especially when a family of non-dominated measures comes into play.
The reason for this gap is the lack of uniqueness, for which it is well-known that the key element is the validity of the comparison principle; its generalization can be regarded the domination principle we are dealing with in the current work.
Indicative for the difficulties arising in the general cases and briefly described above are the restrictive assumptions imposed on the comparison result in \citet[Theorem 5.13]{KaziTani2BSDE_PIDE}, which provides a general probabilistic representation of a fully nonlinear integro-differential equation.
 The aforementioned work uses in turn \citet[Lemma 5.6]{neufeld2016nonlinear}, \emph{i.e.}, both general cases are restricted in bounded solutions and \cite{KaziTani2BSDE_PIDE} essentially degenerates the "semilinear" term.
In other words, one understands that the extension of classical arguments used for proving the comparison principle for an integro-differential operator of increased complexity does not rely simply on tedious computations, but requires further effort for overcoming the additional difficulties appearing because of the arbitrary growth and of the presence of the non-dominated family of measures.
Before we proceed it is interesting to note that in the local case, \emph{i.e.}, when there are no integral parts in the operator and -consequently- no measures appearing therein, one may use the comparison principle associated to the bounded case in order to derive the comparison principle for the arbitrary (polynomial) growth case, see \citet[Theorem C.2.3]{Peng2019book} and the results following it.

Returning to the general case of integro-differential equations it was \citeauthor{Hollender2017} who has meticulously edited the classical arguments, so that the (sub-/ super-) solutions, as well as the test functions, respect the (arbitrary) polynomial growth at infinity which is imposed by the integrability properties of the family of L\'evy measures.
During this painstaking task in \cite{Hollender2017}, the arguments presented in \citet{Peng2019book,neufeld2016nonlinear} have been extended so that, ultimately, multiple non-dominated L\'evy measures can be taken into account in the fully nonlinear case for the (arbitrary) polynomial growth.
Additionally, a couple of minor mistakes occurred in the literature have been properly taken care under the extension of \cite{Hollender2017}; see the comments on the beginning of the pages 90 and 132 in \cite{Hollender2017}.

The current work, which is based on \cite{Hollender2017}, non trivially expands the results to their fullest extend regarding the integrability restrictions. 
The contribution of the current work is twofold.
Firstly, the required integrability condition for the regularity condition is relaxed and is described in terms of a Young function, thus rendering it sharp, \emph{i.e.}, the required integrability and the growth of the solution go hand by hand, without requiring strictly greater (at least quadratic) integrability as in \cite{Hollender2017}.
Secondly, the ``semilinear'' case is taken care, which is meant that wherever a linear operator appears in \cite{Hollender2017} it can be substituted by a semilinear one. 
In other words and in view of the comments in the two previous paragraphs, \citet[Theorem 5.13]{KaziTani2BSDE_PIDE} is now valid without requiring the degeneracy of the semilinear term and the boundedness of the (sub-/super-) solution.
The applications of the general result presented in the current work can be further appreciated once combined with stability results for stochastic differential equations, \emph{e.g.}, \citet{Papapantoleon2023stabilityBSDE,Papapantoleon2019StabilityMartingale}.
Such a combination allows for proving Trotter--Kato-type results, which include among others numerical approximations, by means of probabilistic arguments.
The generality of the arguments is not restricted in the linear or Markovian case.

For the presentation of this work it is assumed that the reader is familiar with the notion and techniques of viscosity solutions' theory, so that we only introduce the definition, the required notation for the statements and their proofs.
The inexperienced reader may consult \citet{Crandall1992} as an introduction in the topic, \citet{BarlesImbert2008} for an overview of the non-local case with respect to a single measure, \citet[Appendix C]{Peng2019book} for the starting point of the general regularity condition and, of course, \citet{Hollender2017} for the general framework which we are going to adapt.
Naturally, we will use results and techniques whose statements cannot be presented in full detail, for the sake of keeping the length of this work as reasonable as possible.
To this end, we will make use of the results presented in \cite{Hollender2017}; for more details on that see \cref{rem:same_visc_framework}.
Regarding the comparison of the two regularity conditions, the one presented here and the one used in \cite{Hollender2017}, the respective comments are presented after \cref{def:regularity_condition}. 

The structure of the current work is the following: in \cref{sec:HJB_Equations} the admissibility and regularity conditions are presented. 
Within the same section it is presented that for an operator to satisfy the regularity condition, it implies that the comparison principle of the associated integro-differential equation is valid.
Additionally, it is presented that a general class of HJB operators satisfies the regularity condition under consideration.
In \cref{sec:YoungFunctions} results associated to Young functions are presented.
These results are the cornerstones for relaxing and sharpening the regularity condition.
Finally, the lengthiest parts of the proofs of the results presented in \cref{sec:HJB_Equations,sec:YoungFunctions} are relegated to \cref{appendix:HJB_Equations_proofs,appendix:YoungFunctions_proofs}, respectively.         

\subsection{Notation}\label{subsec:notation}
In this subsection we will introduce the notation we will use hereinafter.
Since our results extend those of \cite{Hollender2017}, we will keep the same notation as in \cite{Hollender2017}, whenever possible. 
Hence, it will be easier for the reader to consult the respective auxiliary results cited from \cite{Hollender2017}.

The set of of real numbers will be denoted by $\mathbb{R}$.
For a real number $x$ its positive part is denoted by $x^+$.
The set $(0,\infty)$ may alternatively be denoted by $\mathbb{R}_+$.
Let $d\in\mathbb{N}$ be fixed hereinafter.
The domain of the spatial variable will be $\mathbb{R}^d$.
The set of symmetric $(d\times d)-$matrices will be denoted by $\mathbb{S}^{d\times d}$ and it will be endowed with the usual partial order, denoted by $\le$.
With the same symbol it will be denoted the total order of the real numbers.
The matrix whose elements are all equal to $0$ will be denoted by $\mathbf{0}$, while the identity matrix will be denoted by $I$; since the dimension $d$ of the spatial variable is fixed we have suppressed it in the notation.  
The transpose of the matrix $A$ will be denoted by $A^T$, while the trace of a matrix $A$ will be denoted by $\textrm{Tr}(A)$.
Let, additionally, $p\ge 0$ and a normed space $(X, \|\cdot\|)$. 
The usual Euclidean norm will be denoted by $|\cdot|$, without any reference on the dimension of the Euclidean space.
The supremum norm will be denoted by $\|\cdot\|_{\infty}$.
In a normed space, the closed ball centred at $x$ and of radius $r$ will be denoted by $B[x,r]$ and its complement by $B^c[x,r]$. 
A function $f:(X, \|\cdot\|)\to (\mathbb{R},|\cdot|)$ is said to be of (bounded) \emph{$p-$polynomial growth} if, for some constant $C>0$ and for all $x\in X$,
\begin{align*} 
      |f(x)| \leq C (1+ \|x\|^p).
\end{align*} 
If $\mathcal{C}$ is a family of real-valued functions on $(X,\|\cdot\|)$, then $\mathcal{C}_p$ denotes all functions in $\mathcal{C}$ of $p-$polynomial growth.
The \emph{optimal growth constant} for $f\in \mathcal{C}_p$, den by $\|f\|_p$, is defined by 
\begin{align*} 
      \|f\|_p:= \sup_{x\in X} \frac{|f(x)|}{1+\|x\|^p}.
\end{align*} 
For $(X, \rho)$ a metric space, a function $f:(X, \rho)\to (\mathbb{R},|\cdot|)$ is called \emph{upper semicontinuous} if 
\begin{gather*}
      \limsup_{n\to\infty} f(x_n) \leq f(x),
\end{gather*}
for all sequences $(x_n)_{n\in\mathbb{N}}\subset X$ with $\lim_{n\to \infty} x_n = x$.
A function $f:(X, \rho)\to (\mathbb{R},|\cdot|)$ is said called \emph{lower semicontinuous} if $-f$ is upper semicontinuous.
The set of upper, resp. lower, semicontinuous functions defined on $(X,\rho)$ will be denoted by $\textup{USC}(X)$, resp. $\textup{LSC}(X)$.
Also, we will write 
\begin{align*}
      \textup{SC}(X):=\textup{USC}(X)\cup \textup{LSC}(X)
\end{align*}
for the family of all semicontinuous functions on $(X,\rho)$.
The set of continuous functions defined on $\mathbb{R}^d$ will be denoted by $C(\mathbb{R}^d)$.
The set of twice continuously differentiable functions defined on $\mathbb{R}^d$ will be denoted by $C^2(\mathbb{R}^d)$.
For $f\in C^2(\mathbb{R}^d)$, $Df$ denotes the first derivative, while $D^2f$ denotes the second derivative.
In order to use a compact notation in the description of the settings, we may write an $f\in C^2(\mathbb{R}^d)$ as $D^0f$. 
This will allow us to write a condition simultaneously for $f, Df$ and $D^2f$.

Suppose that $\Omega\subset \mathbb{R}^d$, $f:\Omega\to \mathbb{R}$ and $\varphi : \mathbb{R}^d \to [0,\infty)$.
For $\varepsilon>0$, the \emph{supremal convolution} $\triangle^\varepsilon_\varphi [f]: \mathbb{R}^d \to \mathbb{R}\cup\{+\infty\}$ is defined as 
\begin{align*}
x\longmapsto \triangle^\varepsilon_\varphi [f](x) := \sup_{y\in\Omega} \Big(f(y)-\frac{\varphi(x-y)}{\varepsilon} \Big).
\end{align*}
Dually, the \emph{infimal convolution} $\triangledown^\varepsilon_\varphi [f]$ is defined as $\triangledown^\varepsilon_\varphi [f](x):= - \triangle^\varepsilon_\varphi [-f](x)$, for $x\in\mathbb{R}^d$.
For the parabolic case, we will need to use $[0,T]\times \mathbb{R}^d$ as the domain of the functions, for some $T>0$. 
In this case, a function $f:[0,T]\times \mathbb{R}^d\to \mathbb{R}$ is understood to be of $p-$polynomial growth if 
\begin{align*}
\|f\|_p :=\sup_{t\in[0,T]} \|f(t,\cdot)\|_p <\infty.
\end{align*}
Moreover, for $\varphi: \mathbb{R}^d \to [0,\infty)$, $\varepsilon>0$ and $t\in[0,T]$, the \emph{spatial supremal convolution} $\triangle^\varepsilon_\varphi [f(t,\cdot)]: \mathbb{R}^d \to \mathbb{R}^\{+\infty\}$ is defined as 
\begin{align*}
x\longmapsto \triangle^\varepsilon_\varphi [f(t,\cdot)](x) := \sup_{y\in\Omega} \Big(f(t,y)-\frac{\varphi(x-y)}{\varepsilon} \Big).
\end{align*}
Dually can be defined the \emph{spatial infimal convolution} $\triangledown^\varepsilon_\varphi [f(t,\cdot)]$.

Hereinafter, given $p\ge 0$, we will fix  $\varphi_p$ to be a smooth variant of the function
\begin{align*}
     \mathbb{R}^d \ni x\longmapsto |x|^2 \vee |x|^{p\vee 2} \in [0,\infty). 
\end{align*} 
One way to construct the desired smooth function is described in \cite[Lemma 1.14]{Hollender2017}, which is suitable for our purposes.
In other words,  $\varphi_p$ is a \emph{quasidistance}\footnote{See \cite[Section 1.2]{Hollender2017} for more related results and properties.}, \emph{i.e.}, a symmetric function such that it vanishes only at the origin $0$.
For later reference, we provide the properties of the quasidistance $\varphi_p$ which will be occasionally used in the computations:
\begin{gather}
\varphi_p(x+y) \leq 2^{2(p\vee 2)-2} \big( \varphi_p(x) + \varphi_p(y) \big) \label{quasidist:multiplier}
\shortintertext{and}
2^{2-(p\vee 2)} \big(|x|^2\vee |x|^p\big) \le \varphi_p(x) \leq \big( |x|^2 \vee |x|^p\big). \label{quasidist:growth}
\end{gather}
The necessity for choosing such a function stems from its behaviour in the neighbourhood of 0 and at infinity, a behaviour that is bequeathed to the smooth variant given the construction of \cite[Lemma 1.14]{Hollender2017}.
More precisely, the choice of $\varphi_p$ is such that the supremal convolution $\triangle^\varepsilon_{\varphi_p} [f]$ dominates $f$, for every $f\in \textup{USC}_p(\mathbb{R}^d)$; see  \cite[Lemma 1.12 (i)]{Hollender2017}. 
Additionally, \cite[Lemma 1.12 (iii)]{Hollender2017} guarantees that the supremal convolution $\triangle^\varepsilon_\varphi [f]$ is (locally) semiconvex.
Hence, Alexandroff's theorem can be applied, which is the cornerstone in the viscosity solution theory through Jensen's lemma. 

Finally, 
a \emph{modulus of continuity} is a non-decreasing function $\omega:[0,+\infty]\to [0,+\infty]$ such that $\lim_{h\downarrow 0}\omega(h)=\omega(0)=0$. 
Hereinafter, every considered modulus of continuity will be assumed to additionally be continuous and subadditive.
These mild assumptions can be considered innocuous, since we will deal with sets of functions which are uniformly equicontinuous.
Indeed, the modulus of continuity associated to a uniformly continuous function is continuous, which may additionally be chosen to be subadditive.
Hence, in our framework these properties arise naturally.
Moreover, whenever a finite number of moduli of continuity appear in a statement, for notational convenience, all of them will be denoted by the same symbol.
In other words, we do not assume the minimal modulus of continuity of each set, but an increasing majorant (in a neighbourhood of $0$) of the considered moduli.

In the following we may combine the introduced notation, \emph{e.g.}, $\textup{USP}_p$ will denote the set of upper semicontinuous functions with $p-$polynomial growth \emph{etc}.
Moreover, by $C$ will be denoted an arbitrary positive constant.
Using an index on a constant $C$ should be understood as the constant depending on this quantity, \emph{e.g.}, $C_p$ denotes a constant whose value depends on $p$.

\section{Domination Principle for Hamilton--Jacobi--Bellman Equations}\label{sec:HJB_Equations}
\subsection{Admissibility and regularity condition}\label{subsec:admissibility}

We provide the definition of a viscosity (sub-/super-) solution and, afterwards, the desired properties an operator should possess.
\begin{definition}[Viscosity Solutions]\label{def:viscosity_solutions}
Assume that $p\ge 0$, $\Omega\subset \mathbb{R}^d$ is open and that the operator 
$F:\mathbb{R}^d \times \mathbb{R} \times \mathbb{R}^d \times \mathbb{S}^{d\times d} \times C_p^2(\mathbb{R}^d)\rightarrow \mathbb{R}$ is given.
An upper semicontinuous function $u\in\textup{USC}_p(\mathbb{R}^d)$ is a \emph{viscosity subsolution} in $\Omega$ of the nonlocal equation 
\begin{align}\label{equation:E1}
  F(x,u(x),Du(x),D^2u(x),u(\cdot)) = 0,
  \tag{E1}
\end{align}
if, for all $\phi\in C^2_p(\mathbb{R}^d)$ such that $u-\phi$ has a global maximum in $x\in\Omega$, 
\begin{align*}
  F(x,u(x),D\phi(x),D^2\phi(x),\phi(\cdot)) \le 0.
\end{align*}
A lower semicontinuous function $v\in\textup{LSC}_p(\mathbb{R}^d)$ is a \emph{viscosity supersolution} in $\Omega$ of the nonlocal equation \eqref{equation:E1},
if, for all $\phi\in C^2_p(\mathbb{R}^d)$ such that $v-\phi$ has a global minimum in $x\in\Omega$, 
\begin{align*}
  F(x,u(x),D\phi(x),D^2\phi(x),\phi(\cdot)) \ge 0.
\end{align*}
A \emph{viscosity solution} in $\Omega$ of \eqref{equation:E1} is both a viscosity subsolution and supersolution in $\Omega$ of \eqref{equation:E1}.
\end{definition}
The presented definition is \citet[Definition 2.1]{Hollender2017}. 
The next set of assumptions, which is \cite[Remark 2.3]{Hollender2017}, allows to simultaneously consider a general class of operators, instead of considering special cases.
This is the approach followed in the literature already, \emph{e.g.}, \citet{ALVAREZ1996293,BarlesImbert2008,Jakobsen2006}, and it is also adopted in \citet{Hollender2017}, which we follow.

\begin{assumption}[Assumptions on the operator F]\label{assumption:operator_F}
  Suppose that $p\ge 0$ and that $\Omega\subset \mathbb{R}^d$ is open.
  We assume that for each given operator 
  \begin{align*}
    F:\Omega \times \mathbb{R} \times \mathbb{R}^d \times \mathbb{S}^{d\times d} \times C_p^2(\Omega)\rightarrow \mathbb{R}
  \end{align*}
  and every $0<\kappa<1$, there exist operators (the so-called \emph{generalized operators})
  \begin{align*}
    F^\kappa:\Omega \times \mathbb{R} \times \mathbb{R}^d \times \mathbb{S}^{d\times d} \times \textup{SC}_p (\Omega) \times C_p^2 (\Omega) \rightarrow \mathbb{R}
  \end{align*}
  on which we impose the assumptions \ref{ass:operator_consistency} - \ref{ass:operator_monotonicity} to hold for every 
  $x\in\Omega$, $(x_n)_{n\in\mathbb{N}}\subset \Omega$, 
  $r,s\in\mathbb{R}$, 
  $q\in \mathbb{R}^d$, 
  $X,Y\in\mathbb{S}^{d\times d}$, 
  $u,v\in \textup{SC}_p(\Omega)$, 
  $(u_n)_{n\in\mathbb{N}}\subset \textup{SC}_p(\Omega)$, 
  $\phi,\psi\in C^2(\Omega)$ and 
  $(\phi_n)_{n\in\mathbb{N}}\subset C^2(\Omega)$, except otherwise stated:
  \begin{enumerate}[label=\textup{(A\arabic*)}]
    \item\label{ass:operator_consistency} (Consistency) 
      For each $0<\kappa<1$, the operator $F^\kappa$ is a generalization of $F$, in the sense that the equality
      \begin{align*}
        F^\kappa(x,r,q,X,\phi,\phi) = F(x,r,q,X,\phi)
      \end{align*}
      holds for every $\phi\in C^2_p(\Omega)$.    
    \item\label{ass:operator_degen_ellipt} (Degenerate Ellipticity)
      For each $0<\kappa<1$, the operator $F^\kappa$ is \emph{nonlocal degenerate elliptic}, \emph{i.e.},
      \begin{align*}
        F^\kappa(x,r,q,X,u,\phi) \ge F^\kappa(x,r,q,Y,v,\psi)
      \end{align*}
      holds whenever $u-v$ and $\phi - \psi$ have global maxima in $x\in\Omega$ and whenever $X\le Y$. 

    \item\label{ass:operator_transl_invar} (Translation Invariance)
      For each $0<\kappa<1$, the operator $F^\kappa$ is \emph{translation invariant} in its non-local part, \emph{i.e.},
      \begin{align*}
        F^\kappa(x,r,q,X,u+c_1,\phi+c_2) = F^\kappa(x,r,q,X,u,\phi)
      \end{align*}
      holds for all constants $c_1,c_2\in\mathbb{R}$.

    \item\label{ass:operator_continuity} (Continuity)
      For each $0<\kappa<1$, the operator $F^\kappa$ is \emph{continuous} in the following sense:
      if 
      \begin{enumerate}[label=\textbullet]
        \item $\displaystyle \lim_{n\to+\infty} (x_n,r_n,q_n,X_n) =(x,r,q,X)$,
        \item $\displaystyle \lim_{n\to+\infty} D^m\phi_n = D^m\phi$ locally uniformly for all $m\in\{0,1,2\}$
        and
        \item $\displaystyle \lim_{n\to+\infty} u_n = u$ locally uniformly with $u\in C_p(\Omega)$ and $\sup_{n\in\mathbb{N}} \|u_n\|_p<+\infty$,
      \end{enumerate}
      then
      \begin{align*}
         \lim_{n\to+\infty} F^\kappa(x_n,r_n,q_n,X_n,u_n,\phi_n) = F^\kappa(x,r,q,X,u,\phi).
      \end{align*}

    \item\label{ass:operator_monotonicity} (Monotonicity) 
      For each $0<\kappa<1$, the operator $F^\kappa$ is \emph{non-decreasing} in its second argument, \emph{i.e.}, 
      \begin{align*}
        F^\kappa(x,r,q,X,u,\phi) \le F^\kappa(x,s,q,X,u,\phi)
      \end{align*}
      holds for all $r\le s$.
  \end{enumerate}
\end{assumption}
Some remarks are in order:
\begin{remark}\label{rem:same_visc_framework}
  \begin{enumerate}
    \item \label{rem:same_visc_framework_1}
      The framework we have adopted, \emph{i.e.}, \cref{def:viscosity_solutions} and \cref{assumption:operator_F}, is identical to \citet[Definition 2.1, Remak 2.3]{Hollender2017}.
      Therefore, we can use the results, remarks or comments of \citet[Chapters 1 and 2]{Hollender2017}.
      We will do so for every result up to (excluding) \citet[Definition 2.21]{Hollender2017}, which describes a structure condition.
      This structure condition is named  ``regularity condition'' in \cite{Hollender2017} and this term has been also adopted in this work. 

    \item\label{rem:same_visc_framework_2}
      For the convenience of the reader, when we need to make use of any of the aforementioned results, remarks or comments, we will always refer to the respective part in \citet{Hollender2017}.
      We underlined above that we do not make use of \citet[Definition 2.21]{Hollender2017}, since we will weaken it in \cref{def:regularity_condition}. 
      Comments regarding the difference of the two regularity conditions are presented after \cref{def:regularity_condition}. 

    \item\label{rem:same_visc_framework_3}
      An additional comment should be made before we proceed. 
      In the current work we are interested in proving the comparison principle for parabolic type Hamilton--Jacobi--Bellman equations. 
      In particular, the domain $\Omega$ will be identified as $(0,T)\times \mathbb{R}^d$, for some $T>0$ and $d\in\mathbb{N}$.
      Implicitly, we have already used \citet[Lemma 2.6]{Hollender2017} to write the assumptions on $F$.
      Indeed, given the fact that the non-local parts that appear in the Hamilton--Jacobi--Bellman equations under consideration in the current work only depend on values on the open set $(0,T)\times \mathbb{R}^d$, we are eligible to make use of \cite[Lemma 2.6]{Hollender2017}.
      Consequently, we can write the set of assumptions on the operator $F$ over an open set $\Omega$. 
  \end{enumerate}

\end{remark}
Before we present the regularity condition, we will make a pause to introduce some auxiliary notation.
Let $\varepsilon,T,>0$, $\Omega\subset \mathbb{R}^d$ and 
$F:(0,T) \times \mathbb{R}^d\times \mathbb{R} \times \mathbb{R}^d \times \mathbb{S}^{d\times d} \times C_p^2(\mathbb{R}^d)\rightarrow \mathbb{R}$.
Since $p$ was fixed, instead of $\varphi_p$, we will simply write the quasidistance as $\varphi$.

At this point the reader may assume that we are interested in solving on $(0,T)\times \Omega$ the equation
\begin{align*}  
  \partial_t u (t,x) + F(t,x,u(t,x)Du(t,x),D^2u(t,x), u(t,\cdot))=0
\end{align*}  
and that $u\in\textrm{USP}_p$ is a viscosity subsolution in $(0,T)\times \Omega$.
Given \cite[Remark 2.20]{Hollender2017}, we will need to consider an auxiliary problem which is described by the so-called \emph{smudged operator} 
\begin{align*}
\triangle^\varepsilon_{\varphi}[F]: (0,T)\times \mathbb{R}^d\times \mathbb{R}\times\mathbb{R}^d \times\mathbb{S}^{d\times d} \times C^2_p(\mathbb{R}^d) \to \mathbb{R}
\end{align*}
defined by 
\begin{align*}
  \triangle^\varepsilon_{\varphi}[F](t,x,r,q,X,\psi):=
  \inf\big\{ F(t,y,r,q,X,\psi\circ\tau_{x-y})\, \big|\, y\in \Omega \text{ with } \varphi(y-x) \le \varepsilon \delta(x) \big\},
\end{align*}
for $(t,x,r,q,X,\psi)\in (0,T)\times \mathbb{R}^d\times \mathbb{R}\times \mathbb{R}^d\times \mathbb{S}^{d\times d} \times C^2_p(\mathbb{R}^d)$, 
where $\delta(x):= 2^{3(p\vee 2)} \|u\|_p (1+\varphi(x))$, and
$\tau_{h}:\mathbb{R}^d \to \mathbb{R}^d$ is the transition operator by $h\in\mathbb{R}^d$, \emph{i.e.}, $\tau_h(x) = x+h$, for all $x\in\mathbb{R}^d$.
Additionally, for $k\in\mathbb{N}$ we define
\begin{align*}
J_{k,d}:=
\begin{bmatrix}
(k-1) I & -I & -I & \dots& -I & -I\\
 -I &(k-1) I & -I & \dots& -I & -I\\
 \vdots & \vdots & &\ddots &  & \vdots\\
 -I & -I &- I & \dots& (k-1)I & -I\\
 -I & -I &- I & \dots& -I & (k-1)I\\
\end{bmatrix}\in\mathbb{R}^{kd}.
\end{align*}
\begin{definition}[Regularity Condition]\label{def:regularity_condition}
  Suppose that $p>0$ and $T>0$, that $u_i^\varepsilon\in \textup{USC}_p([0,T]\times \mathbb{R}^d)$ is the spatial supremal convolution of some $u_i\in \textup{USC}_p([0,T]\times \mathbb{R}^d)$ for $\varepsilon>0$ and $i\in\{1,\ldots,k\}$ with $k\in\mathbb{N}$, 
  and that\footnote{\label{footnote:simplified_notation_omit_mu_delta}In the introduced penalization function $\phi_{\delta,\gamma}^{\lambda}$, we notationally omit the dependence on $\mu$, since this parameter will be always assumed fixed and of suitable value when we use the penalization function. We will do so hereinafter for any notational dependence on $(\delta,\mu,\gamma,\lambda,\varepsilon)$, \emph{i.e.}, we will only denote the dependence on $(\delta,\gamma,\lambda,\varepsilon)$. For example, the point of maximum $(\bar{t},\bar{x})=(\bar{t}^{(\delta,\mu,\gamma,\lambda,\varepsilon)},\bar{x}^{(\delta,\mu,\gamma,\lambda,\varepsilon)})$ will be simply denoted to depend on $(\delta,\gamma,\lambda,\varepsilon)$.}
  \begin{align*}
    \phi_{\delta,\gamma}^{\lambda}(t,x) := \frac{\delta}{T-t} + \lambda \sum_{i<j} |x_i-x_j|^2 + \gamma e^{\mu t}\sum_{i=1}^k \Upsilon(|x_i|^p)
  \end{align*}%
  for $(t,x):=(t,x_1,\ldots,x_k)\in[0,T)\times\mathbb{R}^{kd}$ with $\gamma,\delta,\mu,\lambda \ge 0$ and $\Upsilon$ a twice continuously differentiable Young function.\footnote{For the definition see \cref{def:Young_fun}.}
  A family of operators 
  \begin{align*}
    G_i:(0,T)\times \mathbb{R}^d \times \mathbb{R} \times \mathbb{R}^d \times \mathbb{S}^{d\times d}\times C^2_p(\mathbb{R}^d)\longrightarrow \mathbb{R}
  \end{align*}
  for $i\in\{1,\ldots,k\}$ satisfies a \emph{regularity condition for $\beta_1,\ldots,\beta_k>0$,} if there exist $C,\overline{C}>0$ and a modulus of continuity $\omega$ such that the following implication holds:

  \emph{If} the function 
  \begin{align*}
    (0,T)\times\mathbb{R}^{kd}\ni (t,x) = (t,x_1,\ldots,x_k)\longmapsto \sum_{i=1}^k\beta_i u_i^\varepsilon(t,x_i) - \phi_{\delta,\gamma}^{\lambda}(t,x),
  \end{align*}
  has a global maximum in 
  $(\bar{t},\bar{x})
  =(\bar{t},\bar{x}_1,\ldots,\bar{x}_k)
  =(\bar{t}^{(\delta,\gamma,\lambda,\varepsilon)},\bar{x}^{(\delta,\gamma,\lambda,\varepsilon)})
  \in(0,T)\times \mathbb{R}^{kd}$ with
  \begin{align*}
    \sup_{(t,x)\in(0,T)\times \mathbb{R}^{kd}} \bigg(\sum_{i=1}^k\beta_i u_i^\varepsilon(t,x_i) - \phi_{\delta,\gamma}^{\lambda}(t,x)\bigg) = \sum_{i=1}^k\beta_i u_i^\varepsilon(\bar{t},\bar{x}_i) - \phi_{\delta,\gamma}^{\lambda}(\bar{t},\bar{x})\ge0
  \end{align*}
  and if $X_i = X_i^{(\delta,\gamma,\lambda,\varepsilon)}\in\mathbb{S}^{d\times d}$ for $i\in\{1,\ldots,k\}$ satisfy
  \begin{align*}
    \numberthis\label{block_diagonal_form}
    &\textup{diag}(X_1,\ldots,X_k)
      \le 4\lambda J_{k,d}\\ 
     &\hspace{2em}   + C\gamma e^{\mu \bar{t}}\textup{diag}\left[
              (\Upsilon''(|\bar{x}_1|^{p}) |\bar{x}_1|^{2p-2}+\Upsilon'(|\bar{x}_1|^{p})|\bar{x}_1|^{p-2})  I,
                \ldots, 
              (\Upsilon''(|\bar{x}_k|^{p}) |\bar{x}_k|^{2p-2}+\Upsilon'(|\bar{x}_k|^{p})|\bar{x}_k|^{p-2})  I
            \right],
  \end{align*}
  \emph{then} the following inequality holds in the ordinary sense
  \begin{align*}
    \begin{multlined}[0.85\textwidth]
    -\sum_{i=1}^k \beta_i \triangle_\varphi^\varepsilon [G^\kappa_i]\big(\bar{t},\bar{x}_i,u_i^\varepsilon(\bar{t},\bar{x}_i),\beta_i^{-1}D_{x_i}\phi_{\delta,\gamma}^{\lambda}(\bar{t},\bar{x}),X_i,u_i^\varepsilon(\bar{t},\cdot), \beta_i^{-1}\phi_{\delta,\gamma}^{\lambda}(\bar{t},\bar{x}_1,\ldots,\bar{x}_{i-1},\cdot,\bar{x}_{i+1},\ldots,\bar{x}_k)\big)\\
    \le \overline{C}\lambda \sum_{i<j}|\bar{x}_i - \bar{x}_j|^2 
          +\overline{C}\gamma e^{\mu \bar{t}}\Big(1 + \sum_{i=1}^k \Upsilon(|\bar{x}_i|^p)\Big) 
          +\varrho_{\gamma,\lambda,\varepsilon,\kappa}
    \end{multlined}      \numberthis\label{ineq:result_regularity_condition}
  \end{align*}
  for a remainder $\varrho_{\gamma,\lambda,\varepsilon,\kappa}$ with
  \begin{align*}
    \limsup_{\gamma \downarrow 0} \limsup_{\lambda\to+\infty}\limsup_{\varepsilon\downarrow 0}\limsup_{\kappa\downarrow 0}\varrho_{\gamma,\lambda,\varepsilon,\kappa}\le 0.
  \end{align*}
\end{definition}
Now that the regularity condition has been presented in terms of an arbitrary Young function, let us compare the presented framework with that adopted in \citet[Definition 2.21]{Hollender2017} and motivate the block diagonal form in the right-hand side of \eqref{block_diagonal_form} in \cref{def:regularity_condition}.
The reason the regularity condition described in \cref{def:regularity_condition} is significantly weaker than \citet[Definition 2.21]{Hollender2017}  will be clarified in the upcoming comments.
More importantly, it will be justified the fact the presented regularity condition is the natural one, since this allows for sharp integrability conditions.

Let us start by observing that we can rewrite\footnote{One has to replace $q(q-1)$ by a constant, which is innocuous in view of the factor $\gamma$ which will finally tend to $0$.} the regularity condition used in \citet[Definition 2.21]{Hollender2017} in terms of \cref{def:regularity_condition} by using the Young function 
\begin{align}
  \Upsilon(y):=y^{\frac{q}{p}}, \text{ for }0\le y, 0< p<q \text{ and } 2\le q.
  \label{YoungFunction_forHollender_RegularityCondition}
\end{align}
Indeed,
\begin{align*}
  \Upsilon'(|z|^p) |z|^{p-2} = C_{p,q} |z|^{q-2} \text{ and }
  \big[\Upsilon''(|z|^p) |z|^{p}\big] |z|^{p-2}= C_{p,q}|z|^{q-2} \text{ for }z\in\mathbb{R}^d.
\end{align*}
We proceed now with the motivation for the block diagonal form in the right-hand side of \eqref{block_diagonal_form} in \cref{def:regularity_condition}.
The conditions on $q$ required in \citet[Definition 2.21]{Hollender2017}, \emph{i.e.}, $q>p$ and $q\ge 2$, dictate the use of\footnote{Here, we foresee the framework of the equation we are interested in, \emph{i.e.}, \cref{def:HJB,assumption:HJB_coefficients_assumption}.} a family of measures which demand the uniform integration of the function $\mathbb{R}^d \ni z\mapsto |z|^q$ on $\{z\in\mathbb{R}^d:|z|>1\}$.
In the case $p\in[0,2)$ one would be content with $q=2$.
In the case $p\ge 2$ one would naturally require the value $q$ to be ``as close to $p$ as possible'', \emph{i.e.}, $q=p+\delta$ for $\delta$ positive and arbitrarily small, so that they can seek for weaker integrability requirements;
this allows us to further assume $q<2p$.
In the case $2\le p<q<2p$, we have for the Young function of \eqref{YoungFunction_forHollender_RegularityCondition} and for $y>0$ that 
\begin{align}\label{example:power_function}
\Upsilon'(y) =\frac{q}{p} y^{\frac{q}{p}-1} \text{ with }\frac{q}{p}-1>0 \text{ and }\Upsilon''(y) = \frac{q}{p} (\frac{q}{p}-1) y^{\frac{q}{p}-2}\text{ with }\frac{q}{p}-2<0.
\end{align}
Going one step further and seeking an upper bound of the Hessian of the function $z\mapsto |z|^q = \Upsilon(|z|^p)$,
one can easily compute that to be of the form $q(q-1)|x|^{q-2}$,
which is the term appearing in the inequality of the block-diagonal matrices of \citet[Definition 2.21]{Hollender2017}.
Rewriting the extracted upper bound in terms of the first two derivatives of the Young function $\Upsilon$ and of the function $z\mapsto |z|^{p-2}$, then one gets
\begin{align*}
  |z|^{q-2} = C_{p,q}\Big[\big(|z|^p \big)^{\frac{q}{p}-2} |z|^{2p-2} + \big(|z|^p\big)^{\frac{q}{p}-1} |z|^{p-2}\Big]
    = C_{p,q} \big[\Upsilon''(|z|^p) |z|^{2p-2} + \Upsilon'(|z|^p)|z|^{p-2}\big],
\end{align*} 
which is the term appearing in the right-hand side of \eqref{block_diagonal_form}. 
Assuming now that $\Upsilon$ is a twice continuously differentiable Young function and $p>0$ then one gets 
\begin{align*}
  \nabla \Upsilon(|z|^p)= p\Upsilon'(|z|^p) |z|^{p-2}z
\end{align*}
and for $i,j\in\{1,\ldots,d\}$, where $\left(\nabla \Upsilon(|z|^p)\right)_j$ is the $j-$element of the gradient,
\begin{align*}
  &\frac{\partial}{\partial z_i}\left(\nabla \Upsilon(|z|^p)\right)_j
    =\frac{\partial}{\partial z_i} (\Upsilon'(|z|^p) |z|^{p-2}z_j) \\
    &\hspace{1em}=  p\left[p\Upsilon'' (|z|^p) |z|^{2p-4} + (p-2)\Upsilon'(|z|^p)|z|^{p-4}\right] z_iz_j
      + p\Upsilon'(|z|^p)|z|^{p-2}\mathds{1}_{\{j\}}(i)\\
    &\hspace{1em}\le \frac{p}{2} \left[p\Upsilon'' (|z|^p) |z|^{2p-2} + |p-2|\Upsilon'(|z|^p)|z|^{p-2}\right]
      + p\Upsilon'(|z|^p)|z|^{p-2}\\
    &\hspace{1em} = C_{p}\big[ \Upsilon'' (|z|^p) |z|^{2p-2} + \Upsilon'(|z|^p)|z|^{p-2}\big],
    \numberthis\label{bound:Young_Second_Derivative}
\end{align*}
for $C>0$.
In other words, in the right-hand side of \eqref{block_diagonal_form} we have simply used the extracted upper bound of the elements of the Hessian of $\Upsilon(|z|^p)$. 
So far we have argued about seeing \cite[Definition 2.21]{Hollender2017} as a special case of \cref{def:regularity_condition}.
In the next paragraphs we will argue about the difficulties that appear when instead of power functions, as in \cite[Definition 2.21]{Hollender2017}, general Young functions are utilized.
These comments will justify the necessity of the additional properties a Young function is required to possess, hence the necessity of the results of \cref{sec:YoungFunctions}.

Let us, now, assume a Young function $\Upsilon\in C^2(0,\infty)$.
At this point a natural question is whether there are Young functions, apart from power functions, such that the term 
$\Upsilon'' (|z|^p) |z|^{2p-2} + \Upsilon'(|z|^p)|z|^{p-2}$ is controllable and possibly having additional convenient properties, for any $p>0$.
Essentially, we need to examine the behaviour of the aforementioned quantity around origin, \emph{e.g.}, it is important not to explode, as well as at infinity, \emph{e.g.}, it respects the integrability properties of the measures appearing in the operator.
Additionally, considering the well-known fact that there are ``plenty'' of Young functions which grow slower than any power function $[0,+\infty) \ni y\mapsto y^\delta$, for $\delta>0$, our aim is to choose among them one which combines all the above properties, otherwise we will be trapped in the need of the auxiliary power $q$.
Making the long story short by avoiding the extensive computations that appear in the proofs, we will ``choose'' the Young function $\Upsilon$ such that 
\begin{align*}
    \lim_{|z|\to 0}\Upsilon''(|z|^p) |z|^{2p-2}
      = \lim_{|z|\to 0}\Upsilon'(|z|^p) |z|^{p-2} 
      =0
\end{align*}
and that the term $\Upsilon'' (|z|^p) |z|^{2p-2} + \Upsilon'(|z|^p)|z|^{p-2}$ is integrable.
In particular, as we have seen in \eqref{example:power_function} the function $[0,+\infty)\ni x \mapsto x\Upsilon''(x)$ will be bounded. 
It is \cref{lemma:UI_Young_improvement} and its corollaries which verify that it is possible to construct such a Young function based on the data of the problem.

Returning to \citet[Definition 2.21]{Hollender2017}, it is because of these points that described in the previous paragraph that it is required the presence of the auxiliary power $q$ with $p<q$ with $q\ge 2$.
At this point we need to underline that in the comments provided after \cite[Theorem 2.22]{Hollender2017}, it is mentioned that the required condition on exponents can be relaxed to $p<q<2$.
Despite this remark in \citet{Hollender2017}, no additional details are provided on the way the ``elaborate construction'' can be achieved. 
However, even in this case the auxiliary exponent $q$ is needed to be strictly greater than $p$ rendering the suggested construction non-optimal.
To the contrary, with the approach in which we use a suitable Young function there is no integrability restriction, since the Young function is constructed based on the (uniform) integrability of the family of measures, see \cref{assumption:HJB_coefficients_assumption}.\ref{HJB:UI} and \cref{prop:UI_Young}, which is identical to \cite[Remark 2.28 (B2)]{Hollender2017}, and on the value of $p>0$.

A last remark about a subtle, innocuous difference between \cite[Definition 2.21]{Hollender2017} and \cref{def:regularity_condition} is the following.  
In the right hand side of \eqref{ineq:result_regularity_condition} we do not assume a term of the form $\omega\big((\sum_{i<j}|\bar{x}_i-\bar{x}_j|)( 1 + \sum_{i=1}^k |\bar{x}_i|^p)\big)$, since it was preferred for simplicity to be included in the remainder $\varrho_{\gamma,\lambda,\varepsilon,\kappa}$.
The reader may verify that a similar term may appear by examining the proofs of the results in \cref{subsubsec_Appendix:preparatory_lemmata_Regularity_Condition}. 

We proceed now to present for the parabolic case the domination principle, the generalization of the comparison principle, which can be seen as an application of the maximum principle.
The result presented below generalizes \cite[Theorem 2.22]{Hollender2017}, thus covering previous results in the existing literature of the integro-differential equation, \emph{e.g.}, \citet[Section 3]{ALVAREZ1996293}, \citet[Section 5]{BarlesImbert2008}, \citet[Section 7]{Hu2021}, \citet[Section 3]{Jakobsen2006} and \citet[Section 4]{PhamOptimalStoppingControlled}.

\begin{theorem}[Domination Principle]\label{theorem:Domination_Principle}
  Suppose that $p> 0$, $T>0$ and $k\in\mathbb{N}$. 
  Moreover, suppose that $u_i\in\textup{USC}_p([0,T]\times \mathbb{R}^d)$ are viscosity solutions in $(0,T)\times \mathbb{R}^d$ of 
  \begin{align*}
    \partial_t u_i(t,x) + G_i\big(t,x,u_i(t,x),Du_i(t,x), D^2u_i(t,x), u_i(t,\cdot)\big) \le 0,
  \end{align*}
  for $i\in\{1,\ldots,k\}$ and operators $(G_i)_{i\in\{1,\ldots,k\}}$ which satisfy \cref{assumption:operator_F} and the regularity condition described in \cref{def:regularity_condition} for $\beta_1,\ldots,\beta_k>0$ and a twice continuously differentiable Young function $\Upsilon$.
  If the initial conditions $u_i(0,\cdot):\mathbb{R}^d \to \mathbb{R}$ are continuous for all $i\in\{1,\ldots,k\}$ and if 
  \begin{align*}
    \sum_{i=1}^k \beta_i u_i(0,x) \le 0, \text{ for every }x\in\mathbb{R}^d,
  \end{align*}
  then 
  \begin{align*}
    \sum_{i=1}^k \beta_i u_i(t,x) \le 0, \text{ for every }(t,x)\in (0,T)\times \mathbb{R}^d.
  \end{align*}
  
\end{theorem}

\begin{proof}
  The proof essentially follows that of \citet[Theorem 2.22]{Hollender2017}.
  The main difference is that in the current proof we use the function $\Upsilon(|\cdot|^p)$, where $\Upsilon$ is the Young function for which the regularity condition holds, instead of $|\cdot|^q$, for $q> p$ and $q\ge 2$, which is used in \citet[Theorem 2.22]{Hollender2017}.
  Regarding the validity of the substitution we mention in the previous line, the crucial remark is that the function $|\cdot|^p$ is finally dominated by $\Upsilon(|\cdot|^p)$.
  Hence, all the properties extracted for the penalization of the conical combination of the supremal convolutions remain true, since the supremal convolutions respect the polynomial order of the subsolutions; see \citet[Lemma 1.12]{Hollender2017}.
  After these remarks, the reader who is familiar with the technicalities of the viscosity solutions theory should be convinced. 
  However, we present the complete proof in \cref{proof:Thm_Domination_Principle}, not only for the unfamiliar reader, but also in order to justify and make use of properties of the families of global maxima that appear in the proof. 
  These properties will be repeatedly used in \cref{HJB:Regularity_Condition}.
\end{proof}
We close the subsection with the most important corollary of the Domination Principle, the Comparison Principle, which guarantees the uniqueness of the solution of the integro-differential equation.
Additionally, through Perron's method, the Comparison Principle also leads to the existence of the solution, given a subsolution and a supersolution.  
\begin{corollary}[Comparison Principle]\label{cor:Comparison_Principle}
Suppose that $G$ is an operator that satisfies \cref{assumption:operator_F} with $p>0$ and $T>0$.
Additionally, suppose that the two operators $G_1$ and $G_2$, which are defined by
\begin{gather*}
G_1(t,x,r,q,X,\phi):= G(t,x,r,q,X,\phi)
\ \text{ and }\
G_2(t,x,r,q,X,\phi):= -G(t,x,-r,-q,-X,-\phi)
\end{gather*}
for $(t,x,r,q,X,\phi)\in(0,T)\times \mathbb{R}^d\times \mathbb{R}\times \mathbb{R}^d \times\mathbb{S}^{d\times d}\times C^2_p(\mathbb{R}^d)$, satisfy a regularity condition as described in \cref{def:regularity_condition} for $\beta_1=\beta_2=1$ for some twice continuously differentiable Young function $\Upsilon$.

If $\underbar{u}\in \textrm{USC}_p([0,T]\times\mathbb{R}^d)$, resp. $\bar{u}\in \textrm{LSC}_p([0,T]\times\mathbb{R}^d)$, is a viscosity subsolution, resp. supersolution, in $(0,T)\times \mathbb{R}^d$ of 
\begin{align*}
\partial_t u (t,x) + G\big(t,x,u(t,x),Du(t,x),D^2u(t,x),u(t,\cdot)\big)=0
\end{align*}
such that $\underbar{u}(0,\cdot), \bar{u}(0,\cdot)$ are continuous with $\underbar{u}(0,x) \leq \bar{u}(0,x)$ for every $x\in\mathbb{R}^d$, then 
\begin{align*}
\underbar{u}(t,x)\leq \bar{u}(t,x) \text{ for every }(t,x)\in(0,T)\times \mathbb{R}^d.
\end{align*}
\end{corollary}

\begin{proof}
It is easily verified that the assumptions of \cref{theorem:Domination_Principle} are satisfied, hence the desired conclusion is indeed true. 
\end{proof}

\begin{remark}\label{rem:Domin_Princ_corollaries}
The presented Domination Principle (\cref{theorem:Domination_Principle}) generalizes \citet[Theorem 2.22]{Hollender2017}, which in turn generalized \citet[Theorem C.2.2]{Peng2019book} and \citet[Theorem 53]{Hu2021}.
The generality of the Domination Principle allows the generator to bequeath its properties, \emph{e.g.}, subadditivity and/or convexity, to the solutions of the associated initial value problems.
The precise statements can be readily adapted from \citet[Corollary 2.25, Corollary 2.26]{Hollender2017} or \citet[Corollary 56, Corollary 57]{Hu2021}.
\end{remark}

\subsection{A prominent example: Hamilton--Jacobi--Bellman Equations}\label{subsec:prominent_example_HJB}

In this subsection we focus on a special class of operators associated to integro-differential equations, namely 
 Hamilton--Jacobi--Bellman Equations (HJB). 
We initially provide the definition of HJB Equations and afterwards we set a framework for the data of the operator, under which the associated HJB operator satisfies the \cref{assumption:operator_F}, as well as the regularity condition described in \cref{def:regularity_condition}. 
\begin{definition}[Hamilton--Jacobi--Bellman Equations]\label{def:HJB}
A nonlinear, parabolic, second-order nonlocal equation in $(0,T)\times\mathbb{R}^d$ with time horizon $T>0$
\begin{align}\label{equation:HJB}
\partial_t u (t,x) + G\left(t,x,u (t,x),Du (t,x),D^2u (t,x),u (t,\cdot)\right)=0
\tag{HJB}
\end{align}
is called \emph{Hamilton--Jacobi--Bellman equation} (or HJB equation  for short) if 
\begin{gather*}
G(t,x,r,q,X,\phi) = G^\kappa (t,x,r,q,X,\phi,\phi)\\ 
G^\kappa(t,x,r,q,X,u,\phi) = \inf_{\alpha\in\mathcal{A}} G_\alpha^\kappa (t,x,r,q,X,u,\phi)
\end{gather*}
for a family $(G_\alpha^\kappa)_{\alpha\in\mathcal{A}}$ of semi-linear operators
\begin{align*}
  G^\kappa_\alpha :(0,T)\times \mathbb{R}^d\times \mathbb{R}\times \mathbb{R}^d\times \mathbb{S}^{d\times d}\times \textup{SC}_p(\mathbb{R}^d) \times C^2(\mathbb{R}^d) \longrightarrow \mathbb{R}
\end{align*}
with $p\ge 1$ and $0<\kappa<1$ of the form
\begin{align*}
G_\alpha^\kappa(t,x,r,q,X,u,\phi) 
  &= -f_\alpha\big(t,x,r,\sigma_\alpha^T(t,x)q ,\mathcal{K}_\alpha^\kappa(t,x,u,\phi)\big) - \mathcal{L}_\alpha(t,x,r,q,X) 
    - \mathcal{I}_\alpha^\kappa(t,x,u,\phi)\\
\mathcal{L}_\alpha(t,x,q,X) 
  & =b_\alpha^T(t,x)q 
    + \textup{Tr}\left(\sigma_\alpha(t,x)\sigma^T(t,x)X\right)\\
\mathcal{I}_\alpha^\kappa(t,x,u,\phi) 
  & =\check{\mathcal{I}}_\alpha^\kappa(t,x,\phi) 
    + \overline{\mathcal{I}}_\alpha^\kappa(t,x,u,\phi) 
    +\hat{\mathcal{I}}_\alpha (t,x,u,\phi),\\
\mathcal{K}_\alpha^\kappa(t,x,u,\phi) 
  & =\check{\mathcal{K}}_\alpha^\kappa(t,x,\phi) 
    + \overline{\mathcal{K}}_\alpha^\kappa(t,x,u) 
    ,
\end{align*}
where $f_\alpha(t,x,r,q,w)\in\mathbb{R}$, $b_\alpha(t,x)\in\mathbb{R}^d$, $\sigma_\alpha(t,x)\in\mathbb{S}^{d\times d}$, 
\begin{align*}
  \check{\mathcal{I}}_\alpha^\kappa(t,x,\phi) 
    &:= \int_{\{z\in\mathbb{R}^d:|z|\le \kappa\}} \left(\phi(x+j_\alpha(t,x,z))-\phi(t,x) - D\phi(x)^Tj_\alpha(t,x,z)\right) m_\alpha(\textup{d}z)   \\
  \overline{\mathcal{I}}_\alpha^\kappa(t,x,u,\phi) 
    &:= \int_{\{z\in\mathbb{R}^d:\kappa< |z|\le 1\}} \left(u(x+j_\alpha(t,x,z))-u(t,x) - D\phi(x)^Tj_\alpha(t,x,z)\right) m_\alpha(\textup{d}z) \\
  \hat{\mathcal{I}}_\alpha (t,x,u,\phi) 
    &:=\int_{\{z\in\mathbb{R}^d:1< |z|\}} \left(u(x+j_\alpha(t,x,z))-u(t,x) - D\phi(x)^Tj_\alpha(t,x,z)\right) m_\alpha(\textup{d}z) 
\end{align*}
and
\begin{align*}
  \check{\mathcal{K}}_\alpha^\kappa(t,x,\phi) 
    &:= \int_{\{z\in\mathbb{R}^d:|z|\le \kappa\}} \left(\phi(x+j_\alpha(t,x,z))-\phi(t,x)\right)\delta_\alpha(t,x,z) m_\alpha(\textup{d}z)   \\
  \overline{\mathcal{K}}_\alpha^\kappa(t,x,u) 
    &:= \int_{\{z\in\mathbb{R}^d:\kappa< |z|\}} \left(u(x+j_\alpha(t,x,z))-u(t,x)\right)\delta_\alpha(t,x,z) m_\alpha(\textup{d}z) 
\end{align*}
with Borel measures $m_\alpha:\mathcal{B}(\mathbb{R}^d)\rightarrow [0,+\infty)$, $j_\alpha(t,x,z)\in\mathbb{R}^d$, $\delta_{\alpha}(t,x,z)\in\mathbb{R}$.
\end{definition}
\begin{remark}\label{rem:after_defHJB}
  \begin{enumerate}
    \item \label{rem:after_defHJB_i}
  We have slightly deviated from the ``classical'' definition of the HJB operator by allowing the presence of the term $-D\phi(x)j_\alpha(t,x,z)$ in the segment $\hat{\mathcal{I}}_\alpha$ of the nonlocal part of the operator.
  In other words, we do not truncate in the domain $\{z\in\mathbb{R}^d:1<|z|\}$ as in the classical case, where the truncation function  $\mathds{1}_{\{z:|z|\le 1\}}$ is used as a factor of the term $D\phi(x)j_\alpha(t,x,z)$.  

  Our choice is motivated from the probabilistic interpretation of the nonlocal operator.
  More precisely, we are interested in the case where the linear part of the operators $G_\alpha^\kappa$ can be associated to integrable martingales, \emph{e.g.} when $j_\alpha(t,x,\cdot)=\textup{Id}$, for $\textup{Id}$ the identity on $\mathbb{R}^d$, and for this reason the value of $p$ has to be greater than or equal to $1$.
  In other words, the linear part corresponds to the infinitesimal generator of a martingale.

  Despite this choice, the way the results and their proofs are presented allows to easily deduce the analogous results of the classical case, \emph{i.e.}, when the truncation function $\mathds{1}_{\{z\in\mathbb{R}^d:|z|\le 1\}}$ is used, which associates to semimartingales.
  In the semimartingale case we may allow $p> 0$, instead of $p\ge 1$ which is dictated by the integrability requirement of the true martingale.
    \item \label{rem:after_defHJB_ii}
  We can easily translate the initial value problem to a terminal value problem, by using the change of time $t\mapsto T-t$. 
  Before we proceed, let us introduce the following notation: if $h$ is a function defined on $[0,T]\times \mathcal{R}$, for some Euclidean space $\mathcal{R}$, then $\tilde{h}(t,\cdot):=h(T-t,\cdot)$, for every $t\in [0,T]$, 
  Now, if $u$ is a viscosity solution of Equation \eqref{equation:HJB} with $u(0,\cdot)=g(\cdot)$, then $\tilde{u}$ is a viscosity solution of the equation
  \begin{align*}
    -\partial_t\tilde{u}(t,x) +\tilde{G}(t,x,\tilde{u}(t,x), D\tilde{u}(t,x), D^2\tilde{u}(t,x), \tilde{u}(t,\cdot))=0, \text{ with }
    \tilde{u}(T,\cdot)=g(\cdot).
   \end{align*}
    \item \label{rem:after_defHJB_iii}
  If the set $\mathcal{A}$ is a singleton, 
  then the Equation \eqref{equation:HJB} with given terminal value can be associated to a decoupled system of Forward-Backward Stochastic Differential Equations (FBSDE). 
  In \citet{barles1997backward} the existence and uniqueness of the solution of the (degenerate) \eqref{equation:HJB} is provided by probabilistic arguments. 
  In this case the integro-differential equation is a semilinear one.
  
  If the set $\mathcal{A}$ is not a singleton, the respective system of FBSDE should be understood under a family of (non-dominated) probability measures. 
  The interested reader may consult \citet{soner2012wellposedness} for the local case and \citet{KaziTani2BSDE_PIDE} for the non-local case.
  In the aforementioned cases, the Forward Stochastic Differential Equation corresponds to the case $b_{\alpha}\equiv 0$, 
  $\sigma_{\alpha}\equiv A \textrm{Id}$, additionally $j_{\alpha}(t,x,\cdot)\equiv \textrm{Id}$ for every $(t,x)\in(0,T)\times \mathbb{R}^d$ for the nonlocal case, for $A$ lying in a suitable set, as well as the L\'evy meaures $\nu$ should lie in a suitable set.
  The pair $(A,\nu)$ should be thought as the parameter $\alpha\in\mathcal{A}$.
  It is interesting to note that these works provide a probabilistic representation of the solution of the associated HJB equation, but they need to be complemented by a comparison result in order to guarantee the uniqueness of the solution of the (integro-)differential equation.
  \end{enumerate}
\end{remark}
Naturally, the framework described in \cref{def:HJB} is too general and without any further assumptions on the coefficients of the semi-linear operators $G_\alpha^\kappa$ there is no hope that the comparison principle holds.
To this end, we will need to impose specific properties.
This is the purpose of \cref{assumption:HJB_coefficients_assumption}.
Afterwards, it will be proven that under \cref{assumption:HJB_coefficients_assumption} the associated Hamilton--Jacobi--Bellman operator satisfies \cref{assumption:operator_F} and the Regularity Condition as described in \cref{def:regularity_condition}.
Hence, in view of \cref{theorem:Domination_Principle} one immediately gets that the comparison principle holds.

\begin{assumption}[on HJB Coefficients]\label{assumption:HJB_coefficients_assumption}
Suppose that, for every $0<\kappa<1$, $(G_\alpha^{\kappa})_{\alpha\in\mathcal{A}}$ is a family of semilinear operators as defined in \cref{def:HJB} with $p\ge 1$ and $T>0$.
Throughout this section we will assume that the following conditions hold:
\begin{enumerate}[label=\textup{(B\arabic*)}]
  \item\label{HJB:coeff_boundedness} (Boundedness) 
    The coefficients of the local part and the semilinear part are bounded, \emph{i.e.},
    \begin{align*}
      \sup_{\alpha\in\mathcal{A}}\left(|b_\alpha(t,x)|+|\sigma_\alpha(t,x)| + |f_\alpha(t,x,r,q,w)| \right)<+\infty
    \end{align*}
    in every $(t,x,r,q,w)\in(0,T)\times \mathbb{R}^d\times \mathbb{R} \times \mathbb{R}^d \times \mathbb{R}$, and the family of measures of the nonlocal part meet
    \begin{align*}
      \sup_{\alpha \in\mathcal{A}} \int_{\mathbb{R}^d}\left(|z|^2\mathds{1}_{\{z:|z|\le 1\}}(z) + |z|^p\mathds{1}_{\{z:|z|>1\}}(z)\right) m_\alpha(\textup{d}z)<+\infty.
    \end{align*}

  \item\label{HJB:UI} (Uniformly integrable measures) The family of measures of the nonlocal part satisfies
    \begin{gather*}
      \lim_{\kappa\downarrow 0}\sup_{\alpha\in\mathcal{A}} \int_{\{z:|z|\le \kappa\}} |z|^2 m_\alpha(\textup{d}z)
      = 0
      \shortintertext{and}
      \lim_{R\to+\infty} \sup_{\alpha\in\mathcal{A}} \int_{\{z:|z|>R\}} |z|^p m_\alpha(\textup{d}z)
      = 0.
    \end{gather*}

  \item\label{HJB:coef_continuity} (Uniform equicontinuity of $(b_\alpha)_{\alpha\in\mathcal{A}}$ and $(\sigma_\alpha)_{\alpha\in\mathcal{A}}$)
    There exists a constant $C>0$ and a modulus of continuity $\omega$ such that 
    \begin{align*}
     \sup_{\alpha\in\mathcal{A}} \left( |\sigma_\alpha(t,x) - \sigma_\alpha(t',x')| + |b_\alpha(t,x) - b_\alpha(t',x')| \right) \le \omega(|t-t'|)+C|x-x'|,
    \end{align*}
    for all $(t,x)\in(0,T)\times \mathbb{R}^d$.
  \item\label{HJB:monotonicity} (Uniform equicontinuity and monotonicity of the non-linearities $(f_{\alpha})_{\alpha\in\mathcal{A}}$))
    The non-linearities $(f_{\alpha})_{\alpha\in\mathcal{A}}$ possess the following properties: 
    \begin{enumerate}
      \item \label{HJB:monotonicity_1}
    for every $K,R >0$ there exist moduli of continuity $\overline{\omega}_{\!{}_R} $ and $ \widetilde{\omega}_{\!{}_K}$ 
    such that
    \begin{align*}
    \begin{multlined}[c][0.9\textwidth]
      \sup_{\alpha \in \mathcal{A}} \big| f_{\alpha}(t,x,r,q, w) -f_{\alpha}(t',x',r',q', w')\big|\\
      \leq \omega(|t-t'|) + \overline{\omega}_{\!{}_R}\big(|x-x'|(1+|q|\vee |q'|)\big) + \widetilde{\omega}_{\!{}_K}(|r-r'|) + C(|q-q'| + |w-w'|),
    \end{multlined}
    \end{align*}
    for all $t,t'\in(0,T)$, $x,x'\in\mathbb{R}^d$ with $|x|,|x'|\le R$, $r,r'\in\mathbb{R}$, $q,q'\in\mathbb{R}^d$ and $w,w'\in\mathbb{R}$.
    \item \label{HJB:monotonicity_2} they are monotone in the second spatial argument and in the last spatial argument in the following sense:
    \begin{align*}
      \sup_{\alpha \in\mathcal{A}}(r-r')(f_{\alpha}(t,x,r,q,w) - f_{\alpha}(t,x,r',q,w))\le 0
    \end{align*}
    and 
    \begin{align*}
      \mathbb{R} \ni w \mapsto f_\alpha (t,x,r,q,w) \in \mathbb{R}\ \text{ is non-decreasing},
    \end{align*}
    for all $t\in(0,T)$, $x,q\in\mathbb{R}^d$ and $r,r'\in\mathbb{R}$.
    \end{enumerate}

  \item\label{HJB:cond_jump} (Continuity and growth of $(j_{\alpha})_{\alpha\in\mathcal{A}}$) The jump-height coefficients of the nonlocal parts $(\mathcal{I}^{\kappa}_{\alpha})_{\alpha\in\mathcal{A}}$ and $(\mathcal{K}^{\kappa}_{\alpha})_{\alpha\in\mathcal{A}}$ possess the following properties:
  \begin{enumerate}
    \item\label{HJB:cond_jump_partA}
    there exist a constant $C>0$ and a modulus of continuity $\omega$
     such that 
    \begin{align*}
      \sup_{\alpha\in\mathcal{A}}  |j_\alpha(t,x,z) - j_\alpha(t',x',z)|  \le |z|\left(\omega(|t-t'|)+C|x-x'|\right)
    \end{align*}
    for all $(t,x),(t',x')\in(0,T)\times \mathbb{R}^d$,
    \item\label{HJB:cond_jump_partB}
      there exists $C>0$ such that
      \begin{align*}
        \sup_{\alpha\in\mathcal{A}}  |j_\alpha(t,x,z)|\le C(1+|x|)|z|
      \end{align*}
      for all $t\in(0,T)$ and $x,z\in\mathbb{R}^d$, and
    \item\label{HJB:cond_jump_partC}
      if $p\in[1,2)$, then additionally for every $x,x'\in\mathbb{R}^d$ and for every $t\in(0,T)$ there exists a constant $C>0$ such that
          \begin{align*}
            \sup_{\alpha\in\mathcal{A}} |j_\alpha(t,x,z) | \le C(1 + |x|^{\frac{p}{2}} ) \,|z|
            \text{ for every $z\in\mathbb{R}^d$ such that $|z|<1$}
          \end{align*} 
          and
          \begin{align*}
            \sup_{\alpha\in\mathcal{A}} |j_\alpha(t,x,z) - j_\alpha(t,x',z)| \le C|x-x'| \,|z|^{\frac{p}{2}}
            \text{ for every $z\in\mathbb{R}^d$ such that $|z|>1$.}
          \end{align*}
  \end{enumerate}

  \item\label{HJB:growth_cond_delta} (Continuity and growth of $(\delta_{\alpha})_{\alpha\in\mathcal{A}}$) 
For every $\alpha\in\mathcal{A}$ there exist a measurable $\ell_\alpha:\mathbb{R}^d \to [0,+\infty)$, a constant $C>0$ and a modulus of continuity $\omega$
  such that 
\begin{enumerate}
  \item\label{HJB:growth_cond_delta_partA} for every $(t,x,z)\in [0,T]\times \mathbb{R}^d\times \mathbb{R}^d$ it holds
    \begin{align*}
      0\le \delta_{\alpha}(t,x,z) \le \ell_\alpha (z), 
    \end{align*}
  \item\label{HJB:growth_cond_delta_partB} it holds
      \begin{align*}    
        \sup_{\alpha\in\mathcal{A}} \max\Big\{ \int_{\mathbb{R}^d} \ell^2_\alpha(z) m_{\alpha}(\textup{d}z), \Vert \ell_\alpha \Vert_{\infty} \Big\} <+\infty,
      \end{align*}
  \item\label{HJB:growth_cond_delta_pre_partC} 
    for all $t,t'\in(0,T)$ and $x,x' \in \mathbb{R}^d$ it holds 
    \begin{align*}
      |\delta_\alpha(t,x,z) - \delta_\alpha(t',x',z)|  \le \ell_\alpha (z)\left(\omega(|t-t'|) + C|x-x'|\right)
    \end{align*}
      and
  \item\label{HJB:growth_cond_delta_partC} 
    for all $t\in(0,T)$ and $x,x' \in \mathbb{R}^d$ it holds 
    \begin{align*}
      |\delta_\alpha(t,x,z) - \delta_\alpha(t,x',z)|  \le C|x-x'| \times \ell_\alpha^2 (z)  \text{ for }z \in B[0,1].
    \end{align*}
as well as
    \begin{align*}
      |\delta_\alpha(t,x,z) - \delta_\alpha(t,x',z)|  \le C|x-x'| \times \ell_\alpha (z) \text{ for }z \in B[0,1]^c.
    \end{align*}
\end{enumerate}
\setcounter{Conditions_counter}{\value{enumi}}
\end{enumerate}
\end{assumption}

  The reader can verify that \Cref{assumption:HJB_coefficients_assumption} extends \citet[Remark 2.28]{Hollender2017}.
%
  We describe the Condition \ref{HJB:UI} as ``uniform integrability'' instead of ``tightness'', as we prefer to interpret these conditions as the property of the family of measures to uniformly integrate the tails of the functions $\mathbb{R}^d\ni z\mapsto |z|^2$ and $\mathbb{R}^d\ni z\mapsto |z|^p$ on the respective domains. 
  Of course, the notion of tightness can be also used.
  Additionally, we have added special behaviour of the jump-height coefficient for the case $p<2$.
  In fact, the same assumptions would be required in \citet[Remark 2.28 (Assumptions on Coefficients)]{Hollender2017} hadn't been assumed the (stronger) square integrability. 
  More details about this point will be presented in \cref{subsec:Remarks_on_Assumptions_Relaxations} which is devoted for these technical remarks.
%
%
%
%
%
%
%
%
%
%
%
%
%
%
%

We proceed to the statement of the first result we promised at the beginning of this subsection.
\begin{lemma}[Admissibility of HJB Equation]\label{lemma:admissibility_HJB}
  Suppose that $G$ is a Hamilton--Jacobi--Bellman operator as described in \cref{def:HJB}.
  If \cref{assumption:HJB_coefficients_assumption} holds, then the corresponding Hamilton--Jacobi--Bellman equation
  \begin{align*}
    \partial_t u(t,x) + G(t,x,u(t,x),Du(t,x),D^2(t,x),u(t,\cdot)) = 0
  \end{align*}
  satisfies \cref{assumption:operator_F}.
\end{lemma}
The proof of \cref{lemma:admissibility_HJB} is presented in \cref{proof:lemma_admissibility_HJB}.
We present now the main result of this subsection.
\begin{proposition}[The HJB operators satisfy the regularity condition]\label{HJB:Regularity_Condition}
Suppose that $G$ is a Hamilton--Jacobi--Bellman operator as described in \cref{def:HJB}.
If \cref{assumption:HJB_coefficients_assumption} holds,
then the two operators
\begin{gather*}
  G_1(t,x,r,q,X,\phi):= G(t,x,r,q,X,\phi)\\
  G_2(t,x,r,q,X,\phi):= -G(t,x,-r,-p,-X,-\phi)
\end{gather*}
for $(t,x,r,q,X)\in(0,T)\times \mathbb{R}^d \times \mathbb{R} \times \mathbb{R}^d \times \mathbb{S}^{d\times d} \times C^2_p(\mathbb{R}^d)$ satisfy the regularity condition as described in \cref{def:regularity_condition} for the Young function of \cref{lemma:UI_Young_improvement} associated to the measures $(m_\alpha)_{\alpha\in\mathcal{A}}$ and for $\beta_1=\beta_2=1$.
\end{proposition}

%
%
%
  The part of the proof of \cref{HJB:Regularity_Condition} corresponding to the linear parts $\mathcal{L}_\alpha$ and $\mathcal{I}^\kappa_\alpha$ essentially follows the arguments of \citet[Proposition 2.33]{Hollender2017}.
  The crucial difference is that we substitute the function $\mathbb{R}^d \ni w\mapsto |w|^q$ ($q\ge 2$ is used in \cite[Proposition 2.33]{Hollender2017}) by the function $\mathbb{R}^d \ni w\mapsto \Upsilon(|w|^p)$, for a Young function determined by \cref{lemma:UI_Young_improvement}.
  Although the change from a power function to a more general Young function seems minor at first sight, the required computational details are more involved and rely on the technical results associated to \cref{lemma:UI_Young_improvement}. 
  Indeed, the homogeneity and multiplicativity of power functions, which are extremely convenient in computations, have to be substituted by the moderate growth and the (finally) submultiplicativity property of Young functions. 
  Additionally, the convenience of having an increasing second derivative when dealing with superquadratic power functions has to be substituted by a Young function having a non monotone, but bounded, second derivative.
  All the aforementioned properties associated to Young functions have to be complemented by additional ones.
  In any case, the big picture is the same as in \cite[Proposition 2.33]{Hollender2017}.
  Naturally, as in the proof of \cref{lemma:admissibility_HJB}, some additional care will be needed for the segment $\hat{\mathcal{I}}_\alpha$, since there is an additional term in the integrand comparing to \cite[Proposition 2.33]{Hollender2017}, and -of course- substantial additional effort will be required for the non-linearities $f_\alpha$.

  In order to present the proof in a readable form, we present the main body of the proof below, but we have relegated the lengthy computations, which are presented in a series of lemmata, to \cref{subsec_Appendix:aux_lemmata_Regularity_Condition}.    
%
%
%
\begin{proof}[Proof of \cref{HJB:Regularity_Condition}]\label{proof:HJB:Regularity_Condition}
  Let us denote by $\Upsilon$ the Young function associated to the measures $(m_\alpha)_{\alpha\in\mathcal{A}}$ 
  by \cref{lemma:UI_Young_improvement}.
  Let  $u$,$-v$ $\in \textup{USC}_p([0,T]\times \mathbb{R}^d)$ and for $\delta$,$\mu$,$\gamma$,
  $\lambda \ge 0$ we define
  \begin{align}\label{def:penalization_function}
  \phi_{\delta,\gamma}^{\lambda}(t,x,y) := \frac{\delta}{T-t} + \lambda|x-y|^2 +\gamma e^{\mu t} \big(\Upsilon(|x|^p) + \Upsilon(|y|^p)\big),
  \end{align}
  with $(t,x,y)\in(0,T)\times \mathbb{R}^d \times \mathbb{R}^d$.
  For later reference, 
  \begin{align}
  D_x\phi_{\delta,\gamma}^{\lambda} (t,x,y)=2\lambda(x-y) + p\gamma e^{\mu t}\Upsilon'(|x|^p) |x|^{p-2} x\label{gradient_penalized_x}
  \shortintertext{and}
  D_y\phi_{\delta,\gamma}^{\lambda} (t,x,y)=2\lambda(y-x) + p\gamma e^{\mu t}\Upsilon'(|y|^p) |y|^{p-2} y\label{gradient_penalized_y}.
  \end{align}
  Moreover, we denote by $u^\varepsilon$ (resp. $v_\varepsilon$) the spatial supremal (resp. infimal) convolution\footnote{See the definition in \cref{subsec:notation} and \citet[Remark 2.20]{Hollender2017} for its properties.} of $u$ (resp. $v$) and 
  \begin{align}\label{def:smudge_operator}
    \triangle^{\varepsilon}_\varphi [G^k](t,x,r,q,X,w,\psi)
    :=\inf_{y:\varphi(y-x)\le \varepsilon \delta(x)} G^k(t,y,r,q,X,w\circ \tau_{x-y},\psi\circ \tau_{x-y}) 
  \end{align}
  for $(t,x,r,q,X,w,\psi)\in (0,T)\times \mathbb{R}^d \times \mathbb{R} \times \mathbb{R}^d \times \mathbb{S}^{d\times d} \times \textup{SC}_p(\mathbb{R}^d) \times C^2(\mathbb{R}^d)$ and $\varepsilon>0$, 
  where $\varphi=\varphi_{p}$ is the quasidistance described in \cref{subsec:notation}, $\tau_h:\mathbb{R}^d\to \mathbb{R}^d$ is the translation by $h\in\mathbb{R}^d$ operator, \emph{i.e.}, $\tau_h(x)=x+h$ for all $x,h\in \mathbb{R}^d$, and 
  \begin{align}\label{def:delta_x_phi}
    \delta(x):= 2^{3(p\vee 2)}C(1+\varphi(x)), \text{ for }x\in\mathbb{R}^d,
  \end{align}
  for $C$ depending only on $\|u\|_p$ and $\|v\|_p$.

  Through the rest of the proof we assume that $(\bar{t},\bar{x},\bar{y})\footnotemark\in(0,T)\times \mathbb{R}^d\times \mathbb{R}^d$ is a global maximum point of $u^\varepsilon -v_\varepsilon - \phi_{\delta,\gamma}^{\lambda}$ with%
  \footnotetext{In order to keep the notation as simple as possible we omit the dependence on $(\delta,\mu,\gamma,\lambda,\varepsilon)$. Moreover, we may assume $\mu$ fixed, but arbitrary, so that we will no mention it in the proof hereinafter.}
  \begin{align*}
    M\footnotemark&:=\sup_{(t,x,y)\in (0,T)\times \mathbb{R}^d\times \mathbb{R}^d} \Big( u^\varepsilon(t,x) -v_\varepsilon(t,y) - \phi_{\delta,\gamma}^{\lambda}(t,x,y) \Big)\\
    & = u^\varepsilon(\bar{t},\bar{x}) -v_\varepsilon(\bar{t},\bar{y}) - \phi_{\delta,\gamma}^{\lambda}(\bar{t},\bar{x},\bar{y}) \ge 0 
  \numberthis\label{supremum_penalized}
  \end{align*}%
  \footnotetext{$M$ depends on $(\delta,\gamma,\lambda,\varepsilon)$.}
  and $X,Y\footnotemark\in\mathbb{S}^{d\times d}$, $C>0$ such that 
    \begin{equation}\label{reg_cond:assumption_ineq}
      \begin{multlined}[0.9\textwidth]
      \begin{bmatrix}
        X & \mathbf{0}\\
        \mathbf{0} & -Y
      \end{bmatrix}
        \le 4 \lambda
      \begin{bmatrix}
        I & -I\\
        -I & I
      \end{bmatrix}\\ 
      + C\gamma e^{\mu t}
      \begin{bmatrix}
        \big[\Upsilon''(|\bar{x}|^{p}) |\bar{x}|^{2p-2}+\Upsilon'(|\bar{x}|^{p}) |\bar{x}|^{p-2} \big]  I & \mathbf{0}\\
        \mathbf{0} & \big[\Upsilon''(|\bar{y}|^{p}) |\bar{y}|^{^{2p-2}}+\Upsilon'(|\bar{y}|^{p}) |\bar{y}|^{p-2}\big]  I
      \end{bmatrix}
      \end{multlined}.
    \end{equation}%
  \footnotetext{$X$ and $Y$ depend on $(\delta,\gamma,\lambda,\varepsilon)$.}
    Our aim is to determine a constant $C>0$ and a modulus of continuity $\omega:[0,+\infty) \to [0,+\infty)$ such that \eqref{ineq:result_regularity_condition} holds, \emph{i.e.}, 
      \begin{align*}
      &\begin{multlined}[0.9\textwidth]
      \triangle_\varphi^\varepsilon [G^\kappa]\big(\bar{t},\bar{y},v_\varepsilon(\bar{t},\bar{y}), -D_{y}\phi_{\delta,\gamma}^{\lambda}(\bar{t},\bar{x},\bar{y}),Y,v_\varepsilon(\bar{t},\cdot), -\phi_{\delta,\gamma}^{\lambda}(\bar{t},\bar{x},\cdot)\big)\\
      -\triangle_\varphi^\varepsilon [G^\kappa]\big(\bar{t},\bar{x},u^\varepsilon(\bar{t},\bar{x}), D_{x}\phi_{\delta,\gamma}^{\lambda}(\bar{t},\bar{x},\bar{y}),X,u^\varepsilon(\bar{t},\cdot), \phi_{\delta,\gamma}^{\lambda}(\bar{t},\cdot,\bar{y})\big)
      \\
      \le
       C\lambda |\bar{x}-\bar{y}|^2 + C \gamma e^{\mu \bar{t}} \big(1 + \Upsilon(|\bar{x}|^p) +\Upsilon(|\bar{y}|^p)\big) 
              +\varrho_{\gamma,\lambda,\varepsilon,\kappa}
      \end{multlined}
      \end{align*}
  for a remainder $\varrho_{\gamma,\lambda,\varepsilon,\kappa}$ with
  \begin{align*}
    \limsup_{\gamma \downarrow 0} \limsup_{\lambda\to+\infty}\limsup_{\varepsilon\downarrow 0}\limsup_{\kappa\downarrow 0}\varrho_{\gamma,\lambda,\varepsilon,\kappa}\le 0.
  \end{align*}
  Given the form of the HJB operator $G$, \emph{i.e.}, $G=\inf_{\alpha\in\mathcal{A}}G^\kappa_\alpha$ for every $0<\kappa<1$, and the properties of infinma and suprema of two sets $A$ and $B$, \emph{i.e.}, 
  \begin{align*}
    \inf A - \inf B \le \sup A - \inf B = \sup (A-B).
  \end{align*}
For 
  \begin{align}
    \mathcal{A}_\varepsilon:=\{ (\alpha,x,y) \in\mathcal{A}\times \mathbb{R}^d\times \mathbb{R}^d : \varphi(x-\bar{x})\le \varepsilon\delta(\bar{x}),\varphi(y-\bar{y})\le \varepsilon\delta(\bar{y})\} \text{ for }\varepsilon>0,
    \label{def:A_epsilon}
  \end{align}
  an upper bound for the left-hand side of the inequality we intend to prove is the following
  \begin{align*}
     &\begin{multlined}[0.85\textwidth]
     \triangle_\varphi^\varepsilon [G^\kappa]\big(
     \bar{t},\bar{y},v_\varepsilon(\bar{t},\bar{y}), -D_{y}\phi_{\delta,\gamma}^{\lambda}(\bar{t},\bar{x},\bar{y}),Y,v_\varepsilon(\bar{t},\cdot), -\phi_{\delta,\gamma}^{\lambda}(\bar{t},\bar{x},\cdot)
     \big)\\
    -\triangle_\varphi^\varepsilon [G^\kappa]\big(
    \bar{t},\bar{x},u^\varepsilon(\bar{t},\bar{x}), D_{x}\phi_{\delta,\gamma}^{\lambda}(\bar{t},\bar{x},\bar{y}),X,u^\varepsilon(\bar{t},\cdot), \phi_{\delta,\gamma}^{\lambda}(\bar{t},\cdot,\bar{y})
    \big)
    \end{multlined}\\
    &\hspace{1em}
    \begin{multlined}[0.95\textwidth]
    \le
    \sup_{(\alpha,x,y)\in \mathcal{A}_\varepsilon} \Big\{ 
    f_\alpha\big(\bar{t},x,u^\varepsilon(\bar{t},\bar{x}), \sigma^T_\alpha(\bar{t},x) D_x\phi_{\delta,\gamma}^{\lambda} (\bar{t},\bar{x},\bar{y}), \mathcal{K}^\kappa_\alpha(\bar{t},x,u^\varepsilon(\bar{t},\cdot)\circ\tau_{\bar{x}-x},\phi_{\delta,\gamma}^{\lambda}(\bar{t},\cdot,\bar{y})\circ\tau_{\bar{x}-x})\big)
    \\
    - f_\alpha\big(\bar{t},y,v_\varepsilon(\bar{t},\bar{y}),-\sigma^T_\alpha(\bar{t},y) D_y\phi_{\delta,\gamma}^{\lambda} (\bar{t},\bar{x},\bar{y}), \mathcal{K}^\kappa_\alpha(\bar{t},y,v_\varepsilon(\bar{t},\cdot)\circ\tau_{\bar{y}-y},-\phi_{\delta,\gamma}^{\lambda}(\bar{t},\bar{x},\cdot)\circ\tau_{\bar{y}-y})\big) 
    \Big\} 
    \end{multlined}
    \\
    &\hspace{1.5em}
    + \sup_{(\alpha,x,y)\in \mathcal{A}_\varepsilon} \big\{ 
    \mathcal{L}_\alpha \big(\bar{t},x, D_x\phi_{\delta,\gamma}^{\lambda} (\bar{t},\bar{x},\bar{y}), X\big)
    - \mathcal{L}_\alpha \big(\bar{t},y, -D_y\phi_{\delta,\gamma}^{\lambda} (\bar{t},\bar{x},\bar{y}), Y\big) \big\}
    \\
    &\hspace{1.5em}
    + \sup_{(\alpha,x,y)\in \mathcal{A}_\varepsilon} \Big\{
    \mathcal{I}^\kappa_\alpha\big(\bar{t},x,u^\varepsilon(\bar{t},\cdot)\circ\tau_{\bar{x}-x},\phi_{\delta,\gamma}^{\lambda}(\bar{t},\cdot,\bar{y})\circ\tau_{\bar{x}-x}\big)
    -\mathcal{I}^\kappa_\alpha\big(\bar{t},y,v_\varepsilon(\bar{t},\cdot)\circ\tau_{\bar{y}-y},-\phi_{\delta,\gamma}^{\lambda}(\bar{t},\bar{x},\cdot)\circ\tau_{\bar{y}-y}\big)
             \Big\}\\
   &\hspace{1em}\begin{multlined}
    \leq 
    \lambda C |\bar{x} - \bar{y}|^2 \times 
    +C \gamma e^{\mu \bar{t}}
    \big[1    + \Upsilon(|\bar{x}|^p) + \Upsilon(|\bar{y}|^p)\big] 
    +\varrho^{1}_{\gamma,\lambda,\varepsilon,\kappa}
    +\varrho^{2}_{\gamma,\lambda,\varepsilon}
    +\varrho^{3}_{\gamma,\lambda,\varepsilon,\kappa}
   \end{multlined}
  \end{align*}
with 
\begin{align*}
    \limsup_{\lambda\to\infty}
    \limsup_{\varepsilon\to 0}
    \limsup_{\kappa\to 0}
      [\varrho^{1}_{\gamma,\lambda,\varepsilon,\kappa}
          +\varrho^{2}_{\gamma,\lambda,\varepsilon}
          +\varrho^3_{\gamma,\lambda,\varepsilon,\kappa}]
      =0,
\end{align*}
  where we used the lemmata \ref{lemma_appendix:regularity_f_nonlinearities}, \ref{lemma_appendix:regularity_L_diffusion_parts} and \ref{lemma_appendix:regularity_I_nonlocal_parts} in order to conclude the last inequality for the desired upper bound.
\end{proof}
\begin{corollary}[Comparison Principle for HJB operators]\label{cor:HJB_Comparison_Principle}
Suppose that $G$ is a Hamilton--Jacobi--Bellman operator as described in \cref{def:HJB} that satisfies \cref{assumption:HJB_coefficients_assumption} with $p\ge 1$ and $T>0$.
If $\underbar{u}\in \textrm{USC}_p([0,T]\times\mathbb{R}^d)$, resp. $\bar{u}\in \textrm{LSC}_p([0,T]\times\mathbb{R}^d)$, is a viscosity subsolution, resp. supersolution, in $(0,T)\times \mathbb{R}^d$ of 
\begin{align*}
\partial_t u (t,x) + G\big(t,x,u(t,x),Du(t,x),D^2u(t,x),u(t,\cdot)\big)=0
\end{align*}
such that $\underbar{u}(0,\cdot), \bar{u}(0,\cdot)$ are continuous with $\underbar{u}(0,x) \leq \bar{u}(0,x)$ for every $x\in\mathbb{R}^d$, then 
\begin{align*}
\underbar{u}(t,x)\leq \bar{u}(t,x) \text{ for every }(t,x)\in(0,T)\times \mathbb{R}^d.
\end{align*}
\end{corollary}
\begin{proof}
This is immediate in view of \cref{HJB:Regularity_Condition} and \cref{cor:Comparison_Principle}. 
\end{proof}
\subsection{Remarks on the assumptions and on possible relaxations}\label{subsec:Remarks_on_Assumptions_Relaxations}

We devote this subsection to comments on the imposed framework of \cref{subsec:prominent_example_HJB}.

We start with a comparison between the framework imposed in the current work and that of \citet{Hollender2017}, which -to the best of the authors' knowledge- provides the most general result in the literature. 
The Conditions \ref{HJB:coeff_boundedness} and \ref{HJB:UI} are identical to \citet[Remark 2.28, (B1)-(B2)]{Hollender2017}.
Also, for the families $(b_\alpha)_{\alpha\in\mathcal{A}}$ and $(\sigma_\alpha)_{\alpha\in\mathcal{A}}$ the Condition \ref{HJB:coef_continuity} are identical to \cite[Remark 2.28, (B3)]{Hollender2017}.
The family $(f_{\alpha})_{\alpha\in\mathcal{A}}$ in the current work incorporates the families $(f_{\alpha})_{\alpha\in\mathcal{A}}$ and $(c_{\alpha})_{\alpha\in\mathcal{A}}$ of \cite{Hollender2017}. 
Moreover, Conditions \ref{HJB:monotonicity} is a genuine generalization of \cite[Remark 2.28, (B3) and (B5)]{Hollender2017}.
Indeed, in the current work the family $(f_{\alpha})_{\alpha\in\mathcal{A}}$ is allowed to be locally uniformly equicontinuous in the first and second spatial variables as described by the family of moduli of continuity $(\overline{\omega}_{\!{}_R})_{R\geq 0}$ and $(\widetilde{\omega}_{\!{}_K})_{K\geq 0}$, see \ref{HJB:monotonicity}.\ref{HJB:monotonicity_1} and compare with \citet[Remark 2.28, (B3)]{Hollender2017} in conjunction with the linearity with respect to the second spatial derivative, and monotone in the second spatial variable, see \ref{HJB:monotonicity}.\ref{HJB:monotonicity_2} and compare with \citet[Remark 2.28, (B5)]{Hollender2017} in conjunction with the linearity with respect to the second spatial derivative.
Additionally, in the current work the family $(f_{\alpha})_{\alpha\in\mathcal{A}}$ allows for dependence on additional spatial variables.
In the current work the family of jump-height coefficients $(j_{\alpha})_{\alpha\in\mathcal{A}}$ when $p\geq 2$ satisfies the exact same conditions as those imposed in \cite{Hollender2017}, see \ref{HJB:cond_jump}.\ref{HJB:cond_jump_partA}-\ref{HJB:cond_jump_partB} and compare with \cite[Remark 2.28, (B3)- (B4)]{Hollender2017}.
When $1\leq p<2$ it was required to strengthen the behaviour of the jump-height coefficients as described in \ref{HJB:cond_jump}.\ref{HJB:cond_jump_partC}. 
This seems to be unavoidable once we require the family of measures to respect an integrability of $p-$polynomial order, as the computations showcased in \cref{lemma:FirstOrderTaylor} and \cref{lemma:SecondOrderTaylor}.
A careful inspection of the proof of \citet[Proposition 2.33]{Hollender2017} will convince the reader that an analogous to \ref{HJB:cond_jump}.\ref{HJB:cond_jump_partC} condition would be required, hadn't been assumed that $q\geq 2$ in \cite{Hollender2017}. 

We devote a separate paragraph for the Condition \ref{HJB:growth_cond_delta}.
Regarding its part \ref{HJB:growth_cond_delta_partA}, it cannot be assumed any dependence on a spatial variable of its upper bound since this would lead to a higher than $p-$order polynomial behaviour; \emph{e.g.}, see examine the derivation of \eqref{lem_app:f_nonlinearities_third_3}.
Regarding its part \ref{HJB:growth_cond_delta_partB}, the condition $\sup_{\alpha}\|\ell_{\alpha}\|_{\infty}<\infty$ may be relaxed as soon as the growth conditions in \ref{HJB:cond_jump} are described in terms of a $q-$polynomial growth, for $q<p$ and the sub-/ super- solutions are of $q-$polynomial growth as well.
Regarding its part \ref{HJB:growth_cond_delta_partC}, the condition on the behaviour in $B[0,1]$ is only imposed in order to have the second order term $|z|^2$ appearing in a neighbourhood of the origin, see after \eqref{lem_app:f_nonlinearities_third_2}. 
This condition can be relaxed, say to $|\delta_\alpha(t,x,z) - \delta_\alpha(t,x',z)|  \le C|x-x'| \times \ell_\alpha^2 (z)  \text{ for }z \in B[0,1]$, either when the family $(m_{\alpha})_{\alpha\in\mathcal{A}}$ has a singularity of the first order around at the origin, or the sub-\ super- solutions are locally Lipschitz, see also \citet[Remark 3.9]{barles1997backward} for a simpler case.
We underline that the property of being locally Lipschitz is bequeathed to the (spatial) supremal convolutions, hence we can use this remark in our results too.

Before we proceed let us remark that the restriction $p\geq 1$ required in \cref{def:HJB} is because of the addition (compared to the classical case) of the compensation term $D\phi(x)^Tj_{\alpha}(t,x,z)$ in the segment $\hat{\mathcal{I}}_{\alpha}$.
We can drop this restriction on $p$ either by dropping the compensation term or by retaining it but imposing stronger conditions on the growth of the jump-height coefficients for $z\in B[0,1]^c$.
In the former case, it is easily verified that all the stated results remain true.
This is indeed true, because of the construction of the $\Upsilon$ function which still possesses the desired properties when $p\in(0,1)$; see \cref{lemma:UI_Young_improvement}.
In the latter case, a growth condition of the form $\sup_{\alpha\in\mathcal{A}}  |j_\alpha(t,x,z)|\le C(1+|x|)|z|^p$, for $z\in B[0,1]^c$, and allowing \ref{HJB:cond_jump}.\ref{HJB:cond_jump_partC} to further hold for $p\in(0,1)$, complies with the integrability properties of the measures; see the last part of the proof of \cref{lemma_appendix:regularity_I_nonlocal_parts}.

We have not dealt with the case $p=0$.
On the one hand, if the family of the measures satisfies \ref{HJB:coeff_boundedness}-\ref{HJB:UI} for some $p>0$, then one can readily adapt and combine the arguments presented in the current work and in the proof of \citet[Proposition 2.33]{Hollender2017} for the case $p=0$, see at the end of page 130 of \cite{Hollender2017}.
On the other hand, it seems plausible to drop the assumption \ref{HJB:UI} on integrals over $\{z:|z|>R\}$, for $R\to\infty$, by adapting and combining the arguments presented in \citet{BarlesImbert2008} and \citet[Chapter 1]{Hollender2017}.

In this paragraph we will comment on the nature of the conditions and how these are related to the respective conditions imposed on (decoupled Forward) Backward Stochastic Differential Equations, hereinafter (F)BSDE.
The literature is indeed vast, so we will indicatively mention some results associated to (F)BSDE with jumps. 
Let us start with the simplest case, \emph{i.e.} when the set $\mathcal{A}$ is a singleton.
Then, \citet[Theorem 3.5]{barles1997backward} is a special case of \cref{cor:Comparison_Principle} as regards the uniqueness of the solution; one may take into account the comments we provided for the case $p>0$ which corresponds to the semimartingale case.
Indeed, one can verify that the conditions imposed in \cite{barles1997backward} are all incorporated in the framework of the current work. 
Possibly, one should comment on how to derive the monotonicity condition \ref{HJB:monotonicity}.\ref{HJB:monotonicity_2} from the uniform Lipschitz condition of $f$ as assumed in \cite[Theorem 2.1, (A.1.ii)]{barles1997backward}.
But, to this end, a well-known manoeuvre in the BSDE theory is to apply It\^o's formula with the help of an exponential (in time) function; \emph{e.g.}, see the comments after \citet[Condition (H3), p. 7]{kruse2015bsdes}.
Let us remain in \citet{kruse2015bsdes,kruse2017lp} and examine whether our framework can be embedded in \cite{kruse2017lp}.
We consider an $f$ satisfying \cref{HJB:monotonicity}.
Due to the dependence of $f$ on $\mathcal{K}$ assumed in \cref{def:HJB}, it can be proven that the \cite[Condition (H3), on p. 6]{kruse2017lp} is satisfied.   
Also, let us assume that the integer-valued measure associates to an $\alpha-$stable L\'evy process with index $\alpha\in(1,2)$, this is a case which is treated in \cite{kruse2017lp}, see end of page 5.
As it is well-known, \emph{e.g.}, see \citet[Proposition 14.5]{sato1999levy}, the integrability condition \ref{HJB:coeff_boundedness} is satisfied for $p=1$ and trivially (as $\mathcal{A}$ is singleton) the integrability condition \ref{HJB:UI}.  
Therefore, we may use the associated L\'evy martingale.
All the above lead us to the following two observations.
Firstly, this is a case which cannot be treated by \citet{Hollender2017}.
Secondly, given \cref{cor:Comparison_Principle}, the associated integral equation has at most one solution. 
So, \emph{if} the BSDE solution given by \citet[Theorem 2]{kruse2017lp} determines a solution of the integral equation, then this is unique.
For the case the set $\mathcal{A}$ is not a singleton, one may address \citet[Section 5]{KaziTani2BSDE_PIDE}.
It is relatively easy to verify that the nonlinearity of an HJB operator $G$ as described in  \cref{def:HJB} and \cref{assumption:HJB_coefficients_assumption} satisfies \cite[Assumption 5.1 (ii)-(iv)]{KaziTani2BSDE_PIDE}, while it satisfies a weaker assumption compared to \cite[Assumption 5.1 (v)]{KaziTani2BSDE_PIDE}.
Since the existence result in general requires a less stringent framework, one should complement it with the additional properties assumed in the current work in order to arrive at the uniqueness. 
Finally, we have not dealt with the case of coupled FBSDEs and we have left it for the future in order to curb the technicalities to the decoupled case.

          


\section{Young functions}\label{sec:YoungFunctions}

We present the results associated to Young functions that are going to play a crucial role in the current work. 
For more details on properties of (moderate) Young functions and its conjugates, the interested reader may consult
\citet{long2013martingale,dellacherie1978probabilities,he1992semimartingale,krasnoselskii1961,HUDZIK1992313,rao2002applications}.
In the following, we will choose from the aforementioned references the results which are most conveniently stated for our needsand we will use them to prove that the constructions in \cref{lemma:UI_Young_improvement} satisfy the desired properties.

\begin{definition}\label{def:Young_fun}
A function $\Upsilon:\mathbb R_+\to\mathbb R_+$ will be called \emph{Young}\index{Young function} if it is non-decreasing, convex and satisfies 
\begin{align*}
\Upsilon(0)=0 \text{ and }\lim_{x\to \infty} \frac{\Upsilon(x)}x=\infty.
\end{align*}
Moreover, the non-negative, non-decreasing and right-continuous function $\upsilon:\mathbb R_+\to \mathbb R_+$ for which 
\begin{align*}
\Upsilon(x)=\int_{[0,x]} \upsilon(z)\,\textup{d} z
\end{align*} 
will be called the \emph{right derivative of $\Upsilon$}\index{Young function!right derivative of a}.
\end{definition}

The existence and the properties of the right derivative of $\Upsilon$ are well-known results of convex analysis, \emph{e.g.} see \citet[Theorem 23.1, Theorem 24.1, Theorem 24.2, Corollary 24.2.1]{rockafellar1970convex}.

\begin{definition}\label{notation:Young-constants}
To a Young function $\Upsilon$ with right derivative $\upsilon$ we associate the constant
\begin{align*}
\widebar{c}_{_\Upsilon}:=\sup_{x>0} \frac{x\upsilon(x)}{\Upsilon(x)}.
\end{align*}
A Young function $\Upsilon$ is said to be \emph{moderate (or of moderate growth)}\index{Young function!moderate} if $\widebar{c}_{_\Upsilon}<+\infty$.
\end{definition}
\begin{remark}
 For every Young function $\Upsilon$, one can easily derive that $\widebar{c}_{_\Upsilon}\ge1.$
\end{remark}
The following characterization of moderate Young functions will be useful in the computations that arise at many proofs of the current work.
\begin{lemma}\label{lem:Young_equiv_moderate}
Let $\Upsilon$ be a Young function. 
Then 
\begin{enumerate}
\item\label{lem:Young_equiv_moderate_i} $\Upsilon$ is a moderate Young function if and only if $\Upsilon(\lambda x)\le \lambda^{\widebar{c}_{_\Upsilon}} \Upsilon(x)$ for every $x>0$ and for every $\lambda > 1$.
\item\label{lem:Young_equiv_moderate_ii} $\Upsilon$ is a moderate Young function if and only if there exists $C_2>0$ such that $\Upsilon(2 x)\le C_2 \Upsilon(x)$ for every $x>0$.
\end{enumerate}
\end{lemma}
\begin{proof}
\begin{enumerate}
  \item[\ref{lem:Young_equiv_moderate_i}] See \cite[Theorem 3.1.1 (c)-(d)]{long2013martingale}.
  \item[\ref{lem:Young_equiv_moderate_ii}] In view of \ref{lem:Young_equiv_moderate_i} it is only left to prove that assuming that there exists $C_2>0$ such that $\Upsilon(2 x)\le C_2 \Upsilon(x)$ for every $x>0$, then $\Upsilon$ is moderate. 
      For this see \citet[Definition 10.32, Lemma 10.33.2)]{he1992semimartingale}.
\qedhere
\end{enumerate}
\end{proof}
%
%
%
%
%
%
%
%

The notion of uniform integrability of random variables has long been understood to be connected with the existence of a Young function, a result known as the ``de La Vallée Poussin Theorem''.
In \citet{meyer1978sur}, the author proved that the Young function may have additional desired properties, \emph{i.e.}, being moderate such that the new family of random variables becomes uniformly integrable.  
Following analogous arguments we can derive the existence of a moderate Young function associated to a set of measures that uniformly integrate a function at infinity in the sense of \ref{HJB:UI}.
In the framework of viscosity solutions and recalling the comments presented after \cref{def:regularity_condition}, we are interested in deriving a moderate Young function which is twice continuously differentiable and behaves ``nicely'' at the origin and at infinity.  
The following constructions show how one can obtain a Young function with all the desired properties.

For notational simplicity, in the following results we do not denote the dependence of the Young function on the family $(m_\alpha)_{\alpha\in\mathcal{A}}$ and on $p>0$.
Also, in order to ease the presentation, the lengthy proofs of the following results have been relegated to \cref{appendix:YoungFunctions_proofs}.
\begin{proposition}\label{prop:UI_Young}
  Let $\{m_\alpha\}_{\alpha\in\mathcal{A}}$ be a family of measures on $(\mathbb{R}^d,\mathcal{B}(\mathbb{R}^d))$ such that for some $p>0$ 
  \begin{align*}
    \sup_{\alpha\in\mathcal{A}} \int_{\{z:|z|\ge1\}}|z|^p m_\alpha(\textup{d}z)<+\infty.
  \end{align*}
  Then, 
  \begin{align*}
    \lim_{R\to\infty}\sup_{\alpha\in\mathcal{A}} \int_{\{z:|z|>R\}} |z|^p m_\alpha(\textup{d}z) =0
  \end{align*}
  if and only if 
  there exists a Young function $\Upsilon$ such that
  \begin{align*}
    \sup_{\alpha\in\mathcal{A}} \int_{\{z:|z|\ge 1\}} \Upsilon(|z|^p)m_\alpha(\textup{d}z)<+\infty.
  \end{align*}
\end{proposition}
\begin{lemma}\label{lemma:UI_Young_improvement}
  Let $\{m_\alpha\}_{\alpha\in\mathcal{A}}$ be a family of measures on $(\mathbb{R}^d,\mathcal{B}(\mathbb{R}^d))$ such that for some $p> 0$ 
    \begin{align*}
      \sup_{\alpha\in\mathcal{A}} \int_{\{z:|z|\ge1\}}|z|^p m_\alpha(\textup{d}z)<+\infty
    \end{align*}
  with
    \begin{align*}
      \lim_{R\to\infty}\sup_{\alpha\in\mathcal{A}} \int_{\{z:|z|>R\}} |z|^p m_\alpha(\textup{d}z) =0.
    \end{align*}
  Then, the Young function $\Upsilon$ associated to the family of measures from \cref{prop:UI_Young} can be chosen to additionally have the following properties:
  \begin{enumerate}
      \item\label{lemma:UI_Young_improvement:moderate_C2} $\Upsilon$ is moderate and twice continuously differentiable on $(0,+\infty)$. 
      \item\label{lemma:UI_Young_improvement:second_derivative_properties} 
      Its second derivative, denoted by $ \Upsilon''$, is positive, bounded and finally non-increasing.
      Moreover, there exists $C_{\Upsilon,p}>0$ such that
      \begin{align}\label{property:second_derivative_faster_than_identity}
        0< x\Upsilon''(x) \le C_{\Upsilon,p}, \text{ for every }x>0.
      \end{align}
      \item \label{lemma:UI_Young_improvement:submultiplicativity}
      $\Upsilon$ is finally \emph{submultiplicative}, \emph{i.e.}, there exists $K_{\Upsilon,p}>0$ such that 
      \begin{align*}
        \Upsilon(xy) \le K_{\Upsilon,p}\Upsilon(x) \Upsilon(y)\text{ for every } x,y\geq 1.
      \end{align*}%
      \item\label{lemma:UI_Young_improvement:p_smaller_than_2} For $p\in(0,2)$ it can additionally hold  
        \begin{align}\label{property:second_derivative_faster_to_0}
          \lim_{x\downarrow 0}\Upsilon'(x^p) x^{p-2} = \lim_{x\downarrow 0}\Upsilon''(x^p) x^{2p-2} = 0.
        \end{align}
      \item\label{lemma:UI_Young_improvement:p_between_1_and_2} 
      For $p\in[1,2)$ it can additionally hold  that the function $(0,+\infty) \ni x \longmapsto \Upsilon'(x^p) x^{p-2}\in (0,\infty)$ is finally decreasing, \emph{i.e.}, there exists $R>1$ such that $(R,\infty)\ni x \longmapsto \Upsilon'(x^p) x^{p-2}$ is decreasing.
%
%
      \item\label{lemma:UI_Young_improvement:p_smaller_than_1} 
      For $p\in(0,1)$ it can additionally hold that the function $(0,+\infty) \ni x \longmapsto \Upsilon'(x^p) x^{p-1}\in (0,\infty)$ is finally decreasing, \emph{i.e.}, there exists $R>1$ such that $(R,\infty)\ni x \longmapsto \Upsilon'(x^p) x^{p-1}$ is decreasing.
   \end{enumerate} 
\end{lemma}
\begin{remark}\label{rem:submultiplicative}
In the literature related to convex functions and Orlicz spaces, a Young function satisfies the $\Delta_2-$condition if it is (finally) moderate, \emph{i.e.}, there exist $x_0, R>0$ such that 
\begin{align*}
  \Upsilon(2x) \le R \Upsilon(x), \text{ for }x\ge x_0.
\end{align*}
Additionally, a Young function  is said to satisfy the $\Delta'-$condition if it satisfies a (finally) submultiplicative property as in \cref{lemma:UI_Young_improvement}.\ref{lemma:UI_Young_improvement:submultiplicativity}; see \citet[Chapter I, Section 5]{krasnoselskii1961}.
It is interesting to note that if a Young function satisfies the $\Delta'-$condition, then it also satisfies the $\Delta_2-$condition; see \cite[Chapter I, Lemma 5.1, p.30]{krasnoselskii1961}.
In our case, \emph{i.e.}, the Young function of \cref{lemma:UI_Young_improvement}, we have that the $\Delta_2-$condition is satisfied for every $x>0$.

\end{remark}
\begin{corollary}\label{corrolary:Young_subadditive_additional_properties}
Let $p>0$ and $\Upsilon$ be the Young function determined by \cref{lemma:UI_Young_improvement}.
Then, $\Upsilon$ is subadditive, \emph{i.e.}, there exists $C_{\Upsilon}>0$ such that 
\begin{align*}
\Upsilon(s +t )\le C_{\Upsilon}[ \Upsilon(s) + \Upsilon(t)], \text{ for all }s,t>0.
\end{align*}
Additionally, it satisfies the following inequalities
\begin{gather}
    \Upsilon'(w^p)w^{p-k}  \le C_{\Upsilon,p} \big[1 + \Upsilon(w^p) w^{-k} \mathds{1}_{[1,\infty)}(w)\big] 
    \text{ for } w>0\text{ and for } k\in\{1,2\}
    \label{corollary_Young:Upsilon_prime_upper_bound}
    \shortintertext{as well as}
    \Upsilon''(w^p) w^{2p-2}\leq C_{\Upsilon,p}\big[ 1 + w^{p-2}\mathds{1}_{[1,\infty)}(w)\big]
    \text{ for } w>0.
    \label{corollary_Young:Upsilon_double_prime_upper_bound}
\end{gather}
\end{corollary}
\begin{proof}
The subadditivity is a direct consequence of the convexity and the moderate property of $\Upsilon$.
Indeed, let $s,t>0$, then
\begin{align*} 
\Upsilon(s+t) \leq 2^{\widebar{c}_{_\Upsilon}-1} [\Upsilon(s) + \Upsilon(t)].
\end{align*} 

Let, now, $p>0$ and $w>0$.
We proceed to prove \eqref{corollary_Young:Upsilon_prime_upper_bound}. 
We treat initially the case $k=2$.
Since $$\lim_{w\downarrow 0}\Upsilon'(w^p) w^{p-2}=0\footnotemark$$%
and the function $(0,\infty)\ni w \longmapsto \Upsilon'(w^p) w^{p-2}$ is continuous,
\footnotetext{If $p\in(0,2)$ we make use of \eqref{property:second_derivative_faster_to_0}.} 
 then there exists $C_{\Upsilon,p}$ such that 
$$\sup_{w\in(0,1]}\Upsilon'(w^p) w^{p-2} \leq C_{\Upsilon,p}.$$%
For $w>1$, we have by $\Upsilon$ being moderate that 
$\Upsilon'(w^p) w^{p-2} \leq \widebar{c}_{_\Upsilon}\Upsilon(w^p) w^{-2}$.
Hence, Inequality \eqref{corollary_Young:Upsilon_prime_upper_bound} is true for $k=2$ for the maximum of the constants derived in the two just presented cases.
The case $k=1$ can be easily proven in view of the arguments presented in the case $k=2$.
Indeed, $$\lim_{w\downarrow 0}\Upsilon'(w^p) w^{p-1}=0\footnotemark$$%
and we proceed analogously.
\footnotetext{If $p\in(0,1)$ we make use of \eqref{property:second_derivative_faster_to_0}.}

Let, now, $p>0$ and $w>0$.
We argue now about the validity of \eqref{corollary_Young:Upsilon_double_prime_upper_bound}. 
In view of \eqref{property:second_derivative_faster_than_identity} and of
$$\lim_{w\downarrow 0}\Upsilon''(w^p) w^{2p-2} = 0,\footnotemark$$%
we may follow the arguments used for \eqref{corollary_Young:Upsilon_prime_upper_bound}.%
\footnotetext{If $p\in(0,1)$ we make use of \eqref{property:second_derivative_faster_to_0}}
\end{proof}
At some points in the course of the computations, instead of the inequalities \eqref{corollary_Young:Upsilon_prime_upper_bound} and \eqref{corollary_Young:Upsilon_double_prime_upper_bound},  it will more convenient to use the simpler bounds presented below. 
\begin{corollary}\label{cor:product_derivatives_negative_powers_bounded}
Let $\Upsilon$ be the Young function determined by \cref{lemma:UI_Young_improvement}. 
Then,
\begin{enumerate}
  \item \label{cor:product_first_derivative_negative_powers_bounded:p_1_2}
  For $p\in[1,2)$, the function $(0,\infty) \ni x \longmapsto \Upsilon'(x^p)x^{p-2}\in(0,\infty)$ is bounded by a constant $C_{\Upsilon,p}>0$.
  \item \label{cor:product_first_derivative_negative_powers_bounded:p_0_1}
  For $p\in(0,1)$, the function $(0,\infty) \ni x \longmapsto \Upsilon'(x^p)x^{p-1}\in(0,\infty)$ is bounded by a constant $C_{\Upsilon,p}>0$.
  \item \label{cor:product_second_derivative_negative_powers_bounded}
  For $p\in(0,2)$, the function $(0,\infty) \ni x \longmapsto \Upsilon''(x^p)x^{2p-2}\in(0,\infty)$ is bounded by a constant $C_{\Upsilon,p}>0$.
\end{enumerate}
\end{corollary}
\begin{proof}
The properties are immediate in view of the continuity of $\Upsilon'$, $\Upsilon''$ and of the functions 
$(0,\infty) \ni x \longmapsto x^p \in(0,\infty)$, 
$(0,\infty) \ni x \longmapsto x^{p-2} \in(0,\infty)$,
$(0,\infty) \ni x \longmapsto x^{p-1} \in(0,\infty)$,
$(0,\infty) \ni x \longmapsto x^{2p-2} \in(0,\infty)$
in conjunction with the properties presented in \ref{lemma:UI_Young_improvement:second_derivative_properties}, \ref{lemma:UI_Young_improvement:p_smaller_than_2}-\ref{lemma:UI_Young_improvement:p_smaller_than_1}  of \cref{lemma:UI_Young_improvement}.
\end{proof}
\begin{corollary}\label{corollary:UI_Young_improvement:StrongerIntegrability}
  Let $\{m_\alpha\}_{\alpha\in\mathcal{A}}$ be a family of finite measures on $(\mathbb{R}^d,\mathcal{B}(\mathbb{R}^d))$ such that for some $p>0$ 
    \begin{align*}
      \sup_{\alpha\in\mathcal{A}} \int_{\{z:|z|\ge1\}}|z|^p m_\alpha(\textup{d}z)<+\infty,
    \end{align*}
  with
    \begin{align}\label{property:limit_for_UI}
      \lim_{R\to\infty}\sup_{\alpha\in\mathcal{A}} \int_{\{z:|z|>R\}} |z|^p m_\alpha(\textup{d}z) =0.
    \end{align}

  Then, the Young function $\Upsilon$ associated to the family of measures from \cref{lemma:UI_Young_improvement} can be chosen to additionally satisfy
      \begin{align*}
        \lim_{R\to\infty}\sup_{\alpha\in\mathcal{A}} \int_{\{z:|z|>R\}} \Upsilon(|z|^p) m_\alpha(\textup{d}z) =0.
      \end{align*}
  In other words, there exists a Young function $\Psi$ such that 
  \begin{align*}
    \sup_{\alpha\in\mathcal{A}} \int_{\{z:|z|\ge 1\}} \Psi(\Upsilon(|z|^p))m_\alpha(\textup{d}z) <+\infty.
  \end{align*}
\end{corollary}

\section{Proofs of the results presented in Section \ref{sec:HJB_Equations}.}\label{appendix:HJB_Equations_proofs}
\subsection{Proof of Theorem \ref{theorem:Domination_Principle}}\label{proof:Thm_Domination_Principle}
  The complete proof will be presented for two reasons.
  It is not only for the convenience of the unfamiliar reader that we provide the complete proof, but also because of the fact that some results proven therein play a crucial role in the arguments used in the proof of other lemmata. 

  The structure of the proof closely follows that of the proof of \citet[Theorem 2.22]{Hollender2017}.
  The first steps will be preparatory and in the last step we will prove the contradiction under the hypothesis that the conclusion of the theorem is not true.
  The crucial difference is that instead of a superquadratic power function, a Young function is used.
  However, given that, for any $p>0$, the function $\mathbb{R}^d \ni x \longmapsto \Upsilon(|x|^p)\in[0,\infty)$ finally dominates $\mathbb{R}^d \ni x \longmapsto |x|^p\in[0,\infty)$, the arguments presented in \cite[Theorem 2.22]{Hollender2017} can be applied verbatim.
\begin{proof}[Proof of \cref{theorem:Domination_Principle}]
  We use the notation introduced in \cref{assumption:operator_F,def:regularity_condition}. 
  We will denote by $B[x,r]$ the closed ball of radius $r>0$ centred at $x$; the dimension of the Euclidean space, in which the ball is considered, will be clear from the context.
  For the convenience of the reader, we remind that the penalization function is defined as follows:
  \begin{align}\label{notation_domination:penalization_function}
    \phi_{\delta,\gamma}^{\lambda} (t,x) := \frac{\delta}{T-t} + \lambda \sum_{j=2}^k \sum_{i=1}^{j-1} |x_i-x_j|^2 +\gamma e^{\mu t} \sum_{i=1}^k \Upsilon(|x_i|^p)
  \end{align}
  for $(t,x)=(t,x_1,\ldots,x_d)\in [0,T)\times \mathbb{R}^{kd}$ and $u_i^\varepsilon$ is the spatial supremal convolution of $u_i$, for $i\in\{1,\ldots,k\}$ and $\varepsilon>0$.
  For the spatial supremal convolution we have used the function $\varphi_p$ as described in \cref{subsec:notation}, which is of $p-$polynomial growth.  
  The reader may have in mind that the parameters $\gamma,\delta,\varepsilon>0$ lie in a neighbourhood of $0$, while $\lambda$ will become arbitrarily large. 

  We assume that the conclusion of the theorem does not hold, \emph{i.e.},
  \begin{align*}
    \sup_{(t,x)\in (0,T)\times \mathbb{R}^d}\sum_{i=1}^k \beta_i u_i(t,x) > 0.
  \end{align*}
  Therefore, by the uppersemicontinuity there exists $(t_0,x_0)\in(0,T)\times \mathbb{R}^d$ such that 
  \begin{align}
  \label{domin_principle:def_supremum_point}
  m_0:= \sum_{i=1}^k \beta_i u_i(t_0,x_0) > 0.
  \end{align}
  At this point we choose the parameter $\mu$ so that $\mu \ge C>0$, where $C$ is the constant appearing in \eqref{ineq:result_regularity_condition}, \emph{i.e.}, the constant for which the conclusion of the regularity condition holds.
  Since we fixed $\mu$, our previous remark that we omit it from the notation on $\phi_{\delta,\gamma}^{\lambda}$ is justified; see \cref{footnote:simplified_notation_omit_mu_delta}.

  \begin{enumerate}[label=\textbullet\, \textit{Step \arabic*}:, itemindent =*]
    \item\label{Domin_principle:Step_1} 
  In this step we will prove that for $\gamma>0$, $\delta\in(0,1]$, $\varepsilon\in(0,\varepsilon_0]$, where $\varepsilon_0>0$, and for $\lambda >0$, there exist $(\bar{t},\bar{x})=\big(\bar{t}(\delta,\gamma,\lambda,\varepsilon),\bar{x}(\delta,\gamma,\lambda,\varepsilon)\big)\footnote{For notational convenience, we will omit the dependence on $\delta,\gamma,\lambda$ and $\varepsilon$ whenever no confusion may arise.}\in (0,T)\times \mathbb{R}^{kd}$ such that 
  \begin{gather}
      \sup_{(t,x)\in (0,T)\times \mathbb{R}^{kd}} \Big(\sum_{i=1}^k \beta_i u_i^\varepsilon (t,x_i) - \phi_{\gamma,\delta}^\lambda (t,x)\Big)
      = \sum_{i=1}^k \beta_i u_i^\varepsilon (\bar{t},\bar{x}_i) - \phi_{\gamma,\delta}^\lambda (\bar{t},\bar{x})
        \label{domin_principle:supremal_conv_maximum}
  \shortintertext{and}
      \sup_{\lambda >0} \sup_{\varepsilon\in(0,\varepsilon_0]} \sum_{i=1}^k |\bar{x}_i|^p<+\infty    
      \label{domin_principle:bounded_norms_of_maxima}.
  \end{gather}
  From \citet[Remark 2.20]{Hollender2017} we know that the supremal convolutions and the functions they approximate are finally of the same polynomial growth, \emph{i.e.}, there exists $\varepsilon_0>0$ such that 
  \begin{align*}
    \max_{i\in\{1,\ldots,k\}}|u_i^\varepsilon(t,x)|\le 2^{p+(p\vee 2)-1}C(1+|x|^p),
  \end{align*}
  for $C>\max_{i\in\{1,\ldots,k\}}\|u_i\|_p$,  $(t,x)\in[0,T)\times \mathbb{R}^d$ and $\varepsilon\in(0,\varepsilon_0]$.
  Hence,
  \begin{align*}
    \sum_{i=1}^k \beta_i u_i^\varepsilon (t,x_i) - \phi_{\gamma,\delta}^\lambda (t,x) \le 2^{p+(p\vee 2)-1}C\sum_{i=1}^k \beta_i (1+|x_i|^p) - \gamma e^{\mu t}\sum_{i=1}^k \Upsilon(|x_i|^p),
  \end{align*}
  for $(t,x)=(t,x_1,\dots,x_k)\in[0,T)\times \mathbb{R}^{kd}$ and $\varepsilon\in(0,\varepsilon_0]$.
  The right-hand side of the last inequality tends to $-\infty$ as $|x|\to +\infty$, for every $\gamma>0$, since $$\lim_{y\to+\infty} \frac{\Upsilon(y)}{y}=+\infty.$$
  In other words, if $\gamma >0$, 
  then there exists $R_\gamma> |(x_0,x_0,\dots,x_0)|\ge 0$ such that
  \begin{align*}
      \sup_{(t,x)\in [0,T)\times B[0,R_{\gamma}]^c} \Big( \sum_{i=1}^k \beta_i u_i^\varepsilon (t,x_i) - \phi_{\gamma,\delta}^\lambda (t,x)\Big)
          &\le \sum_{i=1}^k \beta_i u_i (t_0,x_0) - \Big(\frac{1}{T-t_0} + \gamma e^{\mu t_0}  k \Upsilon(|x_0|^p) \Big)\\
          &\le \sum_{i=1}^k \beta_i u_i^\varepsilon (t_0,x_0) - \phi_{\gamma,\delta}^\lambda\big(t_0,(x_0,\ldots,x_0)\big),
  \end{align*}
  for 
  $\delta\in(0,1]$, $\varepsilon\in(0,\varepsilon_0]$, $\lambda > 0$ and $t_0,x_0$ as assumed in \eqref{domin_principle:def_supremum_point}.
  In the second inequality we used the decreasing, pointwise approximation of an upper-semicontinuous function from its supremal convolution; see \cite[Lemma 1.12]{Hollender2017}. 

  We want to extract an analogous bound for values of $t$ close to $T$ and for $x\in B[0,R_{\gamma}]$, 
  for $\gamma > 0$, $\delta\in(0,1]$, $\varepsilon \in (0,\varepsilon_0]$ and $\lambda >0$.
  The first observation is that the functions $\sum_{i=1}^k \beta_i u_i^{\varepsilon}$, $(t,x) \mapsto \lambda \sum_{j=2}^{k} \sum_{i=1}^{j-1} |x_i-x_j|^2$ and $(t,x) \mapsto \gamma e^{\mu t} \sum_{i=1}^k \Upsilon(|x_i|^p)$ remain upper bounded on $[0,T]\times B[0,R_{\gamma}]$, as upper semicontinuous defined on a compact set.
  The second observation is that the function $(0,T)\ni t\mapsto \frac{1}{T-t}$ becomes arbitrarily large as $t\uparrow T$.
  Hence, for every $\gamma>0$, $\delta\in(0,1]$, there exists $T_{\gamma,\delta} \in (t_0,T)$ such that  
  \begin{align*}
    T_{\gamma,\delta} > (1-\delta) T + \delta t_0
    \Leftrightarrow
    -\frac{\delta}{T-T_{\gamma,\delta}} < -\frac{1}{T-t_0}.
  \end{align*}
  The reader may recall, see \cite[Lemma 1.12]{Hollender2017}, that the supremal convolutions decreasingly converge to the function of interest. 
  Consequently, we can choose $T_{\gamma,\delta} \in (t_0,T)$ such that
  \begin{align*}
  &\begin{multlined}[0.85\textwidth]
    \sup_{\varepsilon\in(0,\varepsilon_0]}\bigg\{
    \sup_{(t,x)\in[t_0,T]\times B[0,R_\gamma]}
    \Big\{
        \sum_{i=1}^k\beta_i u_i^{\varepsilon}(t,x) - e^{\mu t_0} \sum_{i=1}^k \Upsilon(|x_i|^p)
    \Big\}\\
        - \big(\sum_{i=1}^k\beta_i u_i^{\varepsilon}(t_0,(x_0,\dots,x_0)) - e^{\mu t_0}  k \Upsilon(|x_0|^p)\big)
    \bigg\}
   \end{multlined}\\ 
    &\hspace{2em}\le \frac{\delta}{T-T_{\gamma,\delta}} - \frac{1}{T-t_0}.
  \end{align*}
  In other words, for every $\gamma>0$ and $\delta\in (0,1]$, there exists $T_{\gamma,\delta}\in (t_0,T)$ such that 
  \begin{align*}
      \sup_{(t,x)\in [T_{\gamma,\delta},T)\times B[0,R_{\gamma}]} \Big( \sum_{i=1}^k \beta_i u_i^\varepsilon (t,x_i) - \phi_{\gamma,\delta}^\lambda (t,x)\Big)
          &\le \sum_{i=1}^k \beta_i u_i^\varepsilon (t_0,x_0) - \phi_{\gamma,\delta}^\lambda\big(t_0,(x_0,\ldots,x_0)\big), 
  \end{align*}
  for every $\lambda>0$ and $\varepsilon\in(0,\varepsilon_0]$.
  
  The two previous bounds allow us to arrive to the desired conclusion, \emph{i.e.}, for every $\gamma>0$, $\delta \in (0,1]$, $\varepsilon \in (0,\varepsilon_0]$ and $\lambda >0$
  \begin{align}
    \sup_{(t,x)\in [0,T)\times \mathbb{R}^{kd}} \Big( \sum_{i=1}^k \beta_i u_i^\varepsilon (t,x_i) - \phi_{\gamma,\delta}^\lambda (t,x)\Big) 
    &=\sup_{(t,x)\in [0,T_{\gamma,\delta}]\times B[0,R_{\gamma}]} \Big( \sum_{i=1}^k \beta_i u_i^\varepsilon (t,x_i) - \phi_{\gamma,\delta}^\lambda (t,x)\Big) \nonumber\\
    &= \sum_{i=1}^k \beta_i u_i^\varepsilon (\bar{t},\bar{x}_i) - \phi_{\gamma,\delta}^\lambda (\bar{t},\bar{x}),
    \tag{\ref{domin_principle:supremal_conv_maximum}}
  \end{align}
  for some $(\bar{t},\bar{x})=\big(\bar{t}(\delta,\gamma,\lambda,\varepsilon),\bar{x}(\delta,\gamma,\lambda,\varepsilon)\big) \in [0,T_{\gamma,\delta}]\times B[0,R_{\gamma}]$, where we used the compactness of the set $[0,T_{\gamma,\delta}]\times B[0,R_{\gamma}]$ and the upper semicontinuity of the function whose supremum is evaluated.
  In particular, for every $\gamma>0$, $\delta \in(0,1]$ 
  \begin{align}
    \sup_{\lambda >0}\sup_{\varepsilon\in(0,\varepsilon_0]} \sum_{i=1}^k |\bar{x}_i|^p   <+\infty.
  \tag{\ref{domin_principle:bounded_norms_of_maxima}}
  \end{align}

  \item\label{Domin_principle:Step_2} In this step we will prove that, for every $\gamma > 0$, $\delta\in(0,1]$
  \begin{align}
  \liminf_{\lambda\to+\infty}\liminf_{\varepsilon \downarrow 0} \lambda \sum_{j=2}^k\sum_{i=1}^{j-1} | \bar{x}_i-\bar{x}_j|^2=
  \limsup_{\lambda\to+\infty}\limsup_{\varepsilon \downarrow 0} \lambda \sum_{j=2}^k\sum_{i=1}^{j-1} | \bar{x}_i-\bar{x}_j|^2=0.\footnotemark
  \label{domin_principle:liminfima_limsuprema_to_0}
  \end{align}%
  \footnotetext{The left limit is sufficient for the rest of the proof of \cref{theorem:Domination_Principle}. The right limit will be used in the proof of \cref{HJB:Regularity_Condition}. So, we also mention it here in order to avoid later the repetition of the arguments.}
  
  For every $\gamma > 0$, $\delta\in(0,1]$, the sets $[0,T_{\gamma,\delta}]\times B[0,R_{\gamma}]\subset [0,T)\times \mathbb{R}^{kd}$ are compact.
  We will initially consider the family of points $\{(\bar{t},\bar{x}) : {\varepsilon\in(0,\varepsilon_0]}\}$, \emph{i.e.}, we consider the parameter $\varepsilon$ free in $(0,\varepsilon_0]$, while the rest parameters $\gamma,\delta$ and $\lambda$ are assumed arbitrarily fixed.
  Let $(\hat{t},\hat{x}) = \big(\hat{t}(\delta,\gamma,\lambda),\hat{x}(\delta,\gamma,\lambda)\big)$ be a limit point of the aforementioned family, \emph{i.e.},
  \begin{gather*}
    \exists (\varepsilon_n{(\delta,\gamma,\lambda)})_{n\in\mathbb{N}}\subset (0,\varepsilon_0] \text{ with }
    \lim_{n\to+\infty} \varepsilon_n{(\delta,\gamma,\lambda)} = 0 \text{ such that }\\
    (\hat{t},\hat{x}) = (\hat{t}(\delta,\gamma,\lambda),\hat{x}(\delta,\gamma,\lambda)) := \lim_{n\to +\infty} \big(\bar{t}(\delta,\gamma,\lambda,\varepsilon_n),\bar{x}(\delta,\gamma,\lambda,\varepsilon_n)\big) \in [0,T_{\gamma,\delta}]\times B[0,R_{\gamma}].
  \end{gather*}
  The upper semicontinuity of the spatial supremal convolution $u^\varepsilon_i$,\footnote{This property is bequeathed to the spatial supremal convolution because of the upper semicontinuity of $u_i$ and the continuity of $\varphi_p$.} for $i=1,\ldots,k$, and the continuity of $\phi_{\gamma,\delta}^\lambda$ on $[0,T_{\gamma,\delta}]\times B[0,R_{\gamma}]$ imply 
  \begin{align*}
    &\sum_{i=1}^k \beta_i u^\varepsilon_i(\hat{t},\hat{x}_i) - \phi_{\gamma,\delta}^\lambda(\hat{t},\hat{x}) \\
    &\hspace{1 em}\ge \limsup_{n\to + \infty} \Big(\sum_{i=1}^k \beta_i u^\varepsilon_i(\bar{t}(\delta,\gamma,\lambda,\varepsilon_n),\bar{x}_i(\delta,\gamma,\lambda,\varepsilon_n))
     - \phi_{\gamma,\delta}^\lambda(\bar{t}(\delta,\gamma,\lambda,\varepsilon_n),\bar{x}(\delta,\gamma,\lambda,\varepsilon_n)) \Big)\\
    &\hspace{1 em}\ge \limsup_{n\to + \infty} \Big( \sum_{i=1}^k \beta_i u^{\varepsilon_n}_i(\bar{t}(\delta,\gamma,\lambda,\varepsilon_n),\bar{x}_i(\delta,\gamma,\lambda,\varepsilon_n))
     - \phi_{\gamma,\delta}^\lambda(\bar{t}(\delta,\gamma,\lambda,\varepsilon_n),\bar{x}(\delta,\gamma,\lambda,\varepsilon_n))\Big),
     \numberthis\label{domin_principle:limsup_varepsilon}
  \end{align*}
  for every $\gamma,\lambda>0$, $\varepsilon \in (0,\varepsilon_0]$ and $\delta\in(0,1]$.
  In the second inequality we have used the monotonicity of $(0,\varepsilon_0]\ni \varepsilon\mapsto u_i^\varepsilon(t,x)$, for every $(t,x)\in [0,T]\times\mathbb{R}^{d}$.\footnote{See \citet[Lemma 1.12, Remark 2.20]{Hollender2017}.}
  Inequality \eqref{domin_principle:limsup_varepsilon} is the first piece of information we needed in order to accomplish our aim.
  
  The second piece of information will be extracted from \eqref{domin_principle:supremal_conv_maximum}, since we have for $(\bar{t}(\delta,\gamma,\lambda,\varepsilon_n),\bar{x}(\delta,\gamma,\lambda,\varepsilon_n))_{n\in\mathbb{N}}$ 
  \begin{align*}
    &\sum_{i=1}^k \beta_i u^{\varepsilon_n}_i\big(\bar{t}(\delta,\gamma,\lambda,\varepsilon_n),\bar{x}_i(\delta,\gamma,\lambda,\varepsilon_n)\big) 
                  - \phi_{\gamma,\delta}^\lambda\big(\bar{t}(\delta,\gamma,\lambda,\varepsilon_n),\bar{x}(\delta,\gamma,\lambda,\varepsilon_n)\big) \\
    &\hspace{1 em}=\sup_{(t,x)\in [0,T)\times \mathbb{R}^{kd}} \Big(\sum_{i=1}^k \beta_i u^{\varepsilon_n}_i(t,x_i) - \phi_{\gamma,\delta}^\lambda(t,x)\Big) 
     \numberthis\label{domin_principle:supremum_subsequence_varepsilon_for_the_Whole_space}
    \\ 
    &\hspace{1 em}\ge \sum_{i=1}^k \beta_i u^{\varepsilon_n}_i\big(\hat{t}{(\delta,\gamma,\lambda)}, \hat{x}_i{(\delta,\gamma,\lambda)}\big) 
                    - \phi_{\gamma,\delta}^\lambda\big(\hat{t}{(\delta,\gamma,\lambda)}, \hat{x}{(\delta,\gamma,\lambda)}\big),
     \numberthis\label{domin_principle:supremum_subsequence_varepsilon_for_hat_points}
  \end{align*}
  for every  $n\in\mathbb{N},$ $\gamma,\lambda>0$ and $\delta\in(0,1]$.

  %

  Fitting together the two pieces of information and omitting the dependencies for the sake of a readable notation, one deduces
  \begin{gather*}
    \sum_{i=1}^k \beta_i u_i\big( \hat{t}, \hat{x}_i \big) 
          -\phi_{\gamma,\delta}^\lambda \big(  \hat{t}, \hat{x}\big)
    \overset{\text{\eqref{domin_principle:limsup_varepsilon}}}{\underset{\text{\eqref{domin_principle:supremum_subsequence_varepsilon_for_hat_points}}}{=}}
    \lim_{n\to +\infty} \Big( \sum_{i=1}^k \beta_i u_i^{\varepsilon_n}\big( \bar{t}, \bar{x} \big) 
          -\phi_{\gamma,\delta}^\lambda \big( \bar{t}, \bar{x}\big) \Big)
    \overset{\text{\eqref{domin_principle:supremum_subsequence_varepsilon_for_the_Whole_space}}}{\ge}
    \sum_{i=1}^k \beta_i u_i\big( t, x \big) 
          -\phi_{\gamma,\delta}^\lambda \big(  t, x\big),
  \end{gather*}
  for $(t,x)\in[0,T)\times \mathbb{R}^{kd}$, \emph{i.e.}, the supremum is attained at the point $(\hat{t}(\delta,\gamma,\lambda),\hat{x}(\delta,\gamma,\lambda))$, for every $\delta\in(0,1]$ and $\gamma,\lambda>0$.
  This property trivially verifies the requirements of \citet[Corollary 2.17]{Hollender2017}\footnote{This is actually \citet[Lemma 3.1]{Crandall1992}.
   } for the upper semicontinuous function
  \begin{align*}
    [0,T) \times \mathbb{R}^{kd} \ni (t,x) \longmapsto \sum_{i=1}^k \beta_i u_i(t,x_i) - \big(\frac{\delta}{T-t} + \gamma e^{\mu t}\sum_{i=1}^k \Upsilon(|x_i|^p)\big).
  \end{align*}
  Consequently,  we are eligible to apply the aforementioned corollary and make use of the limit
  \begin{align}
      \lim_{\lambda\to +\infty} \lambda \sum_{j=2}^k\sum_{i=1}^{j-1}|\hat{x}_i(\delta,\gamma,\lambda) - \hat{x}_j(\delta,\gamma,\lambda)|^2 =0, 
      \label{domin_principle:corollary_2.17_limit}
  \end{align}
  for every $\delta\in(0,1]$ and $\gamma>0$.
  
  We can conclude the validity of \eqref{domin_principle:liminfima_limsuprema_to_0} because all the above hold for the arbitrary limit point of the family $\{(\bar{t},\bar{x}) : {\varepsilon\in(0,\varepsilon_0]}\}$. 
  Let us consider the case with the limit suprema; the case with the limit infima is proven analogously.
  By definition of the limit supremum, there exists $(\varepsilon_n)_{n\in\mathbb{N}}$ such that 
  \begin{align*}
      \limsup_{\varepsilon\downarrow 0} \sum_{j=2}^k\sum_{i=1}^{j-1} |\bar{x}_i(\delta,\gamma,\lambda,\varepsilon) - \bar{x}_j(\delta,\gamma,\lambda,\varepsilon)|^2
      &=\lim_{n\to\infty} \sum_{j=2}^k\sum_{i=1}^{j-1} |\bar{x}_i(\delta,\gamma,\lambda,\varepsilon_n) - \bar{x}_j(\delta,\gamma,\lambda,\varepsilon_n)|^2\\
      &=\lim_{m\to\infty} \sum_{j=2}^k\sum_{i=1}^{j-1} |\bar{x}_i(\delta,\gamma,\lambda,\varepsilon_{n_m}) - \bar{x}_j(\delta,\gamma,\lambda,\varepsilon_{n_m})|^2,
  \end{align*}
  where we have  passed (possibly successively, at most $k-$many times) to a subsequence $(\varepsilon_{n_m})_{n\in\mathbb{N}}$ so that 
  \begin{align*}
    \lim_{n\to +\infty} \bar{x}_i(\delta,\gamma,\lambda,\varepsilon_{n_m}) = \hat{x}_i(\delta,\gamma,\lambda), \text{ for }i=1,\ldots,k,
  \end{align*}
  for $\delta\in(0,1]$, $\gamma,\lambda >0$, for some $\big(\hat{x}_i(\delta,\gamma,\lambda)\big)_{i=1,\dots,k}$.
  This is possible since for $\delta\in(0,1]$, $\gamma >0$ all the aforementioned points lie in the compact set $B[0,R_\gamma]$.
  In other words,
  \begin{align*}
    0 &\le \limsup_{\lambda\to +\infty} \limsup_{\varepsilon \downarrow 0} \Big(\lambda \sum_{j=2}^k\sum_{i=1}^{j-1}|\bar{x}_i(\delta,\gamma,\lambda,\varepsilon) - \bar{x}_j(\delta,\gamma,\lambda,\varepsilon)|^2\Big)\\
      &= \limsup_{\lambda\to +\infty} \lim_{m\to +\infty} \Big(\lambda \sum_{j=2}^k\sum_{i=1}^{j-1}|\bar{x}_i(\delta,\gamma,\lambda,\varepsilon_{n_m}) - \bar{x}_j(\delta,\gamma,\lambda,\varepsilon_{n_m})|^2\Big)\\
      &= \limsup_{\lambda\to +\infty} \Big(\lambda \sum_{j=2}^k\sum_{i=1}^{j-1}|\hat{x}_i(\delta,\gamma,\lambda) - \hat{x}_j(\delta,\gamma,\lambda)|^2\Big)\overset{\text{\eqref{domin_principle:corollary_2.17_limit}}}{=}0.
  \end{align*}

  \item\label{Domin_principle:Step_3} In this step it will be proven that there exist $\delta_0,\gamma_0>0$ such that for every $\delta\in(0,\delta_0]$ and $\gamma\in(0,\gamma_0]$  it holds $\bar{t}=\bar{t}(\delta,\gamma,\lambda,\varepsilon)>0$, for all but finitely many $\lambda>0$ and $\varepsilon\in(0,\varepsilon_0]$.

  Recall \eqref{domin_principle:def_supremum_point}, \emph{i.e.}, there exists $(t_0,x_0)\in(0,T)\times \mathbb{R}^d$ such that
  \begin{align}
    \tag{\ref{domin_principle:def_supremum_point}}
    m_0:= \sum_{i=1}^k \beta_i u_i(t_0,x_0) > 0.
  \end{align}
  Then, one immediately has that there exist $\delta_0,\gamma_0>0$ such that for every $\delta\in(0,\delta_0]$, $\gamma\in(0,\gamma_0]$
  \begin{align*}
    \sum_{i=1}^k \beta_i u_i(t_0,x_0) - \phi_{\delta,\gamma}^{\lambda} \big(t_0,(x_0,\ldots,x_0)\big) 
    \ge m_0 -\Big(\frac{\delta}{T-t_0} + \gamma e^{\mu t_0} k\Upsilon(|x_0|^p) \Big) \ge \frac{m_0}{2}>0. 
  \end{align*}
  Then, for every $\delta\in(0,\delta_0]$, $\gamma\in(0,\gamma_0]$, $\lambda>0$ and $\varepsilon\in(0,\varepsilon_0]$ we obtain
  \begin{align*}
    \sum_{i=1}^k \beta_i u_i^\varepsilon(\bar{t},\bar{x}_i) - \phi_{\delta,\gamma}^{\lambda} (\bar{t},\bar{x})
      &= \sup_{(t,x)\in [0,T)\times \mathbb{R}^{kd}} \Big( \sum_{i=1}^k \beta_i u_i^\varepsilon (t,x_i) - \phi_{\delta,\gamma}^{\lambda} (t,x) \Big)\\ 
      &\ge \sum_{i=1}^k \beta_i u_i(t_0,x_0) - \phi_{\delta,\gamma}^{\lambda} \big(t_0,(x_0,\ldots,x_0)\big)\ge \frac{m_0}{2} >0,
      \numberthis
      \label{domin_principle:bar_t_preliminary_contradiction}
  \end{align*}
  where in the first inequality we have used that the supremal convolutions monotonically approximate, from above, the desired functions; see \citet[Lemma 1.12]{Hollender2017}.
  
  We proceed by assuming that the statement of \hyperref[Domin_principle:Step_3]{Step 3} is false, \emph{i.e.}, for some $\tilde\delta\in(0,\delta_0]$ and $\tilde\gamma \in (0,\gamma_0]$ there exist $(\tilde\varepsilon_n)_{n\in\mathbb{N}}\subset (0,\varepsilon_0]$ with $\lim_{n\to+\infty}\tilde\varepsilon_n=0$ and $(\tilde\lambda_n)_{n\in\mathbb{N}}\subset(0,+\infty)$ with $\lim_{n\to+\infty}\tilde\lambda_n=+\infty$ such that 
  \begin{align}\label{domin_principle:bar_t_preliminary_contradiction_3}
    \bar{t} = \bar{t}(\tilde\delta,\tilde\gamma,\tilde\lambda_n,\tilde\varepsilon_n) = 0 \text{ for infinitely many }n\in\mathbb{N}.\footnotemark
  \end{align}%
  \footnotetext{Without loss of generality, it may be assumed that every element is equal to zero.}%
  The pointwise, monotone convergence of the supremal convolutions allow for convergence in the supremum norm whenever the approximated function is continuous and the domain of approximation is compact; recall Dini's theorem.
  Having assumed that $u_i(0,\cdot)$ is continuous, for $i=1,\ldots,k$, we have for 
  $$(\bar{t},\bar{x})= \big(\bar{t}(\tilde\delta,\tilde\gamma,\tilde\lambda_n,\tilde\varepsilon_n),\bar{x}(\tilde\delta,\tilde\gamma,\tilde\lambda_n,\tilde\varepsilon_n)\big)\in[0,T_{\tilde{\gamma},\tilde{\delta}}]\times B[0,R_{\tilde{\gamma}}]$$
  that there exists $\tilde n_0=\tilde n_0(\tilde\gamma)\in\mathbb{N}$ such that
  \begin{align*}
    \sum_{i=1}^k \beta_i u_i^{\tilde\varepsilon_n}(\bar{t},\bar{x}_i) - \sum_{i=1}^k \beta_i u_i(\bar{t},\bar{x}_i) 
    \le \sum_{i=1}^k \beta_i \sup_{\{x\in\mathbb{R}^{kd} : |x|\le R_{\tilde\gamma}\}} \big( u_i^{\tilde\varepsilon_n}(0,x) -u_i(0,x)\big) 
    \le \frac{m_0}{4}, \text{ for } n\ge \tilde n_0.
  \end{align*}
  In particular, rearranging the terms of the above inequality 
  \begin{align*}
    &\sum_{i=1}^k \beta_iu_i(\bar{t},\bar{x}_i)
      \overset{\text{\phantom{\eqref{domin_principle:bar_t_preliminary_contradiction}}}}{\ge} 
      \sum_{i=1}^k \beta_iu_i^{\tilde\varepsilon_n}(\bar{t},\bar{x}_i) - \frac{m_0}{4}\\
      &\hspace{1 em}\overset{\text{\phantom{\eqref{domin_principle:bar_t_preliminary_contradiction}}}}{\ge} 
      \Big(\sum_{i=1}^k \beta_iu_i^{\tilde\varepsilon_n}(\bar{t},\bar{x}_i) - \phi_{\delta,\gamma}^{\lambda}(\bar{t},\bar{x}) \Big) - \frac{m_0}{4}
      \overset{\text{\eqref{domin_principle:bar_t_preliminary_contradiction}}}{\ge}\frac{m_0}{4}>0, \text{ for } n\ge \tilde n_0.
      \numberthis\label{domin_principle:bar_t_preliminary_contradiction_2}
  \end{align*}
  From \hyperref[Domin_principle:Step_2]{Step 2} and the compactness of the set $[0,T_{\tilde\gamma,\tilde\delta}]\times B[0,R_{\tilde\gamma}]\subset [0,T)\times \mathbb{R}^{kd}$ we further deduce that there exist subsequences
  $(\tilde\varepsilon_{n_m})_{m\in\mathbb{N}}$ and $(\tilde\lambda_{n_m})_{m\in\mathbb{N}}$ such that 
  \begin{align*}
    (0,\hat{x})=(0,\hat{x}(\tilde{\delta},\tilde{\gamma})):=\lim_{m\to+\infty}\big(\bar{t}(\tilde\delta,\tilde\gamma,\tilde\lambda_{n_m},\tilde\varepsilon_{n_m}),\bar{x}(\tilde\delta,\tilde\gamma,\tilde\lambda_{n_m},\tilde\varepsilon_{n_m})\big),
  \end{align*}
  with $\hat{x}_1=\ldots=\hat{x}_k$.
  Finally, from the upper-semicontinuity we have
  \begin{align*}
    \sum_{i=1}^k \beta_i u_i(0,\hat{x}_1)  
    = \sum_{i=1}^k \beta_i u_i(0,\hat{x}_i)  
    \ge \limsup_{n\to\infty} \sum_{i=1}^k \beta_i u_i(\bar{t},\bar{x}_i)
    \overset{\text{\eqref{domin_principle:bar_t_preliminary_contradiction_2}}}{\ge} \frac{m_0}{4} >0, 
  \end{align*}
  which contradicts the assumption of the statement of the theorem, \emph{i.e.}, $\sum_{i=1}^k \beta_i u_i(0,x)\le 0$ for every $x\in\mathbb{R}^d$.
  The contradiction arose because of the hypothesis \eqref{domin_principle:bar_t_preliminary_contradiction_3}. 
  Hence, this hypothesis has to be false and the statement of \hyperref[Domin_principle:Step_3]{Step 3} is indeed correct.
  \item\label{Domin_principle:Step_4} In this step we will use the information extracted in the previous steps, so that we can prove the validity of the theorem. 
  
  For every $\delta \in(0,\delta_0]$, $\gamma \in(0,\gamma_0]$, $\lambda >0$ and $\varepsilon\in(0,\varepsilon_0]$, there exists $(\bar{t},\bar{x})\in(0,T)\times B[0,R_\gamma]$ such that \eqref{domin_principle:supremal_conv_maximum} holds.
  The parabolic Maximum Principle for Generalized Operators, see \citet[Corollary 2.19 applied for $\gamma=\frac{1}{2k}$]{Hollender2017}, in conjunction with \cite[Remark 2.20]{Hollender2017} imply the existence of 
  \begin{align*}
    b_i=b_i(\delta,\gamma,\lambda,\varepsilon) \in \mathbb{R}\ \text{ and }\ X_i = X_i(\delta,\gamma,\lambda,\varepsilon) \in\mathbb{S}^{d\times d} 
   \ \text{ satisfying }\
    \sum_{i=1}^k \beta_i b_i = \partial_t\phi_{\delta,\gamma}^{\lambda}(\bar{t},\bar{x}) 
  \end{align*}
  and
    \begin{align*}
    \begin{multlined}[0.85\textwidth]
    \textup{diag}\big( X_1,\ldots,X_k\big) \le 4\lambda J_{k,d}\\
    + C\gamma e^{\mu \bar{t}}\textup{diag}
    \left[(\Upsilon''\big(|\bar{x}_1|^{p}) |\bar{x}_1|^{2p-2}+\Upsilon'(|\bar{x}_1|^{p} ) |\bar{x}_1|^{p-2} \big)I,
                \ldots, 
              \big(\Upsilon''(|\bar{x}_k|^{p}) |\bar{x}_k|^{2p-2}+\Upsilon'(|\bar{x}_k|^{p} ) |\bar{x}_k|^{p-2}\big) I \right],   
     \end{multlined} 
    \end{align*}
  such that the following inequalities hold in the ordinary sense\footnote{Recall that $\varphi$ denotes the quasidistance used in the convolution.}
  \begin{align*}
  b_i + \triangle_\varphi^\varepsilon[G_i^{\kappa}]\big(\bar{t},\bar{x}_i,u_i^\varepsilon(\bar{t},\bar{x}_i),p_i,X_i,u_i^\varepsilon(\bar{t},\cdot), \beta_i^{-1}\phi_{\delta,\gamma}^{\lambda}(\bar{t},\bar{x}_1,\ldots,\bar{x}_{i-1},\cdot,\bar{x}_{i+1},\ldots,\bar{x}_k) \big)\le 0
  \end{align*}
  for $i\in\{1,\ldots,k\}$ and $\kappa\in(0,1)$, where $p_i=p_i(\delta,\gamma,\lambda,\varepsilon):=\beta^{-1}D_{x_i}\phi_{\delta,\gamma}^{\lambda}(\bar{t},\bar{x})$ for $i\in\{1,\ldots,k\}$.

  We can now use the regularity condition for the $\beta_1,\ldots,\beta_k>0$ and the Young function $\Upsilon$ in order to derive the inequality
  \begin{align*}
      \begin{multlined}[0.85\textwidth]
      \frac{\delta}{(T-\bar{t})^2} + \mu \gamma e^{\mu \bar{t}} \sum_{i=1}^k \Upsilon(|\bar{x}_i|^p) 
        = \partial_t \phi_{\delta,\gamma}^{\lambda} (\bar{t},\bar{x}) 
        = \sum_{i=1}^k \beta_i b_i
        \\ 
      \le
          \overline{C} \lambda \sum_{i<j}   |\bar{x}_i - \bar{x}_j|^2
          + \overline{C} \gamma e^{\mu \bar{t}} \Big( 1 + \sum_{i=1}^k \Upsilon(|\bar{x}_i|^p)\Big) + \varrho_{\gamma,\lambda,\varepsilon,\kappa}
      \end{multlined}
  \end{align*}

 The reader may recall that the parameter $\mu$ was chosen such that $\mu\ge \overline{C}$, where $\overline{C}>0$ is the constant appearing in the regularity condition; it may now be assumed $\mu = \overline{C}>0$ without loss of generality.
 In total, one derives
  \begin{align*}
      0<&\frac{\delta}{T^2} < \frac{\delta}{(T-\bar{t})^2}
      \le
           \overline{C} \lambda \sum_{i<j}   |\bar{x}_i - \bar{x}_j|^2
          + \overline{C} \gamma e^{\mu \bar{t}} + \varrho_{\gamma,\lambda,\varepsilon,\kappa}
  \end{align*}
  for $\delta \in (0,\delta_0]$, $\gamma \in (0,\gamma_0]$, $\lambda >0$ and $\varepsilon \in (0,\varepsilon_0]$, which in turns leads to a contradiction once we observe that (the order of the limits is crucial) 
  \begin{gather*}
      \limsup_{\gamma \downarrow 0}
      \limsup_{\lambda \to +\infty}
      \limsup_{\varepsilon \downarrow 0}
      \limsup_{\kappa \downarrow 0}
          \Big\{\lambda \sum_{i<j}   |\bar{x}_i - \bar{x}_j|^2 + \gamma e^{\mu \bar{t}} + \varrho_{\gamma,\lambda,\varepsilon,\kappa} \Big\}=0.
    \end{gather*}
\end{enumerate}      

\end{proof}
\subsection{Proof of Lemma \ref{lemma:admissibility_HJB}}\label{proof:lemma_admissibility_HJB}

\begin{proof}
  Given the fact that the operator $G^{\kappa}$ is described as the infimum over $\mathcal{A}$ of the family of operators $(G^{\kappa}_{\alpha})_{\alpha\in\mathcal{A}}$, we will examine each property initially for the operators $G^{\kappa}_{\alpha}$, for each $\alpha\in\mathcal{A}$. 
  In particular, for the linear parts of the generator $G^{\kappa}_{\alpha}$, \emph{i.e.}, $\mathcal{L}_{\alpha}$ and $\mathcal{I}^{\kappa}_{\alpha}$, we will essentially follow \citet[Lemma 2.32]{Hollender2017}, but with a minor difference in the segment $\hat{\mathcal{I}}_\alpha$ of the nonlocal part of the operator $G$.
  Recalling \cref{rem:same_visc_framework}.\ref{rem:same_visc_framework_1}, we only need to verify that the requirements of \cref{assumption:operator_F} are met for the nonlocal part $ \hat{\mathcal{I}}_\alpha$ as well as for the non-linearities $f_\alpha$, for every $\alpha \in \mathcal{A}$. 

  For the consistency requirement \ref{ass:operator_consistency}, we have immediately that it is verified by the definition of the Hamilton--Jacobi--Bellman operator $G$.
  For the degenerate ellipticity requirement \ref{ass:operator_degen_ellipt}, we can follow exactly the same arguments as in \citet[Lemma 2.32]{Hollender2017} for the terms $\mathcal{L}_{\alpha}, \mathcal{I}^{\kappa}_{\alpha}$ , since the presence of the equal terms $D\phi(t,x),D\psi(t,x)$ in the $\hat{\mathcal{I}}_\alpha$ segment of the nonlocal part is innocuous.
  Next, we need to verify that whenever $u-v$ and $\phi-\psi$ have global maxima in $(t,x)$, then 
  \begin{align*}
    -f_{\alpha}(t,x,r,q^T\sigma_{\alpha}(t,x), \mathcal{K}^{\kappa}_{\alpha}(t,x,u,\phi)) 
    \ge     -f_{\alpha}(t,x,r,q^T\sigma_{\alpha}(t,x), \mathcal{K}^{\kappa}_{\alpha}(t,x,v,\psi)).
  \end{align*}
  This is also immediately true in view of the definition of the operator $\mathcal{K}^{\kappa}_{\alpha}$  and the monotonicity property \ref{HJB:monotonicity}.
  Indeed, since $u-v$ and $\phi-\psi$ have global maxima in $(t,x)$, then for every $(t',x')$
  \begin{gather*}
    u(t',x') - u(t,x) \le v(t',x') - v(t,x) \text { and }
    \phi(t',x') - \phi(t,x) \le \psi(t',x') - \psi(t,x).
  \end{gather*}
  Consequently, since $m_{\alpha}$ are (positive) measures and $\delta_{\alpha}\geq 0$ (see \ref{HJB:growth_cond_delta}.\ref{HJB:growth_cond_delta_partA}), we have by means of \ref{HJB:monotonicity}.\ref{HJB:monotonicity_2}
  \begin{align*}
    f_{\alpha} \big( t,x,r,q^T\sigma_{\alpha}(t,x),\mathcal{K}^{\kappa}_{\alpha}(t,x,u,\phi)\big) \le 
    f_{\alpha} \big( t,x,r,q^T\sigma_{\alpha}(t,x),\mathcal{K}^{\kappa}_{\alpha}(t,x,v,\psi)\big),
  \end{align*}
  for every $(t,x,r,q)\in(0,T)\times \mathbb{R}^d\times \mathbb{R}\times \mathbb{R}^d$.

  The translation invariance requirement \ref{ass:operator_transl_invar} is also met.
  Indeed, in the non-local operators appear differences of the inputs, hence the constant vanishes, and the derivative of a constant also vanishes.
  We bounce to the monotonicity requirement \ref{ass:operator_monotonicity}, which is by definition met.
  Indeed, our claim is verified by the monotonicity assumption \ref{HJB:monotonicity} in conjunction with the fact that the dependence in the third variable of the operator $G^{\kappa}$ comes into play only through the non-linearities $(f_{\alpha})_{\alpha\in\mathcal{A}}$. 
  
  The continuity requirement \ref{ass:operator_continuity} is slightly involved.
  Therefore, we will provide all the details for the segment $\hat{\mathcal{I}}_\alpha$ of the nonlocal part $\mathcal{I}^{\kappa}_\alpha$, \emph{i.e.}, the term dealing with the ``large jump'' part which we have altered compared to \cite{Hollender2017}, as well as for the non-linearity $f_{\alpha}$.
  For the former, we essentially combine the arguments presented in \citet[Lemma 2.32, p.132]{Hollender2017} and which deal with the medium and the large jump terms.
  Analogous arguments may be employed to prove the continuity of the part associated to the non-linearities $(f_{\alpha})_{\alpha\in\mathcal{A}}$.
  We will work under the assumptions
  \begin{enumerate}[label=\textbullet]
    \item $\displaystyle \lim_{n\to+\infty} (t_n,x_n,r_n,q_n) =(t,x,r,q)$,\footnote{Since $X$ does not appear in $\hat{\mathcal{I}}_\alpha$ and in $f_{\alpha}$, we omit it in order to lighten the notational burden.}
    \item $\displaystyle \lim_{n\to+\infty} D^m\phi_n = D^m\phi$ locally uniformly for all $m\in\{0,1,2\}$\footnote{In order to unify the notation, we understand $\displaystyle \lim_{n\to+\infty} D^0\phi_n = D^0\phi$ locally uniformly as $\displaystyle \lim_{n\to+\infty} \phi_n = \phi$ locally uniformly.} and
    \item $\displaystyle \lim_{n\to+\infty} u_n = u$ locally uniformly with $u\in C_p([0,T]\times \mathbb{R}^d)$ and $\sup_{n\in\mathbb{N}} \|u_n\|_p<+\infty$.
  \end{enumerate}
  \vspace{0.5em}
  {\tiny $\blacksquare$} We initially provide some auxiliary remarks and estimations for the term $\mathcal{I}^{\kappa}_{\alpha}$:\footnote{Some of them will become handy also for the case of the non-linearity $f_{\alpha}$.}
  \begin{enumerate}
    \item\label{admissibility:aux_comp_1} In view of the growth condition \ref{HJB:cond_jump}.\ref{HJB:cond_jump_partB}, the quantity 
    $\sup_{n\in\mathbb{N}}\sup_{\alpha\in\mathcal{A}}|j_\alpha(t_n,x_n,z)|$ is finite on any compact set $\{z:1<|z|<R\}$.
    Indeed, the sequence $(x_n)_{n\in\mathbb{N}}$ is bounded as convergent sequence and consequently
    \begin{align}\label{bound:jump}
        \sup_{n\in\mathbb{N}}\sup_{\alpha\in\mathcal{A}}|j_\alpha(t_n,x_n,z)| \le C(1+\sup_{n\in\mathbb{N}}|x_n|)R,
        \text{ for }z \text{ s.t. }|z|<R.
    \end{align}


    \item\label{admissibility:aux_comp_2}   
      We provide some auxiliary computations for the absolute difference of the integrands appearing in $\hat{\mathcal{I}}_\alpha$:
      \begin{align*}
        &\begin{multlined}[0.85\textwidth]
        \big| \big(
            u_n(t_n,x_n+j_\alpha(t_n,x_n,z)) - u_n(t_n,x_n) - D\phi_n(t_n,x_n)^T j_\alpha(t_n,x_n,z)
          \big)\\
        -
          \big( u(t,x+j_\alpha(t,x,z)) - u(t,x) - D\phi(t,x)^T j_\alpha(t,x,z) \big) \big|
        \end{multlined}\\
        &\hspace{1em}\begin{multlined}[0.85\textwidth]
        \le  
              \big| u(t,x+j_\alpha(t,x,z)) - u_n(t_n,x_n+j_\alpha(t_n,x_n,z)) \big|
            + \big| u(t,x) - u_n(t_n,x_n) \big|\\
          +\big| D\phi(t,x)j_\alpha(t,x,z) - D\phi_n(t_n,x_n)j_\alpha(t_n,x_n,z) \big|
        \end{multlined}\\
        &\hspace{1em}\begin{multlined}[0.85\textwidth]
        \le  
              \big| u(t,x+j_\alpha(t,x,z)) - u(t_n,x_n+j_\alpha(t_n,x_n,z)) \big| + \big| u(t,x) - u(t_n,x_n) \big|\\
        + \big| u(t_n,x_n+j_\alpha(t_n,x_n,z)) - u_n(t_n,x_n+j_\alpha(t_n,x_n,z)) \big| + \big| u(t_n,x_n) - u_n(t_n,x_n) \big|\\
          +\big| D\phi(t,x) \big| \cdot |j_\alpha(t,x,z) - j_\alpha(t_n,x_n,z)|
          +\big| D\phi(t,x) - D\phi_n(t_n,x_n) \big| \cdot |j_\alpha(t_n,x_n,z)|
        \end{multlined}
      \end{align*} 

    \item\label{admissibility:aux_comp_3} 
    For $R>1$,  we have  on the set $\{z: |z|\le R\}$
    \begin{align*}
      &\sup_{\alpha\in\mathcal{A}, |z|\le R} \Big\{\big| u(t_n,x_n+j_\alpha(t_n,x_n,z)) - u_n(t_n,x_n+j_\alpha(t_n,x_n,z)) \big| + \big| u(t_n,x_n) - u_n(t_n,x_n) \big|\Big\}\\
      &\hspace{1em}
        \le 2\sup\big\{  \big| u(t_n,x_n+\zeta) - u_n(t_n,x_n+\zeta) \big| : \text{ for } \zeta \text{ s.t. } |\zeta|\le C(1+\sup_{n\in\mathbb{N}}|x_n|)R\big\}\\
      &\hspace{1em}  \le 2 \| (u - u_n)\mathds{1}_{K_1}\|_{\infty} \xrightarrow{n\to+\infty}0,  \numberthis\label{upper_bound:limit_3}
    \end{align*}
    for $K_1$ compact set such that
    \begin{align*}
      \{ (t_n,x_n +\zeta): \text{ for }n,\zeta \text{ s.t. } n\in\mathbb{N}\text{ and } |\zeta|\le C(1+\sup_{n\in\mathbb{N}}|x_n|)R\}  \subset K_1.
    \end{align*}

    \item\label{admissibility:aux_comp_4} 
    In view of the continuity condition \ref{HJB:cond_jump}.\ref{HJB:cond_jump_partA}, the quantity $j_\alpha(t_n,x_n,z)$ is in the ball with centre $j_\alpha(t,x,z)$ and of radius $|z|\big( \omega (|t-t_n|) + C|x-x_n|\big)$, \emph{i.e.},
    \begin{align*}
      j_\alpha(t_n,x_n,z) \in B\Big(j_\alpha(t,x,z), |z|\big( \omega (|t-t_n|) + C|x-x_n|\big)\Big).  
    \end{align*}
    In the following we may increase, without harm, the constants $C$ appearing in \ref{HJB:cond_jump} such that it is greater than or equal to 1.
    For $R>1$, we have on the set $\{z: |z|\le R\}$
    \begin{align*}
    \begin{multlined}[c][0.95\displaywidth]
      \sup_{\alpha\in\mathcal{A},|z|\le R}\big| u(t,x+j_\alpha(t,x,z)) - u(t_n,x_n+j_\alpha(t_n,x_n,z)) \big| \\
       \le \sup\big\{ \big| u(t,x+\zeta) - u\big(t_n,x_n+\zeta + z( \omega (|t-t_n|) + C|x-x_n|) \big) \big| : \\
                      \text{ for } \zeta, z \text{ s.t. } |\zeta|\le C(1+|x|)R \text{ and } |z|\le R \big\} =:L_n^{\ref{admissibility:aux_comp_4}},
    \end{multlined}
    \end{align*}
    for $n\in\mathbb{N}$, where we used a notational simplification by using a constant $C$ which takes into account the fact that $\sup_{n\in\mathbb{N}}|x_n|\le C'|x|$, for some $C'>0$.
    At this point, the following observations allow us to conclude the convergence $L_n^{\ref{admissibility:aux_comp_4}}\xrightarrow{n\to+\infty}0$.
    Initially,   
    \begin{align*}
      \lim_{n\to +\infty} \sup_{|z|\le R}\Big|\Big(t_n,x_n+\zeta + z\big( \omega (|t-t_n|) &+ C|x-x_n|\big)\Big) - (t,x+\zeta)\Big| = 0.
    \end{align*}
    Hence, for every $\delta>0$ there exists $\hat{n}=\hat{n}(\delta,R)\in\mathbb{N}$ such that for $n\ge \hat{n}$
    \begin{align*}
        \Big| (t,x+\zeta) - \Big(t_n,x_n+\zeta + z\big( \omega (|t-t_n|) &+ C|x-x_n|\big)\Big)\Big| <\delta.
    \end{align*}
    Additionally, there exists a compact $K_2\subset (0,T)\times \mathbb{R}^d$ such that
    \begin{align*}
    \big\{\big(t_n,x_n+\zeta + z( \omega (|t-t_n|) &+ C|x-x_n|)\big) : \text{ for }n,\zeta,z \text{ s.t. } n\in\mathbb{N}, |\zeta|\le C(1+|x|)R \text{ and } |z|\le R  \big\} \subset K_2.
    \end{align*}
    Recalling that $u$ is a continuous function, hence uniformly continuous on $K_2$, allows us to deduce that for every $\varepsilon>0$ there exists $\delta=\delta(\varepsilon) >0$ such that, for $(s_1,y_1)$,$(s_2,y_2)\in K_2$  if $|(s_1,y_1)-(s_2,y_2)|<\delta$, then
         $|u(s_1,y_1) - u(s_2,y_2)|<\varepsilon$.
    Combining all the above, we have on the set $\{z:|z|\le R\}$ and for every $\varepsilon>0$ that there exists $\hat{n}=\hat{n}(\delta(\varepsilon),R)\in\mathbb{N}$ such that for $n\ge \hat{n}$
    \begin{align*}
        L_n^{\ref{admissibility:aux_comp_4}}\overset{\triangle}{=}\sup\big\{
         \big| u(t,x+\zeta) - u\big(t_n,x_n+\zeta + &|z|[ \omega (|t-t_n|) + C|x-x_n|] \big) \big| : \\
                        &\text{ for } \zeta,z \text{ s.t. } |\zeta| \le C(1+|x|)R \text{ and } |z|\le R \big\} <\varepsilon.
    \end{align*}
    In other words,
    \begin{align}\label{upper_bound:limit_4}
        L_n^{\ref{admissibility:aux_comp_4}}\xrightarrow{n\to+\infty} 0.
    \end{align}

    \item\label{admissibility:aux_comp_5}
    By definition of the $\|\cdot\|_p-$norm, it is true that \footnote{The constants appearing below may change in value, but we keep the same symbol for the sake of simplicity.}
    \begin{align*}
      |u(t,x)| 
      +\sup_{n\in\mathbb{N}}|u(t_n,x_n)|
       +\sup_{n\in\mathbb{N}}|u_n(t_n,x_n)|
       \le C(\|u\|_p + \sup_{n\in\mathbb{N}}\|u_n\|_p) (1+\sup_{n\in\mathbb{N}}|x_n|^p)
    \end{align*}
as well as (using \ref{HJB:cond_jump}.\ref{HJB:cond_jump_partB})
    \begin{align*}
      & 
      \sup_{\alpha\in\mathcal{A}}
      \Big(|u(t,x+j_\alpha(t,x,z))| 
      + \sup_{n\in\mathbb{N}}\big( |u(t_n,x_n+j_\alpha(t_n,x_n,z))|\big)
      + \sup_{n\in\mathbb{N}}\big( |u_n(t_n,x_n+j_\alpha(t_n,x_n,z))|\big) \Big)\\
      &\hspace{1em}   \le  C(\|u\|_p + \sup_{n\in\mathbb{N}}\|u_n\|_p) (1+\sup_{n\in\mathbb{N}}|x_n|^p)(1+|z|^p) 
      = L(1+|z|^p),
    \end{align*}
    where
    \begin{align}\label{def:constant_K}
        L:=C(\sup_{n\in\mathbb{N}} \|u_n\|_p + \|u\|_p) (1 + \sup_{n\in\mathbb{N}} |x_n|^p)<+\infty.
    \end{align}
    In view of \ref{HJB:UI}, 
    for every $\varepsilon>0$ there exists $R>1$ such that
    \begin{align}\label{upper_bound_5}
      \sup_{\alpha\in\mathcal{A}} \int_{\{z:|z|>R\}} 1 m_\alpha(\textup{d}z)
      \le \sup_{\alpha\in\mathcal{A}} \int_{\{z:|z|>R\}} |z|^p m_\alpha(\textup{d}z)
      <\frac{\varepsilon}{8L}.
    \end{align}
    \setcounter{Preliminary_counter}{\value{enumi}}
  \end{enumerate}

\vspace{0.5 em}
After these preparatory results, we return to our aim, \emph{i.e.}, proving that 
\begin{align*}
  \lim_{n\to +\infty} \sup_{\alpha\in\mathcal{A}} \big| \hat{\mathcal{I}}_\alpha(t_n,x_n,u_n,\phi_n) - \hat{\mathcal{I}}_\alpha (t,x,u,\phi)\big| =0.
\end{align*}
To this end, let $\varepsilon>0$ and consider the value $R>1$ determined by \eqref{upper_bound_5}.
Then,
\begin{align*}
  &\sup_{\alpha\in\mathcal{A}} \big| \hat{\mathcal{I}}_\alpha(t_n,x_n,u_n,\phi_n) - \hat{\mathcal{I}}_\alpha (t,x,u,\phi)\big|\\
%
%
  &\underset{\text{\ref{HJB:cond_jump}}}{\overset{\ref{admissibility:aux_comp_2}}{\le} }
    \sup_{\alpha\in\mathcal{A}}
    \int_{\{z:1<|z|\le R\}} 
    \big| u(t,x+j_\alpha(t,x,z)) - u(t_n,x_n+j_\alpha(t_n,x_n,z)) \big| + \big| u(t,x) - u(t_n,x_n) \big|m_\alpha(\textup{d}z)\\
    &\hspace{0.8em}   
    +\sup_{\alpha\in\mathcal{A}}
    \int_{\{z:|z|> R\}} 
    \big| u(t,x+j_\alpha(t,x,z)) - u(t_n,x_n+j_\alpha(t_n,x_n,z)) \big| + \big| u(t,x) - u(t_n,x_n) \big|m_\alpha(\textup{d}z)\\
    &\hspace{0.8em}   
      +\sup_{\alpha\in\mathcal{A}}  
      \int_{\{z:1<|z|\le R\}}\big| u(t_n,x_n+j_\alpha(t_n,x_n,z)) - u_n(t_n,x_n+j_\alpha(t_n,x_n,z)) \big| + \big| u(t_n,x_n) - u_n(t_n,x_n) \big| m_\alpha(\textup{d}z)\\
    &\hspace{0.8em}   
      +\sup_{\alpha\in\mathcal{A}}  
      \int_{\{z:|z|>R\}}\big| u(t_n,x_n+j_\alpha(t_n,x_n,z)) - u_n(t_n,x_n+j_\alpha(t_n,x_n,z)) \big| + \big| u(t_n,x_n) - u_n(t_n,x_n) \big| m_\alpha(\textup{d}z)\\
    &\hspace{0.8em}
      +\big| D\phi(t,x) \big|\big( \omega(|t-t_n|) + C|x-x_n|\big) \sup_{\alpha\in\mathcal{A}} \int_{\{z:|z|>1\}}|z| m_\alpha(\textup{d}z)\\
    &\hspace{0.8em}
      +\big| D\phi(t,x) - D\phi_n(t_n,x_n) \big| C(1+\sup_{n\in\mathbb{N}}|x_n|) \sup_{\alpha\in\mathcal{A}} \int_{\{z:|z|>1\}} |z|  m_\alpha(\textup{d}z)    
\end{align*}
and let us fix an $\varepsilon>0$.
We will consider each summand of the right-hand side of the last inequality. 
\begin{enumerate}[label=\textbullet]
  \item For the first summand, we have
  \begin{align*}
    \begin{multlined}[c][0.9\displaywidth]
    \sup_{\alpha\in\mathcal{A}}
    \int_{\{z:1<|z|\le R\}} 
    \big| u(t,x+j_\alpha(t,x,z)) - u(t_n,x_n+j_\alpha(t_n,x_n,z)) \big| + \big| u(t,x) - u(t_n,x_n) \big|m_\alpha(\textup{d}z)\\
    \hspace{0.2em}\overset{\ref{admissibility:aux_comp_4}}{\le} 
            \big(L_n^{\ref{admissibility:aux_comp_4}} + \big| u(t,x) - u(t_n,x_n) | \big) \sup_{\alpha\in\mathcal{A}}\int_{\{z:1<|z|\le R\}} 1 m_\alpha(\textup{d}z)
            \xrightarrow[n\to+\infty]{\text{\eqref{upper_bound:limit_4},\ref{HJB:coeff_boundedness}}}0,
    \end{multlined}
  \end{align*}
where we also used continuity of $u$.

  \item The sum of the second and the fourth summands can finally be bounded by $\varepsilon$, in view of \eqref{upper_bound_5} and the estimations in \ref{admissibility:aux_comp_5}.

  \item For the third summand, we have
  \begin{align*}
      &\sup_{\alpha\in\mathcal{A}}  
      \int_{\{z:1<|z|\le R\}}\big| u(t_n,x_n+j_\alpha(t_n,x_n,z)) - u_n(t_n,x_n+j_\alpha(t_n,x_n,z)) \big| + \big| u(t_n,x_n) - u_n(t_n,x_n) \big| m_\alpha(\textup{d}z)\\
      &\hspace{1em}\overset{\eqref{upper_bound:limit_3}}{\underset{\ref{HJB:coeff_boundedness}}{\le}} 2\| (u-u_n)\mathds{1}_{K_1}\|_{\infty}
      \sup_{\alpha\in\mathcal{A}}  \int_{\{z:1<|z|\le R\}} 1 m_{\alpha}(\textup{d}z)
      \xrightarrow{n\to+\infty}0.
  \end{align*}
  
  \item For the fifth and sixth summand, in view of the assumption that $p\ge 1$, \ref{HJB:coeff_boundedness},\ref{HJB:cond_jump} and the locally uniform convergence of $(D\phi_n)_{n\in\mathbb{N}}$, we can conclude that they both tend to zero.
\end{enumerate}
In total, we have proven that 
\begin{align*}
\lim_{n\to +\infty} \sup_{\alpha\in\mathcal{A}} \big| \hat{\mathcal{I}}_\alpha(t_n,x_n,u_n,\phi_n) - \hat{\mathcal{I}}_\alpha (t,x,u,\phi)\big| <\varepsilon.
\end{align*}
Since $\varepsilon$ was arbitrary, we have proven the validity of our claim.

\vspace{0.5em}
{\tiny $\blacksquare$} We deal now with the non-linearities $(f_{\alpha})_{\alpha\in\mathcal{A}}$.
Let us also provide some auxiliary computations:
\begin{enumerate}
  \setcounter{enumi}{\value{Preliminary_counter}}
  \item\label{admissibility:aux_comp_6} Let $0<\kappa<1$ and $z\in B[0,\kappa]$. From the Mean Value Theorem we have
  \begin{align*}
    \phi(t,x+j_\alpha(t,x,z))-\phi(t,x) = \int_0^1 j_\alpha(t,x,z)^T D\phi(t,x+\xi j_\alpha(t,x,z)) \textup{d}\xi
  \end{align*} 
  and similarly for $\phi_n$ at the point $(t_n,x_n)$, for $n\in\mathbb{N}$.
  Using \eqref{bound:jump} for $R=1$ and the locally uniform convergence of $(D\phi_n)_{n\in\mathbb{N}}$ to $D\phi$ we have 
  \begin{align*}
    &|D\phi_n(t_n,x_n+\xi j_\alpha(t_n,x_n,z)) - D\phi(t_n,x_n+\xi j_\alpha(t_n,x_n,z))|\\
    &\quad \le \sup_{\zeta\le C(1+ \sup_{n\in\mathbb{N}}|x_n|)}|D\phi_n(t_n,x_n+\zeta) - D\phi(t_n,x_n+\zeta)|
     \le \sup_{(s,y)\in K_3}|D\phi_n(s,y) - D\phi(s,y)|\xrightarrow[]{n\to\infty}0
  \end{align*}
  for some compact $K_3$ such that $(t_n,x_n+\zeta) \in K_3$, for every $n\in\mathbb{N}$, $(t,x+\zeta)\in K_3$ and $\zeta\in B[0,\kappa]$.
  \footnote{Comparing to $K_1$ of \ref{admissibility:aux_comp_3} we have $K_3\subset K_1$.}

  \item\label{admissibility:aux_comp_7} 
  Analogously to \ref{admissibility:aux_comp_4}, there exists compact $K_4$\footnote{Comparing to $K_2$ of \ref{admissibility:aux_comp_4} we have $K_4\subset K_2$.}
   such that
  \begin{align*}
    \big\{(t_n,x_n+\zeta + z\big( \omega (|t-t_n|) &+ C|x-x_n|\big) : \text{ for }n,\zeta,z \text{ s.t. } n\in\mathbb{N}, |\zeta|\le C(1+|x|)\kappa \text{ and } |z|\le \kappa  \big\} \subset K_4.
  \end{align*}
  Then, $D\phi$ is uniformly continuous on $K_4$ and consequently (arguing completely analogously to \ref{admissibility:aux_comp_4})
  \begin{align*}
    \begin{multlined}[c][0.95\displaywidth]
      L_{n}^{\ref{admissibility:aux_comp_7}}:=\sup_{\alpha\in\mathcal{A},|z|\le \kappa}\big| D\phi(t,x+j_\alpha(t,x,z)) - D\phi(t_n,x_n+j_\alpha(t_n,x_n,z)) \big| \\
       \le \sup\big\{ \big| D\phi(t,x+\zeta) - D\phi\big(t_n,x_n+\zeta + z( \omega (|t-t_n|) + C|x-x_n|) \big) \big| : \\
                      \text{ for } \zeta, z \text{ s.t. } |\zeta|\le C(1+|x|)\kappa \text{ and } |z|\le \kappa \big\} \xrightarrow[]{n\to\infty}0.
      \end{multlined}
  \end{align*}

  \item\label{admissibility:aux_comp_8}
  Let $0<\kappa< 1$. 
  The quantity $B_{\kappa}$, defined below, is finite:
  \begin{align*}
    B_{\kappa}&:=\sup_{\alpha\in\mathcal{A},n\in\mathbb{N}}\int_{\{z:|z|\le \kappa\}} \delta_\alpha(t_n,x_n,z) |j_\alpha(t_n,x_n,z)| m_\alpha(\textup{d}z) + \int_{\{z:|z|\le \kappa\}} \delta_\alpha(t,x,z) |j_\alpha(t,x,z)| m_\alpha(\textup{d}z)\\
    &\le C (1+|x| + \sup_{n\in\mathbb{N}}|x_n|) \sup_{\alpha\in\mathcal{A}}\Big\{\int_{\{z:|z|\le \kappa\}} \ell^2_\alpha(z) m_\alpha(\textup{d}z) +\int_{\{z:|z|\le \kappa\}} |z|^2 m_\alpha(\textup{d}z)  \Big\}<\infty,
  \end{align*}
  for some $C>0$, and we used \ref{HJB:cond_jump}.\ref{HJB:cond_jump_partB}, \ref{HJB:growth_cond_delta}.\ref{HJB:growth_cond_delta_partA}-\ref{HJB:growth_cond_delta_partB} and \ref{HJB:coeff_boundedness}.

  \item\label{admissibility:aux_comp_9} 
  Let $0<\kappa< 1$. 
  Using \ref{HJB:cond_jump} and \ref{HJB:growth_cond_delta}.\ref{HJB:growth_cond_delta_pre_partC}
  \begin{align*}
    &\sup_{\alpha\in\mathcal{A}}\int_{\{z:|z|\le \kappa\}} |\delta_\alpha(t,x,z) j_\alpha(t,x,z) - \delta_\alpha(t_n,x_n,z) j_\alpha(t_n,x_n,z)| m_\alpha(\textup{d}z) \\
    &\begin{multlined}[c][0.95\displaywidth]
    \le \sup_{\alpha\in\mathcal{A}}\int_{\{z:|z|\le \kappa\}} |\delta_\alpha(t,x,z) j_\alpha(t,x,z) - \delta_\alpha(t_n,x_n,z) j_\alpha(t,x,z)| m_\alpha(\textup{d}z) \\
    + \sup_{\alpha\in\mathcal{A}}\int_{\{z:|z|\le \kappa\}} |\delta_\alpha(t_n,x_n,z) j_\alpha(t,x,z) - \delta_\alpha(t_n,x_n,z) j_\alpha(t_n,x_n,z)| m_\alpha(\textup{d}z) 
    \end{multlined}\\
    &
    \le \omega(|t-t_n|) + C|x-x_n|) 
    (1+|x| + \sup_{n\in\mathbb{N}}|x_n|) C_\kappa \xrightarrow[]{n\to+\infty}0,
  \end{align*}
  where $C_{\kappa}$ is finite bound of $B_\kappa$ extracted in \ref{admissibility:aux_comp_8}.
  \item\label{admissibility:aux_comp_10} 
  Let $0<\kappa<1$.
  Then, using \ref{admissibility:aux_comp_6} - \ref{admissibility:aux_comp_9} we have
  \begin{align*}
    &\Big|\check{\mathcal{K}}_\alpha^\kappa ( t,x,\phi) - \check{\mathcal{K}}_\alpha^\kappa ( t_n,x_n,\phi_n) \Big|\\
    &\le \int_{\{z:|z|\le \kappa\}} \delta_\alpha(t,x,z) |j_\alpha(t,x,z)| \int_0^1 |D\phi(t,x+\xi j_\alpha(t,x,z) - D\phi(t_n,x_n+\xi j_\alpha(t_n,x_n,z))|\textup{d}\xi m_\alpha(\textup{d}z) \\
    &\quad+ \int_{\{z:|z|\le \kappa\}} |\delta_\alpha(t,x,z) j_\alpha(t,x,z) - \delta_\alpha(t_n,x_n,z) j_\alpha(t_n,x_n,z)| \int_0^1 |D\phi(t_n,x_n+\xi j_\alpha(t_n,x_n,z))|\textup{d}\xi m_\alpha(\textup{d}z) \\
    &\quad+ \int_{\{z:|z|\le \kappa\}} \delta_\alpha(t_n,x_n,z) |j_\alpha(t_n,x_n,z)| \int_0^1 |D\phi - D\phi_n|(t_n,x_n+\xi j_\alpha(t_n,x_n,z))\textup{d}\xi m_\alpha(\textup{d}z) \\
    &\le B_{\kappa}C\big(L_{n}^{\ref{admissibility:aux_comp_7}} + \omega(|t-t_n|) + C|x-x_n| + \sup_{(s,y)\in K_3} |D\phi_n(s,y) - D\phi(s,y)| )\big)
    \xrightarrow[]{n\to+\infty} 0,
  \end{align*}
  where $B_{\kappa}$ was defined in \ref{admissibility:aux_comp_8} and $K_3$ was determined in \ref{admissibility:aux_comp_6}.

  \item\label{admissibility:aux_comp_11}
  Let $0<\kappa < 1$ and $p\in[1,2)$.
  We have the following estimate for every $n\in\mathbb{N}$ and $\alpha\in\mathcal{A}$
  \begin{align*}
    &\big|\big[ u_n(t_n,x_n+j_\alpha(t_n,x_n,z)) - u_n(t_n,x_n)\big]\delta_\alpha(t_n,x_n,z) 
     -   \big[ u(t,x+j_\alpha(t,x,z)) - u(t,x)\big]\delta_\alpha(t,x,z)\big|\\
    &\le \big|\big[ u_n(t_n,x_n+j_\alpha(t_n,x_n,z)) - u_n(t_n,x_n)\big] 
     -   \big[ u(t_n,x_n+j_\alpha(t_n,x_n,z)) - u(t_n,x_n)\big]\big|\delta_\alpha(t_n,x_n,z)\\
    &\quad + \big|\big[ u(t_n,x_n+j_\alpha(t_n,x_n,z)) - u(t_n,x_n)\big] 
     -   \big[ u(t,x+j_\alpha(t,x,z)) - u(t,x)\big]\big|\delta_\alpha(t_n,x_n,z)\\
    &\quad + \big|\big[ u(t,x+j_\alpha(t,x,z)) - u(t,x)\big]\delta_\alpha(t_n,x_n,z) 
     -   \big[ u(t,x+j_\alpha(t,x,z)) - u(t,x)\big]\delta_\alpha(t,x,z)\big|.
  \end{align*}
  Using the condition
  \begin{align*}
    M:=\sup_{\alpha\in\mathcal{A}} \Vert \ell_\alpha \Vert_{\infty} <+\infty,
  \end{align*}
  see \ref{HJB:growth_cond_delta}.\ref{HJB:growth_cond_delta_partB}, we have
  \begin{align*}
    &|\overline{\mathcal{K}}_\alpha^\kappa(t_n,x_n,u_n) - \overline{\mathcal{K}}_\alpha^\kappa(t,x,u)|\\
    &\le M \int_{\{z: \kappa < |z|\}} \Big|\big[ u_n(t_n,x_n+j_\alpha(t_n,x_n,z)) - u_n(t_n,x_n)\big] 
     -   \big[ u(t_n,x_n+j_\alpha(t_n,x_n,z)) - u(t_n,x_n)\big]\Big| m_\alpha(\textup{d}z)\\
    &\quad + M  \int_{\{z: \kappa < |z|\}} \Big|\big[ u(t_n,x_n+j_\alpha(t_n,x_n,z)) - u(t_n,x_n)\big] 
     -   \big[ u(t,x+j_\alpha(t,x,z)) - u(t,x)\big]\Big| m_\alpha(\textup{d}z)\\
    &\quad + M\left(\omega(|t-t_n|) + C|x-x_n|\right) \int_{\{z:\kappa < |z|\}} \Big|\big[ u(t,x+j_\alpha(t,x,z)) - u(t,x)\big]\Big| m_\alpha(\textup{d}z). 
  \end{align*}
  The right-hand side tends to $0$ uniformly in $\alpha\in\mathcal{A}$ in view of the information extracted in \ref{admissibility:aux_comp_1} - \ref{admissibility:aux_comp_5}.\footnote{In \eqref{upper_bound_5}, we need to incorporate in the denominator the quantity $M$, which can be done since $M<+\infty$.}
\end{enumerate}
  After these preparatory results, we return to our aim, i.e., proving that 
  \begin{align*}
    \sup_{\alpha\in\mathcal{A}} \big| f_\alpha\big(t_n,x_n,r_n,\sigma_\alpha^T(t_n,x_n) q_n,\mathcal{K}^\kappa_\alpha(t_n,x_n,\phi_n,u_n)\big)  - f_\alpha\big(t,x,r,\sigma_\alpha^T(t,x)q,\mathcal{K}^\kappa_\alpha(t,x,\phi,u)\big)\big|\xrightarrow[]{n\to+\infty}0.
  \end{align*}
  In view of \ref{HJB:coeff_boundedness}, \ref{HJB:coef_continuity}, \ref{HJB:monotonicity}.\ref{HJB:monotonicity_1} and in conjunction with the boundedness of $(x_n)_{n\in\mathbb{N}}$ and $(r_n)_{n\in\mathbb{N}}$, say $R$ and $K$ the respective bounds, we have
  \begin{align*}
    &\sup_{\alpha\in\mathcal{A}} \big| 
    f_\alpha\big(t_n,x_n,r_n,\sigma_\alpha^T(t_n,x_n) q_n,\mathcal{K}^\kappa_\alpha(t_n,x_n,\phi_n,u_n)\big)  
    - f_\alpha\big(t,x,r,\sigma_\alpha^T(t,x)q,\mathcal{K}^\kappa_\alpha(t,x,\phi,u)\big)\big|\\
    &\quad
    \overset{\phantom{\ref{HJB:coeff_boundedness}}}{\le}
    \omega(|t_n-t|) 
    + \overline{\omega}_{\!{}_R}\big(|x_n-x|(1 + |\sigma_\alpha^T(t_n,x_n) q_n|\vee | \sigma_\alpha^T(t,x) q| )\big) 
    + \widetilde{\omega}_{\!{}_K}(|r_n-r|) \\
    &\hspace{2em}+ C\sup_{\alpha\in\mathcal{A}}|\sigma_\alpha(t_n,x_n) - \sigma_\alpha(t,x)|\cdot |q_n|\\
    &\hspace{2em}+ C(\sup_{\alpha\in\mathcal{A}}|\sigma_\alpha(t,x)|\cdot |q_n - q|
    + \sup_{\alpha\in\mathcal{A}}\big|\mathcal{K}^\kappa_\alpha(t_n,x_n,\phi_n,u_n) - \mathcal{K}^\kappa_\alpha(t,x,\phi,u)\big|)\\
    &\quad
    \overset{\ref{HJB:coeff_boundedness}}{\le}
    \omega(|t_n-t|) 
    + \overline{\omega}_{\!{}_R}\big(C_{x,q} |x_n-x| )\big) 
    + \widetilde{\omega}_{\!{}_K}(|r_n-r|)
    + C_q\sup_{\alpha\in\mathcal{A}}|\sigma_\alpha(t_n,x_n) - \sigma_\alpha(t,x)|\\
    &\hspace{2em}+ C_{\sigma,t,x,q} |q_n - q|
    + C\sup_{\alpha\in\mathcal{A}}\big|\mathcal{K}^\kappa_\alpha(t_n,x_n,\phi_n,u_n) - \mathcal{K}^\kappa_\alpha(t,x,\phi,u)\big|.
  \end{align*}
  Recalling the information extracted in \ref{admissibility:aux_comp_10} - \ref{admissibility:aux_comp_11}, each summand on the right-hand side of the inequality converges to $0$ as $n\to +\infty$. 
\end{proof}
\subsection{Auxiliary lemmata for the proof of Proposition \ref{HJB:Regularity_Condition}}\label{subsec_Appendix:aux_lemmata_Regularity_Condition}

In this section we adopt the notation introduced in \cref{HJB:Regularity_Condition} as well as in the main body of its proof, see p. \pageref{proof:HJB:Regularity_Condition}.

\subsubsection{Preparatory remarks and lemmata}\label{subsubsec_Appendix:auxiliary_lemmata_Regularity_Condition}

For the reader's convenience, in this subsubsection we will recall and present some facts, which appeared in the proof of \cref{theorem:Domination_Principle}, which are going to be helpful in the completion of the proof.
In this way, it will be always clearly justified the validity of the arguments used.
It is underlined that we are, indeed, eligible in using them because we work under the hypotheses of \cref{def:regularity_condition}, as we did in the proof of \cref{theorem:Domination_Principle}. 
We underline that the results presented in this subsubsection hold for $p>0$, \emph{i.e.}, not only for $p\geq 1$.

For ease of later reference, we collect some basic facts and properties in the next remark:
\begin{remark}\label{rem_append:aux_properties}
Following the notation introduced in the proof of \cref{HJB:Regularity_Condition}:
  \begin{enumerate}
    \item\label{def:varepsilon_0} Since $u$,$-v$ $\in \textup{USC}_p([0,T]\times \mathbb{R}^d)$, then for $C_{u,v}:= \max\{\|u\|_p, \|v\|_p\}\geq 0$ it is true that
    \begin{align*} 
      |u(t,x)|, |v(t,x)| \leq C_{u,v} (1+|x|^p) \text{ for }(t,x)\in[0,T]\times\mathbb{R}^d.
    \end{align*}  
    Also, the reader may recall that for given $p>0$, we have fixed the quasidistance $\varphi_p$ (denoted simply by $\varphi$) to be the smooth variant described in \cite[Lemma 1.14]{Hollender2017}, whose multiplier is $\rho =2^{2(p\vee 2)-2}$.
    From \cite[Lemma 1.12 (iii)]{Hollender2017}, we may choose $\varepsilon_0=\varepsilon_0(p,u,v)$ such that 
    \begin{align*}
        0<\varepsilon_0 \leq \frac{1}{2 \rho C_{u,v}}
    \end{align*}
    and, among other properties, it is satisfied
    \begin{align}\label{bound:convolution_p_norms} 
      \|u^\varepsilon\|_p, \|v_\varepsilon\|_p \leq \rho C_{u,v}, \text{ for }\varepsilon\leq \varepsilon_0,
    \end{align}   
    where $u^\varepsilon$, resp. $v_\varepsilon$, is the spatial supremal, resp. infimal, convolution of $u$, resp. $v$.

    Hereinafter, whenever we use $\varepsilon_0$, it will be understood the value determined in \cite[Lemma 1.12 (iii)]{Hollender2017} and depends only on $p$,$u$ and $v$. 
    \item\label{regularity_aux_remark:phi_compact} The quasidistance $\varphi$ is continuous and such that $\{z\in\mathbb{R}^d : \varphi(z)\le c\}$ is compact, for every $c>0$; see \cite[Lemma 1.12, Lemma 1.14]{Hollender2017}. 
    Hence, the set $\{x\in\mathbb{R}^d : (\alpha,x,y)\in\mathcal{A}_{\varepsilon}\}$, resp $\{y\in\mathbb{R}^d : (\alpha,x,y)\in\mathcal{A}_{\varepsilon}\}$, is a compact neighborhood of $\bar{x}$, resp. $\bar{y}$, for every $\varepsilon\ge0$; recall that the definition of $\mathcal{A}_{\varepsilon}$ is given in \eqref{def:A_epsilon}.
    \item\label{regularity_aux_remark:limsup_lambda_varepsilon_barx_bary} 
    Let us assume that we have fixed $\delta \in (0,\delta_0]$, for the $\delta_0>0$ described in \hyperref[Domin_principle:Step_3]{Step 3} of the proof of \cref{theorem:Domination_Principle}. 
    Hereinafter, we will always assume the value of $\delta$ fixed and within the aforemenioned interval.
    Hence, we will omit it from the notational dependence.
    
    The family of points $\{\bar{x}, \bar{y}: \lambda>0,0<\varepsilon\le \varepsilon_0\}$ (for the $\varepsilon_0>0$ described in \ref{def:varepsilon_0}), lies in a closed ball whose radius depends on $\gamma$; the reader may recall \eqref{domin_principle:bounded_norms_of_maxima}.
    Combined with the information in \ref{regularity_aux_remark:phi_compact}, we can safely assume that there exists a closed ball of radius $R=R(\gamma)$ which contains $\{x\in\mathbb{R}^d : (\alpha,x,y)\in\mathcal{A}_{\varepsilon}\}$ and $\{y\in\mathbb{R}^d : (\alpha,x,y)\in\mathcal{A}_{\varepsilon}\}$, for every $\lambda>0$ and $0<\varepsilon\le \varepsilon_0$.
    Moreover, it is reminded that
    \begin{align}
       \limsup_{\lambda\to+\infty}\limsup_{\varepsilon \downarrow 0} \lambda | \bar{x}-\bar{y}|^2=0,
     \tag{\ref{domin_principle:liminfima_limsuprema_to_0}}
     \end{align} 
    which implies 
    \begin{align}
       \limsup_{\lambda\to+\infty} \lambda  | \hat{x}-\hat{y}|^2=0,
       \label{limsup_lambda_limit_point}
     \end{align}
    for the arbitrary limit point $(\hat{t},\hat{x},\hat{y})=\big(\hat{t}(\gamma,\lambda),\hat{x}(\gamma,\lambda),\hat{y}(\gamma,\lambda)\big)$  of $(\bar{t}(\gamma,\lambda,\varepsilon),\bar{x}(\gamma,\lambda,\varepsilon),\bar{y}(\gamma,\lambda,\varepsilon))_{\varepsilon\in(0,\varepsilon_0]}$.
    \item\label{regularity_aux_remark:bounded_barx_bary} 
    The continuity of $\varphi$ and the boundedness of $\{\bar{x}, \bar{y}: 0<\varepsilon\le \varepsilon_0\}$, for fixed $\gamma,\lambda>0$ (for the $\varepsilon_0>0$ describe in \ref{def:varepsilon_0}), implies that $\{\delta(\bar{x}), \delta(\bar{y}): 0<\varepsilon\le \varepsilon_0\}$ is also bounded.
    Hence, we can use that 
    \begin{align}
    \label{limit:varepsilon_delta_to_0}
    \lim_{\varepsilon \downarrow 0} \varepsilon(\delta(\bar{x})\vee \delta(\bar{y}))=0;
    \end{align}
    for the definition of $\delta(\cdot)$ recall \eqref{def:delta_x_phi}.
    \item\label{regularity_aux_remark:norm_to_delta}
    For $(\alpha,x,y)\in\mathcal{A}_{\varepsilon}$, $\varepsilon \in (0,\varepsilon_0]$ it is true that 
    \begin{align}
    |x-\bar{x}| 
      \overset{\text{\eqref{quasidist:growth}}}{\leq} 
      2^{(\frac{p}{2}\vee1)-1} \varphi(x-\bar{x})^{\frac12} 
      \overset{\text{\eqref{def:A_epsilon}}}{\leq} 
      2^{(\frac{p}{2}\vee1)-1}
      \varepsilon^{\frac12} 
      \delta(\bar{x})^{\frac12};
      \label{ineq:norm_to_delta}
    \end{align}
    for the definition of $\delta(\cdot)$ recall \eqref{def:delta_x_phi}.
    These further imply 
    \begin{align}
    \label{limit:norm_x_barx_to_0}
    \lim_{\varepsilon \downarrow 0} (x,y) = (\bar{x},\bar{y})\text{ for }(\alpha,x,y)\in\mathcal{A}_{\varepsilon}.
    \end{align}
    Additionally, we have the bound
    \begin{align}
        |x- y| \leq |\bar{x} - \bar{y}| 
        + 2^{(\frac{p}{2}\vee1)-1}
      \varepsilon^{\frac12} 
      (\delta(\bar{x})^{\frac12} + \delta(\bar{y})^{\frac12}),
    \label{bound:norm_x_y}
    \end{align}
    and, consequently,
    \begin{align}
        \limsup_{\lambda \to\infty} \limsup_{\varepsilon\downarrow 0} |x- y| \overset{\eqref{limsup_lambda_limit_point}}{=}0.
    \label{limit_superiora:norm_x_y}
    \end{align}
    Analogously, we can derive the same bounds for $|y-\bar{y}|$ and limits superiora.
Finally,
\begin{align}\label{ineq:Upsilon_x_to_barx}
    \Upsilon(|x|^p)
    \leq 
    C_{\Upsilon,p} \big[\Upsilon(|x - \bar{x}|^p) + \Upsilon(|\bar{x}|^p)\big]
    \overset{\eqref{ineq:norm_to_delta}}{\leq}
    C_{\Upsilon,p} \big[\Upsilon\big(\varepsilon^{\frac p2} \delta(\bar{x}^{\frac p2})\big) + \Upsilon(|\bar{x}|^p)\big],
\end{align}
using the subadditivity and the moderate growth of $\Upsilon$; see 
\cref{corrolary:Young_subadditive_additional_properties}
and 
\cref{lemma:UI_Young_improvement}.\ref{lemma:UI_Young_improvement:moderate_C2} in conjunction with \cref{lem:Young_equiv_moderate}.
    %
    %
    \item\label{regularity_aux_remark:sigma} 
    Let $(\hat{t},\hat{x},\hat{y})$ be a limit point of $(\bar{t}(\gamma,\lambda,\varepsilon),\bar{x}(\gamma,\lambda,\varepsilon),\bar{y}(\gamma,\lambda,\varepsilon))_{\varepsilon\in(0,\varepsilon_0]}$ and denote by $(\varepsilon_n)_{n\in\mathbb{N}}$ the associated sequence of indices.
    The reader may recall from \hyperref[Domin_principle:Step_3]{Step 3} that $\hat{t}\in (0,T_{\gamma,\delta}]$, for some $T_{\gamma,\delta}\in(0,T)$, where $\delta$ has already been assumed fixed.
    We also underline that the limit point and the sequence $(\varepsilon_n)_{n\in\mathbb{N}}$ depends on $\gamma,\lambda$. 
    For notational simplicity we will denote only the dependence on $\varepsilon$.
    Our claim is that 
    \begin{align}
    \sup_{(\alpha,x): (\alpha,x,y) \in \mathcal{A}_{\varepsilon_n}, n\in\mathbb{N}}|\sigma_{\alpha} (\bar{t}(\varepsilon_n),x)|
    \leq C_{\gamma}
     <\infty.
    \label{regularity_aux_remark:sigma_bound_for_x}
    \end{align}
    Indeed, in view of \ref{HJB:coef_continuity} one gets
    \begin{align*}
    \sup_{(\alpha,x): (\alpha,x,y) \in \mathcal{A}_{\varepsilon_n}, n\in\mathbb{N}}|\sigma_{\alpha} (\bar{t}(\varepsilon_n),x)|
    &\overset{\phantom{\text{\eqref{ineq:norm_to_delta}}}}{\leq}
    \sup_{(\alpha,x): (\alpha,x,y) \in \mathcal{A}_{\varepsilon_n}, n\in\mathbb{N}}
    \big\{
    |\sigma_{\alpha} (\bar{t}(\varepsilon_n),\bar{x}(\varepsilon_n))| 
    + C | x - \bar{x}(\varepsilon_n)|
    \big\}
    \\&\overset{\text{\eqref{ineq:norm_to_delta}}}{\leq}
    \sup_{n\in\mathbb{N}}
    \big\{
    |\sigma_{\alpha} (\bar{t}(\varepsilon_n),\bar{x}(\varepsilon_n))| 
    + C 2^{(\frac{p}{2}\vee 1)-1} \varepsilon_n^{\frac12} \delta(\bar{x}(\varepsilon_n))^{\frac12}
    \big\}.
    \end{align*}
    Now, we have in particular that the sequence $(\bar{x}(\varepsilon_n))_{n\in\mathbb{N}}$ lies within $B[0,R_{\gamma}]$, hence $\delta(\bar{x}(\varepsilon_n))$ is also bounded; see \eqref{def:delta_x_phi}.
    Let us consider the sequence of the first terms of the right-hand side of the last inequality:
    \begin{align*}
    \sup_{n\in\mathbb{N}}
    \big|\sigma_{\alpha} \big(\bar{t}(\varepsilon_n),\bar{x}(\varepsilon_n)\big)\big|
    &
    \overset{\ref{HJB:coef_continuity}}{\leq} 
    \omega(|\bar{t}(\varepsilon_n) - \hat{t}|) + C | \bar{x}(\varepsilon_n) - \hat{x}|
    + \sup_{\alpha\in\mathcal{A}}|\sigma_{\alpha} (\hat{t},\hat{x})|
    \overset{\ref{regularity_aux_remark:limsup_lambda_varepsilon_barx_bary}}{\underset{\ref{HJB:coeff_boundedness}}{\leq}} 
     C_{\gamma}<\infty;
    \end{align*}
    we have argued above that $(\hat{t},\hat{x})$ lies in the domain of $\sigma_{\alpha}$ for every $\alpha\in\mathcal{A}$.
    The combination of the above inequalities proves the claim \eqref{regularity_aux_remark:sigma_bound_for_x}.
  \end{enumerate}  
\end{remark}      

In the next lemmata we present the crucial estimates which are going to be used in the next subsubsections.
We present their proofs separately, so that the body of the proof of \cref{HJB:Regularity_Condition} remains as compact as possible. 
For simplicity, we state the lemmata only for the points $x,\bar{x}$.
Obviously, the analogous results hold for the points $y,\bar{y}$.

Given the preparation we have made in \cref{sec:YoungFunctions}, the proofs of the next lemmata boil down to carefully applying the properties of the Young function $\Upsilon$.
\begin{lemma}\label{lemma:FirstOrderTaylor}
Let $p>0$, $\Upsilon$ be the Young function determined by \cref{lemma:UI_Young_improvement} 
and $(j_{\alpha})_{\alpha\in\mathcal{A}}$ be a family of jump-height coefficients such that for every $(t,x)\in[0,T]\times \mathbb{R}^d$ it satisfies
\begin{align*}
    \sup_{\alpha\in\mathcal{A}} |j_\alpha(t,x,z) | \le C(1 + |x|^{\frac{{p\wedge 2}}{2}} ) \,|z|
    \text{ for every $z\in\mathbb{R}^d$ such that $|z|<1$}.
\end{align*} 
Then, under the framework assumed in this subsubsection, there exists a constant $C_{\Upsilon,p} >0$ such that, for every $x,\bar{x}\in\mathbb{R}^d$, $\bar{t}\in[0,T]$ and $z\in B[0,1]$,
\begin{align}\label{ineq:FirstOrderTaylor}
\Upsilon(|\bar{x} + j_{\alpha}(\bar{t},x,z)|^p) - \Upsilon(|\bar{x}|^p)
\leq  C_{\Upsilon,p}
\big[1 + \Upsilon(|\bar{x}|^p) + \Upsilon\big(\varepsilon^{\frac p2} \delta(\bar{x}^{\frac p2})\big)\big] |z|.
\end{align}  
\end{lemma}
\begin{proof}
The application of Mean Value Theorem yields
\begin{align*}
\Upsilon(|\bar{x} + j_{\alpha}(\bar{t},x,z)|^p) - \Upsilon(|\bar{x}|^p)
    &\leq p\int_0^1 (1-\xi)\Upsilon'(|\bar{x} + \xi j_{\alpha}(\bar{t},x,z)|^p) |\bar{x} + \xi j_{\alpha}(\bar{t},x,z)|^{p-1}  |j_{\alpha}(\bar{t},x,z)| \textup{d}\xi\\
    &\leq p\sup_{\xi \in [0,1]} \big\{\Upsilon'(|\bar{x} + \xi j_{\alpha}(\bar{t},x,z)|^p) |\bar{x} + \xi j_{\alpha}(\bar{t},x,z)|^{p-1}\big\}  |j_{\alpha}(\bar{t},x,z)| 
\end{align*}
\begin{enumerate}
    \item If $p\in(0,1)$, then the supremum is bounded by \cref{cor:product_derivatives_negative_powers_bounded}.\ref{cor:product_first_derivative_negative_powers_bounded:p_0_1}, which in conjunction with the assumption on the jump coefficients yields
        \begin{align*}
            \Upsilon(|\bar{x} + j_{\alpha}(\bar{t},x,z)|^p) - \Upsilon(|\bar{x}|^p)
             &\leq C_{\Upsilon,p} (1 + |x|^{\frac{p}{2}})|z|
            \overset{\eqref{ineq:Upsilon_x_to_barx}}{\leq}
            C_{\Upsilon,p}    \big[ 1 + \Upsilon\big(|\bar{x}|^p\big) + \Upsilon\big(\varepsilon^{\frac{p}{2}}\delta(\bar{x})^{\frac{p}{2}}\big) \big]   |z|.
        \end{align*}
   %
   %
   %
    \item If $p\geq 1$, then the function $(0,\infty)\ni w \longmapsto \Upsilon'(w)|w|^{p-1}\in (0,\infty)$ is non-decreasing.
    Hence, 
        \begin{align*}
            &\sup_{\xi \in [0,1]} \big\{\Upsilon'(|\bar{x} + \xi j_{\alpha}(\bar{t},x,z)|^p) |\bar{x} + \xi j_{\alpha}(\bar{t},x,z)|^{p-1}\big\}  |j_{\alpha}(\bar{t},x,z)|\\
            &\hspace{2.6em}\leq \Upsilon'\big((|\bar{x}| + |j_{\alpha}(\bar{t},x,z)|)^p\big) (|\bar{x}| + |j_{\alpha}(\bar{t},x,z)|)^{p-1}  
                (1+ |x|^{\frac{p\wedge 2}{2}})|z|\\
            &\hspace{1.8em}\overset{|z|\leq 1}{\underset{\ref{HJB:cond_jump}.\ref{HJB:cond_jump_partC}}{\leq}} 
                \Upsilon'\big([1 + |\bar{x}| + |x|^{\frac{p\wedge 2}{2}} ]^p\big) 
                [1 + |\bar{x}| + |x|^{\frac{p\wedge 2}{2}}]^{p-1}  
                (1+ |x|^{\frac{p\wedge 2}{2}})|z|\\
            &\hspace{2.6em}\leq \Upsilon'\big([1 + |\bar{x}| + |x|^{\frac{p\wedge 2}{2}} ]^p\big) 
                [1 + |\bar{x}| + |x|^{\frac{p\wedge 2}{2}}]^{p} |z|\\
            &\hspace{2.6em}\leq \widebar{c}_{_\Upsilon}\Upsilon\big([1 + |\bar{x}| + |x|^{\frac{p\wedge 2}{2}} ]^p\big)|z|\\
            &\hspace{2.5em}\leq C_{\Upsilon,p} \Upsilon\big(3^p(1 + |x|^p + |\bar{x}|^p ) \big)\\
            &\hspace{2.5em}\leq C_{\Upsilon,p} \big[1 + \Upsilon( |\bar{x}|^p ) + \Upsilon(|x|^p)\big]\\
            &\hspace{2em}\overset{\eqref{ineq:Upsilon_x_to_barx}}{\leq}
            C_{\Upsilon,p}    \big[ 1 + \Upsilon\big(|\bar{x}|^p\big) + \Upsilon\big(\varepsilon^{\frac{p}{2}}\delta(\bar{x})^{\frac{p}{2}}\big) \big]   |z|
        \end{align*}
    where in the fourth and fifth inequality it was used the moderate growth of $\Upsilon$, see \cref{lemma:UI_Young_improvement} and \cref{lem:Young_equiv_moderate}, and in the last but one the subadditivity of $\Upsilon$, see \cref{corrolary:Young_subadditive_additional_properties}.
\end{enumerate}

\end{proof}
\begin{lemma}\label{lemma:SecondOrderTaylor}
Let $p>0$, $\Upsilon$ be the Young function determined by \cref{lemma:UI_Young_improvement} and $(j_{\alpha})_{\alpha\in\mathcal{A}}$ be a family of jump coefficients such that for every $(t,x)\in[0,T]\times \mathbb{R}^d$ it satisfies
\begin{align*}
    \sup_{\alpha\in\mathcal{A}} |j_\alpha(t,x,z) | \le C(1 + |x|^{\frac{{p\wedge 2}}{2}} ) \,|z|
    \text{ for every $z\in\mathbb{R}^d$ such that $|z|<1$}.
\end{align*} 
  Then, under the framework assumed in this subsubsection, there exists a constant $C_{\Upsilon,p,d} >0$ such that, for every $x,\bar{x}\in\mathbb{R}^d$, $\bar{t}\in[0,T]$ and $z\in B[0,1]$,
\begin{align}\label{ineq:SecondOrderTaylor}
\begin{multlined}[0.9\textwidth]
\Upsilon(|\bar{x} + j_{\alpha}(\bar{t},x,z)|^p) - \Upsilon(|\bar{x}|^p) 
    - p \Upsilon'(|\bar{x}|^p)|\bar{x}|^{p-2}\bar{x}^T j_{\alpha}(\bar{t},x,z)\\
\leq  
C_{\Upsilon,p,d} \big[1 + \Upsilon(|\bar{x}|^p) + \Upsilon\big(\varepsilon^{\frac{p}{2}} \delta(\bar{x})^{\frac{p}{2}}\big) \big]|z|^2.
\end{multlined}
\end{align}  
\end{lemma}
\begin{proof}
The application of Taylor's theorem for $\hat{x}=\hat{x}(x,\bar{x},\xi) := \bar{x} + \xi j_{\alpha}(\bar{t},x,z)$ with $\xi \in [0,1]$  yields
\begin{align*}
&\Upsilon(|\bar{x} + j_{\alpha}(\bar{t},x,z)|^p) - \Upsilon(|\bar{x}|^p) 
    - p \Upsilon'(|\bar{x}|^p)|\bar{x}|^{p-2}\bar{x}^T j_{\alpha}(\bar{t},x,z)\\
    &\begin{multlined}[0.9\textwidth]
    \hspace{0.7em}\overset{\eqref{bound:Young_Second_Derivative}}{\leq} d^2 C_p 
    \Big\{
    \int_0^1 (1-\xi)\Upsilon''(|\bar{x} + \xi j_{\alpha}(\bar{t},x,z)|^p) |\bar{x} + \xi j_{\alpha}(\bar{t},x,z)|^{2p-2} 
        |j_{\alpha}(\bar{t},x,z)|^2\textup{d}\xi\\
    +\int_0^1 (1-\xi)\Upsilon'(|\bar{x} + \xi j_{\alpha}(\bar{t},x,z)|^p) |\bar{x} + \xi j_{\alpha}(\bar{t},x,z)|^{p-2}  
        |j_{\alpha}(\bar{t},x,z)|^2\textup{d}\xi
    \Big\}
    \end{multlined}\\
&\begin{multlined}
\hspace{1em}\leq  C_{p,d}
\Big[
\sup_{\xi\in[0,1]} \{\Upsilon''(|\bar{x} + \xi j_{\alpha}(\bar{t},x,z)|^p) |\bar{x} + \xi j_{\alpha}(\bar{t},x,z)|^{2p-2}\} |j_{\alpha}(\bar{t},x,z)|^2\\
+ \sup_{\xi\in[0,1]} \{\Upsilon'(|\bar{x} + \xi j_{\alpha}(\bar{t},x,z)|^p) |\bar{x} + \xi j_{\alpha}(\bar{t},x,z)|^{p-2}\}  |j_{\alpha}(\bar{t},x,z)|^2
\Big]
\end{multlined}
\end{align*}  
\begin{enumerate}[label=\textbullet]
  \item If $p\in(0,2)$, \emph{i.e.}, $\frac{p\wedge 2}{2} = \frac{p}{2}$, then the suprema in the square brackets of the right-hand side of the last inequality are bounded by a positive constant $C_{\Upsilon,p}$ (see \cref{cor:product_derivatives_negative_powers_bounded}).
  Also,
  \begin{align*}
    |j_\alpha(\bar{t},x,z)|^2
      \overset{\eqref{ineq:norm_to_delta}}{\underset{\ref{HJB:cond_jump}.\ref{HJB:cond_jump_partC}}{\leq}}
            C_p (1 + |\bar{x}|^p + \varepsilon^{\frac{p}{2}} \delta(\bar{x})^{\frac{p}{2}} ) |z|^2
      \leq  C_{\Upsilon,p} \big[1 + \Upsilon(|\bar{x}|^p) + \Upsilon\big(\varepsilon^{\frac{p}{2}} \delta(\bar{x})^{\frac{p}{2}}\big)\big]|z|^2.
  \end{align*} 
Hence, we conclude the desired inequality \eqref{ineq:SecondOrderTaylor}.
  %
  %
  %
  \item If $p\ge 2$, \emph{i.e.}, $\frac{p\wedge 2}{2} = 1$, then let us initially consider each summand appearing in the square brackets of the last inequality.
  For the former one has
  \begin{align*}
  &\sup_{\xi\in[0,1]} \{\Upsilon''(|\bar{x} + \xi j_{\alpha}(\bar{t},x,z)|^p) |\bar{x} + \xi j_{\alpha}(\bar{t},x,z)|^{2p-2}\} | j_{\alpha}(\bar{t},x,z)|^2\\
  &  \hspace{2.5em}\leq C_{\Upsilon,p} \sup_{\xi\in[0,1]} \{|\bar{x} + \xi j_{\alpha}(\bar{t},x,z)|^{p-2}\} | j_{\alpha}(\bar{t},x,z)|^2\\
  &  \hspace{2.5em}\leq C_{\Upsilon,p}  (|\bar{x} |+| j_{\alpha}(\bar{t},x,z)|)^{p-2} (1+|x|)^2 |z|^2\\
  &  \hspace{1.6em}\overset{|z|\leq 1}{\underset{\ref{HJB:cond_jump}.\ref{HJB:cond_jump_partB}}{\leq}} 
        C_{\Upsilon,p}  [1 + |x| + |\bar{x}|]^{p-2} (1+|x|)^2 |z|^2\\
  &  \hspace{2.5em}\leq C_{\Upsilon,p}  [1 + |x| + |\bar{x}|]^{p} |z|^2,
  \end{align*}
  where in the first inequality we used \cref{lemma:UI_Young_improvement}.\ref{lemma:UI_Young_improvement:second_derivative_properties} and the fact that the function $(0,\infty) \ni x \longmapsto x^{p-2}\in (0,\infty)$ is non-decreasing.
  For the latter one has
\begin{align*}
&\sup_{\xi\in[0,1]} \{\Upsilon'(|\bar{x} + \xi j_{\alpha}(\bar{t},x,z)|^p) |\bar{x} + \xi j_{\alpha}(\bar{t},x,z)|^{p-2}\}  |j_{\alpha}(\bar{t},x,z)|^2\\
& \hspace{2.5em}\leq 
\Upsilon'\big((|\bar{x}| + |j_{\alpha}(\bar{t},x,z)|)^p\big) (|\bar{x}| + |j_{\alpha}(\bar{t},x,z)|)^{p-2}  |j_{\alpha}(\bar{t},x,z)|^2\\
& \hspace{1.6em}\overset{|z|\leq 1}{\underset{\ref{HJB:cond_jump}.\ref{HJB:cond_jump_partB}}{\leq}}
\Upsilon'\big([1 + |\bar{x}| + |x|]^p\big) (1 + |\bar{x}| + |x|)^{p-2}  (1+|x|)^2 |z|^2\\
& \hspace{2.5em}\leq 
\widebar{c}_{_\Upsilon} \Upsilon\big([1 + |\bar{x}| + |x|]^p\big) |z|^2\\
& \hspace{2em}\overset{\eqref{ineq:norm_to_delta}}{\underset{\eqref{ineq:Upsilon_x_to_barx}}{\leq}} 
C_{\Upsilon,p} \big[1 + \Upsilon(|\bar{x}|^p) + \Upsilon\big(\varepsilon^{\frac{p}{2}} \delta(\bar{x})^{\frac{p}{2}}\big) \big]|z|^2
\end{align*}
where in the first inequality it was used the fact that $\Upsilon'$ and the function $(0,\infty) \ni x \longmapsto x^{p-2}\in (0,\infty)$ are non-decreasing.
In the second inequality it was used the non-decreasingness of $\Upsilon'$ and the power function.
In the third inequality it was used the moderate growth of $\Upsilon$, see \cref{lemma:UI_Young_improvement} and \cref{lem:Young_equiv_moderate}.
Hence, we conclude the desired inequality \eqref{ineq:SecondOrderTaylor}.
\end{enumerate}
\end{proof}

\begin{lemma}\label{prep_aux_lemma:bound_through_submultiplicativity}
Let $p>0$ and $\Upsilon$ be the Young function determined by \cref{lemma:UI_Young_improvement}.
Then, under the framework assumed in this subsubsection, for $z\in\mathbb{R}^d\backslash \{0\}$, there exists $C_{\Upsilon,p}\geq1$ such that
\begin{align}
\Upsilon(|\bar{x} + j_\alpha(\bar{t},x,z)|^p) 
    \leq C_{\Upsilon,p} \big[1 + \Upsilon(|\bar{x}|^p)+ \Upsilon\big(\varepsilon^{\frac p2} \delta(\bar{x}^{\frac p2})\big)\big]\big[1+\Upsilon(|z|^p)\big].
    \label{ineq:bound_through_submultiplicativity}
\end{align}
\end{lemma}

\begin{proof}
We will initially use the facts that $\Upsilon$ is a moderate Young function (see \cref{lemma:UI_Young_improvement} and \cref{lem:Young_equiv_moderate}) and that for $p>0$ the function $[0,\infty)\ni s\mapsto s^p$ satisfies for $s,t\ge 0$
\begin{align*}
(s+t)^p \leq 2^p (s^p + t^p).
\end{align*}
In view of the above properties and the subadditivity of $\Upsilon$, see \cref{corrolary:Young_subadditive_additional_properties}, we have
\begin{align*}
\Upsilon(|\bar{x} + j_{\alpha}(\bar{t},x,z)|^p)
 &\le 2^{\overline{c}_{{}_\Upsilon}-1}\big[\Upsilon(|\bar{x}|^p) + \Upsilon(|j_{\alpha}(\bar{t},x,z)|^p)\big]\\
 &\le 2^{\overline{c}_{{}_\Upsilon}-1}\big[\Upsilon(|\bar{x}|^p) + C^{\overline{c}_{{}_\Upsilon}}\Upsilon\big((1+|x|^p)|z|^p\big)\big].
 \numberthis\label{ineq:aux_submultiplicative_1}
\end{align*}
where in the second inequality we used the assumption on the jump coefficients, see \ref{HJB:cond_jump}.\ref{HJB:cond_jump_partC} for the case  $p\in(0,2)$ with $|z|\leq 1$ and \ref{HJB:cond_jump}.\ref{HJB:cond_jump_partB} for the rest cases as well as the moderate property of $\Upsilon$, see \cref{lemma:UI_Young_improvement} and \cref{lem:Young_equiv_moderate}.
Next, we will use that $\Upsilon$ behaves finally as a submultiplicative function.
Indeed, let $K_{\Upsilon,p}>0$ be the constant determined by \cref{lemma:UI_Young_improvement}.\ref{lemma:UI_Young_improvement:submultiplicativity}, which may be assumed strictly greater than $1$.
Then, since $1+|x|^p$, if $|z|>1$ 
\begin{align*}
\Upsilon\big((1+|x|^p)|z|^p\big) 
	&\leq K_{\Upsilon,p} \Upsilon(1+|x|^p)\Upsilon(|z|^p) 
    \leq K_{\Upsilon,p} 2^{\overline{c}_{{}_\Upsilon}-1}\max\{\Upsilon(1),1\}\big(1+\Upsilon(|x|^p)\big)\Upsilon(|z|^p)\\
    &\overset{\eqref{ineq:Upsilon_x_to_barx}}{\leq}
	C_{\Upsilon,p} \big[1 + \Upsilon(|\bar{x}|^p) + \Upsilon\big(\varepsilon^{\frac p2} \delta(\bar{x}^{\frac p2})\big) \big]
	\Upsilon(|z|^p),
\end{align*}
with $C_{\Upsilon,p}\geq1$.
If $|z|\leq 1$, then \eqref{ineq:aux_submultiplicative_1} can be written as follows
\begin{align*}
\Upsilon(|\bar{x} + j_{\alpha}(\bar{t},x,z)|^p)
 &\overset{\eqref{ineq:aux_submultiplicative_1}}{\leq} 
 2^{\overline{c}_{{}_\Upsilon}-1}\big[\Upsilon(|\bar{x}|^p) + C^{\overline{c}_{{}_\Upsilon}}\Upsilon\big(1+|x|^p\big)\big]
 \leq 
 C_{\Upsilon,p} \big[1 + \Upsilon(|\bar{x}|^p) + \Upsilon\big(\varepsilon^{\frac p2} \delta(\bar{x}^{\frac p2})\big) \big],
\end{align*}
with $C_{\Upsilon,p}>1$,
where we used the subadditivity of $\Upsilon$, see \cref{corrolary:Young_subadditive_additional_properties}.
In total, we have the desired inequality \eqref{ineq:bound_through_submultiplicativity}.
\end{proof}
\subsubsection{Auxiliary lemmata for the main body of the proof of Proposition \ref{HJB:Regularity_Condition}}\label{subsubsec_Appendix:preparatory_lemmata_Regularity_Condition}
We are now ready to prove the lemmata mentioned in the proof of \cref{HJB:Regularity_Condition}. 

\begin{lemma}\label{lemma_appendix:regularity_f_nonlinearities}
Under the framework of \cref{HJB:Regularity_Condition} and given the introduced notation \eqref{def:penalization_function}-\eqref{def:delta_x_phi}, \eqref{def:A_epsilon} and assumptions \eqref{supremum_penalized}-\eqref{reg_cond:assumption_ineq} of the proof of \cref{HJB:Regularity_Condition}, it holds 
\begin{align*}
    &\begin{multlined}[c][\textwidth]
    \sup_{(\alpha,x,y)\in \mathcal{A}_{\varepsilon}} \Big\{ 
        f_\alpha(\bar{t},x,u^\varepsilon(\bar{t},\bar{x}), \sigma^T_\alpha(\bar{t},x) D_x\phi_{\delta,\gamma}^\lambda (\bar{t},\bar{x},\bar{y}), \mathcal{K}^\kappa_\alpha(\bar{t},x,u^\varepsilon(\bar{t},\cdot)\circ\tau_{\bar{x}-x},\phi_{\delta,\gamma}^\lambda(\bar{t},\cdot,\bar{y})\circ\tau_{\bar{x}-x}))\\
        - f_\alpha(\bar{t},y,v_\varepsilon(\bar{t},\bar{y}),-\sigma^T_\alpha(\bar{t},y) D_y\phi_{\delta,\gamma}^\lambda (\bar{t},\bar{x},\bar{y}), \mathcal{K}^\kappa_\alpha(\bar{t},y,v_\varepsilon(\bar{t},\cdot)\circ\tau_{\bar{y}-y},\phi_{\delta,\gamma}^\lambda(\bar{t},\bar{x},\cdot)\circ\tau_{\bar{y}-y})) 
        \Big\}  
    \end{multlined}\\
    &\hspace{2em}\le 
    \lambda C |\bar{x}-\bar{y}|^2
    +\gamma e^{\mu \bar{t}} C
        \big[1 + \Upsilon(|\bar{x}|^p)+\Upsilon(|\bar{y}|^p)\big]
    +     \varrho^1_{\gamma,\lambda,\varepsilon,\kappa},
\end{align*}
with 
    \begin{align*}
    \limsup_{\lambda\to +\infty}
    \limsup_{\varepsilon \downarrow 0} 
    \limsup_{\kappa \downarrow 0} 
    \varrho^1_{\gamma,\lambda,\varepsilon,\kappa} = 0,
    \end{align*}
for some non-negative constants which may depend on $p, T, \Upsilon$ and the families $\{\sigma_{\alpha}\}_{\alpha\in\mathcal{A}}$, $\{\ell_{\alpha}\}_{\alpha\in\mathcal{A}}$ and $\{m_{\alpha}\}_{\alpha\in\mathcal{A}}$.
\end{lemma}

\begin{proof}
For the closed ball of radius $R=R(\gamma)$ assumed in \cref{rem_append:aux_properties}.\ref{regularity_aux_remark:limsup_lambda_varepsilon_barx_bary} and for $0< \varepsilon\leq \varepsilon_0$
  \begin{align*}
    &\begin{multlined}[c][\textwidth]
    \sup_{(\alpha,x,y)\in \mathcal{A}_{\varepsilon}} \Big\{ 
        f_\alpha(\bar{t},x,u^\varepsilon(\bar{t},\bar{x}), \sigma^T_\alpha(\bar{t},x) D_x\phi_{\delta,\gamma}^\lambda (\bar{t},\bar{x},\bar{y}), \mathcal{K}^\kappa_\alpha(\bar{t},x,u^\varepsilon(\bar{t},\cdot)\circ\tau_{\bar{x}-x},\phi_{\delta,\gamma}^\lambda(\bar{t},\cdot,\bar{y})\circ\tau_{\bar{x}-x}))\\
        - f_\alpha(\bar{t},y,v_\varepsilon(\bar{t},\bar{y}),-\sigma^T_\alpha(\bar{t},y) D_y\phi_{\delta,\gamma}^\lambda (\bar{t},\bar{x},\bar{y}), \mathcal{K}^\kappa_\alpha(\bar{t},y,v_\varepsilon(\bar{t},\cdot)\circ\tau_{\bar{y}-y},-\phi_{\delta,\gamma}^\lambda(\bar{t},\bar{x},\cdot)\circ\tau_{\bar{y}-y})) 
        \Big\}  
    \end{multlined}\\
    &\begin{multlined}[c][\textwidth]
    =
    \sup_{(\alpha,x,y)\in \mathcal{A}_{\varepsilon}} \Big\{ 
        f_\alpha(\bar{t},x,u^\varepsilon(\bar{t},\bar{x}), \sigma^T_\alpha(\bar{t},x) D_x\phi_{\delta,\gamma}^\lambda (\bar{t},\bar{x},\bar{y}), \mathcal{K}^\kappa_\alpha(\bar{t},x,u^\varepsilon(\bar{t},\cdot)\circ\tau_{\bar{x}-x},\phi_{\delta,\gamma}^\lambda(\bar{t},\cdot,\bar{y})\circ\tau_{\bar{x}-x}))\\
        -f_\alpha(\bar{t},x,v_\varepsilon(\bar{t},\bar{y}), \sigma^T_\alpha(\bar{t},x) D_x\phi_{\delta,\gamma}^\lambda (\bar{t},\bar{x},\bar{y}), \mathcal{K}^\kappa_\alpha(\bar{t},x,u^\varepsilon(\bar{t},\cdot)\circ\tau_{\bar{x}-x},\phi_{\delta,\gamma}^\lambda(\bar{t},\cdot,\bar{y})\circ\tau_{\bar{x}-x}))\\
        +f_\alpha(\bar{t},x,v_\varepsilon(\bar{t},\bar{y}), \sigma^T_\alpha(\bar{t},x) D_x\phi_{\delta,\gamma}^\lambda (\bar{t},\bar{x},\bar{y}), \mathcal{K}^\kappa_\alpha(\bar{t},x,u^\varepsilon(\bar{t},\cdot)\circ\tau_{\bar{x}-x},\phi_{\delta,\gamma}^\lambda(\bar{t},\cdot,\bar{y})\circ\tau_{\bar{x}-x}))\\
        -f_\alpha(\bar{t},y,v_\varepsilon(\bar{t},\bar{y}), \sigma^T_\alpha(\bar{t},x) D_x\phi_{\delta,\gamma}^\lambda (\bar{t},\bar{x},\bar{y}), \mathcal{K}^\kappa_\alpha(\bar{t},x,u^\varepsilon(\bar{t},\cdot)\circ\tau_{\bar{x}-x},\phi_{\delta,\gamma}^\lambda(\bar{t},\cdot,\bar{y})\circ\tau_{\bar{x}-x}))\\
        +f_\alpha(\bar{t},y,v_\varepsilon(\bar{t},\bar{y}), \sigma^T_\alpha(\bar{t},x) D_x\phi_{\delta,\gamma}^\lambda (\bar{t},\bar{x},\bar{y}), \mathcal{K}^\kappa_\alpha(\bar{t},x,u^\varepsilon(\bar{t},\cdot)\circ\tau_{\bar{x}-x},\phi_{\delta,\gamma}^\lambda(\bar{t},\cdot,\bar{y})\circ\tau_{\bar{x}-x}))\\
        - f_\alpha(\bar{t},y,v_\varepsilon(\bar{t},\bar{y}),-\sigma^T_\alpha(\bar{t},y) D_y\phi_{\delta,\gamma}^\lambda (\bar{t},\bar{x},\bar{y}), \mathcal{K}^\kappa_\alpha(\bar{t},y,v_\varepsilon(\bar{t},\cdot)\circ\tau_{\bar{y}-y},-\phi_{\delta,\gamma}^\lambda(\bar{t},\bar{x},\cdot)\circ\tau_{\bar{y}-y})) 
        \Big\}
    \end{multlined}\\
    &\begin{multlined}[c][\textwidth]
    \leq
    \sup_{(\alpha,x,y)\in \mathcal{A}_{\varepsilon}}\Big\{ 
      \overline{\omega}_{\!{}_R}\big(|x-y|(1+|\sigma^T_\alpha(\bar{t},x) D_x\phi_{\delta,\gamma}^\lambda (\bar{t},\bar{x},\bar{y}) |)\big)
      \\ 
      + C |\sigma^T_\alpha(\bar{t},x) D_x\phi_{\delta,\gamma}^\lambda (\bar{t},\bar{x},\bar{y}) +\sigma^T_\alpha(\bar{t},y) D_y\phi_{\delta,\gamma}^\lambda (\bar{t},\bar{x},\bar{y}) |
      \\
      +C \Big[ \mathcal{K}^\kappa_\alpha(\bar{t},x,u^\varepsilon(\bar{t},\cdot)\circ\tau_{\bar{x}-x},\phi_{\delta,\gamma}^\lambda(\bar{t},\cdot,\bar{y})\circ\tau_{\bar{x}-x}) 
          -\mathcal{K}^\kappa_\alpha(\bar{t},y,v_\varepsilon(\bar{t},\cdot)\circ\tau_{\bar{y}-y},-\phi_{\delta,\gamma}^\lambda(\bar{t},\bar{x},\cdot)\circ\tau_{\bar{y}-y})\Big]^+
      \Big\}  
    \end{multlined}
    \\
    &
    \begin{multlined}[0.9\textwidth]
    \leq
    \sup_{(\alpha,x,y)\in \mathcal{A}_{\varepsilon}} 
        \overline{\omega}_{\!{}_R}\big(
        |x-y|(1+|\sigma^T_\alpha(\bar{t},x) D_x\phi_{\delta,\gamma}^\lambda (\bar{t},\bar{x},\bar{y}) |)
        \big)
        \\ 
    + 
    C \sup_{(\alpha,x,y)\in \mathcal{A}_{\varepsilon}}\big\{
        |\sigma^T_\alpha(\bar{t},x) D_x\phi_{\delta,\gamma}^\lambda (\bar{t},\bar{x},\bar{y}) +\sigma^T_\alpha(\bar{t},y) D_y\phi_{\delta,\gamma}^\lambda (\bar{t},\bar{x},\bar{y}) |
        \big\}    
        \\
    +C \sup_{(\alpha,x,y)\in \mathcal{A}_{\varepsilon}}
        \Big[ \mathcal{K}^\kappa_\alpha(\bar{t},x,u^\varepsilon(\bar{t},\cdot)\circ\tau_{\bar{x}-x},\phi_{\delta,\gamma}^\lambda(\bar{t},\cdot,\bar{y})\circ\tau_{\bar{x}-x})\\ 
                        -\mathcal{K}^\kappa_\alpha(\bar{t},y,v_\varepsilon(\bar{t},\cdot)\circ\tau_{\bar{y}-y},-\phi_{\delta,\gamma}^\lambda(\bar{t},\bar{x},\cdot)\circ\tau_{\bar{y}-y})\Big]^+.
    \end{multlined}
        \numberthis\label{HJB:reg_cond_nonlinearities_aux_part_a}
  \end{align*}
In the first inequality it was utilized the fact that $u^{\varepsilon}(\bar{t},\bar{x}) - v_{\varepsilon}(\bar{t},\bar{y}) \ge 0$, see \eqref{supremum_penalized}, in conjunction with the monotonicity assumed in the second spatial variable; see \ref{HJB:monotonicity}.\ref{HJB:monotonicity_2}.
  Additionally, it was used the continuity assumption \ref{HJB:monotonicity}.\ref{HJB:monotonicity_1}, while the positive part appears because of the monotonicity of the non-linearities $f_\alpha$ in their last argument; see \ref{HJB:monotonicity}.\ref{HJB:monotonicity_2}.

  Let us, now, separately consider each summand of \eqref{HJB:reg_cond_nonlinearities_aux_part_a}:

\textbullet\hspace{0.3em}
 For the first summand of \eqref{HJB:reg_cond_nonlinearities_aux_part_a} let us initially deal with the quantity within the modulus $\overline{\omega}_{\!{}_R}$. 
To this end, we have
\begin{align*}
    &    \limsup_{\lambda\to +\infty}   \limsup_{\varepsilon \downarrow 0} 
    \big\{\sup_{(x,y,\alpha)\in\mathcal{A}_\varepsilon}
        |x-y|\big(1+|\sigma^T_\alpha(\bar{t},x) D_x\phi_{\delta,\gamma}^\lambda (\bar{t},\bar{x},\bar{y}) |\big)\big\}\\
    &\hspace{2.3em}\overset{\eqref{gradient_penalized_x}}{\underset{\eqref{gradient_penalized_y}}{\le}}
          \limsup_{\lambda\to +\infty}    \limsup_{\varepsilon \downarrow 0} 
    \big\{\sup_{(x,y,\alpha)\in\mathcal{A}_\varepsilon}|x- y|  
      \big(1+ [2\lambda |\bar{x}-\bar{y} | + p\gamma e^{\mu \bar{t}} \Upsilon'(|\bar{x}|^p) |\bar{x} |^{p-1}]
        \times |\sigma_\alpha(\bar{t},x)|\big)\big\}\\
    &\hspace{2em}
    \overset{\eqref{bound:norm_x_y}}{\underset{\eqref{regularity_aux_remark:sigma_bound_for_x}}{\le}}
           \limsup_{\lambda\to +\infty}
    \limsup_{\varepsilon \downarrow 0} \big(|\bar{x}- \bar{y}| + C_p \varepsilon^{\frac12}\big[\delta(\bar{x})^{\frac12} + \delta(\bar{y})^\frac12\big]\big)  
      \big(1+ [2\lambda |\bar{x}-\bar{y} | + p\gamma e^{\mu \bar{t}} \Upsilon'(|\bar{x}|^p) |\bar{x} |^{p-1}]\times C_{\gamma}]\big)\\
    &\hspace{2em}
    \overset{\phantom{\eqref{bound:norm_x_y}}}{\underset{\phantom{\eqref{regularity_aux_remark:sigma_bound_for_x}}}{\le}}
           \limsup_{\lambda\to +\infty}
    |\hat{x}- \hat{y}|  
      (1+ [2\lambda |\hat{x}-\hat{y} | + p\gamma e^{\mu \hat{t}} \Upsilon'(|\hat{x}|^p) |\hat{x} |^{p-1}]\times C_{\gamma}]\\
    &\hspace{2em}
    \overset{\phantom{\eqref{bound:norm_x_y}}}{\underset{\phantom{\eqref{regularity_aux_remark:sigma_bound_for_x}}}{=}}0,
\end{align*}
where $(\hat{t},\hat{x},\hat{y})$ is the limit point of $(\bar{t},\bar{x},\bar{y})_{\varepsilon\leq \varepsilon_0}$ for which the respective limit supremum is attained.
In the third inequality it was also used \eqref{limit:varepsilon_delta_to_0}.
In order to conclude they were used \eqref{limsup_lambda_limit_point} and the fact that for fixed $\gamma$ the limit points $\hat{x}(\gamma,\lambda),\hat{y}(\gamma,\lambda)$ lie in $B[0,R_{\gamma}]$, see \cref{rem_append:aux_properties}.\ref{regularity_aux_remark:limsup_lambda_varepsilon_barx_bary}, in conjunction with \eqref{corollary_Young:Upsilon_prime_upper_bound}.
In other words, $ (\Upsilon'(|\hat{x}|^p) |\hat{x} |^{p-1})_{\lambda \geq 0}$ remains also bounded.
Hence, for 
    \begin{align*}
    \bar{\varrho}^1_{\gamma,\lambda,\varepsilon}:=
    \sup_{(\alpha,x,y)\in \mathcal{A}_{\varepsilon}} 
        \overline{\omega}_{\!{}_R}\big(
        |x-y|(1+|\sigma^T_\alpha(\bar{t},x) D_x\phi_{\delta,\gamma}^\lambda (\bar{t},\bar{x},\bar{y}) |)
        \big)
    \numberthis\label{lem_app:f_nonlinearities_first}
    \end{align*}
    we have, by using the that every modulus of continuity has been assumed to be increasing and continuous,
    \begin{align*}
    \limsup_{\lambda\to +\infty}
    \limsup_{\varepsilon \downarrow 0} 
    \bar{\varrho}^1_{\gamma,\lambda,\varepsilon} = 0.
    \numberthis\label{limsup_lambda_epsilon_bar_rho_1}
    \end{align*}
%
%
%
%
%
%
%
%
%
%
%
%
%
%
%
%
%
%

\textbullet\hspace{0.3em}
    For the second summand of \eqref{HJB:reg_cond_nonlinearities_aux_part_a}, using \eqref{gradient_penalized_x}, \eqref{gradient_penalized_y} and \ref{HJB:coef_continuity},  we have
    \begin{align*}
      &\sup_{(\alpha,x,y)\in \mathcal{A}_{\varepsilon}}
      \big\{
        |\sigma^T_\alpha(\bar{t},x) D_x\phi_{\delta,\gamma}^\lambda (\bar{t},\bar{x},\bar{y}) +\sigma^T_\alpha(\bar{t},y) D_y\phi_{\delta,\gamma}^\lambda (\bar{t},\bar{x},\bar{y}) | 
        \big\}\\
      &\hspace{1em} \overset{\phantom{\text{\eqref{ineq:norm_to_delta}}}}{\leq}
        \sup_{(\alpha,x,y)\in \mathcal{A}_{\varepsilon}}
        \big\{
          2\lambda C|\bar{x}-\bar{y}|\times |x-y|
          +p\gamma e^{\mu \bar{t}}\big[\Upsilon'(|\bar{x}|^p) |\bar{x}|^{p-2} |\sigma^T_\alpha(\bar{t},x) \bar{x}| 
          +\Upsilon'(|\bar{y}|^p) |\bar{y}|^{p-2} |\sigma^T_\alpha(\bar{t},y) \bar{y}| \big] 
        \big\}\\
      &\hspace{1em}
      \begin{multlined}[0.88\textwidth]
        \overset{\eqref{bound:norm_x_y}}{\underset{\ref{HJB:coef_continuity}}{\leq}}
        2\lambda C |\bar{x}-\bar{y}|^2
        +
        2^{(\frac{p}{2}\vee 1)} 
        \varepsilon^{\frac 12}
        \lambda C |\bar{x}-\bar{y}|
        (
        \delta(\bar{x})^{\frac{1}{2}} + \delta(\bar{y})^{\frac{1}{2}}
        )
        \\
        +
        p\gamma e^{\mu \bar{t}}
        \big\{
          \big[\Upsilon'(|\bar{x}|^p) |\bar{x}|^{p-1} 
          +\Upsilon'(|\bar{y}|^p) |\bar{y}|^{p-1}  \big] 
        \big\}
        \sup_{\alpha\in\mathcal{A}} |\sigma_\alpha(\bar{t},0)|
        \\
        +
        2^{(\frac {p}{2}\vee 1)-1}
        \varepsilon^{\frac12}
        p\gamma e^{\mu \bar{t}}
          \big[\Upsilon'(|\bar{x}|^p) |\bar{x}|^{p-1} \delta(\bar{x}) 
          +\Upsilon'(|\bar{y}|^p) |\bar{y}|^{p-1} \delta(\bar{y}) \big] 
        \\
        +
        p\gamma e^{\mu \bar{t}}
        \big\{
          \big[\Upsilon'(|\bar{x}|^p) |\bar{x}|^{p} 
          +\Upsilon'(|\bar{y}|^p) |\bar{y}|^{p} \big] 
        \big\}
      \end{multlined}\\
      &\hspace{1em}
      \begin{multlined}[0.88\textwidth]
        \overset{\eqref{corollary_Young:Upsilon_prime_upper_bound}}{\underset{\ref{HJB:coeff_boundedness}}{\leq}}
        2\lambda C |\bar{x}-\bar{y}|^2
        +
        \gamma e^{\mu \bar{t}}
        C_{p,T,\Upsilon,\sigma}
          \big[1 + \Upsilon(|\bar{x}|^p) +\Upsilon'(|\bar{y}|^p) \big]
        \\
        +
        \varepsilon^{\frac12}
        2^{(\frac {p}{2}\vee 1)-1} 
        3 C_{\Upsilon}
        p\gamma e^{\mu \bar{t}}
          \big[1 + \Upsilon(|\bar{x}|^p)  \delta(\bar{x}) 
          +\Upsilon(|\bar{y}|^p) \delta(\bar{y}) \big] 
        \\
        +
        2^{(\frac{p}{2}\vee 1)} 
        \varepsilon^{\frac 12}
        \lambda C |\bar{x}-\bar{y}|
        (
        \delta(\bar{x})^{\frac{1}{2}} + \delta(\bar{y})^{\frac{1}{2}}
        )
      \end{multlined}      
      \\
      &\hspace{1em} 
      \overset{\phantom{\text{\eqref{ineq:norm_to_delta}}}}{=}
       \lambda C |\bar{x}-\bar{y}|^2
        +
        \gamma e^{\mu \bar{t}}
        C_{p,T, \Upsilon,\sigma}
          \big\{1 + \Upsilon(|\bar{x}|^p)+\Upsilon(|\bar{y}|^p) \big\}
          +\bar{\varrho}^2_{\gamma,\lambda,\varepsilon},
      \numberthis
      \label{lem_app:f_nonlinearities_second}
    \end{align*}
    where it was defined
    \begin{multline*}
    \bar{\varrho}^2_{\gamma,\lambda,\varepsilon}:=
            \varepsilon^{\frac12}
        \big\{2^{(\frac {p}{2}\vee 1)-1} 
                3 \widebar{c}_{{}_\Upsilon}
                p\gamma e^{\mu \bar{t}}
                  \big[1 + \Upsilon(|\bar{x}|^p)  \delta(\bar{x}) 
                  +\Upsilon(|\bar{y}|^p) \delta(\bar{y}) \big] 
                +
                2^{(\frac{p}{2}\vee 1)} 
                \lambda C |\bar{x}-\bar{y}|
                (
                \delta(\bar{x})^{\frac{1}{2}} + \delta(\bar{y})^{\frac{1}{2}}
                )\big\}.
    \end{multline*}
In the course of the above computations, in order to derive the second inequality it was used \eqref{ineq:norm_to_delta} twice. 
Also, in the third inequality the constant $C_{p,T,\Upsilon,\sigma}$ depends on $p,T,\Upsilon$ and on the family $\sup_{\alpha\in\mathcal{A}} |\sigma_\alpha(\bar{t},0)|$, which is finite;  
    \begin{align*}
    \sup_{\alpha\in\mathcal{A}} |\sigma_\alpha(\bar{t},0)| 
    \leq \sup_{\alpha\in\mathcal{A}} |\sigma_\alpha(\frac{T}{2},0)| + \omega(|\bar{t}-\frac{T}{2}|)
    \leq \sup_{\alpha\in\mathcal{A}} |\sigma_\alpha(\frac{T}{2},0)| + \omega(\frac{T}{2})\overset{\ref{HJB:coeff_boundedness}}{<}\infty.
     \end{align*} 
Finally, in the last inequality  it was used the fact that $\Upsilon$ is a moderate Young function, see \cref{lemma:UI_Young_improvement}.\ref{lemma:UI_Young_improvement:moderate_C2}.

Now, in view of the properties described in \cref{rem_append:aux_properties}.\ref{regularity_aux_remark:limsup_lambda_varepsilon_barx_bary}, which essentially guarantee the boundedness of the curly bracket of $\bar{\varrho}^2_{\gamma,\lambda,\varepsilon}$  on a ball of radius $R_{\gamma}$, 
and the continuity of $\Upsilon$, we have that
    \begin{align}
    \limsup_{\varepsilon \downarrow 0} 
    \bar{\varrho}^2_{\gamma,\lambda,\varepsilon} = 0.
    \label{limsup_lambda_epsilon_bar_rho_2}
    \end{align}

\textbullet\hspace{0.3em}
    For the third summand of \eqref{HJB:reg_cond_nonlinearities_aux_part_a} we have
    \begin{align*}
      &\sup_{(\alpha,x,y) \in \mathcal{A}_{\varepsilon}} 
      \Big[\mathcal{K}^\kappa_\alpha(\bar{t},x,u^\varepsilon(\bar{t},\cdot)\circ\tau_{\bar{x}-x},\phi_{\delta,\gamma}^\lambda(\bar{t},\cdot,\bar{y})\circ\tau_{\bar{x}-x}) 
                  -\mathcal{K}^\kappa_\alpha(\bar{t},y,v_\varepsilon(\bar{t},\cdot)\circ\tau_{\bar{y}-y},-\phi_{\delta,\gamma}^\lambda(\bar{t},\bar{x},\cdot)\circ\tau_{\bar{y}-y})\Big]^+\\
      &\hspace{2em}
      \begin{multlined}[c][0.9\textwidth]
      \le 
      \bar{\varrho}^3_{\gamma,\lambda,\varepsilon,\kappa}
      +\sup_{(\alpha,x,y) \in \mathcal{A}_{\varepsilon}}\Big(
      \overline{\mathcal{K}}^\kappa_\alpha(\bar{t},x,u^\varepsilon(\bar{t},\cdot)\circ\tau_{\bar{x}-x}) 
            -\overline{\mathcal{K}}^\kappa_\alpha(\bar{t},y,v_\varepsilon(\bar{t},\cdot)\circ\tau_{\bar{y}-y})\Big)^+.
      \end{multlined}
      \numberthis\label{lem_app:f_nonlinearities_third_1}
    \end{align*}
    where we defined
    \begin{gather*}
      \bar{\varrho}^3_{\gamma,\lambda,\varepsilon,\kappa}:=
       \sup_{(\alpha,x,y) \in \mathcal{A}_{\varepsilon}}\big| 
      \check{\mathcal{K}}^\kappa_\alpha(\bar{t},x,\phi_{\delta,\gamma}^\lambda(\bar{t},\cdot,\bar{y})\circ\tau_{\bar{x}-x}) \big|
      + \sup_{(\alpha,x,y) \in \mathcal{A}_{\varepsilon}}\big| 
            \check{\mathcal{K}}^\kappa_\alpha(\bar{t},y,-\phi_{\delta,\gamma}^\lambda(\bar{t},\bar{x},\cdot)\circ\tau_{\bar{y}-y})\big|
    \end{gather*}

    We consider initially the term $\bar{\varrho}^3_{\gamma,\lambda,\varepsilon,\kappa}$, \emph{i.e.}, we we will work on the set $\{z: |z|\le \kappa\}$ with $\kappa\le 1$. 
    Given the Mean Value Theorem, we have that 
    \begin{align*}
      &\limsup_{\kappa \downarrow 0} \sup_{(\alpha,x,y)\in\mathcal{A}_{\varepsilon}}
        \big|\check{\mathcal{K}}^\kappa_\alpha(\bar{t},x,\phi_{\delta,\gamma}^\lambda(\bar{t},\cdot,\bar{y})\circ\tau_{\bar{x}-x})\big|\\
      &
      \overset{\phantom{\ref{HJB:cond_jump}.\ref{HJB:cond_jump_partB}}}{\underset{\phantom{\ref{HJB:growth_cond_delta}.\ref{HJB:growth_cond_delta_partA}}}{\leq}}
      \limsup_{\kappa \downarrow 0}\sup_{(\alpha,x,y)\in\mathcal{A}_{\varepsilon}}
        \int_{\{z: |z|\le \kappa\}} 
        \big|\phi_{\delta,\gamma}^\lambda(\bar{t}, \bar{x} + j_\alpha(\bar{t},x,z),\bar{y}) - \phi_{\delta,\gamma}^\lambda(\bar{t}, \bar{x} ,\bar{y})\big|\,\delta_{\alpha}(\bar{t},x,z)\, m_\alpha(\textup{d}z)\\
      &
      \overset{\phantom{\ref{HJB:cond_jump}}}{\underset{\phantom{\ref{HJB:growth_cond_delta}.\ref{HJB:growth_cond_delta_partA}}}{=}}
      \limsup_{\kappa \downarrow 0}\sup_{(\alpha,x,y)\in\mathcal{A}_{\varepsilon}}
      \int_{\{z: |z|\le \kappa\}} \Big| \Big(\int_0^1 D_x\phi_{\delta,\gamma}^\lambda(\bar{t}, \bar{x} + \xi j_\alpha(\bar{t},x,z),\bar{y})\textup{d}\xi\Big)^Tj_\alpha(\bar{t},x,z)\Big| \cdot  \delta_{\alpha}(\bar{t},x,z) m_\alpha(\textup{d}z)\\
      &\begin{multlined}[c][0.95\textwidth]
      \overset{\ref{HJB:cond_jump}}{\underset{\ref{HJB:growth_cond_delta}.\ref{HJB:growth_cond_delta_partA}}{\leq }}
      C\limsup_{\kappa \downarrow 0}
      \sup_{(\alpha,x,y)\in\mathcal{A}_{\varepsilon}} 
      \Big\{\sup_{|\zeta|\le 1, \xi\in[0,1]} 
        \big|D_x\phi_{\delta,\gamma}^\lambda\big(\bar{t}, \bar{x} + C\xi (1+|x|)\zeta,\bar{y}\big)\big| C(1+|x|^{\frac{p\wedge2}{2}}) 
      \Big\} \\
            \times \Big(\sup_{\alpha\in\mathcal{A}}\int_{\{z: |z|\le \kappa\}} \ell_\alpha^2(z) m_\alpha(\textup{d}z)\Big)
            \times \Big(\sup_{\alpha\in\mathcal{A}}\int_{\{z: |z|\le \kappa\}} |z|^2 m_\alpha(\textup{d}z)\Big)
      \end{multlined}\\
      &\hspace{0.2em}
      \overset{\ref{HJB:UI}}{\underset{\ref{HJB:growth_cond_delta}.\ref{HJB:growth_cond_delta_partA}}{=}}
      0
    \end{align*}
    where we used Cauchy--Schwarz inequality and the boundedness of a continuous function on a compact set which contains the bounded\footnote{Use \cref{rem_append:aux_properties}.\ref{regularity_aux_remark:limsup_lambda_varepsilon_barx_bary}.} set
    \begin{align*}
      \{(\bar{t},\bar{x}+C\xi (1+|x|)\zeta,\bar{y}) : \xi \in[0,1], |\zeta|\le 1, x \text{ is s.t. } \varphi(x-\bar{x})\le \varepsilon \delta(\bar{x}), \bar{x},\bar{y} \text{ for }\lambda >0 \text{ and }\varepsilon\le\varepsilon_0 \}
    \end{align*} 
    and the fact that we work on the set $\{z: |z|\le \kappa\}$ with $\kappa\le 1$.
    We have the analogous property for the term $\big|\check{\mathcal{K}}^\kappa_\alpha(\bar{t},y,\phi_{\delta,\gamma}^\lambda(\bar{t},\bar{x},\cdot)\circ\tau_{\bar{y}-y})\big|$.
    Thus, combining the two pieces of information just extracted in the above computations we have
    \begin{align}\label{limsup_lambda_epsilon_bar_rho_3}
      \limsup_{\kappa \downarrow 0}
      \bar{\varrho}^3_{\gamma,\lambda,\varepsilon,\kappa}
      =0
    \end{align}
    We proceed with the second summand of the right-hand side of Inequality \eqref{lem_app:f_nonlinearities_third_1}, which means that we will work on the set $\{ z : \kappa <|z|\}$.
    Because of the (possible) singularity of the family of measure $\{m_\alpha\}_{\alpha\in\mathcal{A}}$ at the origin, we consider separately the cases $\{z: \kappa <|z|\le 1\}$ and $\{z: 1<|z|\}$.

    Let us deal with the integrands of the difference 
    \begin{align*}
    \overline{\mathcal{K}}^\kappa_\alpha(\bar{t},x,u^\varepsilon(\bar{t},\cdot)\circ\tau_{\bar{x}-x}) 
            -\overline{\mathcal{K}}^\kappa_\alpha(\bar{t},y,v_\varepsilon(\bar{t},\cdot)\circ\tau_{\bar{y}-y}),
    \end{align*}
    working on the set $\{z: \kappa <|z|\le 1\}$.
    To this end, one has for $(\alpha, x,y)\in \mathcal{A}_\varepsilon$:
    \begin{align*}
      &\big[u^\varepsilon(\bar{t},\bar{x}+j_\alpha(\bar{t},x,z)) - u^\varepsilon(\bar{t},\bar{x})\big] \delta_\alpha(\bar{t},x,z) 
      -
      \big[v_\varepsilon(\bar{t},\bar{y}+j_\alpha(\bar{t},y,z)) - v_\varepsilon(\bar{t},\bar{y})\big] \delta_\alpha(\bar{t},y,z)\\
      & \begin{multlined}[c][0.9\textwidth]
      \overset{\phantom{\eqref{supremum_penalized}}}{\underset{\phantom{\ref{HJB:growth_cond_delta}.\ref{HJB:growth_cond_delta_partC}}}{=}} 
      \big\{\big[u^\varepsilon(\bar{t},\bar{x}+j_\alpha(\bar{t},x,z)) - u^\varepsilon(\bar{t},\bar{x})\big]   
                - \big[v_\varepsilon(\bar{t},\bar{y}+j_\alpha(\bar{t},y,z)) - v_\varepsilon(\bar{t},\bar{y})\big]
            \big\}\delta_\alpha(\bar{t},x,z)\\
      +\big[v_\varepsilon(\bar{t},\bar{y}+j_\alpha(\bar{t},y,z)) - v_\varepsilon(\bar{t},\bar{y})\big] 
      \big[\delta_\alpha(\bar{t},x,z) - \delta_\alpha(\bar{t},y,z)\big]
      \end{multlined}
      \\
      &\begin{multlined}[c][0.9\textwidth]
      \overset{\eqref{supremum_penalized}}{\underset{{\text{\ref{HJB:growth_cond_delta}.\ref{HJB:growth_cond_delta_partA}}}}{\le}} 
      \big[\phi_{\delta,\gamma}^\lambda (\bar{t},\bar{x}+j_\alpha(\bar{t},x,z),\bar{y}+j_\alpha(\bar{t},y,z)) - \phi_{\delta,\gamma}^\lambda (\bar{t},\bar{x},\bar{y})\big]\ell_{\alpha}(z)\\
      + \big[v_\varepsilon(\bar{t},\bar{y}+j_\alpha(\bar{t},y,z)) - v_\varepsilon(\bar{t},\bar{y})\big] 
      \big[\delta_\alpha(\bar{t},x,z) - \delta_\alpha(\bar{t},y,z)\big],
      \end{multlined}
      \numberthis\label{case_v_locally_Lipschitz}
      \\
      &\begin{multlined}[c][0.9\textwidth]
      \overset{\eqref{def:penalization_function}}{\underset{\text{\ref{HJB:growth_cond_delta}.\ref{HJB:growth_cond_delta_partC}}}{\le}} 
      \lambda\big( |\bar{x} + j_\alpha(\bar{t}, x,z) - (\bar{y} + j_\alpha(\bar{t}, y,z) )|^2 - |\bar{x} -\bar{y}|^2 \big) \ell_{\alpha}(z)\\
        +\gamma e^{\mu \bar{t}}
        \Big[  
        \Upsilon\big( |\bar{x}+j_\alpha(\bar{t},x,z)|^p\big) - \Upsilon\big( |\bar{x}|^p\big) 
        + \Upsilon\big( |\bar{y}+j_\alpha(\bar{t},y,z)|^p\big) - \Upsilon\big( |\bar{y}|^p\big)
        \Big] \ell_{\alpha}(z)\\
        +\sup_{\varepsilon\in(0,\varepsilon_0]}\Vert v_\varepsilon\Vert_p 2^p
        (1 + |\bar{y}|^p + |j_{\alpha}(\bar{t},y,z) |^p )\times  C |x-y|  \ell_{\alpha}^2(z)
      \end{multlined}
      \\
    &\begin{multlined}[c][0.9\textwidth]
    \underset{\phantom{\text{\ref{HJB:growth_cond_delta}.\ref{HJB:growth_cond_delta_partC}}}}{\overset{\text{\eqref{ineq:FirstOrderTaylor}}}{\le}} 
        \lambda\big[
        |j_\alpha(\bar{t}, x,z) - j_\alpha(\bar{t}, y,z)|^2 - 2(\bar{x}-\bar{y})^T(j_\alpha(\bar{t},x,z) - j_\alpha(\bar{t},y,z)) 
      \big] \ell_{\alpha}(z)\\
    +\gamma e^{\mu \bar{t}}
    C_{\Upsilon,p}
    \Big[  
    1 + \Upsilon(|\bar{x}|^p) + \Upsilon(|\bar{y}|^p) 
        + \Upsilon\big(\varepsilon^{\frac p2} \delta(\bar{x}^{\frac p2})\big) 
        + \Upsilon\big(\varepsilon^{\frac p2} \delta(\bar{y}^{\frac p2})\big)
    \Big] 
|z|\ell_{\alpha}(z)\\
    +\Vert v_{\varepsilon_0}\Vert_p 2^p
    (1 + |\bar{y}|^p + |y|^p )\times  C |x-y|  \ell_{\alpha}^2(z)
    \end{multlined}
  \\
    &\begin{multlined}[c][0.9\textwidth]
      \overset{\ref{HJB:cond_jump}.\ref{HJB:cond_jump_partA}}{\underset{\text{\eqref{bound:norm_x_y}}}{\le}} 
        \lambda C |\bar{x}-\bar{y}|^2 \big[|z|^2\ell_{\alpha}(z) + 2|z|\ell_{\alpha}(z)\big] \\
        +\varepsilon \lambda C_p [\delta(\bar{x})  +\delta (\bar{y})] |z|^2\ell_{\alpha}(z)
        + \varepsilon^{\frac12}\lambda  C_p |\bar{x}-\bar{y}| (\delta(\bar{x})^{\frac12} + \delta(\bar{y})^\frac12) 
      |z|  \ell_{\alpha}(z)\\
    +\gamma e^{\mu \bar{t}}
    C_{\Upsilon,p}
    \Big[  
    1 + \Upsilon(|\bar{x}|^p) + \Upsilon(|\bar{y}|^p) 
        + \Upsilon\big(\varepsilon^{\frac p2} \delta(\bar{x}^{\frac p2})\big) 
        + \Upsilon\big(\varepsilon^{\frac p2} \delta(\bar{y}^{\frac p2})\big)
    \Big] 
|z|\ell_{\alpha}(z)\\
    +C_{v,p,\varepsilon_0}
    \big[1 + |\bar{y}|^p + \varepsilon^{\frac{p}{2}} \delta(\bar{y})^\frac{p}{2} \big]\times
    \big[|\bar{x}-\bar{y}| + \varepsilon^{\frac{p}{2}} \big(\delta(\bar{x})^{\frac{p}{2}} + \delta(\bar{y})^{\frac{p}{2}}\big] 
     \ell_{\alpha}^2(z)
    \end{multlined}
  \\
    &\begin{multlined}[c][0.9\textwidth]
      \overset{\eqref{bound:convolution_p_norms}}{\underset{\phantom{\ref{HJB:growth_cond_delta}.\ref{HJB:growth_cond_delta_partB}}}{=}}
      \lambda C |\bar{x}-\bar{y}|^2 \big[|z|^2\ell_{\alpha}(z) + |z|\ell_{\alpha}(z)\big]
      +\gamma e^{\mu \bar{t}} C_{\Upsilon,p} \big[1 + \Upsilon(|\bar{x}|^p) + \Upsilon(|\bar{y}|^p)\big] |z| \ell_{\alpha}(z)
      +\bar{\theta}^1_{\alpha,\gamma,\lambda,\varepsilon}(z)
    \end{multlined}
    \numberthis\label{lem_app:f_nonlinearities_third_2}
    \end{align*}
    where 
    \begin{align}
    \begin{multlined}[0.9\textwidth]
    \bar{\theta}^1_{\alpha,\gamma,\lambda,\varepsilon}(z):=
    \varepsilon \lambda C_p [\delta(\bar{x})  +\delta (\bar{y})] |z|^2\ell_{\alpha}(z)
        +\varepsilon^{\frac12} \lambda C_{p}|\bar{x}-\bar{y}| (\delta(\bar{x})^{\frac12} + \delta(\bar{y})^\frac12)  |z| \ell_{\alpha}(z)\\
        +\gamma e^{\mu \bar{t}}
        C_{\Upsilon,p}\big[
        \Upsilon\big(\varepsilon^{\frac p2} \delta(\bar{x}^{\frac p2})\big) 
        + \Upsilon\big(\varepsilon^{\frac p2} \delta(\bar{y}^{\frac p2})\big)
                \big]
                |z|\ell_{\alpha}(z)\\
        +C_{v,p,\varepsilon_0}(1 + |\bar{y}|^p + \varepsilon^{\frac{p}{2}} \delta(\bar{y})^{\frac{p}{2}})
        \big[
        |\bar{x}-\bar{y}| 
        +2^{(\frac{p}{2}\vee 1)}\varepsilon^{\frac12}\big(\delta(\bar{x})^{\frac12} + \delta(\bar{y})^{\frac12}\big) 
        \big] 
        \times  \ell_{\alpha}^2(z)
    \end{multlined}
    \label{def:bar_theta_1}
    \end{align}
    In the course of computations we used the fact that the (spatial) supremal concloutions converge pointwise and decreasingly in order to be eligilble to write 
    $\sup_{\varepsilon\in(0,\varepsilon_0]} \| u_\varepsilon\|_p = \| u_{\varepsilon_0}\|_p.   $
%
%
%
%
%
%
%
%
%
%
%
%
%
%
%
%
%
%
%
%
%
%
%
%
%
%
%
%
%
%
%

Next, we consider the case $\{z: |z|>1 \}$ and we initially manipulate the respective integrands just like the case $\{z:\kappa<|z|\leq 1\}$. 
Then we proceed as follows:
    \begin{align*}
      &\big[u^\varepsilon(\bar{t},\bar{x}+j_\alpha(\bar{t},x,z)) - u^\varepsilon(\bar{t},\bar{x})\big] \delta_\alpha(\bar{t},x,z) 
      -
      \big[v_\varepsilon(\bar{t},\bar{y}+j_\alpha(\bar{t},y,z)) - v_\varepsilon(\bar{t},\bar{y})\big] \delta_\alpha(\bar{t},y,z)\\
    &\begin{multlined}[c][0.9\textwidth]
      \overset{\eqref{case_v_locally_Lipschitz}}{\underset{\phantom{\ref{HJB:growth_cond_delta}.\ref{HJB:growth_cond_delta_partC}}}{\le}} 
      \lambda\big( |\bar{x} + j_\alpha(\bar{t}, x,z) - (\bar{y} + j_\alpha(\bar{t}, y,z) )|^2 - |\bar{x} -\bar{y}|^2 \big) \ell(z)\\
        +\gamma e^{\mu \bar{t}}
        \Big[  
        \Upsilon\big( |\bar{x}+j_\alpha(\bar{t},x,z)|^p\big) - \Upsilon\big( |\bar{x}|^p\big) 
        + \Upsilon\big( |\bar{y}+j_\alpha(\bar{t},y,z)|^p\big) - \Upsilon\big( |\bar{y}|^p\big)
        \Big] \ell(z)\\
        +C_{v,p,\varepsilon_0} (1 + |\bar{y}|^p + |j_{\alpha}(\bar{t},y,z) |^p ) |x-y| \ell_{\alpha}(z)
    \end{multlined}
      \\
    &\begin{multlined}[c][0.9\textwidth]
      \overset{{\ref{HJB:cond_jump}.\ref{HJB:cond_jump_partC}}}{\underset{\eqref{ineq:bound_through_submultiplicativity}}{\le}} 
        C \lambda |\bar{x}-\bar{y}|^2 \big[|z|^{p\wedge 2} + |z|^{\frac{p \wedge 2}{2}} \big] \sup_{\alpha\in\mathcal{A}}\|\ell_{\alpha}\|_{\infty}\\
        + \varepsilon C_p \lambda (\delta(\bar{x})+\delta(\bar{y}))  \sup_{\alpha\in\mathcal{A}}\|\ell_{\alpha}\|_{\infty} |z|^{p\wedge 2}
        + \varepsilon^{\frac12} \lambda C_p|\bar{x}-\bar{y}| [\delta(\bar{x})^{\frac12} + \delta(\bar{y})^\frac12]  
      \sup_{\alpha\in\mathcal{A}}\|\ell_{\alpha}\|_{\infty}  |z|^\frac{p\wedge 2}{2}\\
        \begin{multlined}[0.75\textwidth]
            +\gamma e^{\mu \bar{t}}
            C_{\Upsilon,p} \big[1 
            +\Upsilon(|\bar{x}|^p)+ \Upsilon\big(\varepsilon^{\frac p2} \delta(\bar{x}^{\frac p2})\big)
            +\Upsilon(|\bar{y}|^p)+ \Upsilon\big(\varepsilon^{\frac p2} \delta(\bar{y}^{\frac p2})\big)\big]
            \big[1+\Upsilon(|z|^p)\big] 
             \sup_{\alpha\in\mathcal{A}}\|\ell_{\alpha}\|_{\infty}
        \end{multlined}\\
        +C_{p,v} \big[|\bar{x} - \bar{y}| + \varepsilon^{\frac12} (\delta(\bar{x})^{\frac12} + \delta(\bar{y})^{\frac12}) \big]
        \big[ 
        \big(1 + |\bar{y}|^p + (1 + |\bar{y}|^p + \varepsilon^{\frac{p}{2}} \delta(\bar{x})^{\frac{p}{2}})|z|^p \big) 
        \big]
       \sup_{\alpha\in\mathcal{A}}\|\ell_{\alpha}\|_{\infty}
    \end{multlined}\\
    &\begin{multlined}[c][0.9\textwidth]
      \overset{\phantom{\ref{HJB:growth_cond_delta}.\ref{HJB:growth_cond_delta_partC}}}{\underset{\phantom{\text{\ref{HJB:growth_cond_delta}\ref{HJB:growth_cond_delta_partC}}}}{\leq}} 
      C_{\ell}  \lambda |\bar{x}-\bar{y}|^2 \times |z|^p 
        +\gamma e^{\mu \bar{t}}C_{\Upsilon,p,\ell}
            \big[1 + \Upsilon(|\bar{x}|^p) +\Upsilon(|\bar{y}|^p)\big]
            [1 + \Upsilon(|z|^p)]
      +\bar{\theta}^2_{\alpha,\gamma,\lambda,\varepsilon}(z)
    \end{multlined}
    \numberthis\label{lem_app:f_nonlinearities_third_3}
    \end{align*}
    where 
\begin{align}
\begin{multlined}[0.85\textwidth]
\bar{\theta}^2_{\alpha,\gamma,\lambda,\varepsilon}(z):=
\varepsilon C_{p,\ell} \lambda (\delta(\bar{x})+\delta(\bar{y}))  |z|^{p}
    +\varepsilon^{\frac12}C_{p,\ell}  \lambda|\bar{x} -\bar{y}| \big(\delta(\bar{x})^{\frac12} + \delta(\bar{y})^{\frac12}\big) \times |z|^p\\
    +\gamma e^{\mu \bar{t}}C_{\Upsilon,p,\ell} \big[ \Upsilon\big(\varepsilon^{\frac p2} \delta(\bar{x}^{\frac p2})\big)
        + \Upsilon\big(\varepsilon^{\frac p2} \delta(\bar{y}^{\frac p2})\big)\big]
        \times \big[1+\Upsilon(|z|^p)\big] \\
    +C_{p,v,\ell} \big[|\bar{x}-\bar{y}| 
    + C_p \varepsilon^{\frac12} \big(\delta(\bar{x})^\frac12  + \delta(\bar{y})^\frac12 \big)\big]
\big[1 + |\bar{y}|^p + \varepsilon^{\frac{p}{2}} \delta(\bar{y})^{\frac{p}{2}} \big] \times[1+ |z|^p].
\end{multlined}
\label{def:bar_theta_2}
\end{align}
In the second inequality of the computations above, it was used \cref{prep_aux_lemma:bound_through_submultiplicativity}.
In the third inequality we distinguished between the cases $p\in[1,2)$ and $p\ge 2$ for the first two summands.
In the former case, for the first summand we used \ref{HJB:cond_jump}.\ref{HJB:cond_jump_partC}, while for the second summand we used that $\frac{p}{2}<p$.
In the later case, we did not need any further assumptions.

Finally, collecting all the terms from \eqref{lem_app:f_nonlinearities_third_2} and \eqref{lem_app:f_nonlinearities_third_3} we obtain 
\begin{align*}
    &\hspace{-1em}\sup_{(\alpha,x,y) \in \mathcal{A}_{\varepsilon}}
    \Big(
    \overline{\mathcal{K}}^\kappa_\alpha(\bar{t},x,u^\varepsilon(\bar{t},\cdot)\circ\tau_{\bar{x}-x}) 
        -\overline{\mathcal{K}}^\kappa_\alpha(\bar{t},y,v_\varepsilon(\bar{t},\cdot)\circ\tau_{\bar{y}-y})\Big)^+\\
    &
    \begin{multlined}[0.95\textwidth]
    \overset{\eqref{lem_app:f_nonlinearities_third_2}}{\underset{\eqref{lem_app:f_nonlinearities_third_3}}{\leq}}
    \sup_{\alpha \in \mathcal{A}}
    \Big\{
        \lambda C_{p,\ell}
            |\bar{x}-\bar{y}|^2 
            \bigg[\int_{\{z: \kappa< |z| \le 1\}}  \big(|z|^2 \ell_{\alpha}(z) + |z| \ell_{\alpha}(z)\big)m_{\alpha} \textup{d}(z)
            + \int_{\{z: 1< |z|\}} |z|^p m_{\alpha} \textup{d}(z)\bigg]\\
            +\gamma e^{\mu \bar{t}}C_{\Upsilon,p}\big[1 + \Upsilon(|\bar{x}|^p) +\Upsilon(|\bar{y}|^p)\big]
            \bigg[\int_{\{z: \kappa< |z| \le 1\}} |z|\ell_{\alpha}(z) m_{\alpha} \textup{d}(z)
            + \int_{\{z: 1< |z|\}}[1 + \Upsilon(|z|^p)] m_{\alpha} \textup{d}(z)\bigg]\\
                +\int_{\{z: \kappa< |z| \le 1\}}\bar{\theta}^1_{\alpha,\gamma,\lambda,\varepsilon}(z) m_{\alpha} \textup{d}(z)
                +\int_{\{z: 1< |z|\}}\bar{\theta}^2_{\alpha,\gamma,\lambda,\varepsilon}(z) m_{\alpha} \textup{d}(z)\Big\}
      \end{multlined}\\
    &
    \begin{multlined}[0.95\textwidth]
    \overset{\phantom{\eqref{lem_app:f_nonlinearities_third_2}}}{\underset{\phantom{\eqref{lem_app:f_nonlinearities_third_3}}}{\leq}}
      \lambda C_{p,\ell} |\bar{x}-\bar{y}|^2 
    \sup_{\alpha \in \mathcal{A}}
    \Big\{
      \int_{\{z:  |z| \le 1\}}  \big[|z|^2 \ell_{\alpha}(z)+ |z| \ell_{\alpha}(z)\big]m_{\alpha} \textup{d}(z)
      +  
      \int_{\{z: 1< |z|\}} |z|^p m_{\alpha} \textup{d}(z) 
    \Big\}
    \\
    \begin{multlined}[0.85\textwidth]
    +\gamma e^{\mu \bar{t}}C_{\Upsilon,p}\big[1 + \Upsilon(|\bar{x}|^p) +\Upsilon(|\bar{y}|^p)\big]\times\\
    \times\sup_{\alpha \in \mathcal{A}}
    \Big\{
      \int_{\{z:  |z| \le 1\}} |z|\ell_{\alpha}(z) m_{\alpha} \textup{d}(z)
      +\int_{\{z: 1< |z|\}}[1 + \Upsilon(|z|^p)] m_{\alpha} \textup{d}(z)
    \Big\}
    \end{multlined}
    \\
      +\sup_{\alpha \in \mathcal{A}}
    \Big\{\int_{\{z: \kappa< |z| \le 1\}}\bar{\theta}^1_{\alpha,\gamma,\lambda,\varepsilon}(z) m_{\alpha} \textup{d}(z)
      +\int_{\{z: 1< |z|\}}\bar{\theta}^2_{\alpha,\gamma,\lambda,\varepsilon}(z) m_{\alpha} \textup{d}(z)\Big\}
    \end{multlined} \\
    &
    \begin{multlined}[0.95\textwidth]
    \overset{{\ref{HJB:UI}}}{\underset{{\ref{HJB:growth_cond_delta}.\ref{HJB:growth_cond_delta_partB}}}{=:}}
    \lambda C_{p,\ell,m} |\bar{x}-\bar{y}|^2 
    + \gamma e^{\mu \bar{t}}C_{p,\Upsilon,m,\ell}
        \big[1 + \Upsilon(|\bar{x}|^p) +\Upsilon(|\bar{y}|^p)\big]
    +\bar{\varrho}_{\gamma,\lambda,\varepsilon}^4,
    \end{multlined} 
    \numberthis\label{lem_app:f_nonlinearities_third}
\end{align*}
where the constants depend on  $p$, $\Upsilon$, on $\sup_{\alpha\in\mathcal{A}}\|\ell_{\alpha}\|_{\infty}$ and on the quantity
\begin{align*}
       \sup_{\alpha\in\mathcal{A}} 
       \int_{\{z:  |z| \le 1\}} |z|^2 + |z|\ell_{\alpha}(z) m_{\alpha} \textup{d}(z)
        +
       \sup_{\alpha\in\mathcal{A}} 
        \int_{\{z: 1< |z|\}}[1 + \Upsilon(|z|^p)] m_{\alpha} \textup{d}(z)
        \overset{\ref{HJB:UI}}{\underset{\ref{HJB:growth_cond_delta}.\ref{HJB:growth_cond_delta_partB}}{<}}\infty.
\end{align*}

Moreover, recalling \eqref{def:bar_theta_1}, \eqref{def:bar_theta_2}, it is immediate verifiable that $\bar{\varrho}^4_{\gamma,\lambda,\varepsilon}$ is finite, as the respective integrals are uniformly bounded.
Therefore, using  \cref{rem_append:aux_properties}.\ref{regularity_aux_remark:limsup_lambda_varepsilon_barx_bary}, \eqref{limit:varepsilon_delta_to_0} and the continuity of $\Upsilon$ with $\Upsilon(0)=0$, we are able to conclude 
\begin{align}
    \label{limsup_lambda_epsilon_bar_rho_4}
    \limsup_{\lambda\to\infty}
    \limsup_{\varepsilon\downarrow 0}
    \bar{\varrho}_{\gamma,\lambda,\varepsilon}^4=0.
\end{align}
In other words, we can rewrite \eqref{HJB:reg_cond_nonlinearities_aux_part_a} with the help of 
\eqref{lem_app:f_nonlinearities_first},   
\eqref{lem_app:f_nonlinearities_second},
\eqref{lem_app:f_nonlinearities_third_1},
 and 
\eqref{lem_app:f_nonlinearities_third} 
as follows
\begin{align*}
  &\begin{multlined}[c][\textwidth]
  \sup_{(\alpha,x,y)\in \mathcal{A}_{\varepsilon}} \Big\{ 
      f_\alpha(\bar{t},x,u^\varepsilon(\bar{t},\bar{x}), \sigma^T_\alpha(\bar{t},x) D_x\phi_{\delta,\gamma}^\lambda (\bar{t},\bar{x},\bar{y}), \mathcal{K}^\kappa_\alpha(\bar{t},x,u^\varepsilon(\bar{t},\cdot)\circ\tau_{\bar{x}-x},\phi_{\delta,\gamma}^\lambda(\bar{t},\cdot,\bar{y})\circ\tau_{\bar{x}-x}))\\
      - f_\alpha(\bar{t},y,v_\varepsilon(\bar{t},\bar{y}),-\sigma^T_\alpha(\bar{t},y) D_y\phi_{\delta,\gamma}^\lambda (\bar{t},\bar{x},\bar{y}), \mathcal{K}^\kappa_\alpha(\bar{t},y,v_\varepsilon(\bar{t},\cdot)\circ\tau_{\bar{y}-y},-\phi_{\delta,\gamma}^\lambda(\bar{t},\bar{x},\cdot)\circ\tau_{\bar{y}-y})) 
      \Big\}  
  \end{multlined}\\
      &\hspace{2em}
  \leq
    \lambda C_{p,\ell,m} |\bar{x}-\bar{y}|^2
    +\gamma e^{\mu \bar{t}} C_{\Upsilon,p,m,\ell}
        \big\{1 + \big[\Upsilon(|\bar{x}|^p)+\Upsilon(|\bar{y}|^p)\big] \big\}
          +\bar{\varrho}^1_{\gamma,\lambda,\varepsilon}
          +\bar{\varrho}^2_{\gamma,\lambda,\varepsilon}  
          +\bar{\varrho}^3_{\gamma,\lambda,\varepsilon,\kappa}
          +\bar{\varrho}^4_{\gamma,\lambda,\varepsilon,\kappa}
\end{align*}
with 
\begin{align*}
  \limsup_{\lambda\to\infty}
  \limsup_{\varepsilon\downarrow0}
  \limsup_{\kappa\downarrow0}
  \big[
    \bar{\varrho}^1_{\gamma,\lambda,\varepsilon}
      +\bar{\varrho}^2_{\gamma,\lambda,\varepsilon}  
      +\bar{\varrho}^3_{\gamma,\lambda,\varepsilon,\kappa}
      +\bar{\varrho}^4_{\gamma,\lambda,\varepsilon,\kappa}
      \big]\leq 0,
\end{align*}
which can be verified in view of 
\eqref{limsup_lambda_epsilon_bar_rho_1},
\eqref{limsup_lambda_epsilon_bar_rho_2},
\eqref{limsup_lambda_epsilon_bar_rho_3},
\eqref{limsup_lambda_epsilon_bar_rho_4}
and
the sub-additivity of the $\limsup$ operation.
\end{proof}

\begin{lemma}\label{lemma_appendix:regularity_L_diffusion_parts}
Under the framework of \cref{HJB:Regularity_Condition} and given the introduced notation \eqref{def:penalization_function}-\eqref{def:delta_x_phi}, \eqref{def:A_epsilon} and assumptions \eqref{supremum_penalized}-\eqref{reg_cond:assumption_ineq} of the proof of \cref{HJB:Regularity_Condition}, it holds 
\begin{align*}
    \begin{multlined}[0.85\textwidth]
    \sup_{(\alpha,x,y)\in \mathcal{A}_\varepsilon} \big\{ 
    \mathcal{L}_\alpha (\bar{t},x, D_x\phi_{\gamma,\delta,\mu}^\lambda (\bar{t},\bar{x},\bar{y}), X)
    - \mathcal{L}_\alpha (\bar{t},y, -D_y\phi_{\gamma,\delta,\mu}^\lambda (\bar{t},\bar{x},\bar{y}), Y) \big\}  
    \\
    \le 
C \lambda |\bar{x} - \bar{y}|^2 
      + \gamma e^{\mu \bar{t}} C  \big[1+ \Upsilon(|\bar{x}|^{p})+\Upsilon(|\bar{y}|^{p}) \big]
      +\varrho^2_{\gamma,\lambda,\varepsilon},
    \end{multlined}
\end{align*}
with 
    \begin{align*}
    \limsup_{\lambda\to +\infty}
    \limsup_{\varepsilon \downarrow 0} 
    \varrho^2_{\gamma,\lambda,\varepsilon} = 0,
    \end{align*}
for some non-negative constants which may depend on $d, p, T, \Upsilon$ and the families $\{b_{\alpha}\}_{\alpha\in\mathcal{A}}$ and $\{\sigma_{\alpha}\}_{\alpha\in\mathcal{A}}$.
\end{lemma}

\begin{proof}
For the closed ball of radius $R=R(\gamma)$ assumed in \cref{rem_append:aux_properties}.\ref{regularity_aux_remark:limsup_lambda_varepsilon_barx_bary} and for $0< \varepsilon\leq \varepsilon_0$
  \begin{align*}
  &\sup_{(\alpha,x,y)\in \mathcal{A}_\varepsilon} \big\{ 
        \mathcal{L}_\alpha (\bar{t},x, D_x\phi_{\gamma,\delta,\mu}^\lambda (\bar{t},\bar{x},\bar{y}), X)
        - \mathcal{L}_\alpha (\bar{t},y, -D_y\phi_{\gamma,\delta,\mu}^\lambda (\bar{t},\bar{x},\bar{y}), Y) \big\}\\
    &\begin{multlined}[0.9\textwidth]
        \hspace{1em}\underset{\phantom{\eqref{gradient_penalized_x},\eqref{gradient_penalized_y}}}{\overset{\phantom{\eqref{supremum_penalized}}}{\le} }
        \sup_{(\alpha,x,y)\in \mathcal{A}_\varepsilon} 
        \big\{ 
        	b_\alpha^T (\bar{t},x) D_x\phi_{\gamma,\delta,\mu}^\lambda(\bar{t},\bar{x},\bar{y}) 
        	- b_\alpha^T (\bar{t},y) (-D_y\phi_{\gamma,\delta,\mu}^\lambda(\bar{t},\bar{x},\bar{y})) 
    	\big\}\\
	+\sup_{(\alpha,x,y)\in \mathcal{A}_\varepsilon} 
        \big\{
        	\textup{Tr}\big[ \sigma_\alpha(\bar{t},x)\sigma_\alpha^T(\bar{t},x) X\big] 
        	- \textup{Tr}\big[ \sigma_\alpha(\bar{t},y)\sigma_\alpha^T(\bar{t},x) Y\big] 
        \big\}
        \end{multlined}
        \\
    &\begin{multlined}[0.9\textwidth]
        \hspace{1em}\overset{\eqref{gradient_penalized_x},\eqref{gradient_penalized_y}}{\underset{\eqref{supremum_penalized}}{\leq} }
        \sup_{(\alpha,x,y)\in \mathcal{A}_\varepsilon} 
        \big\{
	    2\lambda ( b_\alpha (\bar{t},x) - b_\alpha (\bar{t},y))^T(\bar{x}-\bar{y}) 
        \big\}\\     
         +  C_p\sup_{(\alpha,x,y)\in \mathcal{A}_\varepsilon}
        \big\{
        \gamma e^{\mu \bar{t}}
        \big( 
        	\Upsilon'(|\bar{x}|^p) |\bar{x}|^{p-2} b_\alpha^T(\bar{t},x)\bar{x} 
        	+ \Upsilon'(|\bar{y}|^p) |\bar{y}|^{p-2} b_\alpha^T(\bar{t},y)\bar{y}  
    	\big)    
    	\big\}
    	\numberthis\label{HJB:reg_cond_diffusion_aux_part}\\
    +\sup_{(\alpha,x,y)\in \mathcal{A}_\varepsilon} \Big\{ 
    	\sum_{j=1}^d \Big(
    		(\sigma_\alpha(\bar{t},x)^Te_j)^T X (\sigma_\alpha(\bar{t},x)^Te_j) 
    		- (\sigma_\alpha(\bar{t},y)^Te_j)^T Y (\sigma_\alpha(\bar{t},y)^Te_j) 
		\Big) \Big\},
        \end{multlined}
  \end{align*}
  where $(e_j)_{1\leq j\leq d}$ is the standard basis.
  Let us now consider each term of the right-hand side part of \eqref{HJB:reg_cond_diffusion_aux_part}:

\textbullet\hspace{0.3em}
       For the first term
       \begin{align*}
        \sup_{(\alpha,x,y)\in \mathcal{A}_\varepsilon} 
        \big\{
        2\lambda ( b_\alpha (\bar{t},x) - b_\alpha (\bar{t},y))(\bar{x}-\bar{y}) 
        \big\} 
        &\overset{\text{\ref{HJB:coef_continuity}}}{\underset{\eqref{bound:norm_x_y}}{\le}} 
        C_p\lambda 
        \big[
        |\bar{x}-\bar{y}|^2 
        + \varepsilon^{\frac12}|\bar{x}-\bar{y}|\times \big( \delta(\bar{x})^{\frac12} + \delta(\bar{y})^{\frac12} \big)
        \big]\\
        &=:C_p\lambda \big[|\bar{x}-\bar{y}|^2 +\bar{\varrho}^5_{\gamma,\lambda,\varepsilon}
        \numberthis\label{limit_second_term_L}
       \end{align*}
       with
       \begin{align}
       \lim_{\varepsilon\downarrow 0}
       \bar{\varrho}^5_{\gamma,\lambda,\varepsilon}=0.
       \label{limsup_lambda_epsilon_bar_rho_5}
       \end{align}
%
%

\textbullet\hspace{0.3em}
   For the second term, 
   \begin{align*}
    &\sup_{(\alpha,x,y)\in \mathcal{A}_\varepsilon}
	\big\{\gamma e^{\mu \bar{t}}\big(
      	\Upsilon'(|\bar{x}|^p) |\bar{x}|^{p-2} b_\alpha^T(\bar{t},x)\bar{x} 
      	+ \Upsilon'(|\bar{y}|^p) |\bar{y}|^{p-2} b_\alpha^T(\bar{t},y)\bar{y}  
	\big)    \big\}\\
      &\hspace{2em}\underset{\phantom{\eqref{corollary_Young:Upsilon_prime_upper_bound}}}{\overset{\phantom{\ref{HJB:coeff_boundedness}}}{\le}}
        \sup_{(\alpha,x,y)\in \mathcal{A}_\varepsilon}
	\big\{\gamma e^{\mu \bar{t}}\big( 
      	\Upsilon'(|\bar{x}|^p) |\bar{x}|^{p-1}|b_{\alpha}(\bar{t},x)| 
      	+ \Upsilon'(|\bar{y}|^p) |\bar{y}|^{p-1}|b_{\alpha}(\bar{t},y)|  \big)    
	\big\}\\
	&\hspace{2em}\underset{\phantom{\eqref{corollary_Young:Upsilon_prime_upper_bound}}}{\overset{\text{\ref{HJB:coef_continuity}}}{\le}}
	\sup_{(\alpha,x,y)\in \mathcal{A}_\varepsilon}
	C_{T,b}\times \gamma e^{\mu \bar{t}}\big[ 
            \Upsilon'(|\bar{x}|^p) |\bar{x}|^{p-1}(1+|x|) 
            + \Upsilon'(|\bar{y}|^p) |\bar{y}|^{p-1}(1+|y|)  
        \big]\\
	&\hspace{2em}\underset{\phantom{\eqref{corollary_Young:Upsilon_prime_upper_bound}}}{\overset{\text{\eqref{bound:norm_x_y}}}{\le}}
	C_{T,b}\times\gamma e^{\mu \bar{t}}\big[ 
            \Upsilon'(|\bar{x}|^p) |\bar{x}|^{p-1}\big(1+|\bar{x}| + \varepsilon^{\frac12} \delta(\bar{x})^{\frac12}\big) 
            + \Upsilon'(|\bar{y}|^p) |\bar{y}|^{p-1}\big(1+|\bar{y}| + \varepsilon^{\frac12} \delta(\bar{y})^{\frac12}\big) 
        \big]\\
	&
	\begin{multlined}[0.85\textwidth]
	\hspace{2em}\overset{{\eqref{corollary_Young:Upsilon_prime_upper_bound}}}{\leq}
	C_{T,b,\Upsilon} \times \gamma  e^{\mu \bar{t}}\big[1 + \Upsilon(|\bar{x}|^p)  + \Upsilon(|\bar{y}|^p)\big]\\
	+C_{T,b,\Upsilon} \times \varepsilon^{\frac12}\gamma  e^{\mu \bar{t}} 
	\big[ 
            \Upsilon'(|\bar{x}|^p) |\bar{x}|^{p-1} \delta(\bar{x})^{\frac12} 
            + \Upsilon'(|\bar{y}|^p) |\bar{y}|^{p-1} \delta(\bar{y})^{\frac12} 
    \big]
	\end{multlined}\\
	&	\hspace{2em}\overset{{\eqref{corollary_Young:Upsilon_prime_upper_bound}}}{=}
	C_{T,b,\Upsilon} \times \gamma  e^{\mu \bar{t}}\big[1 + \Upsilon(|\bar{x}|^p)  + \Upsilon(|\bar{y}|^p)\big]
	+ \bar{\varrho}^6_{\gamma,\lambda,\varepsilon}
	\numberthis
	\label{limit_third_term_L}
   \end{align*}
   with
   \begin{align}
   \limsup_{\lambda\to\infty}
       \limsup_{\varepsilon\downarrow 0}
       \bar{\varrho}^6_{\gamma,\lambda,\varepsilon}=0.
       \label{limsup_lambda_epsilon_bar_rho_6}
   \end{align}
   In the fourth inequality it was also used the moderate growth of $\Upsilon$, see \cref{lemma:UI_Young_improvement}.\ref{lemma:UI_Young_improvement:moderate_C2} and \cref{lem:Young_equiv_moderate}.
The constant $C_{T,b,\Upsilon}$ depends on the value $T$, the finite numbers $\sup_{\alpha\in\mathcal{A}} |b_\alpha(\frac{T}{2},0)|$ and $\omega(\frac{T}{2})$ and the function $\Upsilon$.

%
%
\textbullet\hspace{0.3em}
   For the third term, we have from Inequality \eqref{reg_cond:assumption_ineq} 
   for every $\xi,\eta\in\mathbb{R}^d$
   \begin{align*}
      &\xi^T X \xi - \eta^T Y\eta\\ 
      &\hspace{1em}
      \le 4\lambda |\xi - \eta|^2
		+ C \gamma e^{\mu \bar{t}} 
	\Big[ 
	    \big(\Upsilon''(|\bar{x}|^{p}) |\bar{x}|^{p}+\Upsilon'(|\bar{x}|^{p} )\big) |\bar{x}|^{p-2} |\xi|^2
    	+(\Upsilon''(|\bar{y}|^{p}) |\bar{y}|^{p}+\Upsilon'(|\bar{y}|^{p} )\big) |\bar{y}|^{p-2} |\eta|^2
    \Big].
   \end{align*}
Therefore, we have for the canonical basis $(e_j)_{j=1,\dots,d}$ of $\mathbb{R}^d$
\begin{align*}
&\sup_{(\alpha,x,y)\in \mathcal{A}_\varepsilon} \Big\{ 
      \textrm{Tr} \Big(\sigma_\alpha(\bar{t},x)\sigma_\alpha(\bar{t},x)^T X \Big) 
    - \textrm{Tr} \Big(\sigma_\alpha(\bar{t},y)\sigma_\alpha(\bar{t},y)^T Y \Big) 
      \Big\}\\
&\overset{\phantom{\eqref{corollary_Young:Upsilon_prime_upper_bound}}}{\underset{\phantom{\eqref{corollary_Young:Upsilon_double_prime_upper_bound}}}{=}} 
\sup_{(\alpha,x,y)\in \mathcal{A}_\varepsilon} \Big\{ 
      \sum_{j=1}^d \Big(
      	(\sigma_\alpha(\bar{t},x)^Te_j)^T X (\sigma_\alpha(\bar{t},x)^Te_j) 
      - (\sigma_\alpha(\bar{t},y)^Te_j)^T Y (\sigma_\alpha(\bar{t},y)^Te_j) 
      \Big) 
      \Big\}\\
&
\begin{multlined}[0.85\textwidth]
\overset{\phantom{\eqref{corollary_Young:Upsilon_prime_upper_bound}}}{\underset{\phantom{\eqref{corollary_Young:Upsilon_double_prime_upper_bound}}}{\le}} 
\sup_{(\alpha,x,y)\in \mathcal{A}_\varepsilon} C_{d}\lambda  |\sigma_{\alpha}(\bar{t},x)-\sigma_{\alpha}(\bar{t},y)|^2 
      + C_p\gamma e^{\mu \bar{t}} 
      \sup_{(\alpha,x,y)\in \mathcal{A}_\varepsilon}\Big\{
        \big[ 
            \Upsilon''(|\bar{x}|^{p}) |\bar{x}|^{p}+\Upsilon'(|\bar{x}|^{p})
        \big] |\bar{x}|^{p-2} |\sigma_\alpha(\bar{t},x)|^2\Big\}\\
      + C_p\gamma e^{\mu \bar{t}} 
      \sup_{(\alpha,x,y)\in \mathcal{A}_\varepsilon}\Big\{
          \big(\Upsilon''(|\bar{y}|^{p}) |\bar{y}|^{p}+\Upsilon'(|\bar{y}|^{p}) \big) |\bar{y}|^{p-2} |\sigma_\alpha(\bar{t},y)|^2
          \big)\Big\}
    \end{multlined} \\
&
\begin{multlined}[0.85\textwidth]
\overset{\ref{HJB:coef_continuity}}{\underset{\eqref{ineq:norm_to_delta}}{\le}} 
   C_d\lambda  \big[|\bar{x}-\bar{y}|^2 
       + \varepsilon \big(\delta(\bar{x})+ \delta(\bar{y})\big)\big]\\
	+C_{p,T,\sigma}\times \gamma e^{\mu \bar{t}}
		\big( \Upsilon''(|\bar{x}|^{p}) |\bar{x}|^{p}+\Upsilon'(|\bar{x}|^{p})\big) |\bar{x}|^{p-2} 
		\big(1 + |\bar{x}|^2 + \varepsilon \delta(\bar{x})\big)\\
  + C_{T,\sigma}\gamma e^{\mu \bar{t}} 
		\big(\Upsilon''(|\bar{y}|^{p}) |\bar{y}|^{p}+\Upsilon'(|\bar{y}|^{p}) \big) |\bar{y}|^{p-2} 
		\big(1+|\bar{y}|^2 + \varepsilon\delta(\bar{y})\big)
\end{multlined}
          \\ 
&
\begin{multlined}[0.85\textwidth]
\overset{\eqref{corollary_Young:Upsilon_prime_upper_bound}}{\underset{\eqref{corollary_Young:Upsilon_double_prime_upper_bound}}{\le}}
	C_d \lambda  |\bar{x}-\bar{y}|^2
	+C_{p,T,\sigma, \Upsilon}\times\gamma e^{\mu \bar{t}}\big[
		1 + \Upsilon(|\bar{x}|^{p}+ \Upsilon(|\bar{y}|^{p})
		\big]\\
	+\varepsilon C_d \lambda 
	\big(\delta(\bar{x}) + \delta(\bar{y})\big)
	+ \varepsilon C_{p,T,\sigma, \Upsilon}\times\gamma e^{\mu \bar{t}}\big[
		(1 + \Upsilon(|\bar{x}|^{p}) \delta(\bar{x}) + (1+ \Upsilon(|\bar{y}|^{p})\delta(\bar{y})
		\big]
\end{multlined}
	\\
&
\overset{\phantom{\eqref{corollary_Young:Upsilon_prime_upper_bound}}}{\underset{\phantom{\eqref{corollary_Young:Upsilon_double_prime_upper_bound}}}{=:}} 
	C_d\lambda  |\bar{x}-\bar{y}|^2
	+C_{p,T,\sigma, \Upsilon}\times\gamma e^{\mu \bar{t}}\big(1+ \Upsilon(|\bar{x}|^{p})+\Upsilon(|\bar{y}|^{p}) \big)
	+ \bar{\varrho}^6_{\gamma,\lambda,\varepsilon}
      \numberthis \label{limit_fourth_term_L}
   \end{align*}
   with
   \begin{align}
   \limsup_{\lambda\to\infty}
       \limsup_{\varepsilon\downarrow 0}
       \bar{\varrho}^7_{\gamma,\lambda,\varepsilon}=0.
       \label{limsup_lambda_epsilon_bar_rho_7}
   \end{align}
Combining the Inequalities 
\eqref{limit_second_term_L}, \eqref{limit_third_term_L} and \eqref{limit_fourth_term_L} we have 
\begin{align*}
  &\sup_{(\alpha,x,y)\in \mathcal{A}_\varepsilon} \big\{ 
    \mathcal{L}_\alpha (\bar{t},x,u^\varepsilon(\bar{t},\bar{x}), D_x\phi_{\gamma,\delta,\mu}^\lambda (\bar{t},\bar{x},\bar{y}), X)
    - \mathcal{L}_\alpha (\bar{t},y,v_\varepsilon(\bar{t},\bar{y}), -D_y\phi_{\gamma,\delta,\mu}^\lambda (\bar{t},\bar{x},\bar{y}), Y) \big\}\\
      &\hspace{5em}\le 
      C_{p,d} \lambda |\bar{x} - \bar{y}|^2 
      + C_{p,T,b,\sigma, \Upsilon} \times \gamma e^{\mu \bar{t}} \big[1+ \Upsilon(|\bar{x}|^{p})+\Upsilon(|\bar{y}|^{p}) \big]
      +\varrho^2_{\gamma,\lambda,\varepsilon,\kappa}
\end{align*}
   with
   \begin{gather*}
       {\varrho}^2_{\gamma,\lambda,\varepsilon}:=
       \bar{\varrho}^5_{\gamma,\lambda,\varepsilon}+\bar{\varrho}^6_{\gamma,\lambda,\varepsilon}+\bar{\varrho}^7_{\gamma,\lambda,\varepsilon}
\shortintertext{and}
   \limsup_{\lambda\to\infty}
       \limsup_{\varepsilon\downarrow 0}
       {\varrho}^2_{\gamma,\lambda,\varepsilon}
       =0.
   \end{gather*}
as required; for the conclusion we used \eqref{limsup_lambda_epsilon_bar_rho_5}, \eqref{limsup_lambda_epsilon_bar_rho_6}, \eqref{limsup_lambda_epsilon_bar_rho_7} and the subadditivity of the $\limsup$ operation.
\end{proof}

\begin{lemma}
\label{lemma_appendix:regularity_I_nonlocal_parts}
Under the framework of \cref{HJB:Regularity_Condition} and given the introduced notation \eqref{def:penalization_function}-\eqref{def:delta_x_phi}, \eqref{def:A_epsilon} and assumptions \eqref{supremum_penalized}-\eqref{reg_cond:assumption_ineq} of the proof of \cref{HJB:Regularity_Condition}, it holds
\begin{align*}
&
\sup_{(\alpha,x,y)\in \mathcal{A}_\varepsilon} \Big\{
    \mathcal{I}^\kappa_\alpha\big(\bar{t},x,u^\varepsilon(\bar{t},\cdot)\circ\tau_{\bar{x}-x},\phi_{\delta,\gamma}^\lambda(\bar{t},\cdot,\bar{y})\circ\tau_{\bar{x}-x}\big)
    -\mathcal{I}^\kappa_\alpha\big(\bar{t},y,v_\varepsilon(\bar{t},\cdot)\circ\tau_{\bar{y}-y},-\phi_{\delta,\gamma}^\lambda(\bar{t},\bar{x},\cdot)\circ\tau_{\bar{y}-y}\big)
             \Big\}
\\
&
\hspace{2em}\leq 
\lambda C |\bar{x} - \bar{y}|^2 \times 
+C \gamma e^{\mu \bar{t}}
  \big[1    + \Upsilon(|\bar{x}|^p) + \Upsilon(|\bar{y}|^p)\big] 
+\varrho^{3}_{\gamma,\lambda,\varepsilon,\kappa}
\end{align*}
with 
\begin{align*}
\limsup_{\lambda\to 0}
\limsup_{\varepsilon\downarrow 0}
\limsup_{\kappa\downarrow 0}
\varrho^{3}_{\gamma,\lambda,\varepsilon,\kappa}=0,
\end{align*}
for some non-negative constants which may depend on $\Upsilon,p$ and the family $\{m_{\alpha}\}_{\alpha\in\mathcal{A}}$.
\end{lemma}

\begin{proof}
We will consider  separately the following cases $\{z: |z|\le \kappa\}$, $\{z:\kappa<|z|\le 1\}$ and $\{ z: |z|>1\}$.
For the closed ball of radius $R=R(\gamma)$ assumed in \cref{rem_append:aux_properties}.\ref{regularity_aux_remark:limsup_lambda_varepsilon_barx_bary} and for $0< \varepsilon\leq \varepsilon_0$
\begin{align*}
&
\sup_{(\alpha,x,y)\in \mathcal{A}_\varepsilon} \Big\{
    \mathcal{I}^\kappa_\alpha\big(\bar{t},x,u^\varepsilon(\bar{t},\cdot)\circ\tau_{\bar{x}-x},\phi_{\delta,\gamma}^\lambda(\bar{t},\cdot,\bar{y})\circ\tau_{\bar{x}-x}\big)
    -\mathcal{I}^\kappa_\alpha\big(\bar{t},y,v_\varepsilon(\bar{t},\cdot)\circ\tau_{\bar{y}-y},-\phi_{\delta,\gamma}^\lambda(\bar{t},\bar{x},\cdot)\circ\tau_{\bar{y}-y}\big)
             \Big\}
\\
&\begin{multlined}[0.85\textwidth]
\hspace{1em}\leq 
\sup_{(\alpha,x,y)\in \mathcal{A}_\varepsilon}
    \big|\check{\mathcal{I}}^\kappa_\alpha\big(\bar{t},x,\phi_{\delta,\gamma}^\lambda(\bar{t},\cdot,\bar{y})\circ\tau_{\bar{x}-x}\big)\big|
+\sup_{(\alpha,x,y)\in \mathcal{A}_\varepsilon}
    \big|\check{\mathcal{I}}^\kappa_\alpha\big(\bar{t},y,-\phi_{\delta,\gamma}^\lambda(\bar{t},\bar{x},\cdot)\circ\tau_{\bar{y}-y}\big)\big|
\numberthis\label{I_nonlocal_term}
\end{multlined}\\
&
\hspace{1em}+
\sup_{(\alpha,x,y)\in \mathcal{A}_\varepsilon} \Big\{
    \overline{\mathcal{I}}^\kappa_\alpha\big(\bar{t},x,u^\varepsilon(\bar{t},\cdot)\circ\tau_{\bar{x}-x},\phi_{\delta,\gamma}^\lambda(\bar{t},\cdot,\bar{y})\circ\tau_{\bar{x}-x}\big)
    -\overline{\mathcal{I}}^\kappa_\alpha\big(\bar{t},y,v_\varepsilon(\bar{t},\cdot)\circ\tau_{\bar{y}-y},-\phi_{\delta,\gamma}^\lambda(\bar{t},\bar{x},\cdot)\circ\tau_{\bar{y}-y}\big)
             \Big\}
\\
&
\hspace{1em}+  
\sup_{(\alpha,x,y)\in \mathcal{A}_\varepsilon} \Big\{
    \hat{\mathcal{I}}^\kappa_\alpha\big(\bar{t},x,u^\varepsilon(\bar{t},\cdot)\circ\tau_{\bar{x}-x},\phi_{\delta,\gamma}^\lambda(\bar{t},\cdot,\bar{y})\circ\tau_{\bar{x}-x}\big)
    -\hat{\mathcal{I}}^\kappa_\alpha\big(\bar{t},y,v_\varepsilon(\bar{t},\cdot)\circ\tau_{\bar{y}-y},-\phi_{\delta,\gamma}^\lambda(\bar{t},\bar{x},\cdot)\circ\tau_{\bar{y}-y}\big)
             \Big\}
\end{align*}

\textbullet\hspace{0.3em}
      Case $\{z: |z|\le \kappa\}$ for $\kappa \in(0,1)$: Applying Taylor's Theorem one has
\begin{align*}
&\sup_{(\alpha,x,y)\in \mathcal{A}_\varepsilon}
\big|\check{\mathcal{I}}^\kappa_\alpha\big(\bar{t},x,\phi_{\delta,\gamma}^\lambda(\bar{t},\cdot,\bar{y})\circ\tau_{\bar{x}-x}\big)\big|\\
&\hspace{1em}
\overset{\phantom{\ref{HJB:cond_jump}}}{\leq} 
\sup_{(\alpha,x,y)\in \mathcal{A}_\varepsilon}
    \Big|
    \int_{\{z:|z|\leq \kappa\}}\int_0^1 (1-\xi) 
        \big| D^2_{xx} \phi_{\delta,\gamma}^\lambda(\bar{t},\bar{x}+\xi j_{\alpha}(\bar{t},x,z),\bar{y})\big|
        \times |j_{\alpha}(\bar{t},x,z)|^2
    \textup{d}\xi\ m_{\alpha}(\textup{d}z)
    \Big|\\
&\hspace{1em}
\begin{multlined}[0.9\textwidth]
\overset{\ref{HJB:cond_jump}}{\underset{\eqref{ineq:norm_to_delta}}{\leq}} 
\sup_{(\alpha,x,y)\in \mathcal{A}_\varepsilon}
    \Big\{\sup_{|\zeta|\leq 1, \xi\in[0,1]}
            \big| D^2_{xx} \phi_{\delta,\gamma}^\lambda(\bar{t},\bar{x}+\xi j_{\alpha}(\bar{t},x,z),\bar{y})\big|
    \Big\}\times \\
        \times C_p(1+ |\bar{x}|^{\frac{p\wedge 2}{2}} + \varepsilon^{\frac{p\wedge 2}{2}} \delta(\bar{x})^{\frac{p\wedge 2}{2}})
        \times 
\sup_{\alpha\in \mathcal{A}}
        \int_{\{z:|z|\leq \kappa\}}|z|^2 m_{\alpha}(\textup{d}z).
\end{multlined}
\end{align*}
In view of the above estimation, we can conclude
\begin{align*}
\limsup_{\kappa\downarrow 0}
\sup_{(\alpha,x,y)\in \mathcal{A}_\varepsilon}
\big|\check{\mathcal{I}}^\kappa_\alpha\big(\bar{t},x,\phi_{\delta,\gamma}^\lambda(\bar{t},\cdot,\bar{y})\circ\tau_{\bar{x}-x}\big)\big| =0,
\end{align*}
where we used the boundedness of a continuous function on a compact set which contains the bounded\footnote{Use \cref{rem_append:aux_properties}.\ref{regularity_aux_remark:limsup_lambda_varepsilon_barx_bary}.} set
\begin{align*}
  \{(\bar{t},\bar{x}+C\xi (1+|x|)\zeta,\bar{y}) : \xi \in[0,1], |\zeta|\le 1, x \text{ is s.t. } \varphi(x-\bar{x})\le \varepsilon \delta(\bar{x}), \bar{x},\bar{y} \text{ for }\lambda >0 \text{ and }\varepsilon\le\varepsilon_0 \}
\end{align*} 
and the fact that we work on the set $\{z: |z|\le \kappa\}$ with $\kappa\le 1$.
Analogously, we have for the term 
\begin{align*}
\limsup_{\kappa\downarrow 0}
\sup_{(\alpha,x,y)\in \mathcal{A}_\varepsilon}
\big|\check{\mathcal{I}}^\kappa_\alpha\big(\bar{t},y,-\phi_{\delta,\gamma}^\lambda(\bar{t},\bar{x},\cdot)\circ\tau_{\bar{y}-y}\big)\big| =0,
\end{align*}
In total, for
\begin{align*}
\begin{multlined}[0.85\textwidth]
\bar{\varrho}^8_{\gamma,\lambda,\varepsilon,\kappa}:= 
\sup_{(\alpha,x,y)\in \mathcal{A}_\varepsilon}
    \big|\check{\mathcal{I}}^\kappa_\alpha\big(\bar{t},x,\phi_{\delta,\gamma}^\lambda(\bar{t},\cdot,\bar{y})\circ\tau_{\bar{x}-x}\big)\big|
+\sup_{(\alpha,x,y)\in \mathcal{A}_\varepsilon}
    \big|\check{\mathcal{I}}^\kappa_\alpha\big(\bar{t},y,-\phi_{\delta,\gamma}^\lambda(\bar{t},\bar{x},\cdot)\circ\tau_{\bar{y}-y}\big)\big|
\end{multlined}
\end{align*}
it holds
\begin{align}\label{limsup_lambda_epsilon_bar_rho_8}
  \limsup_{\kappa \downarrow 0}
  \bar{\varrho}^8_{\gamma,\lambda,\varepsilon,\kappa}
  =0.
\end{align}
%
%
%
%
%
%
%
%
%
%
%
%
%
%
%

\textbullet\hspace{0.3em}
    Case $\{z: \kappa<|z|\le 1\}$  for $\kappa \in(0,1)$: 
    For $(\alpha,x,y)\in\mathcal{A}_\varepsilon$, we initially compute the difference of the integrands the $\overline{\mathcal{I}}^\kappa_\alpha-$parts:
\begin{align*}
&
\begin{multlined}[0.9\textwidth]
u^\varepsilon (\bar{t},\bar{x}+j_\alpha(\bar{t},x,z)) - u^\varepsilon (\bar{t},\bar{x}) - D_x \phi_{\delta,\gamma}^\lambda(\bar{t},\bar{x},\bar{y})j_\alpha(\bar{t},x,z) \\
-\big[ v_\varepsilon (\bar{t},\bar{y}+j_\alpha(\bar{t},y,z)) - u_\varepsilon (\bar{t},\bar{y}) + D_y \phi_{\delta,\gamma}^\lambda(\bar{t},\bar{x},\bar{y})j_\alpha(\bar{t},y,z)\big]
\end{multlined}\\
&
\begin{multlined}[0.8\textwidth]
\overset{\phantom{\eqref{gradient_penalized_x}, \eqref{gradient_penalized_y}}}{\underset{\phantom{\eqref{supremum_penalized}}}{=} } 
\big[ u^\varepsilon (\bar{t},\bar{x}+j_\alpha(\bar{t},x,z)) - v_\varepsilon (\bar{t},\bar{y}+j_\alpha(\bar{t},y,z)) - \phi_{\delta,\gamma}^\lambda \big(\bar{t},\bar{x} + j_\alpha(\bar{t},x,z),\bar{y}+j_\alpha(\bar{t},y,z)\big)\big] \\
\hspace{-14em}
-\big[ u^\varepsilon (\bar{t},\bar{x}) - u_\varepsilon (\bar{t},\bar{y}) - \phi_{\delta,\gamma}^\lambda (\bar{t},\bar{x},\bar{y})\big] \\
+\big[ \phi_{\delta,\gamma}^\lambda \big(\bar{t},\bar{x} + j_\alpha(\bar{t},x,z),\bar{y}+j_\alpha(\bar{t},y,z)\big) - \phi_{\delta,\gamma}^\lambda (\bar{t},\bar{x},\bar{y})\big] \\
-\big[ 
D_y \phi_{\delta,\gamma}^\lambda(\bar{t},\bar{x},\bar{y})^Tj_\alpha(\bar{t},y,z) 
+D_x \phi_{\delta,\gamma}^\lambda(\bar{t},\bar{x},\bar{y})^Tj_\alpha(\bar{t},x,z)\big]
\end{multlined}\\
&\begin{multlined}[0.85\textwidth]
\overset{\eqref{gradient_penalized_x}, \eqref{gradient_penalized_y}}{\underset{\eqref{supremum_penalized}}{\le } } 
\lambda |j_\alpha(\bar{t},x,z) - j_\alpha(\bar{t},y,z) |^2\\
+\gamma e^{\mu \bar{t}} \Big[
    \Upsilon(|\bar{x}+j_\alpha(\bar{t},x,z)|^p) +\Upsilon(|\bar{y}+j_\alpha(\bar{t},y,z)|^p)
    - \Upsilon(|\bar{x}|^p) - \Upsilon(|\bar{y}|^p)\\ 
- p\big[
    \Upsilon'(|\bar{x}|^p) |\bar{x}|^{p-2} \bar{x}^T j_\alpha(\bar{t},x,z)
    +
    \Upsilon'(|\bar{y}|^p) |\bar{y}|^{p-2} \bar{y}^T j_\alpha(\bar{t},y,z)\big]
                        \Big]
\end{multlined}
\numberthis\label{ineq:nonlocal_first_computations}\\
&\begin{multlined}[0.85\textwidth]
\overset{\ref{HJB:cond_jump},\eqref{ineq:norm_to_delta}}{\underset{\eqref{ineq:SecondOrderTaylor}}{\le } } 
\lambda C \big[|\bar{x} - \bar{y}|^2 + \varepsilon \big(\delta(\bar{x})+\delta(\bar{y})\big)\big]\times |z|^2\\
+ \gamma e^{\mu \bar{t}} C_{p,d,\Upsilon}\Big[
1
+ \Upsilon(|\bar{x}|^p) + \Upsilon\big(\varepsilon^{\frac p2} \delta(\bar{x}^{\frac p2})\big) 
+ \Upsilon(|\bar{y}|^p) + \Upsilon\big(\varepsilon^{\frac p2} \delta(\bar{y}^{\frac p2})\big) 
\Big]
|z|^2
\end{multlined}\\
&
\overset{\phantom{\ref{HJB:cond_jump},\eqref{ineq:norm_to_delta}}}{\underset{\phantom{\eqref{ineq:SecondOrderTaylor}}}{=:} } 
\begin{multlined}[0.85\textwidth]
\lambda C |\bar{x} - \bar{y}|^2 |z|^2
+ \gamma e^{\mu \bar{t}} C_{p,d,\Upsilon}\big[ 1 + \Upsilon(|\bar{x}|^p) + \Upsilon(|\bar{y}|^p) \big] |z|^2 
+ \bar{\theta}^3_{\gamma,\lambda,\varepsilon}(z)
\end{multlined}
\numberthis\label{lem_app:I_nonlocal_second}
\end{align*}
with 
\begin{align}\label{def:bar_theta_3}
\bar{\theta}^3_{\gamma,\lambda,\varepsilon}:=
 \big[
\varepsilon \big(\delta(\bar{x})+\delta(\bar{y})\big)
+
\Upsilon\big(\varepsilon^{\frac p2} \delta(\bar{x}^{\frac p2})\big) 
+ \Upsilon\big(\varepsilon^{\frac p2} \delta(\bar{y}^{\frac p2})\big) 
\big]\times |z|^2
\end{align}
Therefore,
\begin{align*}  
&
\sup_{(\alpha,x,y)\in \mathcal{A}_\varepsilon} \Big\{
    \overline{\mathcal{I}}^\kappa_\alpha\big(\bar{t},x,u^\varepsilon(\bar{t},\cdot)\circ\tau_{\bar{x}-x},\phi_{\delta,\gamma}^\lambda(\bar{t},\cdot,\bar{y})\circ\tau_{\bar{x}-x}\big)
    -\overline{\mathcal{I}}^\kappa_\alpha\big(\bar{t},y,v_\varepsilon(\bar{t},\cdot)\circ\tau_{\bar{y}-y},-\phi_{\delta,\gamma}^\lambda(\bar{t},\bar{x},\cdot)\circ\tau_{\bar{y}-y}\big)
             \Big\}
\\
&\begin{multlined}
\hspace{2em}\overset{\eqref{lem_app:I_nonlocal_second}}{\leq}
\Big\{\lambda C |\bar{x} - \bar{y}|^2 
+ \gamma e^{\mu \bar{t}} C_{p,d,\Upsilon}\big[ 1 + \Upsilon(|\bar{x}|^p) + \Upsilon(|\bar{y}|^p) \big]\Big\} 
\sup_{(\alpha,x,y)\in \mathcal{A}_\varepsilon}\int_{\{z:|z|\leq 1\}} |z|^2 m_{\alpha}(\textup{d}z)\\ 
+ \sup_{(\alpha,x,y)\in \mathcal{A}_\varepsilon}\int_{\{z:|z|\leq 1\}}\bar{\theta}^3_{\gamma,\lambda,\varepsilon}(z)m_{\alpha}(\textup{d}z)
\end{multlined}\\
&\begin{multlined}
\hspace{2em}\overset{\phantom{\eqref{lem_app:I_nonlocal_second}}}{=:}
\lambda C_{m} |\bar{x} - \bar{y}|^2 
+ \gamma e^{\mu \bar{t}} C_{p,d,\Upsilon,m}\big[ 1 + \Upsilon(|\bar{x}|^p) + \Upsilon(|\bar{y}|^p) \big]
+\bar{\varrho}_{\gamma,\lambda,\varepsilon}^9
\end{multlined}
\numberthis\label{lem_app:I_nonlocal_second_final}
\end{align*}    
with 
\begin{align}
\label{def:limsup_varrho_9}
\lim_{\varepsilon\downarrow 0} \bar{\varrho}_{\gamma,\lambda,\varepsilon}^9=
\lim_{\varepsilon\downarrow 0} \sup_{(\alpha,x,y)\in \mathcal{A}_\varepsilon}\int_{\{z:|z|\leq 1\}}\bar{\theta}^3_{\gamma,\lambda,\varepsilon}(z)m_{\alpha}(\textup{d}z)
\overset{\eqref{def:bar_theta_3}}{\underset{\ref{HJB:coeff_boundedness}}{=}}
0.
\end{align}

%
%
%
%
%
%
%
%
%

\textbullet\hspace{0.3em}
      Case $\{z: |z|\ge 1\}$:  For $(\alpha,x,y)\in\mathcal{A}_\varepsilon$, we initially compute the difference of the integrands 
    of the $\hat{\mathcal{I}}_\alpha-$parts.
\begin{align*}
  &
  \begin{multlined}[0.9\textwidth]
  u^\varepsilon (\bar{t},\bar{x}+j_\alpha(\bar{t},x,z)) - u^\varepsilon (\bar{t},\bar{x}) - D_x \phi_{\delta,\gamma}^\lambda(\bar{t},\bar{x},\bar{y})j_\alpha(\bar{t},x,z) \\
  -\big[ v_\varepsilon (\bar{t},\bar{y}+j_\alpha(\bar{t},y,z)) - u_\varepsilon (\bar{t},\bar{y}) + D_y \phi_{\delta,\gamma}^\lambda(\bar{t},\bar{x},\bar{y})j_\alpha(\bar{t},y,z)\big] 
  \end{multlined}\\
  &\begin{multlined}[0.85\textwidth]
\hspace{-0.3em}\overset{\eqref{ineq:nonlocal_first_computations}}{\underset{\phantom{\eqref{supremum_penalized}}}{\leq } } 
\lambda |j_\alpha(\bar{t},x,z) - j_\alpha(\bar{t},y,z) |^2\\
+\gamma e^{\mu \bar{t}} \Big[
    \Upsilon(|\bar{x}+j_\alpha(\bar{t},x,z)|^p) +\Upsilon(|\bar{y}+j_\alpha(\bar{t},y,z)|^p)
    - \Upsilon(|\bar{x}|^p) - \Upsilon(|\bar{y}|^p)\\ 
- p\big[
    \Upsilon'(|\bar{x}|^p) |\bar{x}|^{p-2} \bar{x}^T j_\alpha(\bar{t},x,z)
    +
    \Upsilon'(|\bar{y}|^p) |\bar{y}|^{p-2} \bar{y}^T j_\alpha(\bar{t},y,z)\big]
                        \Big]
\end{multlined}\\
  &
  \begin{multlined}[0.85\textwidth]
\hspace{-0.7em}  \overset{\ref{HJB:cond_jump}.\ref{HJB:cond_jump_partC}}{\underset{\eqref{ineq:bound_through_submultiplicativity}}{\leq} } 
  \lambda |x-y|^2 |z|^{p\wedge2}\\
+C_{\Upsilon,p} \times \gamma e^{\mu \bar{t}}
  \Big[1    + \Upsilon(|\bar{x}|^p)+ \Upsilon\big(\varepsilon^{\frac p2} \delta(\bar{x}^{\frac p2})\big)
            + \Upsilon(|\bar{y}|^p)+ \Upsilon\big(\varepsilon^{\frac p2} \delta(\bar{y}^{\frac p2})\big)
            \Big]\big[1+\Upsilon(|z|^p)\big]
  \\ 
+ p\big[\Upsilon'(|\bar{x}|^p) |\bar{x}|^{p-1} (1+|x|)
    +\Upsilon'(|\bar{y}|^p) |\bar{y}|^{p-1} (1+|y|)\big] \times |z|
\end{multlined}
\numberthis\label{remark:about_p_geq_1}\\
  &
  \begin{multlined}[0.85\textwidth]
\hspace{-0.3em}  \overset{\eqref{bound:norm_x_y}}{\underset{\eqref{ineq:norm_to_delta}}{\leq} } 
\lambda C \big[|\bar{x} - \bar{y}|^2 + \varepsilon \big(\delta(\bar{x})+\delta(\bar{y})\big)\big]\times |z|^{p\wedge2}\\
+C_{\Upsilon,p} \times \gamma e^{\mu \bar{t}}
  \Big[1    + \Upsilon(|\bar{x}|^p)+ \Upsilon\big(\varepsilon^{\frac p2} \delta(\bar{x}^{\frac p2})\big)
            + \Upsilon(|\bar{y}|^p)+ \Upsilon\big(\varepsilon^{\frac p2} \delta(\bar{y}^{\frac p2})\big)
            \Big]\big[1+\Upsilon(|z|^p)\big]
  \\ 
+ p\gamma e^{\mu \bar{t}}\big[\Upsilon'(|\bar{x}|^p) |\bar{x}|^{p-1} \big(1+|\bar{x}|+ \varepsilon^{\frac12}\delta(\bar{x})^{\frac12}\big)
    +\Upsilon'(|\bar{y}|^p) |\bar{y}|^{p-1} \big(1+|\bar{y}|+ \varepsilon^{\frac12}\delta(\bar{y})^{\frac12}\big)\big] \times |z|
\end{multlined}\\
  &
  \begin{multlined}[0.85\textwidth]
  \overset{{\eqref{corollary_Young:Upsilon_prime_upper_bound}}}{\underset{}{\leq} } 
\lambda C |\bar{x} - \bar{y}|^2 \times |z|^{p\wedge2}
+C_{\Upsilon,p} \gamma e^{\mu \bar{t}}
  \Big[1    + \Upsilon(|\bar{x}|^p) + \Upsilon(|\bar{y}|^p)
            \Big]\big[1+\Upsilon(|z|^p)\big]
  \\ 
+ C_{\Upsilon,p}\gamma e^{\mu \bar{t}}\big[\Upsilon(|\bar{x}|^p) + \Upsilon(|\bar{y}|^p) \big] \times |z|
+\bar{\theta}^4_{\gamma,\lambda,\varepsilon}(z)
\end{multlined}
  \numberthis\label{lem_app:I_nonlocal_third}
\end{align*}
with 
\begin{align}
\label{def:bar_theta_4}
\begin{multlined}[0.85\textwidth]
\bar{\theta}^4_{\gamma,\lambda,\varepsilon}(z)
:=\varepsilon \big(\delta(\bar{x})+\delta(\bar{y})\big) |z|^{p\wedge2}\\
    +C_{\Upsilon,p}\times \gamma e^{\mu \bar{t}}\big[ \Upsilon\big(\varepsilon^{\frac p2} \delta(\bar{x}^{\frac p2})\big)
    + \Upsilon\big(\varepsilon^{\frac p2} \delta(\bar{y}^{\frac p2})\big)
            \big]\big[1+\Upsilon(|z|^p)\big]\\  
    + p \gamma e^{\mu \bar{t}}\varepsilon^{\frac12}\big[\Upsilon'(|\bar{x}|^p) |\bar{x}|^{p-1} \delta(\bar{x})^{\frac12}
    +\Upsilon'(|\bar{y}|^p) |\bar{y}|^{p-1} \delta(\bar{y})^{\frac12}\big] \times |z|
\end{multlined}.
\end{align}
In \eqref{lem_app:I_nonlocal_third} it was also used the moderate growth of $\Upsilon$, see \cref{lemma:UI_Young_improvement} and \cref{lem:Young_equiv_moderate}.
%
%
%
%
%
%
%
Therefore,
\begin{align*}
&
\sup_{(\alpha,x,y)\in \mathcal{A}_\varepsilon} \Big\{
    \hat{\mathcal{I}}^\kappa_\alpha\big(\bar{t},x,u^\varepsilon(\bar{t},\cdot)\circ\tau_{\bar{x}-x},\phi_{\delta,\gamma}^\lambda(\bar{t},\cdot,\bar{y})\circ\tau_{\bar{x}-x}\big)
    -\hat{\mathcal{I}}^\kappa_\alpha\big(\bar{t},y,v_\varepsilon(\bar{t},\cdot)\circ\tau_{\bar{y}-y},-\phi_{\delta,\gamma}^\lambda(\bar{t},\bar{x},\cdot)\circ\tau_{\bar{y}-y}\big)
             \Big\}
\\
&\begin{multlined}[0.85\textwidth]
\hspace{2em} \overset{\eqref{lem_app:I_nonlocal_third}}\leq 
\lambda C |\bar{x} - \bar{y}|^2 
\sup_{\alpha\in\mathcal{A}}\int_{\{z:1<|z|\}} |z|^{p\wedge2} m_{\alpha}(\textup{d}z)\\
+C_{\Upsilon,p}\times \gamma e^{\mu \bar{t}}
  \Big[1    + \Upsilon(|\bar{x}|^p) + \Upsilon(|\bar{y}|^p)\Big]
\sup_{\alpha\in\mathcal{A}}\int_{\{z:1<|z|\}} \big[1+\Upsilon(|z|^p)\big] m_{\alpha}(\textup{d}z)
  \\ 
+ C_{\Upsilon,p}\times \gamma e^{\mu \bar{t}}\big[\Upsilon(|\bar{x}|^p) + \Upsilon(|\bar{y}|^p) \big] 
\sup_{\alpha\in\mathcal{A}}\int_{\{z:1<|z|\}} |z| m_{\alpha}(\textup{d}z)\\
+\sup_{\alpha\in\mathcal{A}}\int_{\{z:1<|z|\}} \bar{\theta}^4_{\gamma,\lambda,\varepsilon}(z) m_{\alpha}(\textup{d}z)
\end{multlined}\\
&\begin{multlined}[0.85\textwidth]
\hspace{2em} \overset{\ref{HJB:coeff_boundedness}}{\underset{}{\leq}} 
\lambda C_m |\bar{x} - \bar{y}|^2 
+C_{\Upsilon,p,m}  \times\gamma e^{\mu \bar{t}}
  \Big[1    + \Upsilon(|\bar{x}|^p) + \Upsilon(|\bar{y}|^p)\Big] 
+\sup_{\alpha\in\mathcal{A}}\int_{\{z:1<|z|\}} \bar{\theta}^4_{\gamma,\lambda,\varepsilon}(z) m_{\alpha}(\textup{d}z)
\end{multlined}\\
&\begin{multlined}[0.85\textwidth]
\hspace{2em} \overset{\phantom{\ref{HJB:coeff_boundedness}}}{=:} 
\lambda C_m |\bar{x} - \bar{y}|^2  
+C_{\Upsilon,p,m}  \times \gamma e^{\mu \bar{t}}
  \Big[1    + \Upsilon(|\bar{x}|^p) + \Upsilon(|\bar{y}|^p)\Big] 
+\bar{\varrho}^{10}_{\gamma,\lambda,\varepsilon}
\end{multlined}
\numberthis\label{lem_app:I_nonlocal_third_final}
\end{align*}
where the combination of \ref{HJB:UI} and \cref{prop:UI_Young} yields
\begin{align*}
\sup_{\alpha\in\mathcal{A}}\int_{\{z:1<|z|\}} \Upsilon(|z|^p)\, m_{\alpha}(\textup{d}z)<\infty
\end{align*}
and it was also used the assumption $p\geq 1$.
It is only left to observe that 
\begin{align}\label{def:limsup_varrho_10}
\lim_{\varepsilon\downarrow 0}
\bar{\varrho}^{10}_{\gamma,\lambda,\varepsilon}
\overset{\eqref{def:bar_theta_4}} =0.
\end{align}
Collecting the terms \eqref{limsup_lambda_epsilon_bar_rho_8}, \eqref{lem_app:I_nonlocal_second_final}, \eqref{def:limsup_varrho_9}, \eqref{lem_app:I_nonlocal_third_final} and \eqref{def:limsup_varrho_10} we rewrite \eqref{I_nonlocal_term} in the desired form for 
\begin{align*}
\varrho^3_{\gamma,\lambda,\varepsilon,\kappa}=
\bar{\varrho}_{\gamma,\lambda,\varepsilon,\kappa}^8
+\bar{\varrho}_{\gamma,\lambda,\varepsilon}^9
+\bar{\varrho}_{\gamma,\lambda,\varepsilon}^{10}.
\end{align*}

\end{proof}
\section{Proofs of the results presented in Section \ref{sec:YoungFunctions}.}\label{appendix:YoungFunctions_proofs}
\subsection{Proof of Proposition \ref{prop:UI_Young}}

\begin{proof}[Proof of \cref{prop:UI_Young}]
  The equivalence is -\emph{mutatis mutandis}- a restatement of the classical de La Vallée Poussin Theorem. 
  We provide all the details for the validity of the equivalence, since they will be crucial later for the construction of the twice continuously differentiable,  moderate Young function.
  For the equivalence, we follow \citet[Theorem II.22]{dellacherie1978probabilities}.
  We define 
  \begin{align}
    M:= \sup_{\alpha\in\mathcal{A}} \int_{\{z:|z|\ge 1\}} \Upsilon(|z|^p)m_\alpha(\textup{d}z).
    \label{def:M_supremum_of_Young_integrals}
  \end{align}
%
  and we assume initially that $M$ is finite. 
  Let $\varepsilon>0$ and define $K:=\frac{M}{\varepsilon}$.
  Since $\Upsilon$ is a Young function, there exists $R>0$ such that $\displaystyle \frac{\Upsilon(x)}{x}\ge K$, for every $x>R$.
  Hence, we have 
  \begin{align*}
    \sup_{\alpha\in\mathcal{A}} \int_{\{z:|z|>R^{\frac{1}{p}}\}} |z|^p m_\alpha(\textup{d}z)
      \le \sup_{\alpha\in\mathcal{A}} \frac{1}{K}\int_{\{z:|z|>R^{\frac{1}{p}}\}} \Upsilon(|z|^p) m_\alpha(\textup{d}z)
      \le \varepsilon.
  \end{align*}   

  In order to prove the other direction of the equivelence, \emph{i.e.}, the existence of a moderate Young function $\Upsilon$ such that $M$ is finite, we will construct a piecewise constant, non-negative, increasing function $\upsilon$, which will serve as the right-derivative of the desired Young function. 
  The construction will be based on the assumption that 
  \begin{align}
    \lim_{R\to\infty}\sup_{\alpha\in\mathcal{A}} \int_{\{z:|z|>R\}} |z|^p m_\alpha(\textup{d}z) =0.
  \label{Uniform_Integr_Limit_zero}
  \end{align}
  The function $\upsilon$ will initially be chosen to satisfy $\upsilon(x)=\upsilon_n$, for $x\in [n,n+1)$, for $n\in\mathbb{N}\cup\{0\}$, where $0\le \upsilon_n\le \upsilon_{n+1}$ for every $n\in\mathbb{N}\cup\{0\}$.
  We are allowed to choose $\upsilon_0=0$ and we do so. 
  Before we choose the rest values $(\upsilon_n)_{n\in\mathbb{N}}$, the reader may observe that for every $n\in\mathbb{N}$,
  \begin{align*}  
    \Upsilon(x) =  \sum_{k=1}^{n-1} \upsilon_k + \upsilon_n(x-n), \text{ for }x\in[n,n+1),
  \end{align*}  
  which in particular implies 
  \begin{align*}
    \Upsilon(x) \le  \sum_{k=1}^{n} \upsilon_k, \text{ for }x\in[n,n+1).
  \end{align*}
  Let us proceed by determining an upper bound of $M$ in terms of $(\upsilon_n)_{n\in\mathbb{N}}$:
  \begin{align*}
    \int_{\{z:|z|\ge 1\}} \Upsilon(|z|^p) m_\alpha(\textup{d}z) 
    \le \sum_{n=1}^{\infty} (\sum_{k=1}^n \upsilon_k) m_\alpha(\{z:n\le |z|^p<n+1\})
    = \sum_{n=1}^{\infty} \upsilon_n m_\alpha(\{z:|z|^p\ge n\}).
  \end{align*}
  Hence,
  \begin{align}
    M=\sup_{\alpha\in\mathcal{A}} \int_{\{z:|z|\ge 1\}} \Upsilon(|z|^p) m_\alpha(\textup{d}z) 
    \le \sup_{\alpha\in\mathcal{A}} \sum_{n=1}^{\infty} \upsilon_n m_\alpha(\{z:|z|^p\ge n\}).
    \label{ineq:M_upperBound}
  \end{align}
  As a next step, we intend to find a way to control the quantity on the right-hand side of Inequality \eqref{ineq:M_upperBound} for sequences $(\upsilon_n)_{n\in\mathbb{N}}$ with desired properties, namely $(\upsilon_n)_{n\in\mathbb{N}}$ is non-decreasing and non-negative. 
  To this end, let us initially fix a decreasing sequence $(d_n)_{n\in\mathbb{N}}$ of positive numbers which is summable, \emph{e.g.} $d_n:=2^{-n}$ for every ${n\in\mathbb{N}}$. 
  We further define the set 
  \begin{align*}
    \mathcal{N}_n:=\Big\{k\in \mathbb{N}: \sup_{\alpha\in\mathcal{A}} \int_{\{z:|z|^p\ge k\}} |z|^p m_\alpha(\textup{d}z) < d_n\Big\}, \text{ for every }n\in \mathbb{N}.
  \end{align*}
  In view of \eqref{Uniform_Integr_Limit_zero} it is immediate that for each $n\in\mathbb{N}$ the set $\mathcal{N}_n$ contains all but finitely many integers and that $(\mathcal{N}_n)_{n\in\mathbb{N}}$ is a decreasing sequence of sets.
%
%
  Let $n\in\mathbb{N}$ and $c_n \in \mathcal{N}_n$.
  Then, observe that
  \begin{align*}
    \sum_{k=c_n}^{\infty} m_\alpha(\{z:|z|^p \ge k\}) 
    &\le c_n m_\alpha (\{z: |z|^p\ge c_n\}) + \sum_{k=c_n+1}^{\infty} m_\alpha(\{z: |z|^p\ge k\})\\
    &= \sum_{k=c_n}^{\infty} k m_\alpha(\{z: k\le |z|^p <k+1\})\\
    &\le \int_{\{z:|z|^p\ge c_n\}} |z|^p m_\alpha(\textup{d}z)<d_n.
  \end{align*}
In other words, by applying a diagonal argument, we can choose an increasing sequence $(c_n)_{n\in\mathbb{N}}$,\footnote{By increasing we mean $c_{n}<c_{n+1}$ for every $n\in\mathbb{N}$. Hence, $c_n\ge n$ and consequently $\lim_{n\to+\infty} c_n=+\infty$.} such that $c_n\in\mathcal{N}_n$ for every $n\in\mathbb{N}$,  in order to obtain
 \begin{align*}
  \sup_{\alpha \in\mathcal{A}}\sum_{n=1}^{\infty} \sum_{k=c_n}^{\infty} m_\alpha(\{z:|z|^p \ge k\}) <\sum_{n=1}^{\infty} d_n <+\infty.
 \end{align*}
 It is only left to observe that the double series on the left-hand side of the above inequality can be written as a series of the form 
 $\sum_{k=1}^{\infty} \upsilon_k m_\alpha(\{z:|z|^p \ge k\})$, where 
 \begin{align}
   \upsilon_k := |\{n\in\mathbb{N} : c_n\le k\}|,\footnotemark \text{ for every }k\in\mathbb{N};
   \label{def:upsilon_n}
 \end{align}%
 \footnotetext{\label{footnote:comparison_of_upsilon_n}Since we chose the sequence $(c_n)_{n\in\mathbb{N}}$ to be increasing, we immediately have that $\upsilon_n\le \upsilon_{n+1}\le \upsilon_n +1$, for every $n\in\mathbb{N}$, where exactly one of the comparisons is a true equality.}
  $|A|$ denotes the cardinality of the set $A$.
%
%
 The equivalence has been proven. 
\end{proof}	
\begin{remark}\label{property:c_n_choice}
In this remark we will underline a property of $(\upsilon_k)_{k\in\mathbb{N}}$ and argue about the freedom one has in choosing the sequence $(c_n)_{n\in\mathbb{N}}$ in the proof of \cref{prop:UI_Young}.
Both of them are crucial in the construction implemented in the proof of \cref{lemma:UI_Young_improvement}.
To this end, let us follow the notation of the proof of \cref{prop:UI_Young}.
\begin{enumerate}
  \item 
  From the definition \eqref{def:upsilon_n} it is clear that
  \begin{align}\label{property:upsilon_inequalities}
    \upsilon_k \le \upsilon_{k+1} \le \upsilon_{k} + 1, \text{ for every }k\in\mathbb{N},
  \end{align}
  where exactly one of the comparisons is a true equality; hence, the other is a strict inequality. 
  For $k\in\mathbb{N}$, the right-hand equality is attained if, and only if,  $k=c_n$ for some $n\in\mathbb{N}$, recall \eqref{def:upsilon_n}, \emph{i.e.},
  \begin{gather*}
  \upsilon_{c_n-1}<\upsilon_{c_n}=\upsilon_{c_n-1}+1 \text{ for every }n\in\mathbb{N}\backslash\{1\}.
  \numberthis \label{property:upsilon_c_RHSEquality}
  \shortintertext{as well as}
  \upsilon_{k-1}=\upsilon_{k}<\upsilon_{k-1}+1 \text{ for every }k\in \mathbb{N} \backslash \{ c_n : n\in\mathbb{N}\}.
  \end{gather*}
  In other words,
  \begin{align}
  \upsilon_{c_n} = n, \text{ for all }n\in\mathbb{N}.
  \label{property:upsilon_index_cn_is_value}
  \end{align}
  \item
By definition, the set $\mathcal{N}_{n}$ contains all but finitely many integers, for every $n\in\mathbb{N}$, and the sequence $(\mathcal{N}_n)_{n\in\mathbb{N}}$ is decreasing, because the sequence $(d_n)_{n\in\mathbb{N}}$ was chosen decreasing. 
In particular,  for every $n\in\mathbb{N}$, if $m\in\mathcal{N}_n$, then $m+1\in\mathcal{N}_n$.
These properties provide ample freedom in the subsequent choice of the sequence $(c_n)_{n\in\mathbb{N}}$, which in turn determines the sequence $(\upsilon_k)_{k\in\mathbb{N}}$.
Without loss of generality, we may assume that $1\in\mathcal{N}_1$, so that we can choose $c_1$ equal to $1$.\footnote{This can be done by increasing suitably $d_1$, which does not affect the finiteness of the series associated to $(d_n)_{n\in\mathbb{N}}$. Indeed, we can assign to $d_1$ the value $2\sup_{\alpha\in\mathcal{A}}\int_{\{z:|z|\ge 1\}} |z|^p m_\alpha(\textup{d}z)$ so that $\mathcal{N}_1=\mathbb{N}$} 
For every $n\in\mathbb{N}$, a natural rule for choosing $c_{n+1}$, given the determination of $c_n$, is 
\begin{align*} 
c_{n+1}:=\max\{ c_n+1, \min\mathcal{N}_{n+1}\} \text{ for }n\in\mathbb{N}.
\end{align*}
Our later needs dictate to modify the aforementioned rule: for $n\in\mathbb{N}$
\begin{align} 
\label{property:c_doubling}
c_{n+1}:=\max\{ 2c_n, \min\mathcal{N}_{n+1}\}.
\end{align}
\end{enumerate}
\subsection{Proof of Lemma \ref{lemma:UI_Young_improvement}}

      In order to prove \cref{lemma:UI_Young_improvement}, we will suitably choose the right-derivative of the Young function $\Upsilon$, denoted by $\upsilon$, based on the behaviour of the family $\{m_\alpha\}_{\alpha\in\mathcal{A}}$.
      In the next lines we will try to explain intuitively the procedure.
      To this end, we follow the notation introduced in the proof of \cref{prop:UI_Young}.

      In \citet[Lemme, p. 770]{meyer1978sur} it is described how to alter the sequence $(\upsilon_n)_{n\in\mathbb{N}}$ in order to obtain the (seemingly) stronger integrability described in \cref{corollary:UI_Young_improvement:StrongerIntegrability} and to obtain a moderate Young function.
      For the former, it is chosen a new sequence which is dominated by the sequence $(\upsilon_n)_{n\in\mathbb{N}}$, \emph{i.e.}, the right derivative is ``pushed downwards''.
      For the latter, the newly chosen sequence applies to intervals of increasing length, \emph{i.e.}, the right derivative is ``stretched''. 
      Inspired by \citet[Lemme, p. 770]{meyer1978sur}, we will initially follow the arguments presented there, but we will need go beyond \citet{meyer1978sur} in order to choose more conveniently some values. 
      Following the intuition introduced in the previous lines, the idea for proving \cref{lemma:UI_Young_improvement} is to repeatedly alter the right derivative of the Young function by initially ``properly pushing it downwards'', afterwards ``stretching it'', then making it continuous and finally smoothing it.

\begin{proof}[Proof of \cref{lemma:UI_Young_improvement}]
  In view of 
  the comments provided at the beginning of this subsection,
  the current proof should be seen as a continuation of the proof of \cref{prop:UI_Young}.
  To this end, we will assume that the sequence $(\upsilon_n)_{n\in\mathbb{N}}$ is the one described in the second half of the proof of \cref{prop:UI_Young}, where the sequence $(c_n)_{n\in\mathbb{N}}$ obeys to \eqref{property:c_doubling} with $c_1=1$.
  We choose to present the construction of the (right) derivative step by step, so that it is intuitively clear to the reader the final choice of the (right) derivative.
  As a result, we sacrifice the compactness of the proof.   

\vspace{0.5em}
 Following \citet{meyer1978sur}, we ``push downwards'' the right derivative by  passing to a new sequence $(\bar{\upsilon}_n)_{n\in\mathbb{N}}$ which is defined by $\bar{\upsilon}_n:=\vartheta(\upsilon_n)$, for every $n\in\mathbb{N}$, where $\vartheta:[0,\infty) \to [0,\infty)$ is increasing and such that 
 \begin{enumerate}[label=(\alph*)]
 \item\label{vartheta:a} $\vartheta(1)=1$,
 \item\label{vartheta:b} $\lim_{x\to\infty}\frac{\vartheta(x)}{x}=0$ as well as $\lim_{x\to\infty}\vartheta(x)=\infty$ and 
 \item\label{vartheta:c} the function $[1,\infty) \ni x \longmapsto \vartheta(x+1) - \vartheta(x) \in (0,\infty)$ is (strictly) decreasing and bounded by $\frac{1}{2}$,%
 \footnote{In \citet{meyer1978sur}, it was chosen the square root function. 
 The freedom in choosing the function $\vartheta$ allows to endow the Young function $\Upsilon$ with convenient properties. 
 For example, passing to the sequence $(\bar{\upsilon}_n)_{n\in\mathbb{N}}$ will later allow for the property $\lim_{R\to\infty} \sup_{\alpha\in\mathcal{A}} \int_{\{z:|z|>R\}}\Upsilon(|z|^p) m_\alpha(\textup{d}z)=0,$ \emph{i.e.}, the uniform integrability of the function $\Upsilon(|z|^p)$ with respect to the family $(m_\alpha)_{\alpha\in\mathcal{A}}$. 
 One can verify it in the case $\vartheta$ is the square root function, by checking the proof of \cite[Lemme]{meyer1978sur} and modifying the arguments properly, \emph{i.e.}, as we did in the proof of the equivalence of \cref{prop:UI_Young}. 
 For completeness, we provide the full details in the proof of \cref{corollary:UI_Young_improvement:StrongerIntegrability}. 
 Other possible choices for the function $\vartheta$ are $x\longmapsto x^{\frac12(1-\alpha)}$, for $\alpha \in (0,1)$, or $x\longmapsto x^{\frac{1}{2}(2-\alpha)}$, for $\alpha \in [1,2)$. 
 } 
\emph{i.e.},
\begin{align}\label{property:theta_difference} 
\vartheta(x+1) - \vartheta(x)<\frac{1}{2}, \text{ for every }x\in[1,\infty).
\end{align} 
 \end{enumerate}
 In view of \eqref{property:upsilon_inequalities} and the last required property of the function $\vartheta$, we get 
 \begin{gather*}
    \bar{\upsilon}_n \le \bar{\upsilon}_{n+1}< \bar{\upsilon}_n+\frac{1}{2}, \text{ for every }n\in\mathbb{N},
    \shortintertext{as well as}
    \bar{\upsilon}_{c_n}= \vartheta(\upsilon_{c_n}) \overset{\eqref{property:upsilon_index_cn_is_value}}{=} \vartheta(n) \text{ for every }n\in\mathbb{N}.
    \numberthis\label{property:bar_upsilon_index_cn_is_theta_value}
 \end{gather*}
 Then, abusing notation, we define the function $\bar{\upsilon}$
  \begin{gather*}
    \bar{\upsilon}(x) := 
      \begin{cases}
        \frac{1}{2},        &\text{ for }x\in[0,2)\\
        \bar{\upsilon}_n,   &\text{ for }x\in[2^{n},2^{n+1}), n\in\mathbb{N}
      \end{cases}.
  \end{gather*}
 Obviously, the function $\bar{\upsilon}$ is right-continuous and non-decreasing.
 Recalling the description in 
 the comments provided at the beginning of this subsection,
 the above definition amounts to ``stretching'' the right derivative.
 For the sake of unified notation, we define $\bar{\upsilon}_0:=\frac{1}{2}$, for which choice $\bar{\upsilon}_0<\bar{\upsilon}_1=\bar{\upsilon}_0+\frac{1}{2}$.
 Also, $c_0:=0$ and we will write $\mathbb{N}\cup\{0\}$ whenever we include $0$ in the natural numbers.
Given the set $(c_n)_{n\in\mathbb{N}}$ and recalling \eqref{property:upsilon_c_RHSEquality}, we can alternatively write the function $\bar{\upsilon}$ as
 \begin{gather*}
    \bar{\upsilon}(x) := 
      \begin{cases}
        \bar{\upsilon}_{0}(=\frac12),        &\text{ for }x\in[0,2)\\
        \bar{\upsilon}_{c_n}(=\vartheta(n)),   &\text{ for }x\in[2^{c_{n}},2^{c_{n+1}}), \ n\in\mathbb{N}
      \end{cases}.
  \end{gather*}

 We proceed to construct a continuous function which is dominated by $\bar{\upsilon}$.
 To this end, we will linearly interpolate between the right-end points of intervals of constancy of $\bar{\upsilon}$, which are determined by consecutive elements of the sequence $(c_n)_{n\in\mathbb{N}}$ as described above.
 We define now the function
 \begin{align*}
   \hat{\upsilon}(x) :=
    \begin{cases}
      \frac{1}{4}x,                             &\text{for }x\in[0,2)\\
      \bar{\upsilon}_{c_{n-1}} +\frac{ \bar{\upsilon}_{c_{n}}-\bar{\upsilon}_{c_{n-1}}}{2^{c_{n+1}}-2^{c_n}}(x - 2^{c_n}), 
                                                    &\text{for }x\in[2^{c_{n}},2^{c_{n+1}}), \ n\in\mathbb{N}
    \end{cases}.
    \numberthis\label{def:hat_upsilon_piecewise_linear}
 \end{align*}
Some remarks are in order:
\begin{enumerate}[label=(\arabic*)]
  \item The function $\hat{\upsilon}$ is continuous
  and concave on $[0+\infty)$. 
  This can easily be seen, since $\hat{\upsilon}$ is piecewise linear with 
  decreasing slopes on $(0,+\infty)$. 
  Indeed, one can verify that the slopes are decreasing by examining the sequence of nominators and denominators.
  More precisely, we have (by definition) that
  \begin{align*}
       \bar{\upsilon}_{c_{n+1}} \overset{\eqref{property:bar_upsilon_index_cn_is_theta_value}}{=} \vartheta({\upsilon_{c_{n+1}}})=\vartheta(n+1) 
       \overset{\eqref{property:upsilon_c_RHSEquality}}{=}\vartheta({\upsilon_{c_{n}}+1}), \text{ for every }n\in\mathbb{N}.  
  \end{align*}
  Given that the function $[1,+\infty)\ni x\mapsto \vartheta({x+1})-\vartheta({x})$ is decreasing and bounded between $0$ and $\frac{1}{2}$, we have that 
  \begin{align}\label{property:diff_bar_upsilon}
    0<\bar{\upsilon}_{c_{n+2}}-\bar{\upsilon}_{c_{n+1}} 
    < \bar{\upsilon}_{c_{n+1}} - \bar{\upsilon}_{c_{n}}<\frac{1}{2} = \bar{\upsilon}_1- \bar{\upsilon}_0, \text{ for every }n\in\mathbb{N}.
  \end{align}
  Moreover, for every $n\in\mathbb{N}$,
  \begin{gather*}
     2^{c_2}-2^{c_1}  > 2^{c_1} - 2
  \shortintertext{as well as}
    2^{c_{n+2}} - 2^{c_{n+1}} 
      = 2^{c_{n+1}}(2^{c_{n+2}-c_{n+1}} - 1) 
      > 2^{c_{n+1}}(1 - 2^{c_{n}-c_{n+1}})
      =2^{c_{n+1}} - 2^{c_{n}}\geq 2.
  \end{gather*}
  Hence,  the sequence $(\frac{ \bar{\upsilon}_{c_{n}}-\bar{\upsilon}_{c_{n-1}}}{2^{c_{n+1}}-2^{c_n}})_{n\in\mathbb{N}}$ is decreasing, since the nominators decrease and the denominators increase, as $n$ increases.
  Additionally, 
  \begin{align*}
  \frac{ \bar{\upsilon}_1 - \bar{\upsilon}_0}{2^{c_2}-2^{c_1}}\le \frac{1}{4},
  \end{align*}
  which completes the verification of the claim that $\hat{\upsilon}$ has decreasing slopes.

  In total, 
  \begin{align}
  \sup_{n\in\mathbb{N}}\big\{\frac{ \bar{\upsilon}_{c_{n}}-\bar{\upsilon}_{c_{n-1}}}{2^{c_{n+1}}-2^{c_n}}
                          -\frac{ \bar{\upsilon}_{c_{n+1}}-\bar{\upsilon}_{c_{n}}}{2^{c_{n+2}}-2^{c_{n+1}}} \big\}
  \leq \frac{1}{4}.
  \label{upper_bound:slopes}
  \end{align}
\end{enumerate}

\vspace{0.5em}
We proceed to define the derivative of the desired Young function, but we need to treat separately the cases $p\ge 2$ and $p\in(0,2)$:

\vspace{0.5em}
{\tiny $\blacksquare$} $p\ge 2$. We will smooth the function $\hat{\upsilon}$\footnote{The function $\hat{\upsilon}$ is continuous, but not everywhere differentiable. Hence, the associated Young function will not be twice continuously differentiable.}  by standard tools.
More precisely, for $0<\varepsilon<\frac{1}{2^3}$ we use the standard mollifier $\eta_{\varepsilon}$, \emph{i.e.},
\begin{align*}
 \eta_\varepsilon(x):= \frac{C}{\varepsilon} \exp\Big(\frac{1}{\varepsilon^{-2}x^2-1}\Big) \mathds{1}_{[-\varepsilon,\varepsilon]}(x),
\end{align*}
where $C$ is the normalizing constant.
We define
  \begin{align}\label{def:hat_upsilon_varepsilon_smooth}
   \hat{\upsilon}_\varepsilon(x) :=
   \begin{cases}
     \hat{\upsilon}(x), &\text{ for }x\in[0,1]\\
     \displaystyle \int_{[-\varepsilon, \varepsilon]} \hat{\upsilon} (x+y) \eta_\varepsilon(y)\textup{d}y, &\text{ for } x> 1.
   \end{cases}
  \end{align}
Hereinafter, we \emph{fix} an $0<\varepsilon<\frac{1}{2^3}$ such that
\begin{align}\label{property:error_bound}
 0\leq  \hat{\upsilon}(x) - \hat{\upsilon}_\varepsilon(x)<\frac{1}{2^4}, \text{ for every }x\in[0,+\infty).
\end{align}
We present in \cref{comput:error_hat_u} the elementary arguments for the existence of an $\varepsilon$ with the aforementioned property, despite this can be seen as a special case of a general approximation result.
The function $\hat{\upsilon}_\varepsilon$ (for the $\varepsilon$ we fixed) will be the derivative of the desired Young function $\Upsilon$.
Before we prove that $\Upsilon$ possesses the desired properties, some remarks are in order:
\begin{enumerate}[label=(\arabic*)]
\setcounter{enumi}{1}
 \item $\hat{\upsilon}_{\varepsilon}(0)=0$ and $\hat{\upsilon}_{\varepsilon}(x)>0$ for $x>0$. 

 \item\label{property:identical_on_linear} For $x\in(1,2-\varepsilon)\cup(2+\varepsilon,2^{c_2}-\varepsilon) \cup \big(\bigcup_{k\in\mathbb{N}} (2^{c_n}+\varepsilon, 2^{c_{n+1}}-\varepsilon)\big)$ 
 it is true that $\hat{\upsilon}_\varepsilon(x)=\hat{\upsilon}(x)$, because of the linearity of $\hat{\upsilon}$ on the interval $(x-\varepsilon,x+\varepsilon)$.

\item\label{property:concavity_domination}
On $(0,+\infty)$ the function $\hat{\upsilon}_\varepsilon$ is concave, with $\hat{\upsilon}_\varepsilon \le \hat{\upsilon}$.
The concavity of $\hat{\upsilon}_\varepsilon$ is inherited from the concavity of $\hat{\upsilon}$. 
The domination can be easily proven using the piecewise linearity of $\hat{\upsilon}_\varepsilon$, the fact that the slopes are non-increasing and the definition of the mollification.\footnote{In \cref{comput:error_hat_u} we prove that the difference $\hat{\upsilon}(x) -\hat{\upsilon}_\varepsilon(x)$ is non-negative everywhere.}
In particular, 
\begin{align}\label{bound:slopes_of_smooth_upsilon}
\hat{\upsilon}'_{\varepsilon}(2^{c_n}) \le 
\begin{cases}
\frac14, &\text{for }n=1\\
\frac{\bar{\upsilon}_{c_{n-1}}-\bar{\upsilon}_{c_{n-2}}}{2^{c_n}-2^{c_{n-1}}}, &\text{for }n>1
\end{cases}
\quad\text{ and }\quad
\hat{\upsilon}'_{\varepsilon}(2^{c_n}+\varepsilon) = \frac{\bar{\upsilon}_{c_{n}}-\bar{\upsilon}_{c_{n-1}}}{2^{c_{n+1}}-2^{c_{n}}}
.
\end{align}
\item\label{property:second_derivative_finally_decreasing} The function $\hat{\upsilon}_\varepsilon$ is infinitely differentiable.
In particular, its first derivative is positive and non-increasing due to concavity. 
\end{enumerate}
\vspace{0.5em}
The function 
\begin{align*}
  \Upsilon(x) := \int_{[0,x]} \hat{\upsilon}_\varepsilon (u) \textup{d}u, \text{ for }x\ge 0,
\end{align*}
satisfies the claimed properties: 
\begin{enumerate}
  \item[\ref{lemma:UI_Young_improvement:moderate_C2}] Since $\hat{\upsilon}_\varepsilon$ is infinitely differentiable, so is $\Upsilon$. 
  In particular, $\Upsilon$ is twice continuously differentiable.

  \vspace{0.5em}
  We prove, now, that $\Upsilon$ is moderate, \emph{i.e.}, $\bar{c}_{\Upsilon}< +\infty.$
  We will avoid explicit computations and we will make use of the aforementioned properties for $\hat{\upsilon}_{\varepsilon}$ and of known results from the theory of convex functions.
  More precisely, the concavity of $\hat{\upsilon}_{\varepsilon}$ in conjunction with the fact that $\hat{\upsilon}_{\varepsilon}(0)=0$ implies its subadditivity.
  In particular,
  \begin{align*}
  \hat{\upsilon}_{\varepsilon} (2x)\leq 2 \hat{\upsilon}_{\varepsilon}(x), \text{ for every }x>0.
  \numberthis\label{property:hat_upsilon_epsilon_moderate}
  \end{align*}
  In turn, this implies
\begin{align*}
  x\hat{\upsilon}_\varepsilon(x) 
    &\le \int_{[x,2x]}\hat{\upsilon}_\varepsilon(z)\textup{d}z 
      = \int_{[0,2x]}\hat{\upsilon}_\varepsilon(z)\textup{d}z - \int_{[0,x]}\hat{\upsilon}_\varepsilon(z)\textup{d}z\\
    & = \int_{[0,x]}2\hat{\upsilon}_\varepsilon(2w)\textup{d}w - \int_{[0,x]}\hat{\upsilon}_\varepsilon(z)\textup{d}z
      \overset{\eqref{property:hat_upsilon_epsilon_moderate}}{\le} \int_{[0,x]}4\hat{\upsilon}_\varepsilon(w)\textup{d}w - \int_{[0,x]}\hat{\upsilon}_\varepsilon(z)\textup{d}z\\
    &=3\Upsilon(x),
\end{align*}
  where in the first inequality we used that $\hat{\upsilon}_\varepsilon$ is increasing.  
In other words,
\begin{align*}
\sup_{x>0}\frac{x\hat{\upsilon}_\varepsilon(x)}{\Upsilon(x)}\leq 3.
\end{align*}
   \item[\ref{lemma:UI_Young_improvement:second_derivative_properties}] Recall from \ref{property:second_derivative_finally_decreasing} that $\hat{\upsilon}_\varepsilon'=\Upsilon''$ is continuous, positive and non-increasing.
   Moreover, $\hat{\upsilon}_\varepsilon'(x)=\frac{1}{4}$, for $x\in(0,2)$. 
   Therefore, $\Upsilon''$ is bounded.
   
   \vspace{0.5em}
   We proceed to prove \eqref{property:second_derivative_faster_than_identity}.
   The reader may recall that the value of $\varepsilon$ is fixed and smaller that $\frac{1}{2^3}$.
   We will consider the following cases:
   \begin{enumerate}[label=\textbullet]
      \item 
      $x\in[0,2]$: then, immediately $x\Upsilon''(x) \le \frac{1}{2}$.



      \item 
       $x\in(2^{c_n},2^{c_n}+\varepsilon]$, for $n\in\mathbb{N}$:
        then,  recalling that $\varepsilon$ is fixed and lies in $(0,\frac{1}{2^3})$,
        \begin{align*}
          x\Upsilon''(x)
            \overset{\text{\ref{property:second_derivative_finally_decreasing}}}{<} 
            (2^{c_n} + \varepsilon) \Upsilon''(2^{c_n})
            &\overset{\eqref{bound:slopes_of_smooth_upsilon}}{\le} 
            \begin{cases}
            \frac{2 + \varepsilon}{4}, &\text{ for }n=1 \\
            (2^{c_n} + \varepsilon) \frac{\bar{\upsilon}_{c_{n-1}}-\bar{\upsilon}_{c_{n-2}}}{2^{c_n}-2^{c_{n-1}}},  &\text{ for }n>1 
            \end{cases}\\
            &\overset{\eqref{property:diff_bar_upsilon}}{<}
            \begin{cases}
            1, &\text{ for }n=1 \\
             \frac{1}{2(1-2^{c_{n-1}-c_n})} + \frac{1}{2^5}   
            ,&\text{ for }n>1 
            \end{cases}\\
            &\overset{\phantom{\eqref{property:diff_bar_upsilon}}}{\leq}
            2,
        \end{align*}
        since $c_n-c_{n-1}\ge 1$.

     \item 
        $x\in(2^{c_n}+\varepsilon,2^{c_{n+1}}]$, for $n\in\mathbb{N}$:
        then,
        \begin{align*}
          x\Upsilon''(x) 
            \overset{\text{\ref{property:second_derivative_finally_decreasing}}}{<} 
            2^{c_{n+1}}\Upsilon''(2^{c_{n}}+\varepsilon)
            \overset{\eqref{bound:slopes_of_smooth_upsilon}}{=}  
            2^{c_{n+1}}\frac{\bar{\upsilon}_{c_{n}}-\bar{\upsilon}_{c_{n-1}}}{2^{c_{n+1}}-2^{c_{n}}}
            \overset{\eqref{property:diff_bar_upsilon}}{<}\frac{1}{2(1-2^{c_n-c_{n+1}})}\le 1,
        \end{align*}
        since $c_{n+1}-c_n\ge 1$.
    In total, a constant $C_{\Upsilon,p}$ for the case $p\ge 2$ is $2$.
   \end{enumerate}
    \item[\ref{lemma:UI_Young_improvement:submultiplicativity}] 
    In view of \citet[Proposition 1]{HUDZIK1992313},\footnote{In \cite{HUDZIK1992313} the property of being ``finally submultiplicative'' is described by the term ``submultiplicative at infinity''.} it is sufficient to prove the desired property for some $R>1$, instead of $R=1$, \emph{i.e.},
    \begin{align*}
    \Upsilon(xy)\leq C \Upsilon(x) \Upsilon(y)\text{ for every }x,y\geq R.
    \end{align*}
    In order to prove this claim we will detour and we will prove the analogous claim for the Young function, denoted by $\widehat{\Upsilon}$, with  derivative $\hat{\upsilon}$ defined in \eqref{def:hat_upsilon_piecewise_linear}. 
    This detour is indeed valid in view of \cref{lemma:Young_equiv_submult}.
    From \citet[Proposition 2]{HUDZIK1992313}, $\widehat{\Upsilon}$ is finally submultiplicative if and only if $\hat{\upsilon}$ is finally submultiplicative, \emph{i.e.}, there exist $C,R>0$ such that
    \begin{align}\label{ineq:submult_derivative}
    \hat{\upsilon}(xy)\leq C \hat{\upsilon}(x) \hat{\upsilon}(y)\text{ for every }x,y\geq R.
    \end{align}
    We proceed to prove the equivalent statement. 
    Let $x,y\geq2^{c_2}$. 
    Then, there exist $c_{n_x},c_{n_y}$ with ${n_x},{n_y}\geq 2$ such that 
    \begin{align*}
    x\in [2^{c_{n_x}},2^{c_{n_x+1}}) 
    \text{ and }
    y\in [2^{c_{n_y}},2^{c_{n_y+1}}). 
    \end{align*}
    Now, it is sufficient to prove the existence of a $C>0$ such that 
    \begin{align*} 
    \hat{\upsilon}(2^{c_{n_x+1}+c_{n_y+1}}) \leq C \hat{\upsilon}(2^{c_{n_x}}) \hat{\upsilon}(2^{c_{n_y}}). 
    \end{align*} 
    Indeed, the fact that $\hat\upsilon$ is increasing, allows us to conclude \eqref{ineq:submult_derivative}. 
    Assuming that $c_{n_x}\leq c_{n_y}$, one gets
    \begin{align*}
    c_{n_x+1}+c_{n_y+1} \leq 2c_{n_y+1}  \overset{\eqref{property:c_doubling}}{\leq} c_{n_y+2}.
    \end{align*}
    Hence, since $\hat{\upsilon}$ is increasing,
    \begin{align*}
    \hat{\upsilon}(2^{c_{n_x+1}+c_{n_y+1}}) 
    &\overset{\phantom{\eqref{property:theta_difference}}}{\leq} \hat{\upsilon}(2^{c_{n_y+2}})
    \overset{\eqref{def:hat_upsilon_piecewise_linear}}{=}\bar{\upsilon}(c_{n_y+1}) 
    \overset{\eqref{property:bar_upsilon_index_cn_is_theta_value}}{=} \vartheta(n_y+1)\\
    &\overset{\eqref{property:theta_difference}}{<} 1+\vartheta(n_y-1)
    \leq 2\vartheta(n_x-1)\vartheta(n_y-1)\\
    &\overset{\eqref{property:bar_upsilon_index_cn_is_theta_value}}{=} 2\bar{\upsilon}(c_{n_x-1})\bar{\upsilon}(c_{n_y-1})
    \overset{\eqref{def:hat_upsilon_piecewise_linear}}{=} 2\hat{\upsilon}(2^{c_{n_x}})\hat{\upsilon}(2^{c_{n_y}}).
    \end{align*}
    where in the third inequality it was used the fact that $\vartheta(n)\ge1$ for every $n\in\mathbb{N}$.
    Hence, \eqref{ineq:submult_derivative} is true and we conclude in view of the arguments presented above.
\end{enumerate}

{\tiny $\blacksquare$} $p\in(0,2)$. Before we smooth the function $\hat{\upsilon}$, we will need to suitably modify it on $[0,2]$.
We have in mind that \eqref{property:second_derivative_faster_to_0} has to hold. 
To this end, we will choose the derivative $\hat{\upsilon}$ to behave in a neighbourhood of $0$ as a power function, say of power $q$, so that 
\begin{align*}
  \lim_{x\downarrow 0}(x^p)^qx^{p-2}=\lim_{x\downarrow 0}(x^p)^{q-1}x^{2p-2}=0.
\end{align*}
\cref{lem:def_q_for_smoothness} verifies that there is a legitimate choice of $q$ for each $p\in(0,2)$:
\begin{align}\label{def:q}
  q:=\frac{4}{p}-1>0, \text{ for }p\in(0,2).
\end{align}
Observe that for the choice of $q$ it holds $pq+p-2=p(q-1)+2p-2=2$.
Additionally, from \cref{lem:def_q_for_smoothness}, for the point 
\begin{align*}
   x_0(p):=\big(\frac{p}{6(4-p)}\big)^{\frac{p}{4-2p}}\in(0,1)
\end{align*}
and for $(0,\infty)\ni x\overset{f}{\mapsto} x^q$, it holds $f'(x_0)=\frac{1}{6}$.
Moreover, it is true that 
\begin{gather*}
  0<f(x_0) = x_0^q 
    = \Big(\frac{1}{6q}\Big)^{\frac{q}{q-1}} 
    =\Big(\frac{p}{6(4-p)}\Big)^{\frac{4-p}{4-2p}}
    <\frac{p}{6(4-p)}
    < \frac{1}{6},
\end{gather*}
where we used that $\frac{p}{4-p}<1$ and that $\frac{4-p}{4-2p}>1$, as well as that
\begin{gather*}
  0<f(x_0) + \frac{1}{6} (\frac{5}{4}-x_0)
    <\frac{1}{6}+\frac{1}{6} \frac{5}{4} = \frac{1}{6} \frac{9}{4} = \frac{3}{8}<\frac{1}{2}=\bar{\upsilon}_0.
\end{gather*}
The above properties permit us to \emph{re}-define the function $\hat{\upsilon}$ on $[0,2]$ as follows
\begin{align}
  \hat{\upsilon}(x):=
  \begin{cases}
    f(x), &\text{for }x\in[0,x_0]\\
    f(x_0) + \frac{1}{6} (x-x_0), &\text{for }x\in(x_0,\frac{5}{4}]\\ 
    \hat{\upsilon}(\frac{5}{4}) + \frac{4}{3}\big[\bar{\upsilon}_0 - \hat{\upsilon}(\frac{5}{4})\big](x-\frac{5}{4}), &\text{for }x\in(\frac{5}{4},2]\\ 
  \end{cases}.
  \label{def:redefined_hat_upsilon}
\end{align}
The reader may observe that on the interval $[x_0,\infty)$ the redefined function $\hat{\upsilon}$ is continuous and piecewise linear.
In view of the computations provided before \eqref{def:redefined_hat_upsilon}, it can be immediately verified that the redefined function $\hat{\upsilon}$ is positive, increasing, continuous on $[0,+\infty)$ and continuously differentiable on $(0,\frac{5}{4})$.
For the continuity of the derivative $\hat{\upsilon}'$ at $x_0$ we have used that the right and left derivatives of $\hat{\upsilon}$ at $x_0$ are equal.
For $0<\varepsilon<\frac{1}{2^3}$ we define
\begin{align}
\hat{\upsilon}_\varepsilon(x) :=
  \begin{cases}
    \hat{\upsilon}(x), &\text{for } x\in [0,\frac{5}{4}-\varepsilon]\\
     \displaystyle \int_{[-\varepsilon, \varepsilon]} \hat{\upsilon} (x+y) \eta_\varepsilon(y)\textup{d}y, &\text{ for } x> \frac{5}{4}-\varepsilon
  \end{cases}.
  \label{def:redefined_hat_upsilon_smoothed}
\end{align}
The function 
\begin{align*}
  \Upsilon(x) := \int_{[0,x]} \hat{\upsilon}_\varepsilon (u) \textup{d}u, \text{ for }x\ge 0,
\end{align*}
satisfies the claimed properties:\footnote{We will focus our attention on the interval $[0,2^{c_2}]$, since for the rest we can use the analogous arguments presented for the case $p\ge 2$.}
\begin{enumerate}
    \item[\ref{lemma:UI_Young_improvement:moderate_C2}]
    The function $\hat{\upsilon}_\varepsilon$ is continuously differentiable at every point of $(0,x_0]$ and infinitely differentiable at every point of $(x_0,2^{c_2}]$.
    Hence, $\Upsilon$ is twice differentiable.

    \vspace{0.5em}
    In order to prove that $\Upsilon$ is moderate, we cannot work analogously to the case $p\geq 2$, since the function $\hat{\upsilon}_{\varepsilon}$ is only concave on $[x_0,\infty)$.
    Nevertheless, we can use the fact that the Young function constructed for the case $p\geq 2$ is finally submultiplicative in conjunction with the fact that the derivatives for the two cases are identical on $[2+\varepsilon,\infty)$.
    Hence, the Young function with derivative \eqref{def:redefined_hat_upsilon_smoothed} is also finally submultiplicative on the interval $[2^{c_2},\infty)$, as the arguments can be applied verbatim. 
    From \cite[Proposition 1]{HUDZIK1992313} we can assume the property to holds on $[1,\infty)$ and the property \ref{lemma:UI_Young_improvement:submultiplicativity} has been proven.
    This, in turn, implies that $\Upsilon$ is moderate on $2^{c_2}+2$ with constant $2 \Upsilon(2^{c_2}+2)$, see \citet[Chapter I, Lemma 5.1, p.30]{krasnoselskii1961}, where $2$ is the constant appearing in the submultiplicative property.
    From \cite[Chapter I, Theorem 4.1]{krasnoselskii1961} we have that 
    \begin{align*} 
      \sup_{x\geq 2^{c_2}+2} \frac{x\hat{\upsilon}_{\varepsilon}(x)}{\Upsilon(x)} \leq 2 \Upsilon(2^{c_2}+2)<\infty.
    \end{align*} 
    We consider now separately the intervals $(0,x_0]$ and $[x_0,2^{c_2}+2]$.
    For the former interval, immediate computations verify
    \begin{align*}
      \sup_{x\in(0,x_0]}\frac{x\hat{\upsilon}_\varepsilon(x)}{\Upsilon(x)} = q+1<+\infty.
    \end{align*}
    For the latter interval, we use the continuity of $x\mapsto x\hat{\upsilon}_\varepsilon$ and $\Upsilon$, the compactness of the interval under consideration as well as the fact that $0<\inf_{x\in[x_0,2^{c_2}]}\Upsilon(x),\sup_{x\in[x_0,2^{c_2}]}\Upsilon(x)<+\infty$.
    In total,
    \begin{align*}
        \bar{c}_{\Upsilon} =\sup_{x>0} \frac{x\hat{\upsilon}_\varepsilon(x)}{\Upsilon(x)}<+\infty.
    \end{align*}
    \item[\ref{lemma:UI_Young_improvement:second_derivative_properties}] 
    We can readily adapt the arguments we used in the case $p\ge 2$.
    By definition, $\hat{\upsilon}'$ is increasing on $(0,x_0]$, but on $[x_0,\infty)$ is non-increasing.
    It is also clear by its definition that it remains positive and by continuity and the described monotonicity, it remains bounded.
    Finally, the term $x \hat{\upsilon}'(x)$ remains bounded on $[0,2]$; for the interval $[2,\infty]$ we can follow verbatim the arguments of the case $p\geq 2$. 
%
%
%
%
%
%
%
%
%
%
%
    \item[\ref{lemma:UI_Young_improvement:submultiplicativity}] 
    The arguments have been presented in the course of proving that $\Upsilon$ is moderate; see the proof of \ref{lemma:UI_Young_improvement:moderate_C2}.
%
%
%
%
%
%
%
%
%
%
%
    \item[\ref{lemma:UI_Young_improvement:p_smaller_than_2}]
    This property is true because of the choice of $q$; see \eqref{def:q} and the comment following it.
\end{enumerate}
\begin{enumerate}
  \item[\ref{lemma:UI_Young_improvement:p_between_1_and_2}]
  Let $p\in[1,2)$. 
  For the desired property to hold, we may alternatively prove that the function 
  \begin{align*}
  [2^{c_2},\infty) \ni y \overset{h}{\longmapsto} \hat{\upsilon}_\varepsilon(y) y^{1-\frac{2}{p}} \in (0,\infty) \text{ is finally decreasing.}
  \end{align*}
  Let us define the function $[2^{c_2},\infty) \ni y \overset{g}{\longmapsto} y^{1-\frac{2}{p}} \in (0,\infty)$.
  Since $h$ is continuously differentiable, it is sufficient to prove that 
  \begin{align*}
  h'(x)\leq 0 
  \Leftrightarrow \frac{\hat{\upsilon}_{\varepsilon}'(x)}{\hat{\upsilon}_{\varepsilon}(x)} 
  \leq -\frac{g'(x)}{g(x)} = \big| 1-\frac{2}{p} \big| \frac{1}{x}, \text{ for finally every }x>2^{c_2}, 
  \end{align*} 
  since $\hat{\upsilon}_{\varepsilon}(x),g(x)>0$ for $x>0$.
  In order to prove this claim, let us recall that\footnote{We use the concavity of $\hat{\upsilon}_{\varepsilon}$ and $\hat{\upsilon}$}
  \begin{align*}
  \hat{\upsilon}_{\varepsilon}'(x) \overset{\eqref{bound:slopes_of_smooth_upsilon}}{\leq} 
  \begin{cases}
  \hat{\upsilon}'_{-}(2^{c_n})
    =\cfrac{\bar{\upsilon}_{c_{n-1}} - \bar{\upsilon}_{c_{n-2}}}{2^{c_{n}}-2^{c_{n-1}}}
    \overset{\eqref{property:bar_upsilon_index_cn_is_theta_value}}{=}
    \cfrac{\vartheta(n-1) - \vartheta(n-2)}{2^{c_{n}}-2^{c_{n-1}}}, 
      &\text{ for }x\in [2^{c_n},2^{c_n}+\varepsilon)\\
  \hat{\upsilon}'(2^{c_n}+\varepsilon)
    =\cfrac{\bar{\upsilon}_{c_n} - \bar{\upsilon}_{c_{n-1}}}{2^{c_{n+1}}-2^{c_n}}
    \overset{\eqref{property:bar_upsilon_index_cn_is_theta_value}}{=}
    \cfrac{\vartheta(n)- \vartheta(n-1)}{2^{c_{n+1}}-2^{c_n}}, 
      &\text{ for }x\in[2^{c_n}+\varepsilon,2^{c_{n+1}}) 
  \end{cases}, \text{for }n\ge 2,
  \end{align*}
where $\hat{\upsilon}'_{-}$ denotes the left-derivative of $\hat{\upsilon}$.
Moreover, $\hat{\upsilon}_{\varepsilon}$ and $\hat{\upsilon}$ are identical on $[2^{c_n}+\varepsilon, 2^{c_{n+1}}-\varepsilon]$; recall \ref{property:identical_on_linear}.
We distinguish between the following cases: 
\begin{enumerate}[label=\textbullet]
  \item Let $x\in[2^{c_n},2^{c_n}+\varepsilon)$, for $n\geq 2$:

On the one hand,
\begin{align*}
  \frac{\hat{\upsilon}_{\varepsilon}'(x)}{\hat{\upsilon}_{\varepsilon}(x)}
    &\leq \frac{\hat{\upsilon}_{-}'(2^{c_n})}{\hat{\upsilon}_{\varepsilon}(2^{c_n})}    
    \overset{\eqref{property:bound_of_error}}{\leq} \frac{\hat{\upsilon}_{-}'(2^{c_n})}{\hat{\upsilon}(2^{c_n}) -2^{-4}}
    = \frac{\vartheta(n-1) - \vartheta(n-2) }{(\vartheta(n-1) - 2^{-4})(2^{c_n}- 2^{c_{n-1}})}\\
    &\overset{\eqref{property:theta_difference}}{<}
    \frac{1}{2(\vartheta(n-1) - 2^{-4})(2^{c_n}- 2^{c_{n-1}})}.
  \numberthis\label{ineq:bound_aux_frac_1}
\end{align*}  
On the other hand,
\begin{align*}
  -\frac{g'(x)}{g(x)} = \big| 1-\frac{2}{p} \big| \frac{1}{x} 
  \geq \big| 1-\frac{2}{p} \big| \frac{1}{2^{c_{n}}+\varepsilon}. 
  \numberthis\label{ineq:bound_aux_frac_2}
\end{align*}
Comparing the right-hand sides of the inequalities \eqref{ineq:bound_aux_frac_1} and \eqref{ineq:bound_aux_frac_2}, we have 
\begin{align*}
\frac{1}{2(\vartheta(n-1) - 2^{-4})(2^{c_n}- 2^{c_{n-1}})} 
  < \big| 1-\frac{2}{p} \big| \frac{1}{2^{c_{n}}+\varepsilon}
  \Leftrightarrow
\frac{2^{c_{n}}+\varepsilon}{2^{c_n}- 2^{c_{n-1}}} 
  < \big| 1-\frac{2}{p} \big| 2(\vartheta(n-1) - 2^{-4}).
\end{align*}
Now, we observe that the last inequality is certainly true for finally all values of $n$, since $\lim_{n\to\infty}\vartheta(n)=\infty$, $[1,\infty) \ni x\longmapsto \frac{1}{1-2^{-x}} \in (0,2]$, \emph{i.e.}, the function is bounded by $2$, and 
$\mathbb{N}\backslash \{1\} \ni n \longmapsto \cfrac{\varepsilon}{2^{c_n}- 2^{c_{n-1}}} \in (0,2^{-4}]$, \emph{i.e.}, the function is bounded by $2^{-4}$.\footnote{We use the fact that $\varepsilon<\frac{1}{2^3}$ and that $c_{n+1}\geq c_n + 1$ with $c_1=1$; recall also the stronger assumption \eqref{property:c_doubling}.}
In other words, for finally all $n\in\mathbb{N}$,
$h'(x)<0$ for $x\in (2^{c_n}+\varepsilon,2^{c_{n+1}}]$.

  \item Let $x\in(2^{c_n}+\varepsilon,2^{c_{n+1}}]$, for $n\geq 2$:

On the one hand, using that $\hat{\upsilon}(2^{c_n}+\varepsilon) = \hat{\upsilon}_{\varepsilon}(2^{c_n}+\varepsilon)\leq \hat{\upsilon}_{\varepsilon}(x)$,
\begin{align*}
  \frac{\hat{\upsilon}_{\varepsilon}'(x)}{\hat{\upsilon}_{\varepsilon}(x)}
  &\leq 
  \frac{\hat{\upsilon}'(2^{c_n}+\varepsilon)}{\hat{\upsilon}(2^{c_n}+\varepsilon)}
  = \frac{\hat{\upsilon}'(2^{c_n}+\varepsilon)}{\vartheta(n-1) + \hat{\upsilon}'(2^{c_n}+\varepsilon) \varepsilon}
  < \frac{\hat{\upsilon}'(2^{c_n}+\varepsilon)}{\vartheta(n-1) }
  = \frac{\ \cfrac{\vartheta(n) - \vartheta(n-1)}{2^{c_{n+1}}-2^{c_n}}\ }{\vartheta(n-1) }\\
  &= \frac{\vartheta(n) - \vartheta(n-1)}{\vartheta(n-1)}\frac{1}{2^{c_{n+1}}-2^{c_n}}
  < \frac{1}{2\vartheta(n-1)}\frac{1}{2^{c_{n+1}}-2^{c_n}}.
  \numberthis\label{ineq:bound_aux_frac_3}
\end{align*}
On the other hand, 
\begin{align*}
  -\frac{g'(x)}{g(x)} = \big| 1-\frac{2}{p} \big| \frac{1}{x} 
  \geq \big| 1-\frac{2}{p} \big| \frac{1}{2^{c_{n+1}}}. 
  \numberthis\label{ineq:bound_aux_frac_4}
\end{align*}
Comparing the right-hand sides of the inequalities \eqref{ineq:bound_aux_frac_3} and \eqref{ineq:bound_aux_frac_4}, we have 
\begin{align*}
  \frac{1}{2\vartheta(n-1)}\frac{1}{2^{c_{n+1}}-2^{c_n}} 
      < \big| 1-\frac{2}{p} \big| \frac{1}{2^{c_{n+1}}}
  &\Leftrightarrow 
  \frac{2^{c_{n+1}}}{2^{c_{n+1}}-2^{c_n}} 
      < 2\big| 1-\frac{2}{p} \big| \vartheta(n-1)\\
  &\Leftrightarrow 
  \frac{1}{1-2^{c_n - c_{n+1}}} 
      < 2\big| 1-\frac{2}{p} \big| \vartheta(n-1).
\end{align*}
Now, we observe that the last inequality is certainly true for finally all values of $n$, since $\lim_{n\to\infty}\vartheta(n)=\infty$ and the function $[1,\infty) \ni  x\longmapsto \frac{1}{1-2^{-x}}\in (0,2]$, \emph{i.e.}, is bounded by $2$.
In other words, for finally all $n\in\mathbb{N}$,
$h'(x)<0$ for $x\in (2^{c_n}+\varepsilon,2^{c_{n+1}}]$.
\end{enumerate}

  \item[\ref{lemma:UI_Young_improvement:p_smaller_than_1}]
  Let $p\in(0,1)$. 
  We can readily adapt the arguments of the previous case to this one. 
  Essentially, nothing changes apart from the substitution of the term $1-\frac{2}{p}$ from the term $1 - \frac{1}{p}$.
  For the appearing functions, all the properties remain the same, especially for the power function $[2^{c_2},\infty) \ni y \overset{g}{\longmapsto} y^{1-\frac{1}{p}} \in (0,\infty)$ its exponent remains negative.
\end{enumerate}
\end{proof}
\begin{lemma}\label{lemma:Young_equiv_submult}
Let $\Upsilon$, resp. $\widehat{\Upsilon}$, be the Young function associated to the right-derivative $\hat{\upsilon}_{\varepsilon}$ defined in \eqref{def:hat_upsilon_varepsilon_smooth} for $p\geq 2$ or \eqref{def:redefined_hat_upsilon_smoothed} for $p\in(0,2)$, resp. $\hat{\upsilon}$ defined in \eqref{def:hat_upsilon_piecewise_linear}.
Then, the Young functions $\Upsilon$ and $\widehat{\Upsilon}$ are equivalent, \emph{i.e.},
there exist $C_1,C_2, x_0>0$ such that
\begin{align*}
\widehat{\Upsilon}(C_1x) \leq \Upsilon(x) \leq \widehat{\Upsilon}(C_2x) \text{ for every }x\geq x_0.
\end{align*}
Moreover, $\Upsilon$ is finally submultiplicative, \emph{i.e.}, there exist $C,R>0$ such that 
\begin{align*}
\Upsilon (xy) \leq C \Upsilon(x) \Upsilon(y) \text{ for all }x,y>R,
\end{align*}
 if and only if $\widehat{\Upsilon}$ is finally submultiplicative.
\end{lemma}
\begin{proof}
  We start be reminding the fact that the definitions \eqref{def:hat_upsilon_varepsilon_smooth} and \eqref{def:redefined_hat_upsilon_smoothed} are identical on $[2,\infty)$.
  Hence, the following arguments apply for both cases.

  In order to prove that $\Upsilon$ and $\widehat{\Upsilon}$ are equivalent we will compare their derivatives $\hat{\upsilon}_{\varepsilon}$ and $\hat{\upsilon}$.
  More precisely, from \citet[Chapter I, Lemma 3.2, p.18]{krasnoselskii1961} we need to check the existence and the positiveness of the limit of their ratio.
  Indeed,
  \begin{align*}
    \lim_{x\to\infty} \frac{\hat{\upsilon}(x)}{\hat{\upsilon}_{\varepsilon}(x)} = 1,
  \end{align*}
  where we used  $0\leq \hat{\upsilon}(x) - \hat{\upsilon}_{\varepsilon}(x)<\frac{1}{2^4}$ for every $x\geq0$, see \eqref{property:error_bound}, and that $\hat{\upsilon}_{\varepsilon}(x) \uparrow \infty$, as $x\uparrow \infty$.

For the validity of the second statement, the reader may consult the comments following  \citet[Chapter I, Lemma 5.1, p30]{krasnoselskii1961}.
\end{proof}
%
%
%
%
%
%
%
%
%
%
%
%
%
%
%
%
%
%
%
\begin{lemma}\label{comput:error_hat_u}
  There exists $0<\varepsilon<\frac{1}{2^3}$ for which the bound \eqref{property:error_bound} holds.
\end{lemma}
\begin{proof}
Let us define the auxiliary function
\begin{align*}
  g(x):= \hat{\upsilon}(x) - \hat{\upsilon}_\varepsilon(x), \text{for }x>0.
\end{align*}
Let $x_0\in\{2^{c_n}:n\in\mathbb{N}\}$.
Then for any such $x_0$ we denote by $\lambda_l x +\beta_l$ the representation of the affine part of $\hat{\upsilon}$ associated to the subinterval on the left of $x_0$ and by $\lambda_r x +\beta_r$ the representation of the affine part of $\hat{\upsilon}$ associated to the subinterval on the right of $x_0$.
Due to continuity of $\hat{\upsilon}$
\begin{align}\label{equality:lambda_beta}
  \lambda_l x_0 +\beta_l=\lambda_r x_0 +\beta_r
  \Leftrightarrow 
  \beta_l - \beta_r = (\lambda_r - \lambda_l)x_0.
\end{align}
Then, for $0<\varepsilon<\frac{1}{2^3}$
\begin{align*}
    g(x_0)=\hat{\upsilon}(x_0) - \hat{\upsilon}_\varepsilon(x_0) 
      &\overset{\phantom{\eqref{equality:lambda_beta}}}{=} 
      (\lambda_l x_0 +\beta_l) - \int_{[-\varepsilon,0]} [\lambda_l (x_0+y) +\beta_l]\eta_\varepsilon(y) \textup{d}y - \int_{[0,\varepsilon]} [\lambda_r (x_0+y) +\beta_r]\eta_\varepsilon(y) \textup{d}y \\
      &\overset{\phantom{\eqref{equality:lambda_beta}}}{=} 
      \int_{[-\varepsilon,\varepsilon]}[\lambda_l (x_0+y) +\beta_l]\eta_\varepsilon(y)\textup{d}y 
        - \int_{[-\varepsilon,0]} [\lambda_l (x_0+y) +\beta_l]\eta_\varepsilon(y) \textup{d}y \\
      &\hspace{4em}  - \int_{[0,\varepsilon]} [\lambda_r (x_0+y) +\beta_r]\eta_\varepsilon(y) \textup{d}y \\
      &\overset{\phantom{\eqref{equality:lambda_beta}}}{=} 
      \int_{[0,\varepsilon]} [(\lambda_l - \lambda_r) (x_0+y) +(\beta_l - \beta_r)]\eta_\varepsilon(y) \textup{d}y \\
      &\overset{\eqref{equality:lambda_beta}}{=} 
      (\lambda_l - \lambda_r) \int_{[0,\varepsilon]} y\eta_\varepsilon(y) \textup{d}y \ge 0,
\end{align*}
given that $\lambda_l\ge\lambda_r$.
The quantity $g(x_0)$ can vanish only  in the case $\lambda_l=\lambda_r$.
Hence,
\begin{align*}
  0\le \lambda_l - \lambda_r
  \le \max\Big\{
    \frac{1}{4}-\frac{1}{2(2^{c_2}-2)}, 
    \sup_{n\in\mathbb{N}}\{\frac{ \bar{\upsilon}_{c_{n+1}}-\bar{\upsilon}_{c_{n}}}{2^{c_{n+1}}-2^{c_n}}
                          -\frac{ \bar{\upsilon}_{c_{n+2}}-\bar{\upsilon}_{c_{n+1}}}{2^{c_{n+2}}-2^{c_{n+1}}}\}
    \Big\}
  \overset{\eqref{upper_bound:slopes}}{\leq}  \frac{1}{4}
\end{align*}
and given also that 
\begin{align*}
  \lim_{\varepsilon \downarrow 0}\int_{[0,\varepsilon]} y\eta_\varepsilon(y) \textup{d}y=0,
\end{align*}
we can chose $0<\varepsilon<\frac{1}{2^3}$ such that 
\begin{align*}
  \int_{[0,\varepsilon]} y\eta_\varepsilon(y) \textup{d}y \le \frac{1}{2^2}.
 \end{align*}
In total, we have 
\begin{align}\label{property:bound_of_error}
  \sup\big\{\hat{\upsilon}(x_0) - \hat{\upsilon}_\varepsilon(x_0) :x_0\in\{2^{c_n}:n\in\mathbb{N}\} \big\}<\frac{1}{2^4}.
\end{align}
We will proceed to estimate the function $g$ in the rest points of its domain:
\begin{enumerate}[label=\textbullet]
  \item 
For $x\in(1,2-\varepsilon)\cup(2+\varepsilon,2^{c_2}-\varepsilon) \cup \big(\bigcup_{k\in\mathbb{N}} (2^{c_n}+\varepsilon, 2^{c_{n+1}}-\varepsilon)\big)$, which is the set considered in \ref{property:identical_on_linear}, we have that $g(x)=0$. 

  \item 
For $x\in(x_0-\varepsilon,x_0]$ it holds
\begin{align*}
  g(x) = \hat{\upsilon}(x) - \hat{\upsilon}_\varepsilon(x) 
  &\overset{\phantom{\eqref{equality:lambda_beta}}}{=} 
  \lambda_l x + \beta_l 
    - \int_{[-\varepsilon,x_0-x]} [\lambda_l(x+y)+\beta_l]\eta_\varepsilon(y)\textup{d}y
    - \int_{[x_0-x,\varepsilon]} [\lambda_r(x+y)+\beta_r]\eta_\varepsilon(y)\textup{d}y\\
  &\overset{\phantom{\eqref{equality:lambda_beta}}}{=}
   \int_{[x_0-x,\varepsilon]} [(\lambda_l-\lambda_r)(x+y)+(\beta_l-\beta_r)]\eta_\varepsilon(y)\textup{d}y\\
  &\overset{\eqref{equality:lambda_beta}}{=}
   (\lambda_l-\lambda_r) \int_{[x_0-x,\varepsilon]} [(x+y)-x_0 ]\eta_\varepsilon(y)\textup{d}y\ge0,
\end{align*}
since $(x+y)-x_0\ge 0$ for $y\in[x_0-x,\varepsilon]$ and the assumption on $x$.
Moreover, 
\begin{align*}
 g'(x) 
 &= (\lambda_l-\lambda_r) \int_{[x_0-x,\varepsilon]} \eta_\varepsilon(y)\textup{d}y
  + (\lambda_l-\lambda_r) x \eta_\varepsilon(x_0-x) +(-x)\eta_\varepsilon(x_0-x)\\
 & = (\lambda_l-\lambda_r) \int_{[x_0-x,\varepsilon]} \eta_\varepsilon(y)\textup{d}y >0,
\end{align*} 
with $\lim_{x\downarrow x_0-\varepsilon}g(x)=0$ and $\lim_{x\uparrow x_0}g(x)=g(x_0)$. 
In other words, on $(x_0-\varepsilon,x_0]$ the function $g$ attains its maximum at $x_0$. 

  \item For $x\in(x_0,x_0+\varepsilon]$ it holds
\begin{align*}
  g(x)= \hat{\upsilon}(x) - \hat{\upsilon}_\varepsilon(x) 
  &\overset{\phantom{\eqref{equality:lambda_beta}}}{=} 
  \lambda_r x + \beta_r 
    - \int_{[-\varepsilon,x_0-x]} [\lambda_l(x+y)+\beta_l]\eta_\varepsilon(y)\textup{d}y
    - \int_{[x_0-x,\varepsilon]} [\lambda_r(x+y)+\beta_r]\eta_\varepsilon(y)\textup{d}y\\
  &\overset{\phantom{\eqref{equality:lambda_beta}}}{=}
   \int_{[-\varepsilon,x_0-x]} [(\lambda_r-\lambda_l)(x+y)+(\beta_r-\beta_l)]\eta_\varepsilon(y)\textup{d}y\\
  &\overset{\eqref{equality:lambda_beta}}{=}
   (\lambda_r-\lambda_l) \int_{[-\varepsilon,x_0-x]} [(x+y)-x_0 ]\eta_\varepsilon(y)\textup{d}y\ge0,
\end{align*}
since $\lambda_r-\lambda_l\le 0$, $(x+y)-x_0\le 0$ for $y\in[x_0-x,\varepsilon]$ and the assumption on $x$.
Moreover, 
\begin{align*}
 g'(x) 
 &= (\lambda_r-\lambda_l) \int_{[x_0-x,\varepsilon]} \eta_\varepsilon(y)\textup{d}y
  + (\lambda_r-\lambda_l) x \eta_\varepsilon(x_0-x) +(-x)\eta_\varepsilon(x_0-x)\\
 & = (\lambda_r-\lambda_l) \int_{[x_0-x,\varepsilon]} \eta_\varepsilon(y)\textup{d}y <0,
\end{align*} 
with $\lim_{x\downarrow x_0}g(x)=g(x_0)$ and $\lim_{x\uparrow x_0+\varepsilon}g(x)=g(x_0)$. 
In other words, on $[x_0,x_0+\varepsilon]$ the function $g$ attains its maximum at $x_0$. 
\end{enumerate}

Therefore,  the function $g$ is non-negative with  
\begin{align*}
  \sup_{x\ge0} g(x) 
    =\sup_{x_0 \in \{2\}\cup\{2^{c_n}:n\in\mathbb{N}\} } \sup_{x\in (x_0-\varepsilon, x_0+\varepsilon)} |\hat{\upsilon}(x) - \hat{\upsilon}_\varepsilon(x) | 
    =\sup_{x_0 \in \{2\}\cup\{2^{c_n}:n\in\mathbb{N}\} } |\hat{\upsilon}(x_0) - \hat{\upsilon}_\varepsilon(x_0) | 
    <\frac{1}{2^4}.
\end{align*}
\end{proof}
\begin{lemma}\label{lem:def_q_for_smoothness}
For $p\in(0,2)$, we define 
\begin{align*}
  q:=\frac{4}{p}-1 \text{ and }x_0(p):=\Big(\frac{p}{6(4-p)}\Big)^{\frac{p}{4-2p}}. 
\end{align*}
Then, $x_0(p)\in[0,1]$ and for the function $[0,\frac{3}{2}]\ni x \overset{f}{\mapsto} x^q$, it holds $f'(x_0(p))=\frac{1}{6}$.
\end{lemma}
\begin{proof}
It is easily verified that
\begin{align*}
   f'(x_0)=\frac{1}{6}
   \Leftrightarrow q x_0^{q-1}=\frac{1}{6}
   \overset{(q>1)}{\Leftrightarrow} x_0^{q-1}=\frac{1}{6q}
   {\Leftrightarrow} x_0=\Big(\frac{1}{6q}\Big)^{\frac{1}{q-1}} = \big(\frac{p}{6(4-p)}\big)^{\frac{p}{4-2p}}.
\end{align*}
The point $x_0$ lies, indeed, in the domain of $f$.
In order to prove the last claim, we will initially prove that $(0,2)\ni p\mapsto x_0(p)\in(0,+\infty)$ is a decreasing function.
To this end, we will calculate its derivative, which is given by
\begin{align*}
  x_0'(p)&=x_0(p)\big(\frac{p}{4-2p} \ln(\frac{p}{6(4-p)})\big)'.
\end{align*}
Given that $x_0(p)>0$, for every $p\in(0,2)$, we proceed with the second factor only.
Hence,
\begin{align*}
  \big(\frac{p}{4-2p} \ln(\frac{p}{6(4-p)})\big)' 
  = \frac{2\ln(\frac{p}{6(4-p)})(4-p) + 2-p}{2(2-p)^2(4-p)}
\end{align*}
The denominator is positive for every $p\in(0,2)$.
The nominator remains negative in the same region of $p$, since for $p\in(0,2)$ 
\begin{align*}
 2\ln(\frac{p}{6(4-p)})(4-p) + 2-p < 2\ln(\frac{p}{6(4-p)})(4-p) + 4-p
 <(2\ln(\frac{p}{6(4-p)})+1)(4-p)<0.
\end{align*} 

Finally, we compute the limits
\begin{align*}
  \lim_{p\downarrow 0} x_0(p) 
    & = \lim_{p\downarrow 0} \big(\frac{p}{6(4-p)}\big)^{\frac{p}{4-2p}}
      = \lim_{p\downarrow 0} \exp\big(\frac{p}{4-2p}\ln(\frac{p}{6(4-p)})\big)\\
    & = \lim_{z\downarrow 0} \exp\big(z\ln(z)\big) = 1
\end{align*}
and 
\begin{align*}
  \lim_{p\uparrow 2} x_0(p) 
    & = \lim_{p\uparrow 2} \big(\frac{p}{6(4-p)}\big)^{\frac{p}{4-2p}} = 0,
\end{align*}
which in conjunction with the monotonicity proved above, verify our claim.
\end{proof}
\subsection{Proof of Corollary \ref{corollary:UI_Young_improvement:StrongerIntegrability}}

\begin{proof}[Proof of \cref{corollary:UI_Young_improvement:StrongerIntegrability}]
 We will adopt the notation introduced in \cref{lemma:UI_Young_improvement}. 
 The arguments presented below hold for $p>0$.

 \vspace{0.5em}
 For notational convenience we define $\bar{\upsilon}_0=\frac{1}{4},\bar{\upsilon}_{\frac{1}{2}}=\frac{1}{2}$, so that the non-decreasing sequence $(\bar{\upsilon}_n)_{n\in\{0,\frac{1}{2}\}\cup\mathbb{N}}$ associates to the right derivative $\bar{\upsilon}$.
 In view of the introduced notation, we have 
 \begin{align*}
    \overline{\Upsilon}(x)\le 
    \begin{cases}
      \bar{\upsilon}_0 x, &\text{for }x\in[0,1]\\
      \bar{\upsilon}_0 + \bar{\upsilon}_{\frac{1}{2}}(x-1), &\text{for }x\in(1,2]\\
      \bar{\upsilon}_0 + \bar{\upsilon}_{\frac{1}{2}} + \sum_{k=1}^{n-1}\bar{\upsilon}_k + \bar{\upsilon}_{n}(x-n), &\text{for }x\in(n,n+1]
    \end{cases}.
 \end{align*}

 \vspace{0.5em}
 Following the arguments of the proof of \cite[Lemme]{meyer1978sur} and modifying them as in the proof of \cref{prop:UI_Young}, we can prove that for the Young function $\overline{\Upsilon}$ associated to the right-derivative $\bar{\upsilon}$, it holds
 \begin{align}\label{limit:overline_upsilon}
   \lim_{R\to+\infty}\sup_{\alpha\in\mathcal{A}}\int_{{\{z:|z|>R\}}} \overline{\Upsilon}(|z|^p)m_\alpha(\textup{d}z)=0.
 \end{align}

 Indeed, for $c>1$ and $k_0:=\min\{k\in\mathbb{N}  : \bar{\upsilon}_0 + \bar{\upsilon}_{\frac{1}{2}} + \sum_{q=1}^k \bar{\upsilon}_q>c\}$
 \begin{align*}
   &\int_{\{z:\overline{\Upsilon}(|z|^p)>c\}} \overline{\Upsilon}(|z|^p) m_\alpha (\textup{d}z)
    \le \sum_{k=k_0}^\infty (\bar{\upsilon}_0 + \bar{\upsilon}_{\frac{1}{2}} + \sum_{q=1}^k \bar{\upsilon}_q) m_\alpha(\{k< |z|^p \le k+1\})\\
    &\hspace{2em}= (\bar{\upsilon}_0 + \bar{\upsilon}_{\frac{1}{2}} + \sum_{q=1}^{k_0} \bar{\upsilon}_q) m_\alpha(\{|z|>k_0\}) + \sum_{k=k_0+1}^\infty \bar{\upsilon}_k m_\alpha(\{|z|^p>k\}).
    \numberthis \label{ineq:rhs_bar_upsilon}
 \end{align*}
 Let us choose an arbitrarily small $\delta>0$.
 Then, there exists $n_0$ such that for $n\ge n_0\ge 1$ it holds $ \frac{\bar{\upsilon}_n}{\upsilon_n}<\delta$.
 Bounding the right part of Inequality \eqref{ineq:rhs_bar_upsilon}, we have for $k_0>n_0$\footnote{The reader may observe that $\lim_{c\to+\infty}k_0 = +\infty$.}
 \begin{align*}
   &(\bar{\upsilon}_0 + \bar{\upsilon}_{\frac{1}{2}} + \sum_{q=1}^{k_0} \bar{\upsilon}_q) m_\alpha(\{k_0 < |z|^p\le k_0+1\}) + \sum_{k=k_0+1}^\infty \bar{\upsilon}_k m_\alpha(\{k< |z|^p \le k+1\})\\
    &\hspace{2em}\le (\bar{\upsilon}_0 + \bar{\upsilon}_{\frac{1}{2}} + \sum_{q=1}^{n_0} \bar{\upsilon}_q) m_\alpha(\{|z|^p>k_0\}) 
                      + \delta \sum_{q=n_0+1}^{k_0} \upsilon_q m_\alpha(\{|z|^p>k_0\})
                      + \delta \sum_{k=k_0+1}^\infty \upsilon_k m_\alpha(\{|z|^p>k\})\\
    &\hspace{2em}\le (\bar{\upsilon}_0 + \bar{\upsilon}_{\frac{1}{2}} + \sum_{q=1}^{n_0} \bar{\upsilon}_q) m_\alpha(\{|z|^p>k_0\}) 
                      +\delta \sum_{k=n_0+1}^\infty {\upsilon}_k m_\alpha(\{|z|^p>k\}).
    \numberthis \label{ineq:rhs_limit_k0}
 \end{align*}
 From the proofs of \cref{prop:UI_Young} and \cref{lemma:UI_Young_improvement} we have for the sequence $(\upsilon_n)_{n\in\mathbb{N}}$
  \begin{align*}
    S:=\sup_{\alpha\in\mathcal{A}} \sum_{n=1}^{\infty}\upsilon_n m_\alpha(\{z:|z|^p>k\}) <+\infty. 
  \end{align*}
  In total, in view of \eqref{property:limit_for_UI}, which implies 
  \begin{align*}
    \lim_{k\to+\infty} \sup_{\alpha\in\mathcal{A}} m_\alpha(\{z:|z|^p>k\}) = 0,
  \end{align*}
  we have from \eqref{ineq:rhs_bar_upsilon}, \eqref{ineq:rhs_limit_k0} and the above limit
  \begin{align*}
    &\lim_{R\to+\infty}\sup_{\alpha\in\mathcal{A}}\int_{{\{z:|z|>R\}}} \overline{\Upsilon}(|z|^p)m_\alpha(\textup{d}z) 
    \le (\bar{\upsilon}_0 + \bar{\upsilon}_{\frac{1}{2}} + \sum_{q=1}^{n_0} \bar{\upsilon}_q) \lim_{k\to+\infty} \sup_{\alpha\in\mathcal{A}} m_\alpha(\{z:|z|^p>k\}) + \delta S=\delta S. 
  \end{align*}
  Since $\delta$ was chosen arbitrarily small, we have that the limit on the left hand side of the last inequality is indeed $0$.

 \vspace{0.5em}
 In view of the validity of \eqref{limit:overline_upsilon} and recalling \cref{prop:UI_Young}, there exists a Young function $\Psi$ (with right derivative $\psi$) such that 
 \begin{align*}
    \sup_{\alpha\in\mathcal{A}}\int_{\{z:|z|\ge 1\}} \Psi\big(\overline{\Upsilon}(|z|^p)\big)m_\alpha (\textup{d}z)<+\infty.
 \end{align*}
 Now, observing that $\hat{\upsilon}_\varepsilon(w)<\bar{\upsilon}(w)$, for every $w\ge 2$, we can extract a point $w_0>2$ such that
  \begin{align*}
    \Upsilon(w)<\overline{\Upsilon}(w) \text{ for every }w\ge w_0. 
  \end{align*}  
 Therefore, we get
  \begin{align*}
    &\sup_{\alpha\in\mathcal{A}}\int_{\{z:|z|\ge 1\}} \Psi\big(\Upsilon(|z|^p)\big)m_\alpha (\textup{d}z)\\
    &\hspace{1em}=\sup_{\alpha\in\mathcal{A}}\Big\{ \int_{\{z:1\le |z|^p\le w_0\}} \Psi\big(\Upsilon(|z|^p)\big)m_\alpha (\textup{d}z) + 
                                                    \int_{\{z:|z|^p> w_0\}} \Psi\big(\Upsilon(|z|^p)\big)m_\alpha (\textup{d}z)\Big\}\\
    &\hspace{1em}\le\sup_{\alpha\in\mathcal{A}}\Big\{ \int_{\{z:1\le |z|^p\le w_0\}} \psi\big(\Upsilon(w_0)\big) \Upsilon(|z|^p) m_\alpha (\textup{d}z) + 
                                                      \int_{\{z:|z|^p> w_0\}} \Psi\big(\overline{\Upsilon}(|z|^p)\big)m_\alpha (\textup{d}z)\Big\}\\
    &\hspace{1em}\le  \psi\big(\Upsilon(w_0)\big)  
                      \sup_{\alpha\in\mathcal{A}}  \int_{\{z:|z|\ge 1\}} \Upsilon(|z|^p) m_\alpha (\textup{d}z) + 
                      \sup_{\alpha\in\mathcal{A}}  \int_{\{z:|z|\ge 1\}} \Psi\big(\overline{\Upsilon}(|z|^p)\big)m_\alpha (\textup{d}z)\Big\}
    <+\infty.
  \end{align*}  
  where we have used that $\Psi(x)\le x\psi(x)$ and $\psi$ is non-decreasing.
\end{proof}

\bibliographystyle{abbrvnat}
\small
\bibliography{bibliographyDylan}

\end{document}